\documentclass[11pt,reqno]{amsart}
\usepackage[letterpaper,margin=1in]{geometry}
\usepackage[T1]{fontenc}
\usepackage[utf8]{inputenc}
\usepackage{lmodern,microtype}
\usepackage{amsmath,amssymb,amsthm,mathtools,mathrsfs}
\usepackage{enumitem,booktabs,tabularx,array,longtable,xurl,needspace,float}
\usepackage[section]{placeins}
\usepackage[table]{xcolor}
\usepackage{tikz,tikz-cd}
\usetikzlibrary{arrows.meta,positioning,calc,decorations.pathreplacing,fit}
\definecolor{Ink}{HTML}{20384A}
\definecolor{Teal}{HTML}{147D85}
\definecolor{Amber}{HTML}{926018}
\definecolor{Blue}{HTML}{285F91}
\definecolor{Purple}{HTML}{76528D}
\definecolor{Coral}{HTML}{B54640}
\definecolor{Pale}{HTML}{F0F5F6}
\usepackage[colorlinks=true,linkcolor=Blue,citecolor=Teal,urlcolor=Blue]{hyperref}
\hypersetup{pdftitle={Counterexamples to Haiman's conjectures in every dimension greater than two},pdfauthor={Eren Ercan},bookmarksnumbered=true,pdfstartview=FitH,pdfkeywords={Haiman conjectures, isospectral Hilbert scheme, polygraph, diagonal harmonics}}
\setlist{itemsep=3pt,topsep=5pt,leftmargin=2.2em}
\allowdisplaybreaks[2]
\theoremstyle{plain}
\newtheorem{theorem}{Theorem}[section]
\newtheorem{maintheorem}{Theorem}

\newtheorem{lemma}[theorem]{Lemma}
\newtheorem{proposition}[theorem]{Proposition}
\newtheorem{corollary}[theorem]{Corollary}
\theoremstyle{definition}

\newtheorem{problem}[theorem]{Problem}
\theoremstyle{remark}
\newtheorem{remark}[theorem]{Remark}

\newtheorem*{theorem*}{Theorem}
\newtheorem*{lemma*}{Lemma}
\newtheorem*{proposition*}{Proposition}
\newtheorem*{corollary*}{Corollary}
\providecommand{\Q}{\mathbb Q}
\providecommand{\C}{\mathbb C}

\providecommand{\N}{\mathbb N}

\providecommand{\A}{\mathbb A}

\providecommand{\Sym}{\operatorname{Sym}}
\providecommand{\Spec}{\operatorname{Spec}}

\providecommand{\Hom}{\operatorname{Hom}}

\providecommand{\Aut}{\operatorname{Aut}}
\providecommand{\Ext}{\operatorname{Ext}}
\providecommand{\Tor}{\operatorname{Tor}}

\providecommand{\sgn}{\operatorname{sgn}}
\providecommand{\Hilb}{\operatorname{Hilb}}

\providecommand{\Tr}{\operatorname{Tr}}
\providecommand{\Gr}{\operatorname{Gr}}

\providecommand{\ord}{\operatorname{ord}}

\numberwithin{equation}{section}
\newcommand{\resultref}[1]{\hyperref[#1]{\textbf{Theorem~\ref*{#1}}}}
\newcommand{\proofref}[1]{\hyperref[#1]{p.~\pageref*{#1}}}
\newcommand{\readinglink}[2]{\hyperref[#1]{\textcolor{Blue}{#2}}}

\title[Counterexamples to Haiman's conjectures]
{Counterexamples to Haiman's conjectures in every dimension greater than two}
\author{Eren Ercan}
\date{30 November 2026}
\begin{document}
\begin{abstract}
We disprove Haiman's polygraph freeness conjecture and the normality
assertion for principal isospectral Hilbert schemes in every ambient
dimension greater than two. Punctual divisors at $(n,d)=(96,3)$,
$(20,4)$ and $(10,5)$ give counterexamples at all greater lengths
and dimensions. The ordered covers fail $R_1$, seminormality and
the Du Bois property along these divisors. An invariant Koszul
complex transfers this obstruction to marked polygraphs, whose
eventual projective dimensions at the three lengths are $94$,
$18$ and $8$.

The alternating pairing has the same ambient boundary. Every
minimum-degree determinant in the plane has apolar dimension $n!$.
In each fixed dimension $d\geq3$, the derivative span of all
minimum-degree alternants has exponentially vanishing dimension
relative to $n!$, and is stable under every positive-order
generalized polarization. This conclusion extends to a growing
degree band; in dimension three, a band reaching $101/100$ of the
minimum degree has exponential rate below $49/50$.

We describe transverse orbit rings and canonical modules, and determine
the conductors when the principal component is regular at the divisor. At any homogeneous punctual divisor of
length $n\geq3$, the normalization defect satisfies
$\delta\geq\mu n!h/2$, where $\mu$ is the least alternating degree
and $h$ is the unordered reconstruction degree. At length ten the
tangent cone has $10!$ reduced lines, embedding dimension $46$
and conductor exponent $13$. We also compute invariant Tor
algebras and prove that their exterior products persist in
marked polygraphs.
\end{abstract}
\maketitle
% Initialize the math Roman font before its first shipped use.
\begingroup\fontencoding{OT1}\fontfamily{lmr}\fontsize{9}{11}\selectfont\endgroup
\enlargethispage{\baselineskip}

\Needspace{10\baselineskip}
\section{Introduction}\label{sec:introduction}
\subsection*{State of the problem}
A polygraph records every way of choosing $\ell$ marked points,
with repetition, from an ordered configuration of $n$ base points.
Over a characteristic-zero field, its reduced scheme is
\[
 Z(n,\ell;d)=\bigcup_{f:[\ell]\to[n]}
 \{(x_1,\ldots,x_n,z_1,\ldots,z_\ell):z_a=x_{f(a)}\}
 \subset(\A^d)^{n+\ell}.
\]
Put
\[
 S=\C[x_{ij}:1\le i\le n,\ 1\le j\le d],\quad
 P_j=\C[x_{1j},\ldots,x_{nj}],\quad
 R_{n,\ell,d}=\C[Z(n,\ell;d)].
\]
The subring $P_j$ is one \emph{coordinate alphabet}: it retains
one coordinate of each base point. Haiman proved that the planar
polygraph ring is free over either alphabet. This theorem is a
central ingredient in his proof of Macdonald positivity
\cite{gg-HaimanPolygraph}.

Let $H_{n,d}$ be the reduced principal component of $\Hilb^n(\A^d)$,
which parametrizes limits of $n$ distinct points, and let
\[
 X_{n,d}=
 \bigl(H_{n,d}\times_{\Sym^n(\A^d)}(\A^d)^n\bigr)_{\mathrm{red}}.
\]
The isospectral cover $X_{n,d}$ records an ordering of the
Hilbert--Chow cycle. Haiman's Conjecture~5.2.1 predicts that these
covers are normal and Cohen--Macaulay in every ambient dimension.
His Conjecture~5.2.2 predicts that every $R_{n,\ell,d}$ is free
over $P_j$ \cite[\S5.2]{gg-HaimanPolygraph}. Haiman called the first
rather speculative, and wrote that failure of the second would surprise him.

Lee proved that $H_{9,8}$ is not Cohen--Macaulay at $\mathfrak m^2$
\cite[Theorem~A]{dim:Lee}. Since Cohen--Macaulayness descends through
$H_{n,d}=X_{n,d}/S_n$, this refutes the Cohen--Macaulay assertion of
Conjecture~5.2.1. We settle the normality and polygraph assertions
in every ambient dimension greater than two.

\begin{table}[!htbp]
\centering
\caption{Principal results. The theorem labels and proof locations are linked.
The bounds on $n$ give counterexamples in the indicated ranges; the sharp
least length in the third row concerns failure of $R_1$.}
\label{tab:principal-results}
\small
\setlength{\tabcolsep}{5pt}
\renewcommand{\arraystretch}{1.2}
\begin{tabularx}{\textwidth}{@{}>{\raggedright\arraybackslash}p{.22\textwidth}>{\raggedright\arraybackslash}X >{\raggedright\arraybackslash}p{.16\textwidth}@{}}
\toprule
\rowcolor{Pale}
\textbf{Result} & \textbf{Conclusion and scope} & \textbf{Statement / proof}\\
\midrule
Haiman's conjectures &
Divisorial nonnormality and eventual polygraph nonfreeness for
$d=3$, $n\geq96$; $d=4$, $n\geq20$; and $d\geq5$, $n\geq10$. &
\resultref{fh:j:main}\newline
\proofref{dim:all-dimensional-counterexamples}\\
\addlinespace[3pt]
Alternating pairing &
Planar apolar dimension $n!$; an exponential deficit in every fixed
$d\geq3$, including all minimum-degree alternants and a band above them. &
\resultref{dim:intro-pairing}\newline
\proofref{dim:apolar}\\
\addlinespace[3pt]
Least divisorial failure &
For $n\leq10$, failure of $R_1$ occurs exactly when
$n=10$ and $d\geq5$. &
\resultref{new:sharp}\newline
\proofref{new:small}\\
\addlinespace[3pt]
Transverse local ring &
At $(10,5)$: $10!$ branches, embedding dimension $46$, conductor
$t^{13}\widetilde B$, and
$23\,738\,952\leq\delta\leq29\,746\,080$. &
\resultref{new:quantitative}\newline
\proofref{new:orbit}\\
\addlinespace[3pt]
Canonical module &
At every homogeneous punctual divisor of length $n\geq3$:
$\delta\geq\mu n!h/2$ and a canonical module for arbitrary conductor;
$\mu$ is the minimum alternating degree and $h$ the reconstruction degree. &
\resultref{dim:general-canonical-module}\newline
\proofref{dim:graded-defect-section}\\
\addlinespace[3pt]
Characters of alternating quotients &
For point Hilbert function $(1,d,r)$, every hook survives and only
Young diagrams with at most $r$ boxes outside the first row and column occur. &
\resultref{next:intro-apolar}\newline
\proofref{next:uniform-apolar}\\
\addlinespace[3pt]
Invariant Tor products &
An exterior algebra on $n-2$ invariant classes at the orbit cones;
nonzero top products transfer to the marked polygraphs. &
\resultref{next:intro-tor}\newline
\proofref{next:invariant-koszul}\\
\addlinespace[3pt]
Isospectral local rings &
At fixed length, every completion on $X_{n,d}$ comes from one multiplication
chart up to formal variables. Its vertex has maximal Cohen--Macaulay
defect, at least $\lfloor n/9\rfloor$. &
\resultref{fh:j:cone}\newline
\hyperref[fh:j:defect-growth]{Corollary~\ref*{fh:j:defect-growth}}\\
\bottomrule
\end{tabularx}
\end{table}

\subsection*{Counterexamples and the ambient boundary}

\Needspace{20\baselineskip}
\begin{maintheorem}[Counterexamples in every dimension above two]
\label{fh:j:main}
\enlargethispage{\baselineskip}
Set
\[
 N(3)=96,\qquad N(4)=20,\qquad N(d)=10\quad(d\geq5).
\]
For every $d\geq3$ and $n\geq N(d)$, the following assertions hold.
\begin{enumerate}[label=\textup{(\roman*)}]
\item $X_{n,d}$ has a divisor along which it fails $R_1$,
seminormality and the Du Bois property. Its general cycle has
one point of multiplicity $N(d)$ and $n-N(d)$ simple points.
For $d\geq5$ this divisor also gives failure of Gorensteinness.
\item Projection to every coordinate alphabet has a fibre of
dimension at least $n(d-1)+N(d)-2$, while the generic fibre has
dimension $n(d-1)$.
\item For each seed $m\in\{10,20,96\}$ there is a mark constant
$L_m$, independent of extra points, extra ambient coordinates
and the chosen alphabet, such that all $\ell\geq L_m$ give
nonflat polygraph rings in the corresponding range.
At the seed length $m$, their projective and flat dimensions are
exactly $m-2$, with invariant Tor in every degree $1,\ldots,m-2$.
\item If $J=S\cdot S^{\sgn}$ is the alternating ideal, all
sufficiently large powers $J^r$ are nonflat over every $P_j$, and
\begin{equation}\label{fh:j:ideal}
 J^r\subsetneq\overline{J^r}\subseteq\bigcap_{i<k}I_{ik}^{\,r},
 \qquad I_{ik}=(x_{i1}-x_{k1},\ldots,x_{id}-x_{kd}).
\end{equation}
At the seed lengths these powers have projective dimension $m-2$.
\end{enumerate}
\end{maintheorem}

Combined with Haiman's planar theorem and the elementary
one-dimensional calculation, this gives the sharp ambient boundary
\[
\begin{split}
 X_{n,d}\text{ is normal for every }n
     &\quad\Longleftrightarrow\quad d\leq2,\\
 k[Z(n,\ell;d)]\text{ is free over every alphabet for all }n,\ell
     &\quad\Longleftrightarrow\quad d\leq2.
\end{split}
\]
Among lengths at least three, $96$ and $20$ are the first lengths
admitting a homogeneous full punctual divisor in dimensions three
and four, respectively. Every such divisor automatically has maximal
point Hilbert function. Homogeneous algebras with that Hilbert function
are parametrized by $\Gr(r,N_s)$, where $s$ is the top degree,
$N_s=\binom{d+s-1}{d-1}$ and $r$ is the last Hilbert-function value.
A homogeneous punctual divisor can therefore exist only if
\[
 r(N_s-r)\geq nd-d-1.
\]
At the three seeds this condition reads
\[
\begin{aligned}
 (n,d)=(10,5):&\quad\dim\Gr(4,15)=44=5\cdot10-6,\\
 (n,d)=(20,4):&\quad\dim\Gr(5,20)=75=4\cdot20-5,\\
 (n,d)=(96,3):&\quad\dim\Gr(12,36)=288\geq3\cdot96-4=284,
\end{aligned}
\]
and Proposition~\ref{d3:homogeneous-optimal} shows that it fails at every
smaller length $n\geq3$ in dimensions three and four. The double-point
divisor at $n=2$ satisfies the same condition, with cyclic inertia $S_2$.

\begin{problem}\label{prob:least-length}
Determine the least $n$ for which $X_{n,3}$, respectively $X_{n,4}$,
is not normal. In particular, is $X_{n,3}$ normal for every $n\leq95$,
and is $X_{n,4}$ normal for every $n\leq19$?
\end{problem}

The ten-point construction makes the normality obstruction visible.
Ga\l\k{a}zka--Keneshlou--\v{S}ivic prove that every homogeneous algebra
with Hilbert function $(1,5,4)$ is smoothable
\cite[Proposition~4.2]{new:GKS}. The homogeneous locus is
$\Gr(4,15)$, so its translates form a divisor $D\subset H_{10,5}$:
\[
 \dim D=\dim\Gr(4,15)+5=4(15-4)+5=49=\dim H_{10,5}-1.
\]
Above its generic point, the ordered cover has inertia $S_{10}$.
If the corresponding height-one local ring were normal, it would be
a discrete valuation ring. In characteristic zero, its inertia acts
faithfully on the one-dimensional cotangent space and is therefore
cyclic. Since $S_{10}$ is noncyclic, $X_{10,5}$ fails $R_1$.

Projection sends $D$ into the one-dimensional small diagonal of
the ten-dimensional coordinate alphabet. Thus
\[
 \dim(\text{fibre over the small diagonal})\geq49-1=48
 >40=50-10=\dim(\text{generic fibre}).
\]
The alternating Rees charts and marked alternating summands carry
this dimension jump to the polygraph rings. Section~\ref{new:polygraph-homology}
identifies the invariant Tor classes that survive this passage.

\begin{figure}[!htbp]
\centering
\begin{minipage}[t]{.48\textwidth}
\centering\small
\textcolor{Ink}{\textbf{(a) Projection to one coordinate alphabet}}
\medskip
\[
\begin{tikzcd}[column sep=8mm,row sep=9mm]
 D_X \arrow[r,hook] \arrow[d] & X_{10,5}\arrow[d,"\pi_j"]\\
 \Delta_j\simeq\A^1\arrow[r,hook] & \A^{10}
\end{tikzcd}
\]
\[
\begin{array}{c|cc}
 &D_X&X_{10,5}\\ \hline
 \dim&49&50\\
 \dim\pi_j(\,\cdot\,)&1&10\\ \hline
 \text{fibre dimension}&\textcolor{Teal}{48}&40
\end{array}
\]
\end{minipage}\hfill
\begin{minipage}[t]{.48\textwidth}
\centering\small
\textcolor{Ink}{\textbf{(b) Inertia at the generic punctual point}}
\medskip
\[
 I_\xi=S_{10},\qquad F=k(\xi).
\]
\[
\begin{gathered}
 \mathcal O_{X,\xi}\text{ normal}\\
 \Downarrow\\[-2pt]
 \mathcal O_{X,\xi}\text{ a discrete valuation ring}\\
 \Downarrow\\[-2pt]
 S_{10}\hookrightarrow
 \Aut_F(\mathfrak m_\xi/\mathfrak m_\xi^2)=F^\times
\end{gathered}
\]
\textcolor{Coral}{\(S_{10}\) is noncyclic, so this injection is impossible.}
\end{minipage}
\caption{The ten-point divisor gives both obstructions.
Here $D_X$ is the reduced inverse image in $X_{10,5}$ of the
translated $(1,5,4)$ locus, and $\xi$ is its generic point.
In (a), subtracting the image dimensions gives an excess of eight
fibre dimensions. In (b), tame inertia at a normal height-one point
acts faithfully on a one-dimensional cotangent space and must be
cyclic. The polygraph argument transfers the excess dimension
to invariant Tor in Section~\ref{new:polygraph-homology}.}
\label{fig:normality-overview}
\end{figure}

The new divisors come from a single collision map. If degree $s-1$
polynomials evaluate injectively on $n$ points and degree $s$
polynomials evaluate surjectively, let $W$ be the degree-$s$ leading
forms of their vanishing ideal. Scalar contraction gives
$(W)+\mathfrak m^{s+1}$. Translations and dilation leave $W$ fixed.
A differential of rank $nd-d-1$ therefore proves that the
translated limits form a divisor. Nonzero minors modulo a prime certify this rank over $\Q$ for every length $20$ through $26$
in dimension four, and $96$ through $105$ in dimension three.
At each such length the eventual polygraph projective dimension
attains its general upper bound $n-2$.

The collision image also determines the entire homogeneous smoothable
locus, which is irreducible. For Hilbert functions $(1,4,10,r)$,
$r=5,\ldots,11$, its codimensions are
\[
 0,\ 5,\ 8,\ 9,\ 8,\ 5,\ 0.
\]
Every algebra at either endpoint is smoothable; at every intervening
Hilbert function, a general algebra is nonsmoothable.

\subsection*{Alternating quotients and diagonal harmonics}
The alternating pairing has a corresponding dimensional distinction.
Write $\mu(n,d)$ for the minimum total degree of an alternating
polynomial in $n$ blocks of $d$ variables. A minimum-degree evaluation
determinant uses all monomials below one degree and a subspace in that
last degree. Let $G$ be its apolar algebra, namely the quotient of the
constant-coefficient differential operators by those that annihilate
the determinant.

\begin{maintheorem}[The alternating pairing and the minimum-degree component]
\label{dim:intro-pairing}
For $d=1,2$, every such determinant has $\dim G=n!$.
For each fixed $d\geq3$, uniformly over every last-degree subspace,
\[
 \limsup_{n\to\infty}\left(\frac{\dim G}{n!}\right)^{1/n}
 \leq\rho_d:=\frac{e^{d+1}d!}{2^d(d+1)^d}<1.
\]
The same bound holds for the derivative span of all minimum-degree
alternants together, and for the derivative span of every alternating
polynomial of degree at most $\mu(n,d)+o(n^{1+1/d})$.
These spans are stable under all generalized polarizations
$\sum_i z_{ia}f(\partial_{i1},\ldots,\partial_{id})$ of positive
differential order. The minimum-degree span lies in diagonal harmonics.

In dimension three, let $\mathcal D_n$ be the derivative span of
all alternating polynomials of degree at most
$\lfloor101\mu(n,3)/100\rfloor$. Then
\[
 \limsup_{n\to\infty}
 \left(\frac{\dim\mathcal D_n}{n!}\right)^{1/n}<\frac{49}{50}.
\]
\end{maintheorem}

The explicit constants in the first bound are
$\rho_3=3e^4/256<0.640$ and $\rho_4=3e^5/1250<0.357$.
The point-algebra relations sharpen them to less than $0.630068$
and $0.352215$ for the minimum-degree component. These bounds concern
the same alternating multiplication pairing that is perfect in the
plane. At the fat points $k[x_1,\ldots,x_4]/\mathfrak m^5$ and
$k[x_1,x_2,x_3]/\mathfrak m^8$, integer coefficient sums show that
the pairing loses more than $98.7819\%$ and $79.1989\%$ of its generic
factorial rank. Here $G$ denotes the image of the alternating pairing on the
isospectral fibre.

\subsection*{Transverse local rings and canonical modules}
The local geometry has a uniform description even when the principal
Hilbert component is singular. At any homogeneous punctual divisor
of length $n\geq3$,
radial normalization gives $b=e=n!h$, where $h$ is the finite degree
of reconstruction of the unordered configuration up to translation
and dilation. An invariant rational uniformizer and sign trace give
\[
 \delta\geq\frac{\mu}{2}\,n!h,\qquad
 \operatorname{length}(K/R)=2\delta-\mu n!h,
\]
where $K$ is an explicit canonical fractional module and $K/R$ has
perfect sign-valued duality. An invariant parameter of weight one
belongs to the local ring. Its reduction identifies the image
of the canonical pairing with the determinant apolar algebra $G$,
so the stronger estimate is
\[
 2\delta-\mu n!h\geq n!h-\dim G\geq n!(h-1).
\]
These statements allow singular principal components and arbitrary
conductor. In particular,
\[
 \delta_{20,4}\geq47\,441\,589\,159\,444\,480\,000,
 \qquad \delta_{96,3}\geq231\cdot96!.
\]
At the twenty-point divisor the entire negative tangent is
concentrated in weight one. After removing translations it has
dimension $21$. The ordered tangent cone is consequently a
reduced union of exactly $20!h$ lines in embedding dimension at most
$97$. The reconstruction degree satisfies $h=1$ exactly when the
principal component is regular at the generic divisor. At length
$20$ it is bounded by $1\leq h\leq2^{30}$.

\begin{problem}\label{prob:reconstruction}
Determine the reconstruction degrees at lengths $20$ and $96$.
Is $h=1$, so that the principal component is regular at the generic
point of each divisor?
\end{problem}

\begin{maintheorem}[The smallest divisorial counterexample]
\label{new:sharp}
For $1\le n\le10$ and every $d\ge1$,
\[
 X_{n,d}\text{ satisfies }R_1
 \quad\Longleftrightarrow\quad n\le9\text{ or }d\le4.
\]
Thus ten is the least length admitting a failure of $R_1$,
and five is the least ambient dimension admitting that failure
at length ten.
\end{maintheorem}

The codimension-one singularity admits an explicit algebraic
description. Let $B$ be the local ring above the generic point
of the centered $(1,5,4)$ divisor, with residue field
$F=\C(\Gr(4,15))$. Its normalization is $E[t]_{(t)}$, where
$E/F$ is an $S_{10}$-Galois extension of degree $10!$.
The residue coordinates $v_{ij}=x_{ij}/t$ describe the velocities
of the ordered points under dilation. If $F_r\subset E$ consists
of polynomial functions of degree at most $r$ in these velocities,
then
\begin{equation}\label{new:intro-cone}
 B\simeq R_{R_+},\qquad
 R=\bigoplus_{r\ge0}t^rF_r\subset E[t],\qquad
 \operatorname{gr}_{\mathfrak m_B}B\simeq R.
\end{equation}
Thus the tangent cone itself determines the local ring at its vertex.
After extending $F$ to a separable closure, the projective point set
of this cone is the $10!$ orderings of one generic configuration.

\begin{maintheorem}[The transverse singularity]
\label{new:quantitative}
For the preceding local ring $B$, with normalization $C$,
\[
 b(B)=e(B)=10!=3\,628\,800,\qquad
 \operatorname{edim}B=46,\qquad
 \mathfrak c_{B/C}=t^{13}C.
\]
Its tangent cone is geometrically reduced with $10!$ distinct lines.
Writing $\delta(B)=\operatorname{length}_B(C/B)$, one has
\[
 23\,738\,952\le\delta(B)\le29\,746\,080.
\]
Its seminormalization is $F+tC$, and its local ring is not Gorenstein.
\end{maintheorem}

The exponent $13$ comes from the least possible degree of an
alternant in ten blocks of five variables:
\[
 0+5\cdot1+4\cdot2=13.
\]
A degree-controlled interpolation basis proves that every function
on the ordered orbit occurs by this same degree. The resulting
sign trace pairing makes the failure of Gorensteinness explicit.
The conjugate partitions $(5,5)$ and $(2^5)$ both have Specht
dimension $42$. The first has no copies in $F_4$, while the second
has at most $15$ copies in $F_8$. Gorenstein duality would require
these two multiplicities to sum to $42$.

\subsection*{Characters and invariant exterior products}
The orbit algebra also has a canonical Gorenstein quotient.
For a generic configuration let $A=R/(t)$, let $\mu$ be the
minimum alternating degree, and let $\lambda$ be its nonzero
sign-valued functional in degree $\mu$. Define
\[
 G=A/\{a:\lambda(ab)=0\text{ for every }b\in A\}.
\]
The dimension of $A$ is $n!$, and $G$ has sign socle in degree
$\mu$. Its size records how much of the orbit algebra is
generated by differentiating the minimal alternating polynomial.

\begin{maintheorem}[The alternating quotient]
\label{next:intro-apolar}
Let $n=d+1+r$ generic points in $\A^d$ satisfy
$n\le\binom{d+2}{2}$, so their point algebra has Hilbert
function $(1,d,r)$ and $\mu=n-1+r$.
For a partition $\alpha\vdash n$, put
$q(\alpha)=n+1-\alpha_1-\ell(\alpha)$.
Every hook component survives in $G$, while every constituent
with $q(\alpha)>r$ vanishes. In particular,
\[
 \binom{2n-2}{n-1}\le\dim G
 \le\sum_{q(\alpha)\le r}(f^\alpha)^2
 =O_r(n^{2r-1/2}4^n).
\]
Consequently, for every fixed $r$,
\[
 \lim_{n\to\infty}(\dim G)^{1/n}=4,
 \qquad \frac{\dim G}{n!}\longrightarrow0.
\]
At the generic fifteen-point divisor,
\[
 \dim G\le108\,621\,259\,564,
 \qquad\dim A=1\,307\,674\,368\,000.
\]
\end{maintheorem}

The first-occurrence bound for a Specht module is
$n-\alpha_1$. Sign duality in $G$ imposes the complementary
upper bound $\mu-(n-\ell(\alpha))$. Their compatibility is
exactly $q(\alpha)\le r$. The integer $q(\alpha)$ counts the
boxes outside the first row and first column of the Young
diagram. This gives a uniform representation-theoretic
restriction, while the hook lower bound follows from the
exterior-covariant calculation.

The same pairing determines an explicit canonical overmodule
\[
 K=\bigoplus_{s\ge0}t^sF_{\mu-1-s}^{\perp},
 \qquad\operatorname{length}(K/R)=2\delta-\mu n!.
\]
The quotient $K/R$ has perfect sign-valued pairings between
degrees $s$ and $\mu-1-s$. Missing constituents persist
through specified intervals of these degrees. At length
fifteen their contribution gives
\[
 \delta\ge12\,356\,860\,349\,568.
\]
The fifteen-point statement uses generic regularity of the
principal Hilbert component along the Erman--Velasco divisor,
which we establish from the negative tangent space.
A second calculation gives a Grassmannian divisor at length
twelve. The local invariants are
\[
\begin{array}{c|c|r|c|c}
 (n,d)&\text{point Hilbert function}&b=e&\operatorname{edim}&\mu\\\hline
 (10,5)&(1,5,4)&3\,628\,800&46&13\\
 (12,5)&(1,5,6)&479\,001\,600&56&17\\
 (15,11)&(1,11,3)&1\,307\,674\,368\,000&155&17
\end{array}
\]
In each case the geometric tangent cone is a reduced union
of the stated number of lines, and the homogeneous-limit
map reconstructs a generic centered unordered configuration
up to dilation.

\begin{maintheorem}[Invariant products in the derived fibre]
\label{next:intro-tor}
For the orbit cone of a generic configuration with point
Hilbert function $(1,d,r)$, let $K_j$ be the Koszul complex
on the $n-1$ centered differences in one coordinate alphabet.
Then
\[
 H_*(K_j)^{S_n}\simeq F[t]/(t^2)\otimes_F
 \Lambda_F(\eta_1,\ldots,\eta_{n-2}).
\]
Every $\eta_i$ has homological degree one; $d-1$ of them have
internal degree two and $r$ have internal degree three.
The invariant differential graded algebra is formal.

At the three divisorial seeds, these classes transfer to
the marked polygraph rings. For sufficiently many marks,
their Tor algebras over the localized coordinate alphabet
contain exterior algebras on $8$, $10$ and $13$ invariant
degree-one generators, respectively, with nonzero top products.
The projective and flat dimensions at the seed lengths are
exactly $8$, $10$ and $13$.
\end{maintheorem}

For a separated collision with $a$ ten-point blocks,
$b$ twelve-point blocks and $c$ fifteen-point blocks, the
local projective dimension and the maximum nonzero
product length are
\[
 8a+10b+13c.
\]
Uniform seed constants suffice once the number of marks is
at least $aL_{10}+bL_{12}+cL_{15}$, after enlarging them
to carry the products. The ambient dimension is at least
five, and at least eleven when $c>0$.
The ten-point class in internal degree twenty-one is
$t$ times the product of the eight invariant degree-one
classes. The full invariant homology has Hilbert series
$(1+z)(1+sz^2)^4(1+sz^3)^4$.

For each $n\geq3$, a centered multiplication chart represents every
completed local ring of $X_{n,d}$, in every ambient dimension, up to formal variables
(Theorem~\ref{fh:j:cone}). Its square-zero vertex attains the maximum
Cohen--Macaulay defect at that length. The defect is superadditive,
and Lee's length-nine example gives the lower bound
$\lfloor n/9\rfloor$ (Corollary~\ref{fh:j:defect-growth}).
The same corollary gives the sublinear defect criterion.

The sharpness assertion in Theorem~\ref{new:sharp} concerns $R_1$.
The transverse local rings along the constructed divisors are
one-dimensional domains and hence Cohen--Macaulay.
The normality and polygraph conclusions have the all-dimensional
ranges stated above.
The constants $L_m$ are independent of extra points and ambient
coordinates; their existence follows from the finite-generation argument
in Section~\ref{new:polygraph-homology}.
\Needspace{10\baselineskip}

\subsection*{Reading the proofs}
The short normality argument is in
\readinglink{new:ten}{Section~\ref*{new:ten}}.
For the three- and four-dimensional counterexamples, proceed to
\readinglink{dim:four-section}{Sections~\ref*{dim:four-section}}
and~\readinglink{dim:three-boundary}{\ref*{dim:three-boundary}};
the transfer to polygraphs is in
\readinglink{new:polygraph-homology}{Section~\ref*{new:polygraph-homology}}.
The local-ring calculations begin in
\readinglink{new:radial-section}{Section~\ref*{new:radial-section}},
and the alternating-quotient bounds begin in
\readinglink{next:uniform-apolar}{Section~\ref*{next:uniform-apolar}}.

\tableofcontents
\clearpage
\section{A punctual divisor and the normality obstruction}
The reduced ordered cover is integral, and its finite quotient is
\[
 X_{n,d}\longrightarrow H_{n,d}=X_{n,d}/S_n.
\]
To see this, take the closure of the ordered distinct-point locus
in the defining fibre product. Its finite image is closed and
contains the dense configuration locus, so it maps onto $H_{n,d}$.
Each geometric fibre of the fibre product is one $S_n$-orbit.
The invariant closure therefore contains every fibre and equals
$X_{n,d}$ after reduction. Reynolds averaging gives the quotient
identity. Over the configuration locus this is an $S_n$-torsor,
so the function-field extension is Galois with group $S_n$.

\begin{lemma}\label{fh:j:inertia}
At a normal codimension-one point of $X_{n,d}$, the geometric
support multiplicities are all one, or consist of one two and
otherwise ones.
\end{lemma}
\begin{proof}
For support multiplicities $\lambda_1,\ldots,\lambda_s$, the
geometric inertia group is
\[
 I=S_{\lambda_1}\times\cdots\times S_{\lambda_s}.
\]
A normal height-one local ring is a discrete valuation ring $B$.
Let $\pi$ be a uniformizer. The action on its cotangent line is
\[
 \chi:I\longrightarrow\kappa(B)^\times,
 \qquad g\longmapsto g(\pi)/\pi\pmod\pi.
\]
If $\chi(g)=1$, then $g$ acts trivially on the associated graded
ring of $B$. If $g$ has order $r$ and
$g(a)-a\in\mathfrak m_B^v\setminus\mathfrak m_B^{v+1}$, then
\[
 0=g^r(a)-a
   =\sum_{j=0}^{r-1}g^j\bigl(g(a)-a\bigr)
   \equiv r\bigl(g(a)-a\bigr)\pmod{\mathfrak m_B^{v+1}},
\]
a contradiction in residue characteristic zero. Thus $\chi$ is
injective, so $I$ is cyclic. A product of symmetric groups is cyclic only
in the two cases stated. This is the tame inertia argument
\cite[Lemma~15.114.5]{gg-StacksInertia}.
\end{proof}

\begin{corollary}\label{fh:j:support-criterion}
If a divisor in $H_{n,d}$ parametrizes generically schemes with
a support point of multiplicity at least three, or with at least
two double points, then $X_{n,d}$ fails $R_1$ above that divisor.
\end{corollary}
\begin{proof}
The finite map $X_{n,d}\to H_{n,d}$ preserves the codimension
of these points. Their inertia contains $S_3$ or $S_2\times S_2$,
respectively, so Lemma~\ref{fh:j:inertia} applies.
\end{proof}

\subsection{The ten-point divisor and its propagation}
\label{new:ten}
\begin{lemma}\label{new:divisor}
The closure of the translated homogeneous $(1,5,4)$ locus is an
irreducible divisor $D\subset H_{10,5}$. Its generic point is smooth
on the full Hilbert scheme, and its general member is supported at
one point.
\end{lemma}
\begin{proof}
A homogeneous quotient with Hilbert function $(1,5,4)$ is determined
by a four-dimensional quadratic inverse system. Thus the punctual
family is the Grassmannian $\Gr(4,\Sym^2\C^5{}^*)$ of dimension
$4(15-4)=44$. Every such quotient is smoothable by
\cite[Proposition~4.2]{new:GKS}. The proof in that source constructs
a point of this family which is smooth on $\Hilb^{10}(\A^5)$;
an explicit rational witness is also given in
\cite[Lemma~3.1 and Appendix~A]{new:Michalek}.
Smoothness is open, so the generic point of the irreducible family
is smooth. The trace barycenter recovers the support and splits off translations:
\[
 p_j=\frac{\Tr(M_{x_j})}{10},\qquad
 \dim D=4(15-4)+5=49=\dim H_{10,5}-1.
\]
This proves the remaining assertions.
\end{proof}

At the generic point of $D$, the trace formula in the proof gives
a unique $F$-rational support point $p$. The only ordered cycle
above it is $(p,\ldots,p)$, so its inertia is $S_{10}$.
Corollary~\ref{fh:j:support-criterion} proves failure of $R_1$.
The following smooth comparisons carry the same transverse
singularity through the range $n\geq10$, $d\geq5$ in Theorem~\ref{fh:j:main}.

An extra ambient coordinate is an arbitrary function on the
length-ten quotient algebra. Adjoining $d-5$ coordinates is
therefore a graph bundle of rank $10(d-5)$. Adding $n-10$
separated simple points gives an \'etale local product of relative
dimension $d(n-10)$. The ordered values of each extra function
lift the graph bundle to the isospectral cover. The separated
support factors lift because the corresponding idempotents
lift uniquely in infinitesimal deformations. Consequently the
induced divisor $D_{n,d}$ has dimension
\[
 49+10(d-5)+d(n-10)=nd-1.
\]
The local rings differ from the ten-point transverse ring by
smooth base changes. Normality, seminormality, the Du Bois
property and Gorensteinness are preserved and reflected under
these comparisons. Thus the failures proved below at the seed
give Theorem~\ref{fh:j:main}(i) for $n\geq10$ and $d\geq5$.

The general cycle of $D_{n,d}$ has $n-9$ support points.
Its image in one alphabet has dimension at most $n-9$.
Some fibre of the projection therefore has dimension at least
\[
 (nd-1)-(n-9)=n(d-1)+8.
\]
When the selected coordinates are pairwise distinct, the other
$d-1$ coordinates are arbitrary, so the generic fibre dimension
is $n(d-1)$. This proves Theorem~\ref{fh:j:main}(ii) in this range.

\Needspace{10\baselineskip}
\section{Normalization along a homogeneous punctual divisor}
\label{new:radial-section}

Coordinate dilation determines the normalized valuation of a
pointwise fixed divisor on a birational model of affine space.
This identifies its residue field and removes the loss from the
general cyclic-inertia bound.

\begin{theorem}[The radial valuation]\label{new:radial}
Let $k$ be a field of characteristic zero, let $V=\mathbb A_k^N$
carry scalar dilation, and let $Y\to V$ be an equivariant birational
morphism from an integral finite-type $k$-scheme. Suppose that $q$
is a codimension-one point of $Y$, that its closure is fixed
pointwise by dilation, and that it maps to the origin of $V$.
Put $B=\mathcal O_{Y,q}$ and $F=\kappa(q)$.

The normalization $C$ of $B$ is a local DVR. Its valuation on
$k(V)$ is the ordinary order at the origin:
\[
 v(f)=\min\{r:f_r\ne0\},\qquad f=\sum_r f_r\in k[V]\setminus\{0\},
\]
where $f_r$ is homogeneous of degree $r$; the formula extends to
rational functions by subtraction. Its residue field is
\[
 E=k(\mathbb P(V)).
\]
In particular, $E/F$ is finite and
\[
 b(B)=e(B)=[E:F].
\]
Thus there is one arithmetic normalization point, and its residue
extension accounts for all geometric normalization branches.
\end{theorem}

\begin{proof}
The normalization is finite because $B$ is essentially of finite
type over a field. Dilation lifts to normalization. Since the
torus is connected, it fixes each point of the finite set over
$q$. Its action on $F$ is trivial, and its action on every finite
residue extension $E_i/F$ is trivial as well: in characteristic
zero that extension is separable, so its automorphism scheme is
finite. Consequently the cotangent line of each normalized DVR
is an algebraic one-dimensional representation of the split
torus over $E_i$.

Let $v_i$ be the normalized valuation of such a DVR, and let
$m_i\in\mathbb Z$ be its cotangent character. We use the convention
that a coordinate function has dilation weight one. Every nonzero
coordinate $x_j$ vanishes at $q$, so its initial form gives
\[
 1=m_i v_i(x_j),\qquad v_i(x_j)>0.
\]
Both factors on the right are integers. Hence $m_i=1$ and
$v_i(x_j)=1$.

For any nonzero homogeneous $f$ of degree $r$, its initial form
lies in the $v_i(f)$-th cotangent power. The residue field has
weight zero, so comparison of its two weights gives
\[
 r=\operatorname{wt}(\operatorname{in}_{v_i}f)=m_i v_i(f)=v_i(f).
\]
This calculation applies to $f$ itself, including any cancellation
among its monomials: birationality makes $f$ nonzero in $k(V)$.
The nonzero homogeneous parts of a polynomial now have distinct
valuations, and hence
\[
 v_i\!\left(\sum_r f_r\right)=\min_{f_r\ne0}r,
 \qquad v_i(f/g)=v_i(f)-v_i(g).
\]
Thus all $v_i$ are the same valuation of $k(V)$. Distinct
normalization points give distinct valuation rings in this field,
so $C$ is local. Its residue field consists of ratios of homogeneous
polynomials of equal degree and is $k(\mathbb P(V))$.

Every coordinate has valuation one, so
$\mathfrak m_B C=\mathfrak m_C$. The multiplicity formula for finite
birational normalization in dimension one gives
\[
 e_{\mathfrak m_B}(B)
 =[E:F]e_{\mathfrak m_BC}(C)=[E:F].
\]
The residue extension is separable, and hence its degree is also
the geometric branch number.
\end{proof}

\begin{corollary}[Uniqueness of the fixed divisor]
\label{new:radial-unique}
If $Y$ in Theorem~\ref{new:radial} is separated, there is at most
one prime divisor over the origin that is fixed pointwise by
scalar dilation.
\end{corollary}

\begin{proof}
Every such divisor is the center on $Y$ of the radial valuation
identified in Theorem~\ref{new:radial}. A valuation has at most
one center on a separated scheme.
\end{proof}

\begin{corollary}[Finite groups and projective configurations]
\label{new:radial-group}
In Theorem~\ref{new:radial}, suppose that a finite group $G$ acts
linearly and faithfully on $V$, acts compatibly on $Y$, commutes
with dilation, and fixes the closure of $q$ pointwise. Let
\[
 Z=G\cap k^\times I_V.
\]
The inertia of the normalized branch is exactly $Z$. If
$A=B^G$ and
\[
 h=[k(\mathbb P(V))^G:F],
\]
then the normalization of $A$ has one arithmetic point and $h$
geometric branches, and
\[
 b(B)=e(B)=|G/Z|h.
\]
The ramification index of the normalized cover is $|Z|$.
\end{corollary}

\begin{proof}
Theorem~\ref{new:radial} identifies the unique normalized residue
field as $E=k(\mathbb P(V))$. The decomposition group is all of
$G$, and its inertia is the kernel of the action on $E$. A linear
transformation acts trivially on the projective function field
precisely when it is scalar. Thus the inertia is $Z$, and
$[E:E^G]=|G/Z|$.

Since $G$ fixes the original divisor pointwise, $F\subset E^G$.
Normalization commutes with taking finite-group invariants here:
the normalization of $A$ in $\operatorname{Frac}(B)^G$ is $C^G$.
Indeed, $B$ is integral over $A$, so any element integral over
$A$ lies in $C$, while any invariant element of $C$ is integral
over $A$. Exactness of invariants in characteristic zero identifies
the residue field of $C^G$ with $E^G$. This proves the branch formula,
and tame inertia identifies its ramification index as $|Z|$.
\end{proof}

For the Hilbert-scheme applications we take $k=\C$.

\begin{corollary}[Homogeneous punctual divisors]
\label{new:homogeneous-branches}
Let $n\ge3$, and let $D\subset H_{n,d}$ be an irreducible divisor
whose general point is a translate of a homogeneous punctual
scheme. Let $A$ be the local ring at its generic point and $B$ the
local ring at the unique point of the reduced ordered cover above it.
Then
\[
 \boxed{b(B)=e(B)=n!\,b(A).}
\]
Both normalizations have one arithmetic point, and the normalized
ordered cover is unramified there. After removing translation
parameters, their residue fields are respectively
\[
 k\bigl(\mathbb P(\operatorname{Std}_n\otimes k^d)\bigr),
 \qquad
 k\bigl(\mathbb P(\operatorname{Std}_n\otimes k^d)\bigr)^{S_n}.
\]
In particular $b(B)=e(B)\ge n!$. Equality holds if and only if
$A$ is geometrically unibranch; regularity of $A$ is sufficient.
\end{corollary}

\begin{proof}
The trace barycenter splits off translations. On the centered
ordered configuration space the coordinates form
$V=\operatorname{Std}_n\otimes k^d$, and dilation is scalar.
The centered ordered cover maps birationally to $V$, because
distinct ordered points determine their Hilbert point. The
centered part of $D$ is fixed pointwise by dilation. Over a
punctual cycle the ordered fibre has one point with the same
residue field, so $S_n$ fixes the divisor upstairs pointwise.
The quotient property gives $A=B^{S_n}$.

For $n\ge3$, the standard representation has no nonidentity
scalar permutation. It is faithful, and any scalar commutes
with every permutation; the underlying permutation would
therefore lie in the trivial center of $S_n$.
Corollary~\ref{new:radial-group} now applies with $Z=1$.
Adjoining the fixed translation coordinates gives the stated
local invariants on the original, uncentered schemes.
\end{proof}

Since the centered ordered cover is separated,
Corollary~\ref{new:radial-unique} gives at most one homogeneous
punctual divisor on the centered principal component for each
pair $(n,d)$. The constructions in this paper establish existence
in the stated ranges.

At the homogeneous $(1,5,4)$ divisor, the smoothness theorem of
Ga\l\k{a}zka--Keneshlou--\v{S}ivic, as used in
Lemma~\ref{new:divisor}, makes $A$ a DVR. Thus the preceding
corollary gives
\[
 b(B)=e(B)=10!=3\,628\,800.
\]
This deduction uses the geometric smoothness statement and the
radial valuation. In particular, the primitive normal dilation
character also follows from the unramified normalized cover.

\begin{remark}[Normal characters and ramification]
\label{new:weight-divisibility}
The initial-form argument gives a useful general restriction.
Suppose a torus-fixed divisor has a DVR downstairs with nonzero
cotangent weight $w$, and a finite cover carries a compatible
action of the same torus. On each normalized branch the torus
fixes the finite residue extension, so its cotangent line has
an integer weight $m$. For uniformizers $t$ downstairs and $\pi$ upstairs,
\[
 t=a\pi^e,\quad a\in C^\times
 \quad\Longrightarrow\quad w=em
 \quad\Longrightarrow\quad e\mid |w|.
\]
Thus a primitive normal character forces unramified normalization
whenever the torus action lifts to the cover.
\end{remark}
\subsection{Branch and defect bounds for a finite group action}
Normalization separates the point into branches. The inertia
group permutes these branches, and each branch stabilizer acts
faithfully on a cotangent line. The resulting bound applies to
any finite group action on a one-dimensional local domain.

\begin{proposition}\label{fh:j:branch-bound}
Let $B$ be a one-dimensional Noetherian local domain of residue
characteristic zero with finite normalization. Suppose that a
finite group $I$ acts faithfully on $B$ and trivially on
$F=\kappa(B)$. Set
\[
 c(I)=\max\{\ord(g):g\in I\}.
\]
If $b(B)$ is the number of geometric branches, then
\begin{equation}\label{fh:j:general-branch-bound}
 b(B)\geq\frac{|I|}{c(I)},
 \qquad e(B)\geq b(B).
\end{equation}
\end{proposition}
\begin{proof}
Let $C$ be the normalization and $\overline F$ an algebraic
closure of $F$. The geometric branches form the finite $I$-set
\[
 \coprod_{\mathfrak n\in\operatorname{Max}C}
 \Hom_F\bigl(\kappa(\mathfrak n),\overline F\bigr).
\]
The stabilizer $C_\alpha$ of a point in this set fixes
$\mathfrak n$ and its residue field. It acts faithfully on the
discrete valuation ring $C_{\mathfrak n}$, so the character in
Lemma~\ref{fh:j:inertia} makes $C_\alpha$ cyclic. Hence
\[
 |I\cdot\alpha|=[I:C_\alpha]\geq |I|/c(I).
\]
For the multiplicity assertion, write $v_i$ for the
valuations of $C$ and $F_i$ for its residue fields over
$F=\kappa(B)$. A principal minimal reduction $(x)$ of
$\mathfrak m_B$ gives
\[
 e(B)=\operatorname{length}_B(C/xC)
     =\sum_i v_i(x)[F_i:F]
     \geq\sum_i[F_i:F]=b(B).
\]
Here $C/B$ has finite length, so it does not change multiplicity.
\end{proof}

\begin{corollary}\label{fh:j:group-defect}
Under the hypotheses of Proposition~\ref{fh:j:branch-bound}, let
$C$ be the normalization of $B$, and put
$e=\operatorname{edim}B$ and
$\delta=\operatorname{length}_B(C/B)$. Then
\[
 \delta\geq\max_{k\geq1}\left\{
       \frac{k|I|}{c(I)}-\binom{e+k-1}{k-1}\right\}.
\]
If $B$ is seminormal, then $e\geq |I|/c(I)$.
\end{corollary}
\begin{proof}
Completion preserves $e$, $\delta$ and $b(B)$.
Choose a coefficient field $F\subset B$ after completion; then
$C=\prod_iF_i[[t_i]]$ with $\sum_i[F_i:F]=b(B)$.
Write $\mathfrak n=\operatorname{rad}C$. Since
$\mathfrak m_B^k\subseteq\mathfrak n^k$, the image of $B$ in
$C/\mathfrak n^k$ has dimension at most
$\dim_F B/\mathfrak m_B^k$. Thus, for $k\geq1$,
\begin{equation}\label{fh:j:defect-estimate}
 \begin{gathered}
 \dim_F C/\mathfrak n^k=k\,b(B),\qquad
 \dim_F B/\mathfrak m_B^k\leq\binom{e+k-1}{k-1},\\
 \delta\geq k\,b(B)-\binom{e+k-1}{k-1}.
 \end{gathered}
\end{equation}
Apply Proposition~\ref{fh:j:branch-bound}. If $B$ is seminormal,
the conductor calculation~\eqref{fh:j:edim} gives $e=b(B)$,
which proves the second assertion.
\end{proof}

For a seminormal one-dimensional local domain with finite
normalization and separable residue extensions,
\begin{equation}\label{fh:j:edim}
 \operatorname{edim}B=\sum_i[F_i:F].
\end{equation}
Here is the argument. If $B\ne C$, the conductor
$\mathfrak c=\{c\in C:cC\subseteq B\}$ is nonzero and proper.
Seminormality makes it radical in $C$: if $z^r\in\mathfrak c$,
then $(zu)^m\in B$ for every $u\in C$ and all $m\geq r$;
the square-and-cube criterion gives $zu\in B$. Thus
\[
 \mathfrak c=\mathfrak m_B=\operatorname{rad}C,
 \qquad
 \mathfrak m_B^2=(\operatorname{rad}C)^2,
 \qquad
 \dim_F\mathfrak m_B/\mathfrak m_B^2=\sum_i[F_i:F].
\]
The first equality also follows inside $B$, since a nonzero
proper radical ideal in a one-dimensional local domain is its
maximal ideal. The normal case is the discrete valuation case.

\Needspace{10\baselineskip}
\section{A defect bound without a regularity hypothesis}
\label{dim:graded-defect-section}

The radial valuation also controls the normalization defect when the
principal Hilbert component is singular along the homogeneous divisor.
For $n$ points in $d$ alphabets, write $\mu(n,d)$ for the least total
degree of a nonzero alternating polynomial. If
\[
 N_{s-1}=\binom{d+s-1}{d}<n\leq\binom{d+s}{d},
 \qquad r=n-N_{s-1},
\]
then
\begin{equation}\label{dim:minimum-alternant}
 \mu(n,d)=\sum_{j=0}^{s-1}j\binom{d+j-1}{d-1}+sr.
\end{equation}
Indeed, alternation of a monomial is a determinant of $n$ distinct
monomials, and the minimum is obtained by taking the $n$ monomials of
least total degree. Centering the coordinates preserves this minimum,
since translations preserve the alternating representation and its
total-degree filtration.

\begin{theorem}\label{dim:general-sign-defect}
Let $D$ be a homogeneous punctual divisor of the principal component
$H_{n,d}$ over a field of characteristic zero, with $n\geq3$. Let $A$
and $B$ be the generic local rings of $D$ and its reduced ordered cover,
and put $h=b(A)$. Then
\begin{equation}\label{dim:general-sign-bound}
 \boxed{\quad
 \delta(B)\ \geq\ \frac{\mu(n,d)}2\,b(B)
 =\frac{\mu(n,d)}2\,n!h.
 \quad}
\end{equation}
The bound holds for singular principal components and arbitrary
reconstruction degree, with no standard grading assumption.
\end{theorem}

\begin{proof}
Split off the trace barycenter, and let $V=\operatorname{Std}_n\otimes
k^d$ be the centered ordered configuration space. Put $F=\kappa(B)$.
The radial valuation theorem identifies the unique arithmetic
normalization point, its residue field
\[
 E=k(\mathbb P(V)),\qquad [E:F]=b(B)=n!h,
\]
and the valuation of every nonzero homogeneous polynomial with its
degree.

The reduced ordered cover is the blowup of the alternating ideal $J$.
Choose a homogeneous alternant $a$ of minimum degree $\mu=\mu(n,d)$.
The choice of chart matters here: a minimum-degree alternant generates
$JB$ at the divisor. Indeed, if $JB=(g)$, then
\[
 v(g)=\min_{0\ne c\in J}v(c)=\mu=v(a),\qquad
 a/g\in B,\qquad v(a/g)=0.
\]
The normalization has a unique point over the maximal ideal of $B$,
so $a/g\in B^\times$. Hence the Rees chart with denominator $a$
contains the point under consideration.

The coordinate ring of this chart is generated by the centered
coordinates and homogeneous ratios $c/a$, where $c$ is an alternant.
Its grading is nonnegative: coordinates have weight one and $c/a$ has
weight $\deg(c)-\mu\geq0$. The ratios are invariant under $S_n$.
Every nonzero homogeneous element has valuation equal to its weight,
by the radial valuation formula for rational functions. Therefore the
prime of the fixed divisor is precisely the ideal of positive-degree
elements. Localize the degree-zero domain to its fraction field $F$.
This produces a finitely generated graded domain
\[
 R=\bigoplus_{j\geq0}R_j,\qquad R_0=F,\qquad
 B=R_{R_+},\qquad \dim R=1.
\]
Thus $F\subset R$ is a coefficient field in this presentation.

Choose a nonzero homogeneous $S_n$-invariant rational function $u$ of
degree one. For example, the quotient of a nonzero cubic power sum by
a nonzero quadratic power sum in one centered alphabet has this
property. It has valuation one. Set
\[
 V_j=u^{-j}R_j\subseteq E.
\]
The choice of $u$ need not belong to $R$. Its invariance ensures that
each $V_j$ is an $S_n$-subrepresentation of $E$.

The normalization of $R$ is $E[u]$. To see finiteness directly, each
homogeneous generator of $R$ is of the form $\alpha u^j$ with
$\alpha\in E$. Since $E/F$ is finite, all such coefficients are
integral over $F$. A nonzero degree-one coordinate is $\beta u$ with
$\beta\in E^\times$, so $u$ is integral over $R$ after adjoining
$E$. Thus $E[u]$ is finite and integral over $R$, and it is normal in
the common fraction field. After inverting this degree-one coordinate,
the degree-zero ring is a finite $F$-domain and hence a field; its
fraction field is $E$. The localization is therefore
$E[u,u^{-1}]$. It follows that $E[u]/R$ is supported at the vertex,
that $V_j=E$ for all sufficiently large $j$, and that
\begin{equation}\label{dim:graded-delta}
 \delta(B)=\sum_{j\geq0}\bigl([E:F]-\dim_F V_j\bigr).
\end{equation}
This also identifies the valuation-associated graded ring of $B$
with $R$: a denominator outside $R_+$ has nonzero degree-zero term,
so it does not change the least homogeneous part except by an element
of $F^\times$.

There is no sign representation in $R_j$, and hence in $V_j$, for
$j<\mu$. Indeed, before imposing the chart equations, every element is
a sum of coordinate polynomials multiplied by products of invariant
chart generators. Alternating averaging affects only the coordinate
polynomial, whose degree is at most the total weight. Every nonzero
alternating coordinate polynomial has degree at least $\mu$.
Exactness of finite-group invariants in characteristic zero passes
this statement through the quotient and the degree-zero localization.

Choose any nonzero alternating element $\Delta\in E$ and define
\[
 \lambda(z)=\operatorname{Tr}_{E/F}(z/\Delta).
\]
The pairing $(x,y)\mapsto\lambda(xy)$ is nondegenerate because $E/F$
is finite separable. Moreover $\lambda$ transforms by the sign
character. Since $V_iV_j\subseteq V_{i+j}$, sign vanishing implies
\[
 \lambda(V_iV_j)=0\qquad\text{whenever }i+j<\mu.
\]
Pairing complementary degrees gives
\[
 V_{\mu-1-i}\subseteq V_i^\perp
 \quad\Longrightarrow\quad
 \dim_FV_i+\dim_FV_{\mu-1-i}\leq b,
 \qquad b=[E:F],\quad 0\leq i<\mu.
\]
Thus the absence of low-degree sign vectors bounds the sum of all
low-degree dimensions at once:
\[
 2\sum_{i=0}^{\mu-1}\dim_FV_i\leq\mu b,
 \qquad
 \delta(B)\geq\sum_{i=0}^{\mu-1}(b-\dim_FV_i)
             \geq\frac{\mu b}{2}.
\]
This is \eqref{dim:general-sign-bound}.
\end{proof}

For the four-dimensional length-$20$ divisor, $\mu(20,4)=39$, so
\[
 \delta(B)\geq\frac{39}{2}\,20!\,h
 \geq47\,441\,589\,159\,444\,480\,000.
\]
At $(n,d)=(96,3)$ the least alternating degree is $462$, giving
$\delta(B)\geq231\cdot96!\,h$.

\begin{theorem}[The canonical module with arbitrary conductor]
\label{dim:general-canonical-module}
Use the graded presentation and the invariant uniformizer from the
proof of Theorem~\ref{dim:general-sign-defect}. Put $V_j=0$ for $j<0$,
let $\perp$ denote orthogonal complement for
$(x,y)\mapsto\lambda(xy)$ on $E$, and set
\begin{equation}\label{dim:all-integer-canonical}
 K=\bigoplus_{j\in\mathbb Z}
 u^j V_{\mu-1-j}^{\perp}\ \subseteq E[u,u^{-1}].
\end{equation}
Then $K$ is a canonical fractional module of $R$, contains $R$, and
$M=K/R$ has finite length. There are perfect sign-valued pairings
\[
 M_j\times M_{\mu-1-j}\longrightarrow F,
 \qquad ([u^jx],[u^{\mu-1-j}y])\longmapsto\lambda(xy),
\]
and
\begin{equation}\label{dim:general-canonical-defect}
 \boxed{\quad \operatorname{length}_R(M)=2\delta(B)-\mu b(B).
 \quad}
\end{equation}
Formula~\eqref{dim:all-integer-canonical} includes negative degrees;
these can occur when the conductor degree exceeds $\mu$.
The canonical identification uses the trivializing differential
$u^{-\mu}du/\Delta$, of dilation weight $1-\mu$ and sign character;
these are the corresponding grading shift and equivariant twist if
one retains the natural grading and action on dualizing differentials.
\end{theorem}

\begin{proof}
Define a residue pairing on Laurent polynomials by
\[
 \langle f,g\rangle
 =\operatorname{Tr}_{E/F}\operatorname{Res}_{u=0}
       \left(\frac{fg}{\Delta}u^{-\mu}\,du\right).
\]
Its degree-$j$ and degree-$(\mu-1-j)$ coefficient pairing is precisely
$\lambda(xy)$. Thus the annihilator of $R$ under this residue pairing
is the module in \eqref{dim:all-integer-canonical}.

To identify this annihilator with a canonical module, choose a
nonzero degree-one coordinate $z=\beta u\in R$, with $\beta\in E^\times$.
The domain $R$ is finite and torsion free over $F[z]$, so
\[
 \omega_R=\operatorname{Hom}_{F[z]}(R,F[z])\,dz.
\]
The separable function-field trace realizes $\omega_R$ as the fractional
differentials $\omega$ for which
\[
 \operatorname{Tr}_{E(u)/F(z)}(r\omega)\in F[z]\,dz
 \qquad(r\in R).
\]
Since $R$ agrees with its normalization away from the vertex, a
homogeneous such differential is Laurent. For a Laurent polynomial $q$,
\[
 q\in F[z]
 \quad\Longleftrightarrow\quad
 \operatorname{Res}_{z=0}(z^a q\,dz)=0\quad(a\geq0),
\]
because these residues detect every negative coefficient of $q$.
Using $z^aR\subseteq R$ and the compatibility of residues with finite
separable trace, the condition on $\omega$ becomes
\[
 \operatorname{Tr}_{E/F}\operatorname{Res}_{u=0}(r\omega)=0
 \qquad(r\in R).
\]
Trivializing by $u^{-\mu}du/\Delta$ therefore identifies $\omega_R$
with $K$.

The multiplication rule $V_iV_j\subseteq V_{i+j}$ shows directly that
$K$ is an $R$-module. Sign vanishing below $\mu$ gives $R\subseteq K$.
Since $V_j=E$ for large $j$, the module $K$ is bounded below and agrees
with $E[u]$ in large positive degrees. Hence $K/R$ has finite length.

For $t=\mu-1-j$, the two graded pieces of this quotient are
\[
 M_j=V_t^\perp/V_j,\qquad M_t=V_j^\perp/V_t.
\]
The nondegenerate pairing on $E$ descends to a perfect pairing on
these quotients. It transforms by the sign character because $u$ is
invariant and $\lambda$ is alternating.

Finally choose $c\geq\mu$ such that $V_j=E$ for all $j\geq c$. Every
nonzero graded piece of $M$ lies between $\mu-c$ and $c-1$. Summing
its dimensions gives
\begin{align*}
 \operatorname{length}_R(M)
 &=\sum_{j=\mu-c}^{c-1}
   \left(b-\dim_FV_{\mu-1-j}-\dim_FV_j\right)\\
 &=(2c-\mu)b-2\sum_{j=0}^{c-1}\dim_FV_j
 =2\delta(B)-\mu b.
\end{align*}
This proves \eqref{dim:general-canonical-defect}.
\end{proof}

\begin{corollary}\label{dim:reconstruction-canonical-defect}
In the preceding setting, let $G=S_n$ and $h=[E^G:F]$. Then
\[
 \dim_F(R_\mu)^{\mathrm{sgn}}=1,
 \qquad
 \dim_F M_{-1}^{G}=\dim_F M_0^{G}
 =\dim_F M_{\mu-1}^{\mathrm{sgn}}
 =\dim_F M_\mu^{\mathrm{sgn}}=h-1.
\]
Consequently $h>1$ forces the conductor degree to exceed $\mu$ and
forces a nonzero negative degree in the canonical fractional module.
Moreover
\[
 \delta(B)\geq\frac{\mu}{2}\,n!h+2(h-1).
\]
\end{corollary}

\begin{proof}
In weight $\mu$, a sign element of the alternating chart must be a
minimum-degree coordinate alternant multiplied by an invariant of
degree zero. Every such alternant $a'$ satisfies $a'/a\in F$, so
the sign space is the one-dimensional space $Fa$. The normal basis
theorem for $E/E^G$ gives both $\dim_FE^G=h$ and
$\dim_FE^{\mathrm{sgn}}=h$. The sign-trace pairing restricts to a
perfect pairing between these two spaces.

Since $V_{\mu-1}$ has no sign component, its orthogonal complement
contains all of $E^G$. Quotienting by $V_0=F$ gives
$\dim M_0^G=h-1$. Since $V_\mu$ has a one-dimensional sign component,
its orthogonal complement has invariant dimension $h-1$, proving the
same assertion for $M_{-1}=K_{-1}$. The perfect pairings of
Theorem~\ref{dim:general-canonical-module} give the two sign dimensions.
For $n\geq3$ one has $\mu\geq2$, so the four indicated degrees are
distinct. Thus $\operatorname{length}(M)\geq4(h-1)$, and
\eqref{dim:general-canonical-defect} gives the claimed defect bound.
\end{proof}

\subsection{An invariant parameter and the alternating quotient}
\label{dim:actual-alternating-quotient}

Retain \(F,E,R,K,M\) from
Theorems~\ref{dim:general-sign-defect} and
\ref{dim:general-canonical-module}, and put \(b=[E:F]=n!h\).
Choosing the invariant uniformizer inside \(R\) turns reduction
modulo that parameter into a comparison between the canonical
defect and the determinant apolar algebra. This comparison uses
the local ring's own grading. The configuration-orbit filtration
will enter later, under the regularity hypothesis in
Section~\ref{new:orbit}.

\begin{lemma}[The Hilbert function at a homogeneous divisor]
\label{dim:automatic-maximal-hilbert}
Every homogeneous punctual divisor in this setting has the
maximal point Hilbert function at its generic point.
\end{lemma}

\begin{proof}
As in the proof of Theorem~\ref{dim:general-sign-defect}, every
nonzero minimum-degree monomial alternant \(a\) is a local
generator of \(JB\).  The ordered evaluation maps of the
tautological algebra \(\mathcal T\) give a determinant map
\[
 \rho^*\det\mathcal T\longrightarrow\mathcal O_X.
\]
Its image is \(J\mathcal O_X\), since the monomial determinants
generate \(J\).  This identifies the invertible source with
that invertible ideal.  Therefore a generator \(a\) corresponds
to a nonvanishing determinant section of \(\mathcal T\), and
the monomials in its columns form a basis at the divisor point.
Those columns consist of all monomials below degree \(s\) and
\(r\) monomials of degree \(s\). Homogeneity makes this a graded
basis, so the point algebra has Hilbert function
\[
 (\dim S_0,\ldots,\dim S_{s-1},r,0,0,\ldots).
\]
\end{proof}

\begin{lemma}[An invariant parameter of weight one]
\label{dim:degree-one-parameter}
For \(d\geq2\) and \(n\geq3\), the space of centered alternating
polynomials of degree \(\mu(n,d)+1\) is nonzero.  Consequently
\(R\) contains an invariant homogeneous element \(u\) of weight
and radial valuation one.
\end{lemma}

\begin{proof}
Put \(M_j=\binom{d+j-1}{d-1}\), choose \(s\) so that
\(n=\sum_{j<s}M_j+r\), and write \(m=M_s\).
There are \(\binom mr\) minimum-degree monomial alternants.
Increasing their degree by one permits precisely two choices:
replace one degree-\(s\) monomial by one of degree \(s+1\), or
remove one degree-\((s-1)\) monomial and include one additional
degree-\(s\) monomial.  Before centering, the dimension in
degree \(\mu+1\) is therefore
\[
 M_{s+1}\binom m{r-1}+M_{s-1}\binom m{r+1}.
\]
The full polynomial ring is the tensor product of the centered
polynomial ring and the polynomial ring in the \(d\) mean
coordinates.  Since the first alternating degree is \(\mu\),
the centered alternating dimension in degree \(\mu+1\) is
\begin{equation}\label{dim:centered-next-alternants}
 M_{s+1}\binom m{r-1}+M_{s-1}\binom m{r+1}
       -d\binom mr.
\end{equation}
If \(s\geq2\), then \(m\geq3\) and
\(M_{s-1},M_{s+1}\geq d\).  Dividing the first two terms by
\(\binom mr\) and using
\[
 \frac r{m-r+1}+\frac{m-r}{r+1}
       \geq\frac{2m}{m+2}>1
\]
proves positivity.  If \(s=1\), then \(m=d\) and
\(r=n-1\geq2\).  The ratio of
\eqref{dim:centered-next-alternants} to \(\binom dr\) is
\[
 \frac{r(d+1)(d+2)(r-1)}{2(d-r+1)(r+1)}>0.
\]
Choose a nonzero centered alternant \(a_+\) of degree
\(\mu+1\).  On the minimal-alternant chart with denominator
\(a\), the ratio \(u=a_+/a\) belongs to \(R\), is invariant,
and has weight and valuation one.
\end{proof}

\begin{corollary}\label{dim:unibranch-regular}
The principal component is regular at the generic homogeneous
divisor if and only if \(h=1\).  In that case its local ring is
\(F[u]_{(u)}\), and the ordered graded local model is
\(R=F[u,x_{ij}]\).
\end{corollary}

\begin{proof}
The invariant normalization is \(E^{S_n}[u]\), and
\([E^{S_n}:F]=h\).  If \(h=1\), then
\(F[u]\subseteq R^{S_n}\subseteq F[u]\), so equality holds.
All alternating-ratio chart generators are invariant, which
also gives \(R=F[u,x_{ij}]\).  Conversely, a regular
one-dimensional principal local ring is normal, and its
normalization has the same residue field.  Thus \(h=1\).
\end{proof}

\begin{theorem}[The canonical defect and the alternating quotient]
\label{dim:actual-canonical-apolar}
Assume \(d\geq2\) and \(n\geq3\).  Let \(G_Q\) be the
determinant apolar algebra of the generic homogeneous point and
put \(g=\dim_FG_Q\).  Then
\begin{equation}\label{dim:actual-apolar-defect}
 \boxed{\quad
 \operatorname{length}_R M=2\delta(B)-\mu b\ \geq\ b-g,
 \qquad \delta(B)\geq\frac{(\mu+1)b-g}{2}.
 \quad}
\end{equation}
In particular \(\operatorname{length}_R M\geq n!(h-1)\).
These conclusions hold for every reconstruction degree \(h\),
including divisors where the principal component is singular.
\end{theorem}

\begin{proof}
Use the parameter \(u\in R_1^{S_n}\) from the lemma.
The spaces \(V_j=u^{-j}R_j\) are now nested.  Since
\(R\subset E[u]\) is closed under multiplication by \(u\),
it is a finite free \(F[u]\)-module of rank \(b\).
Consequently \(A_u=R/uR\) has dimension \(b\).

Choose \(\Delta=a/u^\mu\) and
\(\lambda(z)=\operatorname{Tr}_{E/F}(z/\Delta)\).
The canonical module of
Theorem~\ref{dim:general-canonical-module} can be written
\[
 K=u^\mu\{q\in E(u):\lambda(qR)\subset F[u]\}.
\]
The module in braces is
\(\operatorname{Hom}_{F[u]}(R,F[u])\) under the trace
pairing.  Therefore reduction modulo \(u\) identifies
\(K/uK\) with the shifted dual of \(A_u\).  The map induced
by \(R\subset K\) is the top-degree multiplication pairing
\begin{equation}\label{dim:actual-reduced-pairing}
 A_u\otimes A_u\longrightarrow F,
 \qquad (x,y)\longmapsto[u^\mu]\lambda(xy).
\end{equation}
This argument retains every integer degree in \(K\), including
negative degrees when the conductor exceeds \(\mu\).

Tensoring \(0\to R\to K\to M\to0\) with \(F[u]/(u)\) gives
\begin{equation}\label{dim:actual-reduced-kernel}
 \ker(A_u\longrightarrow K/uK)\simeq(0:_M u)(-1).
\end{equation}
Consequently
\[
 \dim_F\ker(A_u\longrightarrow K/uK)
 =\dim_F(0:_M u)\leq\operatorname{length}_R M.
\]

Every positive-degree invariant element of \(R\) annihilates
\eqref{dim:actual-reduced-pairing}: its only possible partner
has degree less than \(\mu\), where the sign representation
is absent.  In particular all positive-degree alternating-ratio
chart generators lie in its radical.  Thus the nondegenerate
quotient is generated by the ordered point coordinates alone.

The restricted functional is the exterior multiplication
functional of the homogeneous point algebra, up to a nonzero
scalar.  For a monomial in the ordered rows of total degree
\(\mu\), let \(a'\) be its alternating determinant.  Its
alternating projection is \(a'/n!\), with the fixed ordering
convention.  The ratio \(a'/a\in F\) is the corresponding
determinant coordinate of the tautological point algebra.
Hence
\[
 \lambda\!\left(\frac{a'}{n!u^\mu}\right)
       =\frac b{n!}\frac{a'}a.
\]
Since \(b/n!\ne0\), the nondegenerate quotient is \(G_Q\).
Together with \eqref{dim:actual-reduced-kernel}, this gives
\[
 b-g=\dim_F\ker(A_u\longrightarrow K/uK)
 \leq\operatorname{length}_R M=2\delta(B)-\mu b,
\]
which is \eqref{dim:actual-apolar-defect}.

Finally the universal alternating multiplication map has rank
\(n!\) on the reduced-configuration locus.  Its specialization
at the point under consideration has rank \(g\leq n!\).
Thus \(b-g\geq n!(h-1)\).
\end{proof}

The source bound of Proposition~\ref{dim:source-bound} gives
the concrete consequence
\[
 \delta(B)\geq
 \frac{(\mu+1)b-\min\{n!,U(n,d)\}}2.
\]
Here \(Q\) is the generic homogeneous point of the divisor.
The fat-point examples of lengths \(70\) and \(120\) concern
the alternating pairing at individual points.

\section{Seven homogeneous punctual divisors in four dimensions}
\label{dim:four-section}

Write $S=k[x_1,x_2,x_3,x_4]$. For $r=5,\ldots,11$, put $n=15+r$
and consider homogeneous ideals
\[
 I_W=(W)+\mathfrak m^4,\qquad
 W\subseteq S_3,\qquad \dim W=20-r.
\]
They have Hilbert function $(1,4,10,r)$.

\begin{proposition}\label{dim:collision-criterion}
Let $S=k[x_1,\ldots,x_d]$, let $n,d,s$ satisfy
$\binom{d+s-1}{d}<n\leq\binom{d+s}{d}$, and put
$r=n-\binom{d+s-1}{d}$. For a configuration of $n$ distinct points for
which evaluation is injective on polynomials of degree at most $s-1$
and surjective on polynomials of degree at most $s$, its scalar
collision has ideal
\[
 I_W=(W)+\mathfrak m^{s+1},\qquad
 W\in\operatorname{Gr}\!\left(\dim S_s-r,S_s\right).
\]
If the differential of the resulting rational collision map has rank
$nd-d-1$ at one characteristic-zero configuration, then the closure of
its translated image is a homogeneous punctual divisor of $H_{n,d}$.
If additionally $r(\dim S_s-r)=nd-d-1$, every homogeneous algebra of
this Hilbert function is smoothable.
\end{proposition}

\begin{proof}
Let $Z$ be such a configuration and let $I(Z)$ be its vanishing
ideal. The evaluation hypotheses give
\[
 I(Z)\cap S_{\leq s-1}=0,\qquad
 \dim(I(Z)\cap S_{\leq s})=\dim S_s-r.
\]
Taking the degree-$s$ part is therefore injective on the second
space; denote its image by $W$. The degree filtration has associated
graded algebra $S/((W)+\mathfrak m^{s+1})$, and its Rees family
is the scalar collision of $Z$.

Let $U$ be this open configuration locus and $c$ its collision map. Translations and dilation preserve $W$, and their $d+1$
tangent directions are independent when the points are distinct.
Consequently
\[
 \operatorname{rank}(dc)\leq\dim\overline{c(U)}\leq nd-d-1.
\]
A nonzero minor of size $nd-d-1$ forces equality. The trace barycenter
recovers the support of a translated limit, so translation adds $d$
parameters:
\[
 \dim\overline{\{\text{translated collision limits}\}}
 =(nd-d-1)+d=\dim H_{n,d}-1.
\]
If the Grassmannian has dimension $nd-d-1$, the collision image is
dense in it. The principal Hilbert component is closed, so it then
contains the entire Grassmannian.
\end{proof}

\begin{corollary}[The entire homogeneous smoothable locus]
\label{dim:entire-smoothable-locus}
For any very compressed Hilbert function
\[
 \left(\dim S_0,\ldots,\dim S_{s-1},r\right),
\]
the homogeneous smoothable locus in
$\operatorname{Gr}(\dim S_s-r,S_s)$ equals the closure of the scalar
collision image of configurations with this polynomial Hilbert
function. In particular this entire smoothable locus is irreducible.
\end{corollary}

\begin{proof}
At a homogeneous ideal, choose a border basis consisting of all
monomials below degree $s$ and $r$ monomials of degree $s$.
If $c_{m,b}$ is the coefficient of a basis monomial $b$ in the
reduction of a border monomial $m$, then
\[
 \operatorname{wt}(c_{m,b})=\deg m-\deg b\geq0.
\]
Thus setting the positive-weight coefficients to zero defines a
retraction $\pi$ from the border-basis chart to its Grassmannian
chart. On distinct-point configurations, $\pi$ is the scalar
collision map.

Let $I$ be a homogeneous smoothable point and let $U$ be the
distinct-point locus in its chart. Since $I\in\overline U$ and
$\pi(I)=I$, continuity gives
\[
 I=\pi(I)\in\pi(\overline U)\subseteq\overline{\pi(U)}.
\]
Conversely, $\pi(U)$ and its closure lie in the closed principal
component. Hence the homogeneous smoothable locus is the closure
of the collision image, which is irreducible because the
configuration space is irreducible.
\end{proof}

\begin{theorem}\label{dim:four-divisors}
For every $n=20,21,\ldots,26$, the principal component $H_{n,4}$ has a
homogeneous punctual divisor with Hilbert function $(1,4,10,n-15)$.
The entire homogeneous smoothable locus is irreducible, and its
dimension and codimension in the homogeneous Grassmannian are:
\[
\begin{array}{c|rrrrrrr}
 n&20&21&22&23&24&25&26\\ \hline
 r&5&6&7&8&9&10&11\\
 \dim(\text{smoothable locus})&75&79&83&87&91&95&99\\
 \dim\operatorname{Gr}(r,20)&75&84&91&96&99&100&99\\
 \operatorname{codim}&0&5&8&9&8&5&0
\end{array}
\]
In particular every homogeneous algebra with Hilbert function
$(1,4,10,5)$ or $(1,4,10,11)$ is smoothable.
For each of the five intermediate top entries $r=6,\ldots,10$, a
general homogeneous algebra is nonsmoothable.

At the generic point of each divisor, if $A_n$ is the principal
Hilbert local ring and $B_n$ its reduced ordered-cover local ring,
then, for the finite reconstruction degree $h_n=b(A_n)$,
\[
 b(B_n)=e(B_n)=n!h_n,\qquad
 \delta(B_n)\geq\frac{3n-21}{2}\,n!h_n.
\]
The numbers $h_n$ are not determined by these rank certificates.
\end{theorem}

\begin{proof}
The certificates give rational point configurations and their collision
differentials. All evaluation denominators are units at $32003$ and
$65521$, and reduction at either prime gives rank $4n-5$. Thus
\[
 4n-5=\operatorname{rank}_{\mathbb F_p}(dc)
 \leq\operatorname{rank}_{\mathbb Q}(dc)\leq4n-5.
\]
The last inequality comes from the five translation and dilation
directions. Therefore the rational rank is $4n-5$ at each length.

Proposition~\ref{dim:collision-criterion} gives the divisors, and
Corollary~\ref{dim:entire-smoothable-locus} identifies their centered
images with the entire smoothable loci. The radial valuation
theorem gives the branch and multiplicity identities without a
regularity hypothesis on $A_n$. Finally
$\mu(n,4)=4+20+3(n-15)=3n-21$, so
Theorem~\ref{dim:general-sign-defect} gives the defect bound.
\end{proof}

The dimension inequality factors as
\[
 r(20-r)\geq4(15+r)-5
 \quad\Longleftrightarrow\quad (r-5)(r-11)\leq0.
\]
Thus the seven lengths exhaust the possibilities for this very
compressed cubic construction. The minimum-length question for homogeneous punctual divisors is
examined later in the paper.

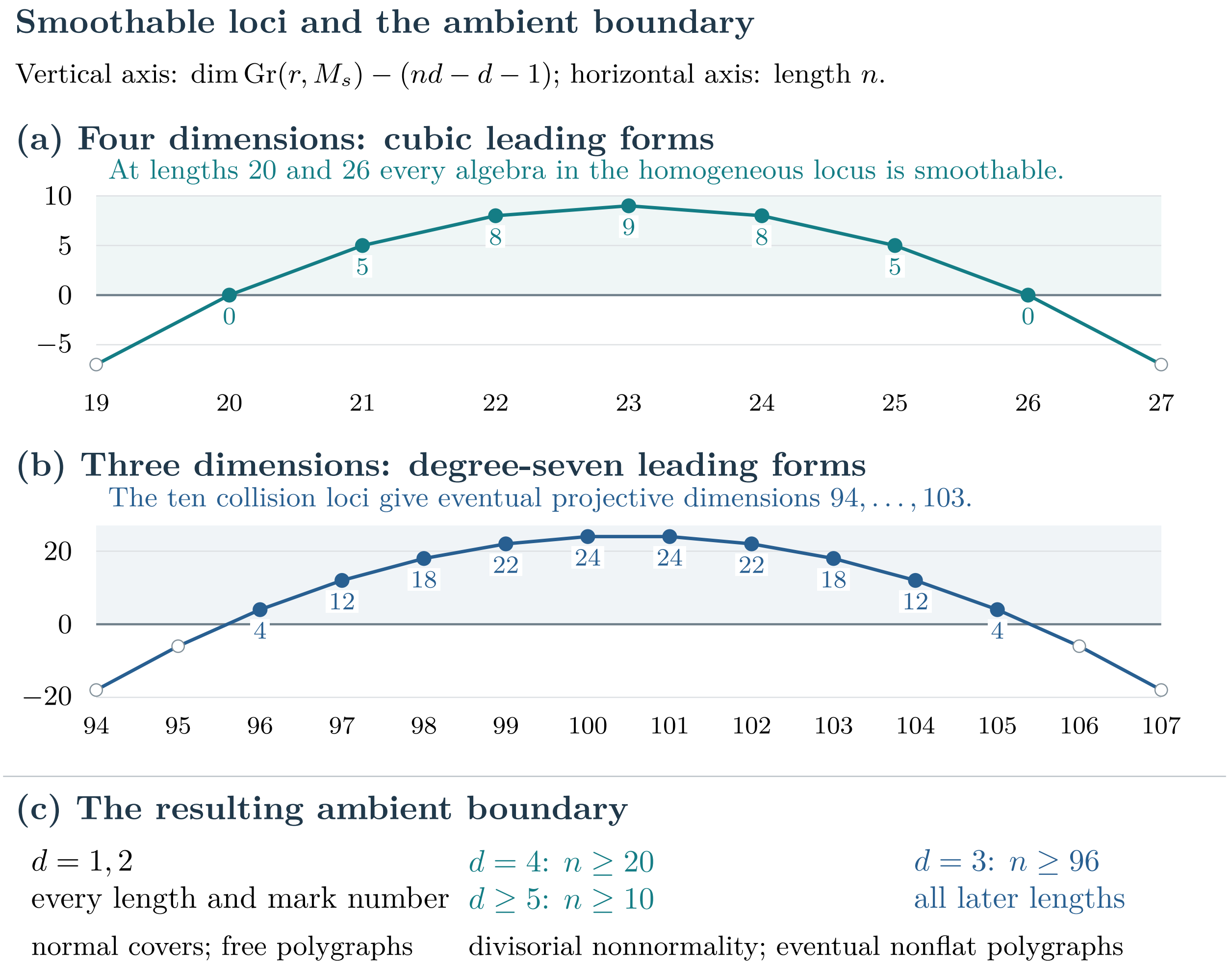
\begin{figure}[!htbp]
\centering
\resizebox{\linewidth}{!}{%
\begin{tikzpicture}[x=.92cm,y=.78cm,every node/.style={font=\small}]
\node[anchor=west,font=\bfseries,Ink] at (0,13.4)
 {Smoothable loci and the ambient boundary};
\node[anchor=west,font=\footnotesize] at (0,12.65)
 {Vertical axis: $\dim\operatorname{Gr}(r,M_s)-(nd-d-1)$; horizontal axis: length $n$.};
\node[anchor=west,font=\bfseries,Ink] at (0,11.7)
 {(a) Four dimensions: cubic leading forms};

\fill[Teal!7] (1.15,9.4556) rectangle (14.3,10.9000);
\draw[Ink!15] (1.15,8.7333)--(14.3,8.7333);\node[anchor=east,font=\footnotesize] at (1.0,8.7333) {$ -5 $};
\draw[Ink!15] (1.15,9.4556)--(14.3,9.4556);\node[anchor=east,font=\footnotesize] at (1.0,9.4556) {$ 0 $};
\draw[Ink!15] (1.15,10.1778)--(14.3,10.1778);\node[anchor=east,font=\footnotesize] at (1.0,10.1778) {$ 5 $};
\draw[Ink!15] (1.15,10.9000)--(14.3,10.9000);\node[anchor=east,font=\footnotesize] at (1.0,10.9000) {$ 10 $};
\draw[Ink!65,line width=.7pt] (1.15,9.4556)--(14.3,9.4556);
\draw[Teal,line width=1.1pt] (1.15000,8.44444) -- (2.79375,9.45556) -- (4.43750,10.17778) -- (6.08125,10.61111) -- (7.72500,10.75556) -- (9.36875,10.61111) -- (11.01250,10.17778) -- (12.65625,9.45556) -- (14.30000,8.44444);
\draw[Ink!55,fill=white] (1.1500,8.4444) circle (2.0pt);
\node[font=\scriptsize,below] at (1.1500,8.1800) {$ 19 $};
\fill[Teal] (2.7938,9.4556) circle (2.4pt);
\node[font=\scriptsize,below] at (2.7938,8.1800) {$ 20 $};
\fill[Teal] (4.4375,10.1778) circle (2.4pt);
\node[font=\scriptsize,below] at (4.4375,8.1800) {$ 21 $};
\fill[Teal] (6.0813,10.6111) circle (2.4pt);
\node[font=\scriptsize,below] at (6.0813,8.1800) {$ 22 $};
\fill[Teal] (7.7250,10.7556) circle (2.4pt);
\node[font=\scriptsize,below] at (7.7250,8.1800) {$ 23 $};
\fill[Teal] (9.3688,10.6111) circle (2.4pt);
\node[font=\scriptsize,below] at (9.3688,8.1800) {$ 24 $};
\fill[Teal] (11.0125,10.1778) circle (2.4pt);
\node[font=\scriptsize,below] at (11.0125,8.1800) {$ 25 $};
\fill[Teal] (12.6562,9.4556) circle (2.4pt);
\node[font=\scriptsize,below] at (12.6562,8.1800) {$ 26 $};
\draw[Ink!55,fill=white] (14.3000,8.4444) circle (2.0pt);
\node[font=\scriptsize,below] at (14.3000,8.1800) {$ 27 $};
\node[below=3pt,font=\scriptsize,text=Teal,fill=white,inner sep=1pt] at (2.7938,9.4556) {$ 0 $};
\node[below=3pt,font=\scriptsize,text=Teal,fill=white,inner sep=1pt] at (4.4375,10.1778) {$ 5 $};
\node[below=3pt,font=\scriptsize,text=Teal,fill=white,inner sep=1pt] at (6.0813,10.6111) {$ 8 $};
\node[below=3pt,font=\scriptsize,text=Teal,fill=white,inner sep=1pt] at (7.725,10.7556) {$ 9 $};
\node[below=3pt,font=\scriptsize,text=Teal,fill=white,inner sep=1pt] at (9.3688,10.6111) {$ 8 $};
\node[below=3pt,font=\scriptsize,text=Teal,fill=white,inner sep=1pt] at (11.0125,10.1778) {$ 5 $};
\node[below=3pt,font=\scriptsize,text=Teal,fill=white,inner sep=1pt] at (12.6562,9.4556) {$ 0 $};
\node[anchor=west,font=\footnotesize,Teal] at (1.15,11.25) {At lengths $20$ and $26$ every algebra in the homogeneous locus is smoothable.};
\node[anchor=west,font=\bfseries,Ink] at (0,6.95) {(b) Three dimensions: degree-seven leading forms};
\fill[Blue!7] (1.15,4.6638) rectangle (14.3,6.1000);
\draw[Ink!15] (1.15,3.6000)--(14.3,3.6000);\node[anchor=east,font=\footnotesize] at (1.0,3.6000) {$ -20 $};
\draw[Ink!15] (1.15,4.6638)--(14.3,4.6638);\node[anchor=east,font=\footnotesize] at (1.0,4.6638) {$ 0 $};
\draw[Ink!15] (1.15,5.7277)--(14.3,5.7277);\node[anchor=east,font=\footnotesize] at (1.0,5.7277) {$ 20 $};
\draw[Ink!65,line width=.7pt] (1.15,4.6638)--(14.3,4.6638);
\draw[Blue,line width=1.1pt] (1.15000,3.70638) -- (2.16154,4.34468) -- (3.17308,4.87660) -- (4.18462,5.30213) -- (5.19615,5.62128) -- (6.20769,5.83404) -- (7.21923,5.94043) -- (8.23077,5.94043) -- (9.24231,5.83404) -- (10.25385,5.62128) -- (11.26538,5.30213) -- (12.27692,4.87660) -- (13.28846,4.34468) -- (14.30000,3.70638);
\draw[Ink!55,fill=white] (1.1500,3.7064) circle (2.0pt);
\node[font=\scriptsize,below] at (1.1500,3.4800) {$ 94 $};
\draw[Ink!55,fill=white] (2.1615,4.3447) circle (2.0pt);
\node[font=\scriptsize,below] at (2.1615,3.4800) {$ 95 $};
\fill[Blue] (3.1731,4.8766) circle (2.4pt);
\node[font=\scriptsize,below] at (3.1731,3.4800) {$ 96 $};
\fill[Blue] (4.1846,5.3021) circle (2.4pt);
\node[font=\scriptsize,below] at (4.1846,3.4800) {$ 97 $};
\fill[Blue] (5.1962,5.6213) circle (2.4pt);
\node[font=\scriptsize,below] at (5.1962,3.4800) {$ 98 $};
\fill[Blue] (6.2077,5.8340) circle (2.4pt);
\node[font=\scriptsize,below] at (6.2077,3.4800) {$ 99 $};
\fill[Blue] (7.2192,5.9404) circle (2.4pt);
\node[font=\scriptsize,below] at (7.2192,3.4800) {$ 100 $};
\fill[Blue] (8.2308,5.9404) circle (2.4pt);
\node[font=\scriptsize,below] at (8.2308,3.4800) {$ 101 $};
\fill[Blue] (9.2423,5.8340) circle (2.4pt);
\node[font=\scriptsize,below] at (9.2423,3.4800) {$ 102 $};
\fill[Blue] (10.2538,5.6213) circle (2.4pt);
\node[font=\scriptsize,below] at (10.2538,3.4800) {$ 103 $};
\fill[Blue] (11.2654,5.3021) circle (2.4pt);
\node[font=\scriptsize,below] at (11.2654,3.4800) {$ 104 $};
\fill[Blue] (12.2769,4.8766) circle (2.4pt);
\node[font=\scriptsize,below] at (12.2769,3.4800) {$ 105 $};
\draw[Ink!55,fill=white] (13.2885,4.3447) circle (2.0pt);
\node[font=\scriptsize,below] at (13.2885,3.4800) {$ 106 $};
\draw[Ink!55,fill=white] (14.3000,3.7064) circle (2.0pt);
\node[font=\scriptsize,below] at (14.3000,3.4800) {$ 107 $};
\node[below=3pt,font=\scriptsize,text=Blue,fill=white,inner sep=1pt] at (3.1731,4.8766) {$ 4 $};
\node[below=3pt,font=\scriptsize,text=Blue,fill=white,inner sep=1pt] at (4.1846,5.3021) {$ 12 $};
\node[below=3pt,font=\scriptsize,text=Blue,fill=white,inner sep=1pt] at (5.1962,5.6213) {$ 18 $};
\node[below=3pt,font=\scriptsize,text=Blue,fill=white,inner sep=1pt] at (6.2077,5.834) {$ 22 $};
\node[below=3pt,font=\scriptsize,text=Blue,fill=white,inner sep=1pt] at (7.2192,5.9404) {$ 24 $};
\node[below=3pt,font=\scriptsize,text=Blue,fill=white,inner sep=1pt] at (8.2308,5.9404) {$ 24 $};
\node[below=3pt,font=\scriptsize,text=Blue,fill=white,inner sep=1pt] at (9.2423,5.834) {$ 22 $};
\node[below=3pt,font=\scriptsize,text=Blue,fill=white,inner sep=1pt] at (10.2538,5.6213) {$ 18 $};
\node[below=3pt,font=\scriptsize,text=Blue,fill=white,inner sep=1pt] at (11.2654,5.3021) {$ 12 $};
\node[below=3pt,font=\scriptsize,text=Blue,fill=white,inner sep=1pt] at (12.2769,4.8766) {$ 4 $};
\node[anchor=west,font=\footnotesize,Blue] at (1.15,6.48) {The ten collision loci give eventual projective dimensions $94,\ldots,103$.};
\draw[Ink!30] (0,2.45)--(15.1,2.45);
\node[anchor=west,font=\bfseries,Ink] at (0,1.95) {(c) The resulting ambient boundary};
\node[anchor=west,align=left] at (.2,.92)
 {$d=1,2$\\every length and mark number};
\node[anchor=west,align=left,Teal] at (5.6,.92)
 {$d=4$: $n\geq20$\\$d\geq5$: $n\geq10$};
\node[anchor=west,align=left,Blue] at (11.1,.92)
 {$d=3$: $n\geq96$\\all later lengths};
\node[anchor=west,font=\footnotesize] at (.2,-.05) {normal covers; free polygraphs};
\node[anchor=west,font=\footnotesize] at (5.6,-.05) {divisorial nonnormality; eventual nonflat polygraphs};
\end{tikzpicture}}
\caption{Homogeneous collision loci. At each filled point, a nonzero
maximal minor gives differential rank $nd-d-1$, so the displayed ordinate
is the codimension of the homogeneous smoothable locus. Hollow points fail
the dimension inequality. In four dimensions, every algebra at either
endpoint is smoothable, while the intermediate loci have codimensions
$5,8,9,8,5$. In three dimensions, the construction attains all ten lengths
$96,\ldots,105$. Among lengths at least three, $96$ and $20$ are the first
homogeneous punctual-divisor lengths in dimensions three and four.
Panel (c) combines these constructions with Haiman's planar theorem.}
\label{dim:collision-figure}
\end{figure}

\subsection{The generic twenty-point tangent cone}

\begin{lemma}\label{dim:degree-one-cone}
Let $C$ be a finitely generated nonnegatively graded algebra over a
field $F$, with $C_0=F$. If the positive cotangent space
$C_+/C_+^2$ is concentrated in degree one, then $C$ is generated by
$C_1$. Every one-dimensional graded domain quotient $R$ of $C$ is
therefore a standard graded domain. Writing
$L=\operatorname{Frac}(R)_0$, one has
\[
 \operatorname{Proj}(R)=\operatorname{Spec}(L),\qquad
 \overline R=L[t],\qquad
 \operatorname{gr}_{\mathfrak m}(R_{R_+})\simeq R.
\]
Here $t$ has degree one. If $L/F$ is separable, the geometric tangent
cone consists of $[L:F]$ distinct reduced lines.
\end{lemma}

\begin{proof}
For $j>1$, the cotangent hypothesis gives
\[
 C_j=\sum_{a=1}^{j-1}C_aC_{j-a}.
\]
Induction proves generation by $C_1$, and the same conclusion holds
in every graded quotient. For a one-dimensional domain $R$, the
projective scheme is integral and zero-dimensional, hence is
$\operatorname{Spec}(L)$ for a finite extension $L/F$. If $t\in R_1$
is nonzero, then
\[
 R[t^{-1}]=L[t,t^{-1}],\qquad \overline R=L[t],\qquad
 (R_+)^a=\bigoplus_{j\geq a}R_j.
\]
The last identity identifies the tangent cone with $R$.
When $L/F$ is separable, its projective point splits into $[L:F]$
distinct geometric points, whose affine cones are the asserted
reduced lines.
\end{proof}

\begin{theorem}\label{dim:twenty-cone}
At the generic homogeneous punctual divisor of $H_{20,4}$, remove the
four translation parameters and put $F=k(\operatorname{Gr}(5,20))$.
The reduced principal local ring is the vertex local ring of a
standard graded one-dimensional $F$-domain generated by at most
$21$ elements of degree one. Its reduced ordered cover has the same
property with at most $97$ generators. Consequently its tangent cone
is reduced and consists geometrically of exactly $20!h_{20}$ lines,
and
\[
 \operatorname{edim}(B_{20})\leq97.
\]
These conclusions permit the full Hilbert scheme to have other
components or nonreduced structure along the divisor.
\end{theorem}

\begin{proof}
The certificate exhibits a homogeneous algebra $A=S/I$ of Hilbert
function $(1,4,10,5)$ for which $S_1\otimes I_3\to S_4$ is surjective.
Thus $I$ has $15$ cubic generators and $60-35=25$ linear syzygies.
Their multiples span all $94$ degree-five syzygies.
Since $\operatorname{reg}(I)\leq4$, the first syzygies can only have
degrees four and five, giving a presentation
\[
 S(-4)^{25}\longrightarrow S(-3)^{15}\longrightarrow I\longrightarrow0.
\]

Applying $\operatorname{Hom}(-,A)$ gives negative-degree coefficient
maps with the following dimensions and ranks:
\[
\begin{array}{c|rrr}
 \text{degree of homomorphism}&-1&-2&-3\\ \hline
 \text{source dimension}&150&60&15\\
 \text{target dimension}&125&250&100\\
 \text{rank}&125&60&15\\
 \text{kernel dimension}&25&0&0
\end{array}
\]
There are no lower negative homomorphisms because the ideal is
generated in degree three. Nonzero modular minors give the same rank lower bounds over
$\mathbb Q$. The matrix sizes force equality in degree $-1$,
and all these rank conditions hold on a nonempty Grassmannian open set.

Use the border-basis chart consisting of all monomials of degree at
most two and five monomials of degree three. Every border-basis
coefficient has nonnegative dilation weight. The chart therefore
retracts to its homogeneous fixed locus. Pass to its fibre over the
generic point of the Grassmannian. Its positive cotangent space is
the dual of the negative tangent space just computed. Trace
barycentering splits off four independent translation directions in
weight one. The entire centered fibre, with its possible nilpotents
and extra components, consequently has positive cotangent
space concentrated in weight one and dimension $21$.

Lemma~\ref{dim:degree-one-cone} applies before taking the reduced
principal component. The principal quotient is an integral graded
algebra of dimension $80-4-75=1$, hence a standard graded domain with
at most $21$ generators. The reduced ordered cover is obtained by
adjoining the centered ordered coordinates. They have weight one and
span a space of dimension $4(20-1)=76$. It is therefore a standard
graded domain with at most $21+76=97$ generators. The radial theorem
identifies its normalization degree as $20!h_{20}$. The lemma now
gives its reduced tangent cone and the embedding-dimension bound.
\end{proof}

\subsection{Equations for the filtered deformations}

The same source makes the reconstruction problem finite and explicit.
Let $W_1,\ldots,W_{15}$ be its cubic generators, and let
$\ell_1,\ldots,\ell_{25}$ be a basis of their linear syzygies. Write a
filtered deformation as
\[
 f_a=W_a+\phi_a+\psi_a+\chi_a,
 \qquad \phi_a\in S_2,\quad\psi_a\in S_1,\quad\chi_a\in k.
\]
The linear tangent condition requires $\ell_s\phi\in W$. Write
$\ell_s\phi=C_s(\phi)W$, where $C_s$ is linear in the tangent
parameters. Then the full filtered-flatness equations are
\begin{equation}\label{dim:pbw-equations}
 \ell_s\psi=C_s(\phi)\phi,\qquad
 \ell_s\chi=C_s(\phi)\psi,\qquad
 C_s(\phi)\chi=0,
 \quad 1\leq s\leq25.
\end{equation}
The equations are the homogeneous components of
\[
 (\ell_s-C_s(\phi))(W+\phi+\psi+\chi)=0,
\]
using $\ell_sW=0$ and $\ell_s\phi=C_s(\phi)W$.
Thus every generating syzygy lifts. Reducing any relation by these
lifts lowers its filtration degree, so the initial ideal remains
$(W)$. Conversely, filtered flatness requires these syzygy lifts
and therefore the displayed equations.

The vanishing of the degree $-2$ and $-3$ tangent kernels makes
$\psi$ and $\chi$ unique polynomial functions of $\phi$ whenever they
exist. After a linear translation slice, this gives a presentation in
$21$ variables with equations of degree at most four. In the
certificate modulo $32003$, the quadratic obstruction map has rank
$154$ on the $231$ quadratic monomials. The coefficient choices use
invertible linear-algebra pivots, so the quadratic relation space
has generic characteristic-zero dimension at least $154$.
The next argument uses this bound to control the degree of the
principal component inside the filtered deformation scheme.

\begin{corollary}\label{dim:reconstruction-bound}
The generic reconstruction degree of the twenty-point collision map
satisfies
\[
 1\leq h_{20}\leq2^{30}=1\,073\,741\,824.
\]
\end{corollary}

\begin{proof}
The vanishing and generation conditions used in
\eqref{dim:pbw-equations} hold on a nonempty open subset of the
Grassmannian, so they give the generic centered fibre by
equations of degrees at most four in $21$ degree-one variables. The
principal Hilbert component is a component of the full Hilbert scheme
and dominates this Grassmannian. Its prime therefore remains a
minimal prime after localization to the generic point. Its
one-dimensional centered cone consequently gives an isolated
zero-dimensional component $\operatorname{Spec}(L)$ of the
projectivized fibre in $\mathbb P^{20}_F$, of degree $[L:F]=h_{20}$.
Let $Q$ be the ideal generated by the quadratic relations. Its
quotient has degree-two dimension at most $231-154=77$.
If a component had affine dimension $d\geq12$, linear Noether
normalization would give
\[
 F[y_1,\ldots,y_d]\hookrightarrow F[z_1,\ldots,z_{21}]/\mathfrak p,
 \qquad
 78\leq\binom{d+1}{2}\leq77,
\]
a contradiction. Hence $\operatorname{ht}Q\geq10$.
Over the infinite field $F$, ten generic linear combinations of
its quadrics form a regular sequence. Their complete intersection
in $\mathbb P^{20}_F$ has dimension ten and degree $2^{10}$.
B\'ezout's inequality bounds the total degree of the isolated
components cut out by the remaining equations by
\[
 h_{20}\leq2^{10}\,4^{10}=2^{30}.
\]
The isolated-component form of the inequality allows the full
deformation fibre to have positive-dimensional components.
\end{proof}

\begin{corollary}\label{dim:four-uniform-pd}
For each seed length $m\in\{20,\ldots,26\}$, the eventual projective
and flat dimensions of the marked polygraph over any coordinate
alphabet are exactly $m-2$, in every ambient dimension at least four.
Choose mark thresholds $L_m$ for these seven seeds and put
$L=\max_m L_m$. For every $n\geq80$, set $b=\lceil n/26\rceil$.
Then for every $\ell\geq bL$,
\[
 \operatorname{pd}_{P_j}A_\ell\geq
 n-2\left\lceil\frac n{26}\right\rceil,
 \qquad
 \operatorname{fd}_{P_j}A_\ell\geq
 n-2\left\lceil\frac n{26}\right\rceil.
\]
At a partition collision into $b$ blocks of sizes between $20$ and
$26$, the localized dimensions equal the displayed lower bound.
Invariant Tor under the product of the block symmetric groups is
nonzero in every degree from one through this value.
\end{corollary}

\begin{proof}
Apply the punctual-divisor invariant Tor theorem to the seven seeds.
For $n\geq80$, the integer intervals $[20b,26b]$, $b\geq4$, overlap
or are adjacent, and $b=\lceil n/26\rceil$ satisfies $20b\leq n$.
Thus $n$ is the sum of $b$ integers in $[20,26]$. The separated-block
theorem gives the local dimension
$\sum(m_i-2)=n-2b$ and the invariant classes in every lower degree.
Ambient projection and its zero section give the same mark constants
in every larger ambient dimension.
\end{proof}

\section{The three-dimensional divisor and the ambient boundary}
\label{dim:three-boundary}

\begin{theorem}\label{dim:three-divisor}
The principal Hilbert component $H_{96,3}$ has a homogeneous punctual
divisor whose general point has Hilbert function
\[
 (1,3,6,10,15,21,28,12).
\]
Its centered collision image has dimension $284$ and codimension
four in $\operatorname{Gr}(12,36)$. At its generic ordered point,
with unordered reconstruction degree $h_{96}$,
\[
 b=e=96!h_{96},\qquad
 \delta\geq231\cdot96!h_{96}.
\]
The local ring fails $R_1$, seminormality and the Du Bois property.
\end{theorem}

\begin{proof}
For the certified rational configuration, evaluation is injective
through degree six and surjective through degree seven. At the
prime $p=2147483647$, every evaluation and Grassmannian-chart
denominator is invertible, and a $284\times284$ collision minor has
determinant $1915860956\ne0$. The translation and dilation directions
give the reverse rank bound, so
\[
 284\leq\operatorname{rank}_{\mathbb Q}(dc)
 \leq3\cdot96-4=284.
\]
The accompanying certificate records the matrix and its four kernel
directions. Proposition~\ref{dim:collision-criterion} gives the
divisor, while $\dim\operatorname{Gr}(12,36)=12\cdot24=288$
gives codimension four.

The radial theorem and Theorem~\ref{dim:general-sign-defect} give the
branch, multiplicity and defect assertions, since the minimum
alternating degree is $462$. For the seminormality claim, use a
monomial-basis chart of the Hilbert scheme. Its point algebra is
represented by $d$ multiplication matrices of size $n$; adjoining
the ordered point coordinates gives an affine embedding with at most
$dn^2+nd$ coordinates. Localizing at a divisor of dimension $nd-1$
therefore gives
\[
 \operatorname{edim}(B)\leq dn^2+1.
\]
At $(96,3)$,
\[
 \operatorname{edim}(B)\leq27649<96!\leq b(B),
\]
so \eqref{fh:j:edim} rules out seminormality. Du Bois singularities are seminormal
\cite[Remark~4.8]{gg-SchwedeDuBois}. Normality already fails because
the full noncyclic inertia group $S_{96}$ fixes the original divisor.
\end{proof}

\begin{corollary}\label{dim:all-dimensional-counterexamples}
For $d=3$, $n\geq96$, and for $d=4$, $n\geq20$, the principal
isospectral cover fails $R_1$, seminormality and the Du Bois property
along a divisor. Projection to each coordinate alphabet has a fibre
of dimension at least
\[
 n(d-1)+94\quad(d=3),\qquad n(d-1)+18\quad(d=4).
\]
All sufficiently highly marked polygraphs and all sufficiently large
powers of the alternating ideal are nonflat over each alphabet.
The mark constants can be chosen independently of the added points,
added ambient coordinates and chosen alphabet. At lengths $96$ and
$20$ the eventual projective and flat dimensions are exactly $94$
and $18$, respectively, with invariant Tor in every positive degree
up to these values.
\end{corollary}

\begin{proof}
The twenty-point cover is not seminormal because its branch number
is at least $20!$, while Theorem~\ref{dim:twenty-cone} bounds its
embedding dimension by $97$. The preceding theorem handles length
$96$. The graph-bundle construction and disjoint-support products
used for length ten preserve these local failures and extend them
to the stated ranges. A divisor with one block of size $m$ and the
remaining points reduced has dimension $nd-1$. Its image in the
selected alphabet has dimension at most $n-m+1$, so an alphabet
fibre has dimension at least
\[
 nd-1-(n-m+1)=n(d-1)+m-2.
\]
The invariant chart theorem and Theorem~\ref{d3:punctual-invariant-tor}
transfer the obstruction to the polygraph rings and to powers of the
alternating ideal. Their upper bound is $m-2$ at the seed length.
Ambient projection, its zero section and the separated-support
retractions give the stated uniformity.
\end{proof}

\begin{corollary}[The sharp ambient-dimension boundary]
\label{dim:ambient-boundary}
For a positive integer $d$, the following conditions are equivalent:
\begin{enumerate}[label=\textup{(\roman*)}]
\item $d\leq2$;
\item $X_{n,d}$ is normal for every $n$;
\item $k[Z(n,\ell;d)]$ is free over each coordinate alphabet for
every $n$ and $\ell$.
\end{enumerate}
\end{corollary}

\begin{proof}
Haiman proves both assertions for $d=2$. For $d=1$,
$\operatorname{Hilb}^n(\mathbb A^1)=\operatorname{Sym}^n(\mathbb A^1)$,
so $X_{n,1}=\mathbb A^n$. The one-dimensional polygraph has coordinate
ring
\[
 k[x_1,\ldots,x_n,z_1,\ldots,z_\ell]\Big/
 \left(\prod_{i=1}^n(z_a-x_i):1\leq a\leq\ell\right).
\]
The monic relations make this a free module of rank $n^\ell$ over
$k[x_1,\ldots,x_n]$, with basis the monomials whose individual
$z_a$-exponents are less than $n$. The generic fibre is reduced;
freeness embeds the algebra into that generic fibre and hence the
whole algebra is reduced. It is therefore the indicated polygraph
ring. For every $d\geq3$, the preceding counterexamples contradict
each of (ii) and (iii).
\end{proof}

The new four-dimensional divisors also meet known nonreduced
components. Jelisiejew proves that the very compressed loci with
Hilbert functions $(1,4,10,r)$ for $r=6,7,8,9$ give generically
nonreduced elementary components \cite[Theorem~1.1]{dim:Jelisiejew}.
Our principal-component collision loci lie in these homogeneous
loci with codimensions $5,8,9,8$, respectively. Thus each of these
elementary components meets the principal component along a
punctual divisor. The length-$25$ calculation has codimension five
in the $r=10$ homogeneous locus and does not settle its generic
nonreducedness.

The local constructions also extend by \'etale coordinates to any
smooth variety of the indicated dimension over a characteristic-zero
field. In particular, a smooth variety of dimension at least three
has nonnormal principal isospectral covers for all lengths in the
corresponding range.

\section{The transverse orbit cone}
\label{new:orbit}

For the geometric seed we work over \(k=\mathbb C\).
Let
\[
 D^0=\operatorname{Gr}(4,15)\subset H^0
\]
be the centered homogeneous divisor constructed in~\ref{new:divisor}.
Here \(H^0\) is the trace-centered principal component of
\(\operatorname{Hilb}^{10}(\mathbb A^5)\). Write
\[
 F=k(D^0),\qquad
 A_0=\mathcal O_{H^0,\eta_{D^0}},
\]
and let \(B\) be the local ring of the reduced ordered cover at its unique
point above \(\eta_{D^0}\). The ring \(A_0\) is a discrete valuation ring.
Theorem~\ref{new:radial} identifies the normalization of \(B\) with the
radial valuation ring of the centered ordered configuration space.

\subsection{Coordinates and the graded local model}

Put \(W=\operatorname{Std}_{10}\otimes k^5\), where
\(\operatorname{Std}_{10}\) is the nine-dimensional standard representation
of \(S_{10}\). If \(x_{ij}\) denotes the \(j\)-th coordinate of the
\(i\)-th centered point, then
\[
 \sum_{i=1}^{10}x_{ij}=0\quad(1\le j\le5),\qquad
 L=k(W),\qquad K=L^{S_{10}}=k(H^0).
\]
The coordinates have dilation weight one. Choose an \(S_{10}\)-invariant
homogeneous rational function \(t\) of weight one. For example, with
\(p_a=\sum_i x_{i1}^{\,a}\), one may take
\[
 t=p_5/p_4.
\]
The radial valuation gives \(v(t)=1\). Thus \(t\) is a uniformizer of
the normalization DVR \(C\). Its invariant subring \(C^{S_{10}}\)
is integral over \(A_0\), with fraction field \(K\). Since \(A_0\)
is normal, \(C^{S_{10}}=A_0\), and \(t\) is also a uniformizer of
\(A_0\). Set
\[
 v_{ij}=x_{ij}/t,\qquad E=k(v_{ij}).
\]
The \(v_{ij}\) are dilation invariant and identify \(E\) with the
radial residue field \(k(\mathbb P(W))\). Exactness of finite-group
invariants in characteristic zero identifies the residue field of
\(A_0=C^{S_{10}}\) with \(E^{S_{10}}\). Since \(L=E(t)\), taking
invariants gives \(K=E^{S_{10}}(t)\); its valuation is the ordinary
\(t\)-adic valuation. Consequently
\begin{equation}
 \label{new:radial-fields}
 E=k(\mathbb P(W)),\qquad L=E(t),\qquad
 F=E^{S_{10}},\qquad A_0=F[t]_{(t)}.
\end{equation}
In particular \(E/F\) is an \(S_{10}\)-Galois field extension of degree
\(10!\). Its \(10!\) geometric embeddings label the orderings of the
generic configuration.

For \(r\ge0\), let \(F_r\subset E\) be the \(F\)-span of the monomials of
total degree at most \(r\) in the \(v_{ij}\), and set \(F_{-1}=0\).
These spaces form the polynomial filtration:
\[
 F_0=F,\qquad F_rF_s\subset F_{r+s},\qquad
 \bigcup_{r\ge0}F_r=E.
\]
Writing \(V=F_1\), one has \(F_r=V^r\), because \(1\in V\).

\begin{theorem}[The transverse local ring]
\label{new:local-cone}
With the preceding notation, define the standard graded algebra
\[
 R=\bigoplus_{r\ge0}t^rF_r
   =F[t,tv_{ij}]\subset E[t].
\]
Then
\begin{equation}
 \label{new:cone-identifications}
 B=R_{R_+},\qquad
 \operatorname{gr}_{\mathfrak m_B}B\simeq R,\qquad
 \widehat B\simeq\prod_{r\ge0}t^rF_r.
\end{equation}
The normalization of \(R\) is \(E[t]\), and the normalization of \(B\)
is \(C=E[t]_{(t)}\). Moreover
\begin{equation}
 \label{new:degree-one-representation}
 F_1\simeq{\bf1}\oplus\operatorname{Std}_{10}^{\oplus5},
 \qquad \operatorname{edim}B=46.
\end{equation}
After extending \(F\) to a separable closure, the tangent cone is the
reduced union of \(10!\) distinct lines with projective directions
\[
 [1:v_{\sigma(i)j}]\in\mathbb P^{45},
 \qquad \sigma\in S_{10}.
\]
\end{theorem}

\begin{proof}
The ordered cover is finite over the principal Hilbert component and is
generated by the ordered point coordinates. Its fibre above the generic
punctual point has one point, with residue field \(F\). Consequently its
finite algebra over \(A_0\) is local, and
\[
 B=A_0[x_{ij}]
   =F[t]_{(t)}[tv_{ij}].
\]
A monic equation for each \(v_{ij}\) homogenizes to an integral
equation for \(tv_{ij}\) over \(F[t]\). It also shows that
\(tv_{ij}\) is nilpotent modulo \(tR\). Thus \(R\) is finite over
\(F[t]\), with \(R_+\) its unique prime above \((t)\), proving
\(B=R_{R_+}\).

Inverting \(t\) recovers every \(v_{ij}\), so
\[
 R[t^{-1}]=E[t,t^{-1}],\qquad \operatorname{Frac}(R)=E(t).
\]
The finite integral extension \(R\subset E[t]\) is therefore its
normalization. Since \(R\) is generated in degree one,
\[
 (R_+)^a=\bigoplus_{r\geq a}t^rF_r.
\]
Localization preserves the initial homogeneous pieces, giving the
associated graded ring and completion in~\eqref{new:cone-identifications}.

For degree one, pass to a splitting field \(\overline F\) of the
torsor and write \(U=\operatorname{Std}_{10}\otimes\overline F\).
The five centered coordinate columns span a five-dimensional subspace
\(W'\subset U\). The matrix-coefficient map of the irreducible
representation \(U\) is injective, hence so is its restriction
\[
 U^*\otimes W'\longrightarrow\overline F^{S_{10}},
 \qquad \alpha\otimes w\longmapsto[g\mapsto\alpha(gw)].
\]
Its image is the span of the ordered coordinate functions and has
dimension \(9\cdot5=45\). Since this image has no invariant vectors,
adjoining the constant function gives
\(F_1\simeq{\bf1}\oplus\operatorname{Std}_{10}^{\oplus5}\).
The representation and dimension descend to \(F\).

Under the splitting \(E\otimes_F\overline F\simeq\overline F^{S_{10}}\),
the normalization becomes a product of polynomial rings, one for each
ordering. The homomorphism from the polynomial ring on the degree-one
generators evaluates those generators as \(t\) times
\((1,v_{\sigma(i)j})\). Its kernel is the intersection of the corresponding
line ideals. The directions are distinct because the points of the
configuration are distinct. This identifies the geometric tangent cone.
\end{proof}

The homogeneous-limit map also acquires a birational description:
\begin{equation}
 \label{new:grassmannian-quotient}
 k(\operatorname{Gr}(4,15))
 =k(\mathbb P(\operatorname{Std}_{10}\otimes k^5))^{S_{10}}.
\end{equation}
Thus a general homogeneous \((1,5,4)\) algebra determines an unordered
centered configuration up to dilation. The next statement gives an
open neighborhood on which this description is a family.

\begin{corollary}[An open line-bundle model]
\label{new:line-bundle}
There is a nonempty open set \(U\subset D^0\) and an equivariant open
immersion
\[
 U\times\mathbb A^1\longrightarrow H^0
\]
that identifies \(U\times\{0\}\) with \(U\). Dilation fixes \(U\) and
acts with weight one on \(\mathbb A^1\). After shrinking \(U\), the
normalization of the ordered cover over this open set is
\[
 \widetilde U\times\mathbb A^1
 \longrightarrow U\times\mathbb A^1,
\]
where \(\widetilde U\to U\) is a finite \'{e}tale \(S_{10}\)-torsor.
The ordered cover itself is defined by the relative subalgebra generated
by \(t\) and the \(tv_{ij}\).
\end{corollary}

\begin{proof}
Let \(T=E^{S_9}\), where \(S_9\) fixes the first label. This is the
degree-ten point algebra over \(F\). Its polynomial filtration in the
five coordinates \(v_{1j}\) has dimensions \(1,6,10\) in degrees
zero, one, and two. One can see this directly from the tautological
family over \(A_0\): choose a basis of its special fibre consisting of
the unit, five linear monomials, and four quadratic monomials.
These lift to a basis over the DVR. Dividing their generic values by
the corresponding powers \(1,t,t^2\) gives a basis of \(T\) with the
stated filtered degrees.

The algebra
\[
 \mathcal T=\bigoplus_{r\ge0}t^rT_{\le r}\subset T[t]
\]
is therefore free of rank ten over \(F[t]\), with basis degrees
\(0,1^5,2^4\). It is generated by \(t\) and the five functions
\(tv_{1j}\), so it defines a flat family of subschemes of
\(\mathbb A^5\). At \(t=0\) it is the generic homogeneous
\((1,5,4)\) algebra, and at every nonzero \(t\) it is reduced.

Spread the finite multiplication table, its filtered basis, and its
generators over a sufficiently small affine open \(U\subset D^0\).
Shrink \(U\) so that the point algebra and its splitting algebra are
finite \'{e}tale, the basis remains a basis, and the special-fibre
classifying map is the identity of \(U\). This produces an equivariant
morphism
\[
 \Phi:U\times\mathbb A^1\longrightarrow H^0.
\]
It induces the function-field isomorphism \(k(H^0)=F(t)\).
Along the zero section its differential is the identity in the divisor
directions and is nonzero in the normal direction, since \(t\) is a
uniformizer. After another shrinking of \(U\), the morphism is
\'{e}tale along the entire zero section.

The non\'{e}tale locus is closed and dilation invariant. Since
\[
 (u,0)\in\overline{\mathbb G_m\cdot(u,t)}\qquad(t\ne0),
\]
any nonempty such locus would meet the zero section. Thus \(\Phi\)
is \'{e}tale everywhere. Its image lies in the smooth locus of
\(H^0\), by construction at \(t=0\) and because the other fibres
are reduced configurations. A separated birational \'{e}tale morphism
to a normal variety is an open immersion.

Finally spread \(E/F\) as a finite \'{e}tale \(S_{10}\)-torsor over
\(U\). The algebra generated by \(t,tv_{ij}\) is the ordered-coordinate
algebra on this open set, and its normalization is the asserted product.
\end{proof}

\subsection{Interpolation and the conductor}

Polynomial interpolation determines when the filtration reaches
all functions on the orderings. That degree is also the conductor
exponent of the orbit cone.

Let \(n\ge2\), \(d\ge1\), let \(F\) be a characteristic-zero field,
and let \(E\) be a finite
\'{e}tale \(S_n\)-torsor algebra of rank \(n!\). After a separable
extension, identify \(E\) with the functions on all orderings of a
configuration \(p_1,\ldots,p_n\). The elements \(v_{ij}\) are the
ordered coordinate functions. As above, \(F_r\) is their polynomial
filtration and
\[
 R=\bigoplus_{r\ge0}t^rF_r.
\]
All rank conditions below may be checked over a separable closure.

\begin{lemma}
\label{new:interpolation-open}
A general configuration of \(n\) points in \(\mathbb A^d\) satisfies
\begin{equation}
 \label{new:maximal-interpolation}
 \operatorname{rank}\left(
 F[z_1,\ldots,z_d]_{\le j}\longrightarrow F^S\right)
 =\min\left\{|S|,\binom{d+j}{j}\right\}
\end{equation}
for every subset \(S\) of the configuration and every \(j\ge0\).
\end{lemma}

\begin{proof}
For a fixed degree \(j\), point evaluations span the dual of the space
of polynomials of degree at most \(j\). Otherwise a nonzero polynomial
would vanish at every point of affine space. Hence one can choose any
number up to \(\binom{d+j}{j}\) of independent evaluations, and add
further distinct points if needed. For each fixed subset, the rank
condition is therefore nonempty and open. Only finitely many subsets
and degrees up to the first \(j\) with
\(\binom{d+j}{j}\ge n\) need be considered. Their open sets have
nonempty intersection in the irreducible configuration space.
\end{proof}

\begin{theorem}[Conductor of a permutation-orbit cone]
\label{new:orbit-conductor}
Assume~\eqref{new:maximal-interpolation} for every subset. Define
\[
 N_j=\binom{d+j}{j},\qquad N_{-1}=0,\qquad
 \nu_m=\min\{j:N_j\ge m\},\qquad
 \mu=\mu(n,d)=\sum_{m=1}^n\nu_m.
\]
Then
\[
 F_\mu=E,\qquad F_{\mu-1}\ne E.
\]
The normalization of \(R\) is \(E[t]\), and its conductor is exactly
\(t^\mu E[t]\). At the homogeneous vertex,
\begin{equation}
 \label{new:delta-filtration}
 \delta=\sum_{r=0}^{\mu-1}(n!-\dim_FF_r).
\end{equation}
Put
\begin{equation}
 \label{new:newton-polynomial}
 Q_{n,d}(q)=
 \prod_{m=1}^n
 \sum_{j\ge0}
 \bigl(\min(m,N_j)-\min(m,N_{j-1})\bigr)q^j.
\end{equation}
There is a basis of the orbit algebra whose elements admit polynomial
representatives with degree bounds encoded by \(Q_{n,d}\). In particular
\begin{equation}
 \label{new:newton-bounds}
 \dim_FF_r\ge\sum_{j=0}^r[q^j]Q_{n,d}(q),
 \qquad
 \delta\le Q_{n,d}'(1).
\end{equation}
\end{theorem}

\begin{proof}
Work first over a splitting field. Order any \(m\)-element subset as
\(p_1,\ldots,p_m\). By~\eqref{new:maximal-interpolation}, for each
\(j\) there is a polynomial \(b_j\) of degree \(\nu_j\) such that
\[
 b_j(p_i)=0\quad(i<j),\qquad b_j(p_j)\ne0.
\]
Indeed, evaluations on the first \(j\) points have full rank \(j\) in
degree \(\nu_j\). In smaller degree, the evaluations on the preceding
points already exhaust the polynomial space, so the chosen degree is
minimal.

Induct on \(m\). For each possible first value \(p_j\), multiply
\(b_j\) in the first slot by a recursively constructed basis for the
remaining \(m-1\) points. Order the evaluation rows and columns by
this first value. If \(M_j\) is the evaluation matrix for the remaining
set, the resulting block matrix \(M\) satisfies
\[
 M_{ij}=0\ (i<j),\qquad
 M_{jj}=b_j(p_j)M_j,\qquad
 \det M=\prod_{j=1}^m b_j(p_j)^{(m-1)!}\det M_j\ne0.
\]
Thus these polynomial functions form a basis of the \(m!\) orderings.

If \(Q_m\) records its degrees, this construction gives
\[
 Q_m(q)=Q_{m-1}(q)\sum_{j=1}^m q^{\nu_j},\qquad Q_0=1.
\]
The number of indices \(a\le m\) with \(\nu_a=j\) is
\(\min(m,N_j)-\min(m,N_{j-1})\), so this is
~\eqref{new:newton-polynomial}. Its largest degree is
\(\sum_m\nu_m=\mu\), and \(Q_{n,d}(1)=n!\).
Thus \(F_\mu=E\). The functions of basis degree at most \(r\)
are independent members of \(F_r\), proving the first inequality in
~\eqref{new:newton-bounds}. These rank assertions descend to \(F\);
choosing a basis along the resulting filtration over \(F\) gives
the same degree bounds.

To see that degree \(\mu\) is necessary, alternate a monomial in the
\(n\) blocks of \(d\) variables. Every alternating polynomial is a
linear combination of determinants
\[
 \det\bigl(m_j(x_i)\bigr)_{1\le i,j\le n},
\]
where the \(m_j\) are distinct monomials in \(d\) variables. A repeated
monomial gives a zero determinant. The smallest total degree of
\(n\) distinct monomials is exactly
\(\sum_{m=1}^n\nu_m=\mu\). Hence the sign component of the
polynomial space is zero below degree \(\mu\).
The torsor algebra contains the sign representation, so
\(F_{\mu-1}\ne E\).

The normalization is \(E[t]\), by the same finite-integrality
argument as in Theorem~\ref{new:local-cone}. Its conductor is graded, and its degree-\(r\) coefficient space is
\[
 \{a\in E:aE\subseteq F_r\}.
\]
Over a splitting field, let \(e_\sigma\) be the indicator of one
ordering. If \(a\ne0\), choose \(\sigma\) with \(a(\sigma)\ne0\).
Then
\[
 aE\subseteq F_r
 \quad\Longrightarrow\quad
 a e_\sigma=a(\sigma)e_\sigma\in F_r
 \quad\Longrightarrow\quad e_\sigma\in F_r.
\]
The \(S_n\)-translates of \(e_\sigma\) span \(E\), so \(F_r=E\).
Thus the conductor starts in degree \(\mu\) and equals \(t^\mu E[t]\).

Finally \(E[t]/R\) is supported at the vertex and has degree-\(r\)
part \(E/F_r\). This proves~\eqref{new:delta-filtration}.
The basis-degree bound gives
\[
 \sum_{r=0}^{\mu-1}(n!-\dim F_r)
 \le\sum_{j\ge0}j[q^j]Q_{n,d}(q)=Q_{n,d}'(1).
\]
\end{proof}

The polynomial \(Q_{n,d}\) records a degree-controlled basis. Further
linear combinations may occur in smaller degree, so its coefficients
are used through the inequalities in~\eqref{new:newton-bounds}.

\begin{corollary}[The ten-point conductor]
\label{new:ten-conductor}
For the local ring of Theorem~\ref{new:local-cone},
\[
 F_{13}=E,\qquad F_{12}\ne E,\qquad
 \mathfrak c_{B\subset C}=t^{13}C.
\]
Moreover
\begin{equation}
 \label{new:ten-delta-upper}
 \delta(B)\le29\,746\,080.
\end{equation}
\end{corollary}

\begin{proof}
By~\eqref{new:grassmannian-quotient}, the generic point of \(D^0\)
corresponds to a generic centered configuration modulo dilation.
The interpolation conditions of Lemma~\ref{new:interpolation-open}
are preserved by translation and nonzero dilation, so they hold here.
For ten points in five-space, the degrees of ten smallest monomials are
\[
 0,\quad1,1,1,1,1,\quad2,2,2,2.
\]
Their sum is thirteen. The theorem gives the conductor, and its
degree polynomial specializes to
\[
 Q_{10,5}(q)=
 \prod_{m=1}^{10}
 \left(1+\min(5,m-1)q+\max(0,m-6)q^2\right).
\]
Since \(Q_{10,5}(1)=10!\),
\[
 \begin{split}
 Q_{10,5}'(1)
 &=10!\left(
   \sum_{m=2}^{6}\frac{m-1}{m}
   +\sum_{m=7}^{10}\frac{2m-7}{m}\right)\\
 &=10!\,\frac{2951}{360}
 =29\,746\,080.
 \end{split}
\]
\end{proof}

Figure~\ref{new:orbit-figure} illustrates the cone construction on a
six-ordering example. Its evaluation matrix makes visible the passage
from low-degree functions to the full algebra of branch directions.

\begin{figure}[!htbp]
\centering
\begin{tikzpicture}[x=1cm,y=1cm,font=\small]
\path[use as bounding box] (0,-1.25) rectangle (16.2,8.0);
\node[anchor=west,font=\bfseries\small,Ink] at (0,7.65) {(a) Six branches at one vertex};
\node[anchor=west,font=\bfseries\small,Ink] at (6.35,7.65) {(b) Evaluation distinguishes the six orderings};
\draw[Ink!20] (5.95,1.35)--(5.95,7.2);
\draw[Ink!30,fill=Teal!3] plot[domain=0:360,samples=80,smooth,variable=\a]
 ({2.7+1.94*cos(\a)+.35*sin(\a)},{5.55+.56*sin(\a)}) -- cycle;
\draw[Ink!25,densely dashed] (2.7,2.25)--(2.7,6.55);
\node[left,Ink!70] at (.45,5.55) {$1$};

\draw[Teal,line width=1.1pt] (2.7,2.3)--(0.760,5.550);
\draw[Teal!45,line width=.7pt] (2.7,2.3)--(2.933,1.910);
\fill[Teal] (0.760,5.550) circle (2.5pt);
\node[Teal,fill=white,inner sep=.5pt] at (0.760,5.750) {$p_1$};
\draw[Blue,line width=1.1pt] (2.7,2.3)--(2.033,6.035);
\draw[Blue!45,line width=.7pt] (2.7,2.3)--(2.780,1.852);
\fill[Blue] (2.033,6.035) circle (2.5pt);
\node[Blue,fill=white,inner sep=.5pt] at (2.033,6.235) {$p_2$};
\draw[Amber,line width=1.1pt] (2.7,2.3)--(1.427,5.065);
\draw[Amber!45,line width=.7pt] (2.7,2.3)--(2.853,1.968);
\fill[Amber] (1.427,5.065) circle (2.5pt);
\node[Amber,fill=white,inner sep=.5pt] at (1.427,4.835) {$p_3$};
\draw[Purple,line width=1.1pt] (2.7,2.3)--(3.973,6.035);
\draw[Purple!45,line width=.7pt] (2.7,2.3)--(2.547,1.852);
\fill[Purple] (3.973,6.035) circle (2.5pt);
\node[Purple,fill=white,inner sep=.5pt] at (3.973,6.235) {$p_4$};
\draw[Ink,line width=1.1pt] (2.7,2.3)--(3.367,5.065);
\draw[Ink!45,line width=.7pt] (2.7,2.3)--(2.620,1.968);
\fill[Ink] (3.367,5.065) circle (2.5pt);
\node[Ink,fill=white,inner sep=.5pt] at (3.367,4.835) {$p_5$};
\draw[Coral,line width=1.1pt] (2.7,2.3)--(4.640,5.550);
\draw[Coral!45,line width=.7pt] (2.7,2.3)--(2.467,1.910);
\fill[Coral] (4.640,5.550) circle (2.5pt);
\node[Coral,fill=white,inner sep=.5pt] at (4.640,5.750) {$p_6$};
\fill[Ink] (2.7,2.3) circle (2.7pt);
\node[below=5pt,align=center] at (2.7,2.15) {$t=0$: one vertex\\six branches in the normalization};
\node[align=center,font=\footnotesize] at (2.7,6.9) {$p_i=(1,a_i,b_i)$,\quad $(a_i,b_i,c_i)\in S_3\cdot(-1,0,1)$};
\node[align=center] at (2.7,.25) {$x^2+xy+y^2=t^2$\\$xy(x+y)=0$};

\fill[Blue!6] (8.660,3.480) rectangle (9.340,6.570);
\node at (9.000,6.890) {$1$};
\fill[Teal!9] (9.380,3.480) rectangle (10.060,6.570);
\node at (9.720,6.890) {$x$};
\fill[Teal!9] (10.100,3.480) rectangle (10.780,6.570);
\node at (10.440,6.890) {$y$};
\fill[Blue!6] (10.820,3.480) rectangle (11.500,6.570);
\node at (11.160,6.890) {$y^2$};
\fill[Blue!6] (11.540,3.480) rectangle (12.220,6.570);
\node at (11.880,6.890) {$xy$};
\fill[Amber!14] (12.260,3.480) rectangle (12.940,6.570);
\node at (12.600,6.890) {$\Delta/2$};
\node[anchor=east,font=\footnotesize] at (8.47,6.250) {$p_1:\;(-1,0,1)$};
\node[Teal] at (9.000,6.250) {$ 1 $};
\node[Amber] at (9.720,6.250) {$ -1 $};
\node[Ink] at (10.440,6.250) {$ 0 $};
\node[Ink] at (11.160,6.250) {$ 0 $};
\node[Ink] at (11.880,6.250) {$ 0 $};
\node[Amber] at (12.600,6.250) {$ -1 $};
\node[anchor=east,font=\footnotesize] at (8.47,5.750) {$p_2:\;(-1,1,0)$};
\node[Teal] at (9.000,5.750) {$ 1 $};
\node[Amber] at (9.720,5.750) {$ -1 $};
\node[Teal] at (10.440,5.750) {$ 1 $};
\node[Teal] at (11.160,5.750) {$ 1 $};
\node[Amber] at (11.880,5.750) {$ -1 $};
\node[Teal] at (12.600,5.750) {$ 1 $};
\node[anchor=east,font=\footnotesize] at (8.47,5.250) {$p_3:\;(0,-1,1)$};
\node[Teal] at (9.000,5.250) {$ 1 $};
\node[Ink] at (9.720,5.250) {$ 0 $};
\node[Amber] at (10.440,5.250) {$ -1 $};
\node[Teal] at (11.160,5.250) {$ 1 $};
\node[Ink] at (11.880,5.250) {$ 0 $};
\node[Teal] at (12.600,5.250) {$ 1 $};
\node[anchor=east,font=\footnotesize] at (8.47,4.750) {$p_4:\;(0,1,-1)$};
\node[Teal] at (9.000,4.750) {$ 1 $};
\node[Ink] at (9.720,4.750) {$ 0 $};
\node[Teal] at (10.440,4.750) {$ 1 $};
\node[Teal] at (11.160,4.750) {$ 1 $};
\node[Ink] at (11.880,4.750) {$ 0 $};
\node[Amber] at (12.600,4.750) {$ -1 $};
\node[anchor=east,font=\footnotesize] at (8.47,4.250) {$p_5:\;(1,-1,0)$};
\node[Teal] at (9.000,4.250) {$ 1 $};
\node[Teal] at (9.720,4.250) {$ 1 $};
\node[Amber] at (10.440,4.250) {$ -1 $};
\node[Teal] at (11.160,4.250) {$ 1 $};
\node[Amber] at (11.880,4.250) {$ -1 $};
\node[Amber] at (12.600,4.250) {$ -1 $};
\node[anchor=east,font=\footnotesize] at (8.47,3.750) {$p_6:\;(1,0,-1)$};
\node[Teal] at (9.000,3.750) {$ 1 $};
\node[Teal] at (9.720,3.750) {$ 1 $};
\node[Ink] at (10.440,3.750) {$ 0 $};
\node[Ink] at (11.160,3.750) {$ 0 $};
\node[Ink] at (11.880,3.750) {$ 0 $};
\node[Teal] at (12.600,3.750) {$ 1 $};
\node[anchor=west,font=\footnotesize] at (6.4,3.05) {$\Delta=(x-y)(2x+y)(x+2y)$};
\node[anchor=west,font=\footnotesize] at (6.4,2.55) {Degree $r$};
\node[anchor=west,font=\footnotesize] at (6.4,1.98) {$\dim F_r$};
\node[anchor=west,font=\footnotesize] at (6.4,1.41) {$6-\dim F_r$};

\node at (10.0,2.55) {$ 0 $};\node[Teal] at (10.0,1.98) {$ 1 $};\node[Amber] at (10.0,1.41) {$ 5 $};
\fill[Teal!65] (9.640,.76) rectangle ++(.105,.24);
\fill[Amber!22] (9.760,.76) rectangle ++(.105,.24);
\fill[Amber!22] (9.880,.76) rectangle ++(.105,.24);
\fill[Amber!22] (10.000,.76) rectangle ++(.105,.24);
\fill[Amber!22] (10.120,.76) rectangle ++(.105,.24);
\fill[Amber!22] (10.240,.76) rectangle ++(.105,.24);
\node at (11.45,2.55) {$ 1 $};\node[Teal] at (11.45,1.98) {$ 3 $};\node[Amber] at (11.45,1.41) {$ 3 $};
\fill[Teal!65] (11.090,.76) rectangle ++(.105,.24);
\fill[Teal!65] (11.210,.76) rectangle ++(.105,.24);
\fill[Teal!65] (11.330,.76) rectangle ++(.105,.24);
\fill[Amber!22] (11.450,.76) rectangle ++(.105,.24);
\fill[Amber!22] (11.570,.76) rectangle ++(.105,.24);
\fill[Amber!22] (11.690,.76) rectangle ++(.105,.24);
\node at (12.9,2.55) {$ 2 $};\node[Teal] at (12.9,1.98) {$ 5 $};\node[Amber] at (12.9,1.41) {$ 1 $};
\fill[Teal!65] (12.540,.76) rectangle ++(.105,.24);
\fill[Teal!65] (12.660,.76) rectangle ++(.105,.24);
\fill[Teal!65] (12.780,.76) rectangle ++(.105,.24);
\fill[Teal!65] (12.900,.76) rectangle ++(.105,.24);
\fill[Teal!65] (13.020,.76) rectangle ++(.105,.24);
\fill[Amber!22] (13.140,.76) rectangle ++(.105,.24);
\node at (14.35,2.55) {$ 3 $};\node[Teal] at (14.35,1.98) {$ 6 $};\node[Amber] at (14.35,1.41) {$ 0 $};
\fill[Teal!65] (13.990,.76) rectangle ++(.105,.24);
\fill[Teal!65] (14.110,.76) rectangle ++(.105,.24);
\fill[Teal!65] (14.230,.76) rectangle ++(.105,.24);
\fill[Teal!65] (14.350,.76) rectangle ++(.105,.24);
\fill[Teal!65] (14.470,.76) rectangle ++(.105,.24);
\fill[Teal!65] (14.590,.76) rectangle ++(.105,.24);
\draw[Ink!30] (6.4,.47)--(15.6,.47);
\node[anchor=west,align=left] at (6.4,-.15) {$\mathfrak c=t^3\widetilde R$,\qquad
 $\delta=5+3+1=9$\\[3pt]
 $\operatorname{Hilb}(R;q)=\dfrac{1+2q+2q^2+q^3}{1-q}$};
\end{tikzpicture}
\caption{A six-ordering permutation-orbit cone. Its lines are
$(t,x,y)=t(1,a_i,b_i)$; panel (a) shows their real positive halves
and the section $t=1$. In (b), adjoining columns of degrees $0,1,2,3$
gives ranks $1,3,5,6$. The alternating column completes the interpolation
basis. The deficits from the six-dimensional normalization are $5,3,1$,
so $\delta=9$ and the conductor exponent is $3$. This example illustrates
the degree filtration and conductor calculation in
Theorem~\ref{new:orbit-conductor}.}
\label{new:orbit-figure}
\end{figure}
\FloatBarrier

\subsection{Seminormalization and the Artinian reduction}

\begin{proposition}
\label{new:seminormal-ring}
For the ten-point transverse local ring,
\[
 B^{\mathrm{sn}}=F+tC,\qquad
 \widehat{B^{\mathrm{sn}}}=F+tE[[t]],
\]
and
\begin{equation}
 \label{new:seminormal-defect}
 \operatorname{length}_B(B^{\mathrm{sn}}/B)
 =\delta(B)-(10!-1).
\end{equation}
In particular \(B\) is not seminormal.
\end{proposition}

\begin{proof}
The ring \(F+tC\) is the inverse image of \(F\subset E=C/tC\).
It has the same residue field and fraction field as \(B\), and the
extension is subintegral. It is the largest intermediate ring with
these properties: any subintegral intermediate ring must still have
residue field \(F\).

For the square-and-cube criterion, suppose \(z^2,z^3\in F+tC\).
Then \(v(z)\geq0\), so \(z\in C\), and its residue \(\alpha\in E\)
satisfies
\[
 \alpha^2,\alpha^3\in F
 \quad\Longrightarrow\quad
 \alpha=0\ \text{or}\ \alpha=\alpha^3/\alpha^2\in F.
\]
Thus \(z\in F+tC\), proving seminormality and identifying
\(B^{\mathrm{sn}}\).

The residue map gives an exact sequence
\[
 0\longrightarrow B^{\mathrm{sn}}/B\longrightarrow C/B
 \longrightarrow E/F\longrightarrow0,
\]
so its lengths give~\eqref{new:seminormal-defect}.
The conductors of \(B^{\mathrm{sn}}\) and \(B\) have exponents
one and thirteen, respectively, hence \(B\ne B^{\mathrm{sn}}\).
\end{proof}

For subsequent use, put
\begin{equation}
 \label{new:artinian-reduction}
 A=R/(t)=\bigoplus_{r\ge0}F_r/F_{r-1},\qquad
 h_r=\dim_FF_r-\dim_FF_{r-1}.
\end{equation}
This associated graded coordinate algebra is the construction used in
orbit harmonics; see~\cite[Section~2.2]{new:Rhoades}. In a polynomial
presentation \(E=P/I\), it is \(A\simeq P/I^{\mathrm{top}}\), where
\(I^{\mathrm{top}}\) is generated by the highest-degree homogeneous
parts of the elements of \(I\).
The element \(t\) is a nonzerodivisor, and \(A\) has length \(10!\).
In the ten-point case its top nonzero degree is thirteen,
\(h_0=1\), and \(h_1=45\). More generally, for the orbit cone of
Theorem~\ref{new:orbit-conductor},
\begin{equation}
 \label{new:delta-h-vector}
 \sum_r h_r=n!,\qquad
 \delta=\sum_r r h_r.
\end{equation}
The second identity follows by summing the tails of the Hilbert
numerator in~\eqref{new:delta-filtration}. It will connect the
normalization defect with Gorenstein duality.

\Needspace{10\baselineskip}
\section{Further homogeneous punctual divisors}
\label{new:fifteen}

Generic regularity along a homogeneous punctual divisor identifies
its transverse local ring with the orbit cone. The negative-degree
Hilbert tangent proves this regularity for the Erman--Velasco divisor
at length fifteen and for a second Grassmannian divisor in five-space.

\begin{lemma}\label{new:fifteen-quantitative}
The centered smoothable locus $D_{15}^0$ of homogeneous
algebras with Hilbert function $(1,11,3)$ is irreducible of
dimension $153$. Its translates form a divisor
$D_{15}\subset H_{15,11}$ of dimension $164$.
Projection to any coordinate alphabet has a fibre of
dimension at least $163$, exceeding the generic dimension
$150$ by at least $13$.
\end{lemma}

\begin{proof}
Erman--Velasco construct the centered locus from homogeneous
limits of reduced configurations and prove its dimension
$153$ \cite[Propositions~3.1--3.2, Lemma~3.4 and proof of
Proposition~1.6]{gg-ErmanVelasco}. After placing twelve points
at $0,e_1,\ldots,e_{11}$, the remaining points $a^{(j)}$
give the inverse-system quadrics
\[
 q_j=\sum_i a_i^{(j)}y_i^2-
       \left(\sum_i a_i^{(j)}y_i\right)^2,
 \qquad1\le j\le3,
\]
with ambient linear changes of coordinates allowed.
The trace barycenter splits off eleven translations, giving
dimension $164$ in the $165$-dimensional principal component.
The divisor maps to the small diagonal in a coordinate
alphabet, so fixing the common coordinate leaves dimension
$163$. The generic projection fibre has dimension
$15(11-1)=150$.
\end{proof}

\subsection{Regularity of the principal component at the fifteen-point divisor}
\label{next:fifteen-regularity-section}

The full Hilbert scheme has intersecting components along the
Erman--Velasco divisor. The negative tangent calculation below
determines the generic local ring of its reduced principal component.

\begin{lemma}[A one-dimensional attracting tangent]
\label{next:attracting-criterion}
Let $D$ be an irreducible locus of centered homogeneous quotients
with Hilbert function $(1,d,e)$, of codimension one in the reduced
principal component $H^0_{1+d+e,d}$. Suppose that at one point of
$D$ the negative-degree Hilbert tangent has dimension at most
$d+1$ in degree $-1$ and is zero in all lower degrees.
Then the principal component is regular at the generic point of $D$.
\end{lemma}

\begin{proof}
A quotient-basis chart with degree distribution $(1,d,e)$ has
nonnegative dilation weights. Its weight-zero projection is the
homogeneous-limit map, and the tangent to its fibre at a homogeneous
point is the negative-degree part of $\operatorname{Hom}_S(I,S/I)$.
The trace barycenter splits off the $d$ translation directions,
all in degree $-1$. Thus the centered attracting fibre has tangent
dimension at most one at the specified point. Upper semicontinuity
gives the same bound at the generic point of $D$.

Pass to $F=k(D)$. Since the homogeneous-limit map restricts to the
identity on $D$, its principal fibre has dimension one and vertex
local ring $\mathcal O_{H^0,\eta_D}$. This is a quotient of the full
centered attracting-fibre local ring, so
\[
 1=\dim\mathcal O_{H^0,\eta_D}
 \leq\operatorname{edim}\mathcal O_{H^0,\eta_D}\leq1.
\]
Hence the principal component is regular at $\eta_D$.
\end{proof}

\begin{theorem}
\label{next:fifteen-regularity}
The reduced principal component $H_{15,11}$ is regular at the
generic point of the Erman--Velasco divisor $D_{15}$.  The homogeneous-limit map is birational:
\[
 \mathbb P(\operatorname{Std}_{15}\otimes k^{11})/S_{15}
 \dashrightarrow D_{15}^0.
\]
\end{theorem}

\begin{proof}
Write $A=k\cdot1\oplus V\oplus W$, where $\dim V=11$,
$\dim W=3$, and the internal degrees are $0,1,2$.
For the three rows of the matrix
\[
 \begin{pmatrix}
 -2&5&-4&1&1&0&-3&-2&0&-4&-1\\
 -3&-1&1&2&0&-3&5&4&5&5&-4\\
 -3&-1&2&1&-1&-2&-3&-1&-1&1&-3
 \end{pmatrix},
\]
put $Q_a=\operatorname{diag}(a)-aa^\mathsf T$ and define
\[
 v_i v_j=\sum_{a=1}^3(Q_a)_{ij}w_a,
 \qquad VW=W^2=0.
\]
The three forms are independent: their columns indexed by
$(1,1),(1,2),(1,3)$ have determinant $414$.
The algebra is the homogeneous
limit of $0,e_1,\ldots,e_{11}$ and the three displayed points, by
Erman--Velasco's explicit formula
\cite[Proposition~3.2]{gg-ErmanVelasco}; it therefore belongs to
the smoothable homogeneous locus.

The negative tangent is computed by linearizing associativity. For $s>0$, let $C_s$ consist of symmetric maps
$\phi:A\times A\to A$ of internal degree $-s$ with
$\phi(1,a)=0$. The defining linear map is
\[
 (d\phi)(a,b,c)=
 \mu(\phi(a,b),c)+\phi(\mu(a,b),c)
 -\mu(a,\phi(b,c))-\phi(a,\mu(b,c)).
\]
The allowed changes of quotient basis fix $1,V$; in degrees
$-1,-2$ they are respectively
$\operatorname{Hom}(W,V)$ and $\operatorname{Hom}(W,k)$.
Their coboundary maps are injective because
$\operatorname{Sym}^2V\to W$ is surjective. The framed algebra chart
is a torsor over the Hilbert chart under these basis changes; hence
\[
 \dim T_{-s}=\dim C_s-\operatorname{rank}(d|_{C_s})
                 -\dim(\text{basis changes}).
\]

All entries of $d$ are integers. Elimination over
$\mathbb F_{32003}$ gives the following rank lower bounds, each
certified by a nonzero minor. The accompanying certificate contains
the row labels, column labels and coefficient matrices.
\[
\begin{array}{c|r|r|r|r}
 s&\dim C_s&\operatorname{rank}(d)\text{ at least}
   &\text{basis-change dimension}&
     \dim T_{-s}\text{ at most}\\ \hline
 1&825&780&33&12\\
 2&447&444&3&0\\
 3&99&99&0&0\\
 4&6&6&0&0
\end{array}
\]
The nonzero modular minors give the same bounds over $\mathbb Q$.
There are no lower-degree cochains. Every nonunit triple with first
index less than the third, and every output of the required degree,
is included; commutativity makes the reversed triples redundant.
Thus the table includes all negative tangent directions, including
constant terms in products.

Erman--Velasco's dimension calculation makes the centered
smoothable homogeneous locus a divisor. Applying
Lemma~\ref{next:attracting-criterion} proves regularity.
The radial residue-field formula then gives $h=1$ and the
birational assertion.
\end{proof}

\begin{corollary}
\label{next:fifteen-exact-cone}
At the generic point of the ordered cover above $D_{15}$, let
$B$ be the local ring and $F$ its residue field. There is an
$S_{15}$-Galois extension $E/F$ of degree $15!$ and a dilation
parameter $t$ such that
\[
 R=\bigoplus_{r\ge0}t^rF_r,\qquad B=R_{R_+},
 \qquad F_r=\langle\text{products of at most $r$ normalized
 ordered coordinates}\rangle_F.
\]
Consequently
\[
 \begin{gathered}
 b(B)=e(B)=15!=1\,307\,674\,368\,000,\\
 \operatorname{edim}B=155,\qquad
 (B:\overline B)=t^{17}\overline B,\\
 \delta(B)\ge\frac{17\cdot15!}{2}
 =11\,115\,232\,128\,000.
 \end{gathered}
\]
Its associated graded ring is the displayed graded ring, and
its geometric tangent cone is the reduced union of $15!$
distinct lines.
\end{corollary}

\begin{proof}
The principal local ring is $F[t]_{(t)}$ by the preceding theorem.
Adjoining the ordered point coordinates gives the orbit-cone
construction of Theorem~\ref{new:local-cone}.
The degree-one representation is
$\mathbf1\oplus\operatorname{Std}_{15}^{\oplus11}$, of dimension
$1+14\cdot11=155$. For a generic fifteen-point configuration
in eleven-space, the least alternating degree is
$0+11\cdot1+3\cdot2=17$.
The general interpolation and conductor theorem therefore
gives the stated conductor; sign duality gives the defect bound.
\end{proof}

\subsection{A second Grassmannian divisor in five-space}
\label{next:twelve-divisor-section}

Zhao proves smoothability of the $(1,5,6)$ locus
\cite{next:Zhao}. The tangent calculation below gives a smooth
point of the full Hilbert scheme and hence determines the
transverse orbit cone along the generic homogeneous divisor.

\begin{theorem}
\label{next:twelve-smoothability}
Every homogeneous algebra with Hilbert function $(1,5,6)$ over a
field of characteristic zero is geometrically smoothable. The
corresponding Grassmannian $\operatorname{Gr}(6,15)$ contains a
nonempty open set of smooth points of the full Hilbert scheme
$\operatorname{Hilb}^{12}(\mathbb A^5)$. Its translates form a
divisor $D_{12}$ in the reduced principal component $H_{12,5}$.
The principal component is regular at the generic point of this
divisor, and
\[
 \mathbb P(\operatorname{Std}_{12}\otimes k^5)/S_{12}
 \dashrightarrow D_{12}^0=\operatorname{Gr}(6,15)
\]
is birational.
\end{theorem}

\begin{proof}
The following six rows, together with $0,e_1,\ldots,e_5$, give
an explicit twelve-point source:
\[
 \begin{pmatrix}
 -2&-1&-2&-1&-1\\
 -2&-3&-3&-1&-2\\
 1&-4&3&-2&5\\
 -3&0&-2&4&5\\
 4&1&5&4&0\\
 2&3&5&5&5
 \end{pmatrix}.
\]
For each row $a$, put $Q_a=\operatorname{diag}(a)-aa^{\mathsf T}$.
These six symmetric forms are independent: their columns
indexed by $(1,1),(1,2),(1,3),(1,4),(1,5),(2,2)$ have
determinant $-1\,038\,336$. On
$A=k\cdot1\oplus V\oplus W$, with dimensions $1,5,6$, take
\[
 v_i v_j=\sum_{a=1}^{6}(Q_a)_{ij}w_a,
 \qquad VW=W^2=0.
\]
This is the homogeneous limit of the displayed reduced
configuration, by the explicit formula in
\cite[Proposition~3.2]{gg-ErmanVelasco}.

Use the normalized symmetric cochains and the linearized
associativity map from the proof of
Theorem~\ref{next:fifteen-regularity}. The negative-weight
basis changes fix $1,V$ and are $\operatorname{Hom}(W,V)$ in
degree $-1$ and $\operatorname{Hom}(W,k)$ in degree $-2$.
Their coboundaries are injective. The following table records
nonzero minors over $\mathbb F_{32003}$ of the integer
coefficient matrices:
\[
\begin{array}{c|r|r|r|r|r}
 s&\text{matrix size}&\dim C_s&\operatorname{rank}(d)\text{ at least}
 &\text{basis changes}&\dim T_{-s}\text{ at most}\\ \hline
 1&300\times255&255&219&30&6\\
 2&1510\times291&291&285&6&0\\
 3&1407\times135&135&135&0&0\\
 4&131\times21&21&21&0&0
\end{array}
\]
The certificate records the selected minors over $\mathbb F_{32003}$;
the matrix ranks also agree over $\mathbb F_{65521}$.
Nonvanishing at the first prime proves the displayed
characteristic-zero bounds.

The degree-zero Hilbert tangent is the tangent of
$\operatorname{Gr}(6,15)$, hence has dimension $6(15-6)=54$.
Indeed, a degree-zero map $I\to A$ is an arbitrary map
$I_2\to A_2$, since $A_j=0$ for $j\ge3$.
Every positive-degree map $I\to A$ vanishes. Since the point lies
on the principal component,
\[
 60=\dim H_{12,5}\leq\dim T_{[I]}\operatorname{Hilb}^{12}(\mathbb A^5)
 \leq54+6=60.
\]
Thus it is a smooth point of the full Hilbert scheme.

Near this point the principal component is the unique
irreducible component. Its intersection with the irreducible
Grassmannian therefore contains a nonempty open subset.
Because the principal component is closed, it contains the
entire Grassmannian. This proves smoothability of every
homogeneous quotient in the Grassmannian, and the smooth locus
of the full Hilbert scheme meets it in a nonempty open set.
Translations add five dimensions, so $\dim D_{12}=54+5=59$,
one less than $\dim H_{12,5}=60$.

The full Hilbert scheme, and hence its principal component,
is regular at the generic point of $D_{12}$. The radial
residue-field formula gives degree one for the displayed
homogeneous-limit map and proves birationality.
\end{proof}

\begin{corollary}
\label{next:twelve-exact-cone}
At the generic point of the ordered cover above the centered
divisor $D_{12}^0$, let $B$ be the local ring, with residue
field $F$. There is an $S_{12}$-Galois extension $E/F$ of
degree $12!$ and normalized ordered coordinate functions
$v_{ij}\in E$ such that, for
\[
 F_r=\langle\text{products of at most $r$ of the $v_{ij}$}\rangle_F,
 \qquad R=\bigoplus_{r\ge0}t^rF_r,
\]
one has $B=R_{R_+}$. In particular,
\[
 \begin{gathered}
 b(B)=e(B)=12!=479\,001\,600,\qquad
 \operatorname{edim}B=56,\\
 (B:\overline B)=t^{17}\overline B,\qquad
 \delta(B)\ge\frac{17\cdot12!}{2}=4\,071\,513\,600.
 \end{gathered}
\]
The associated graded ring of $B$ is $R$, and its geometric
tangent cone is the reduced union of $12!$ distinct lines.
\end{corollary}

\begin{proof}
Generic regularity gives the orbit-cone construction
of Theorem~\ref{new:local-cone}. Its degree-one representation is
$\mathbf1\oplus\operatorname{Std}_{12}^{\oplus5}$, so its
embedding dimension is $1+11\cdot5=56$. The least alternating
degree of a generic twelve-point configuration in five-space
is $0+5\cdot1+6\cdot2=17$. The general interpolation and
conductor theorem gives conductor exponent $17$, and sign
duality gives the stated lower bound on $\delta$.
\end{proof}

\Needspace{10\baselineskip}
\section{Gorenstein duality and its failure}
\label{new:gorenstein}

The sign-trace pairing pairs conjugate Specht components in
complementary degrees. For a Gorenstein orbit cone these pairings
must be perfect. At the ten-point divisor, two multiplicities
already contradict that requirement.

\subsection{The sign-trace pairing}

Retain the notation of Theorem~\ref{new:orbit-conductor}, and write
\(b=n!\). Choose a nonzero minimum-degree alternant
\(\Delta\in E\). Such an element exists because \(F_\mu=E\), while
every lower filtration piece has zero sign component. It is a unit,
since its values on the ordered
configuration differ only by sign. Define
\[
 \lambda:E\longrightarrow F,\qquad
 \lambda(a)=\operatorname{Tr}_{E/F}(a/\Delta),
 \qquad
 \langle a,b\rangle=\lambda(ab).
\]
The trace form of a finite \'{e}tale algebra is nondegenerate, so this
pairing is nondegenerate. Its equivariance is
\begin{equation}
 \label{new:sign-pairing-equivariance}
 \langle ga,gb\rangle=\operatorname{sgn}(g)\langle a,b\rangle.
\end{equation}

\begin{proposition}[Sign duality]
\label{new:sign-duality}
For every \(0\le r\le\mu-1\),
\begin{equation}
 \label{new:filtration-orthogonality}
 F_{\mu-1-r}\subset F_r^\perp,\qquad
 \dim F_r+\dim F_{\mu-1-r}\le b.
\end{equation}
The defect satisfies
\begin{equation}
 \label{new:orthogonal-defect-identity}
 2\delta-\mu b
 =\sum_{r=0}^{\mu-1}
   \dim_F\bigl(F_r^\perp/F_{\mu-1-r}\bigr).
\end{equation}
In particular \(\delta\ge\mu b/2\), with equality if and only if
the local cone ring is Gorenstein.
\end{proposition}

\begin{proof}
The filtration and sign equivariance give
\[
 F_rF_{\mu-1-r}\subset F_{\mu-1},\qquad
 \lambda(F_{\mu-1})=0
 \quad\Longrightarrow\quad F_{\mu-1-r}\subset F_r^\perp.
\]
Since \(\dim F_r^\perp=b-\dim F_r\),
\[
 \sum_{r=0}^{\mu-1}\dim(F_r^\perp/F_{\mu-1-r})
 =\mu b-2\sum_{r=0}^{\mu-1}\dim F_r
 =2\delta-\mu b,
\]
where the last equality is~\eqref{new:delta-filtration}.

Suppose equality \(\delta=\mu b/2\) holds. Every quotient on the
right of~\eqref{new:orthogonal-defect-identity} vanishes, so
\[
 F_r^\perp=F_{\mu-1-r}
\]
for all \(r\). The functional \(\lambda\) therefore induces
pairings
\[
 (F_r/F_{r-1})\otimes
 (F_{\mu-r}/F_{\mu-r-1})\longrightarrow F.
\]
The preimage in \(F_r\) of each left kernel is
\(F_r\cap F_{\mu-r}^\perp=F_{r-1}\), so the quotient
pairings are perfect. Thus
\[
 A=R/(t)\text{ is Gorenstein with socle in degree }\mu
 \quad\Longrightarrow\quad R_{R_+}\text{ is Gorenstein},
\]
because \(t\) is a nonzerodivisor.

Conversely, if \(R_{R_+}\) is Gorenstein, then \(A\) is a graded
Artinian Gorenstein algebra. Its top nonzero degree is \(\mu\),
so its Hilbert numerator is palindromic:
\(h_r=h_{\mu-r}\). Equation~\eqref{new:delta-h-vector} gives
\[
 \delta=\sum_r r h_r
 =\frac{\mu}{2}\sum_rh_r
 =\frac{\mu b}{2}.
\]
\end{proof}

Let \(S^\alpha\) denote the Specht representation indexed by a partition
\(\alpha\vdash n\), and let \(f^\alpha=\dim S^\alpha\).
Write \(a_\alpha(r)\) for its multiplicity in \(F_r\). Since \(E\)
is the regular representation, its multiplicity in \(E\) is
\(f^\alpha\). The sign pairing identifies the dual of the
\(\alpha\)-isotypic component with the \(\alpha'\)-isotypic
component, where \(\alpha'\) is the conjugate partition. Consequently
~\eqref{new:filtration-orthogonality} gives
\begin{equation}
 \label{new:paired-multiplicities}
 a_\alpha(r)+a_{\alpha'}(\mu-1-r)\le f^\alpha.
\end{equation}
If the cone is Gorenstein, equality holds in every one of these
inequalities.

\subsection{Multiplicity at the first possible degree}

Let \(m_\alpha(r)\) be the multiplicity of \(S^\alpha\) in the
polynomials of total degree at most \(r\) on
\(\operatorname{Std}_n\otimes F^d\). Coordinate evaluation is
\(S_n\)-equivariant and has image \(F_r\), so
\begin{equation}
 \label{new:source-capacity}
 a_\alpha(r)\le m_\alpha(r),\qquad
 a_\alpha(r)\le f^\alpha.
\end{equation}
The first degree in which a Specht module can occur determines
the multiplicities needed for both conclusions.

\begin{lemma}
\label{new:first-degree}
Let \(s=n-\alpha_1\), and let
\(\overline\alpha=(\alpha_2,\alpha_3,\ldots)\).
Then \(m_\alpha(r)=0\) for \(r<s\), and
\begin{equation}
 \label{new:first-degree-formula}
 m_\alpha(s)=\dim\mathbf S_{\overline\alpha}(F^d).
\end{equation}
The right side is allowed to vanish.
\end{lemma}

\begin{proof}
Begin with the uncentered polynomial ring in \(n\) blocks of
\(d\) variables. A degree-\(r\) monomial occupies at most \(r\)
row indices. The subspace occupying exactly \(j\) rows is induced
from a representation of \(S_j\), with the other \(n-j\) indices
acting trivially:
\[
 \operatorname{Ind}_{S_{n-j}\times S_j}^{S_n}
 ({\bf1}\boxtimes M).
\]
By Pieri's rule, every constituent has first row at least \(n-j\).
Indeed, the horizontal strip of length \(n-j\) contributes at most
one box in each column. Thus \(S^\alpha\) cannot occur when
\(r<n-\alpha_1\).

At \(r=s=n-\alpha_1\), the only possible contribution has exactly
\(s\) occupied rows, and every occupied row has degree one. This
multilinear subspace is
\[
 \operatorname{Ind}_{S_{n-s}\times S_s}^{S_n}
 \bigl({\bf1}\boxtimes(F^d)^{\otimes s}\bigr).
\]
Schur--Weyl duality decomposes its second factor as
\[
 (F^d)^{\otimes s}
 \simeq\bigoplus_{\beta\vdash s}
       \mathbf S_\beta(F^d)\otimes S^\beta.
\]
Pieri's rule requires \(\alpha/\beta\) to be a horizontal strip
of length \(n-s=\alpha_1\). There are exactly \(\alpha_1\)
columns, so one box must be removed from every column. This forces
\(\beta=\overline\alpha\), proving~\eqref{new:first-degree-formula}
for the uncentered ring.

The uncentered permutation representation is
\(\operatorname{Std}_n\oplus{\bf1}\). Its polynomial ring is therefore
the tensor product of the centered polynomial ring and a polynomial
ring in \(d\) invariant variables. Since there are no lower-degree
copies of \(S^\alpha\), those invariant variables do not change its
multiplicity in degree \(s\). The same formula holds after centering.
\end{proof}

\begin{theorem}[Failure of Gorensteinness]
\label{new:not-gorenstein}
The transverse ring \(B\) at the generic ten-point divisor is not
Gorenstein. The same failure occurs along the divisors obtained by
the smooth graph and separated-support comparisons for every
\(n\ge10\) and \(d\ge5\).
\end{theorem}

\begin{proof}
Take
\[
 \alpha=(5,5),\qquad
 \alpha'=(2,2,2,2,2).
\]
Both Specht representations have dimension
\[
 f^\alpha=f^{\alpha'}=\frac1{6}\binom{10}{5}=42.
\]
For \(\alpha\), the first possible degree is five, so
\(a_\alpha(4)=0\). For \(\alpha'\), it is eight, and
Lemma~\ref{new:first-degree} gives
\[
 a_{\alpha'}(8)
 \le m_{\alpha'}(8)
 =\dim\mathbf S_{(2,2,2,2)}(F^5)=15.
\]
The last Schur representation is
\(\det^2\otimes\operatorname{Sym}^2((F^5)^*)\).

The conductor exponent is thirteen, so Gorensteinness would force
\[
 42=a_\alpha(4)+a_{\alpha'}(8)\leq0+15<42
\]
by~\eqref{new:paired-multiplicities}. Hence \(B\) is not Gorenstein.
Gorensteinness ascends and descends under the smooth local products
of Theorem~\ref{fh:j:main}, proving the assertion for every
\(n\geq10\), \(d\geq5\).
\end{proof}

\begin{figure}[!htbp]
\centering
\begin{tikzpicture}[x=1cm,y=1cm,font=\small]
\path[use as bounding box] (0,-1.5) rectangle (16.2,10.4);
\node[anchor=west,font=\bfseries\small,Ink] at (0,10.0) {(a) The minimum alternating degree is thirteen};

\draw[Ink!30,fill=Blue!7] (0.3,8.85) rectangle ++(1.42,.70);\node at (1.01,9.2) {$ 1 $};\node[font=\footnotesize,Ink!70] at (1.01,8.57) {$\deg=0$};
\draw[Ink!30,fill=Teal!12] (1.8,8.85) rectangle ++(1.42,.70);\node at (2.51,9.2) {$ a $};\node[font=\footnotesize,Ink!70] at (2.51,8.57) {$\deg=1$};
\draw[Ink!30,fill=Teal!12] (3.3,8.85) rectangle ++(1.42,.70);\node at (4.01,9.2) {$ b $};\node[font=\footnotesize,Ink!70] at (4.01,8.57) {$\deg=1$};
\draw[Ink!30,fill=Teal!12] (4.8,8.85) rectangle ++(1.42,.70);\node at (5.51,9.2) {$ c $};\node[font=\footnotesize,Ink!70] at (5.51,8.57) {$\deg=1$};
\draw[Ink!30,fill=Teal!12] (6.3,8.85) rectangle ++(1.42,.70);\node at (7.01,9.2) {$ d $};\node[font=\footnotesize,Ink!70] at (7.01,8.57) {$\deg=1$};
\draw[Ink!30,fill=Teal!12] (7.8,8.85) rectangle ++(1.42,.70);\node at (8.51,9.2) {$ e $};\node[font=\footnotesize,Ink!70] at (8.51,8.57) {$\deg=1$};
\draw[Ink!30,fill=Amber!16] (9.3,8.85) rectangle ++(1.42,.70);\node at (10.010000000000002,9.2) {$ a^2 $};\node[font=\footnotesize,Ink!70] at (10.010000000000002,8.57) {$\deg=2$};
\draw[Ink!30,fill=Amber!16] (10.8,8.85) rectangle ++(1.42,.70);\node at (11.510000000000002,9.2) {$ ab $};\node[font=\footnotesize,Ink!70] at (11.510000000000002,8.57) {$\deg=2$};
\draw[Ink!30,fill=Amber!16] (12.3,8.85) rectangle ++(1.42,.70);\node at (13.010000000000002,9.2) {$ ac $};\node[font=\footnotesize,Ink!70] at (13.010000000000002,8.57) {$\deg=2$};
\draw[Ink!30,fill=Amber!16] (13.8,8.85) rectangle ++(1.42,.70);\node at (14.510000000000002,9.2) {$ ad $};\node[font=\footnotesize,Ink!70] at (14.510000000000002,8.57) {$\deg=2$};
\node[anchor=west] at (.3,7.96) {$\deg\det(m_j(v_i))=0+5\cdot1+4\cdot2=13$};
\node[anchor=east,Teal] at (15.2,7.96) {$F_{12}^{\sgn}=0$,\qquad $F_{13}=E$};
\draw[Ink!25] (0,7.5)--(16,7.5);
\node[anchor=west,font=\bfseries\small,Ink] at (0,7.07) {(b) The conjugate partitions expose the failure of Gorenstein duality};

\draw[Ink!65,fill=Teal!10,line width=.35pt] (0.450,6.450) rectangle ++(0.43,-0.43);
\draw[Ink!65,fill=Teal!10,line width=.35pt] (0.880,6.450) rectangle ++(0.43,-0.43);
\draw[Ink!65,fill=Teal!10,line width=.35pt] (1.310,6.450) rectangle ++(0.43,-0.43);
\draw[Ink!65,fill=Teal!10,line width=.35pt] (1.740,6.450) rectangle ++(0.43,-0.43);
\draw[Ink!65,fill=Teal!10,line width=.35pt] (2.170,6.450) rectangle ++(0.43,-0.43);
\draw[Ink!65,fill=Teal!10,line width=.35pt] (0.450,6.020) rectangle ++(0.43,-0.43);
\draw[Ink!65,fill=Teal!10,line width=.35pt] (0.880,6.020) rectangle ++(0.43,-0.43);
\draw[Ink!65,fill=Teal!10,line width=.35pt] (1.310,6.020) rectangle ++(0.43,-0.43);
\draw[Ink!65,fill=Teal!10,line width=.35pt] (1.740,6.020) rectangle ++(0.43,-0.43);
\draw[Ink!65,fill=Teal!10,line width=.35pt] (2.170,6.020) rectangle ++(0.43,-0.43);
\draw[Ink!65,fill=Teal!10,line width=.35pt] (5.150,6.450) rectangle ++(0.43,-0.43);
\draw[Ink!65,fill=Teal!10,line width=.35pt] (5.580,6.450) rectangle ++(0.43,-0.43);
\draw[Ink!65,fill=Teal!10,line width=.35pt] (5.150,6.020) rectangle ++(0.43,-0.43);
\draw[Ink!65,fill=Teal!10,line width=.35pt] (5.580,6.020) rectangle ++(0.43,-0.43);
\draw[Ink!65,fill=Teal!10,line width=.35pt] (5.150,5.590) rectangle ++(0.43,-0.43);
\draw[Ink!65,fill=Teal!10,line width=.35pt] (5.580,5.590) rectangle ++(0.43,-0.43);
\draw[Ink!65,fill=Teal!10,line width=.35pt] (5.150,5.160) rectangle ++(0.43,-0.43);
\draw[Ink!65,fill=Teal!10,line width=.35pt] (5.580,5.160) rectangle ++(0.43,-0.43);
\draw[Ink!65,fill=Teal!10,line width=.35pt] (5.150,4.730) rectangle ++(0.43,-0.43);
\draw[Ink!65,fill=Teal!10,line width=.35pt] (5.580,4.730) rectangle ++(0.43,-0.43);
\node[anchor=west] at (.45,5.18) {$\lambda=(5,5)$};
\node[anchor=west] at (.45,4.65) {$\dim S^\lambda=42$};
\node[anchor=west,Amber] at (.45,4.12) {$a_\lambda(4)=0$};
\node[anchor=west] at (6.35,6.18) {$\lambda'=(2^5)$};
\node[anchor=west] at (6.35,5.65) {$\dim S^{\lambda'}=42$};
\node[anchor=west,Teal] at (6.35,5.12) {$a_{\lambda'}(8)\le15$};
\node[anchor=west,font=\footnotesize,align=left] at (6.35,4.40)
 {$15=\dim\mathbb S_{(2^4)}(\C^5)$\\$\phantom{15}=\dim\Sym^2(\C^5)$};
\draw[-{Latex[length=2mm]},Ink!55] (3.05,6.0)--(4.65,6.0);
\node[above,font=\footnotesize,Ink!70] at (3.85,6.0) {conjugate};
\draw[Ink!30,fill=Amber!5,rounded corners=3pt] (10.7,3.9) rectangle (16,6.5);
\node[align=center,text width=4.9cm] at (13.35,5.85)
 {Gorenstein duality would give};
\node[align=center] at (13.35,5.15)
 {$a_\lambda(4)+a_{\lambda'}(8)=42$};
\node[Amber,font=\large] at (13.35,4.45) {$0+15<42$};
\draw[Ink!25] (0,3.42)--(16,3.42);
\node[anchor=west,font=\bfseries\small,Ink] at (0,2.98) {(c) An orthogonal quotient in the middle degree};

\draw[Ink!65,fill=Teal!10,line width=.35pt] (0.400,2.400) rectangle ++(0.39,-0.39);
\draw[Ink!65,fill=Teal!10,line width=.35pt] (0.790,2.400) rectangle ++(0.39,-0.39);
\draw[Ink!65,fill=Teal!10,line width=.35pt] (1.180,2.400) rectangle ++(0.39,-0.39);
\draw[Ink!65,fill=Teal!10,line width=.35pt] (1.570,2.400) rectangle ++(0.39,-0.39);
\draw[Ink!65,fill=Teal!10,line width=.35pt] (0.400,2.010) rectangle ++(0.39,-0.39);
\draw[Ink!65,fill=Teal!10,line width=.35pt] (0.790,2.010) rectangle ++(0.39,-0.39);
\draw[Ink!65,fill=Teal!10,line width=.35pt] (1.180,2.010) rectangle ++(0.39,-0.39);
\draw[Ink!65,fill=Teal!10,line width=.35pt] (0.400,1.620) rectangle ++(0.39,-0.39);
\draw[Ink!65,fill=Teal!10,line width=.35pt] (0.790,1.620) rectangle ++(0.39,-0.39);
\draw[Ink!65,fill=Teal!10,line width=.35pt] (0.400,1.230) rectangle ++(0.39,-0.39);
\node[anchor=west] at (2.35,2.12) {$\alpha=(4,3,2,1)=\alpha'$};
\node[anchor=west] at (2.35,1.58) {$\dim S^\alpha=768$};
\node[anchor=west] at (2.35,1.04) {$a_\alpha(6)\le280$};

\draw[white,line width=.9pt,fill=Teal!20] (7.300,1.1) rectangle ++(3.135,.95);\node at (8.868,1.58) {$ \le280 $};
\draw[white,line width=.9pt,fill=Amber!35] (10.435,1.1) rectangle ++(2.329,.95);\node at (11.600,1.58) {$ \ge208 $};
\draw[white,line width=.9pt,fill=Blue!15] (12.765,1.1) rectangle ++(3.135,.95);\node at (14.332,1.58) {$ \le280 $};
\node[anchor=west,font=\footnotesize,Ink!75] at (7.3,2.43) {Multiplicities among the $768$ copies of $S^\alpha$ in $E$};
\node[font=\footnotesize] at (8.87,.74) {$(F_6)_\alpha$};
\node[font=\footnotesize,Amber] at (11.6,.74) {$(F_6^\perp/F_6)_\alpha$};
\node[font=\footnotesize] at (14.3,.74) {paired with $(F_6)_\alpha$};
\node[anchor=west] at (.4,-.35)
 {$\dim(F_6^\perp/F_6)\ge208\cdot768=159\,744$};
\node[anchor=west,font=\footnotesize,Ink!80] at (.4,-.93)
 {The middle segment is a lower bound; a smaller $F_6$ enlarges this quotient.};
\end{tikzpicture}
\caption{Failure of Gorenstein duality in the ten-point orbit cone.
The monomials in (a) give an alternant of minimum degree $13$, and
interpolation reaches every ordering in that degree. The sign-trace
pairing therefore gives $F_r\perp F_{12-r}$.
Writing $a_\lambda(r)$ for the multiplicity of $S^\lambda$ in $F_r$,
Gorenstein duality would require
$a_{(5,5)}(4)+a_{(2^5)}(8)=42$; panel (b) gives the contradictory
upper bound $15$. Panel (c) quantifies the defect on a self-conjugate
isotypic component. The bounds follow from the hook-length and
hook-content formulas in Section~\ref{new:gorenstein}.}
\label{new:gorenstein-figure}
\end{figure}
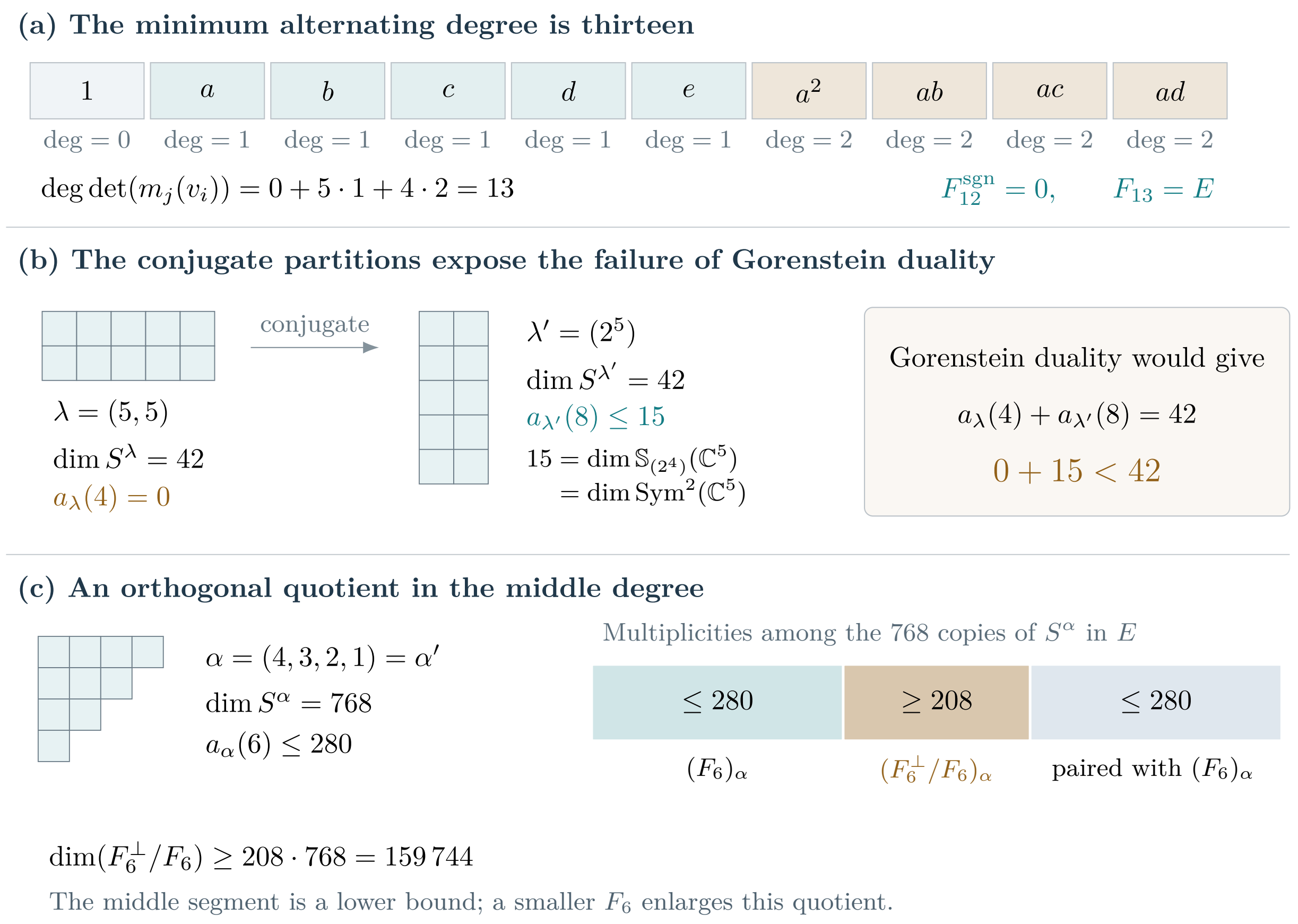

Figure~\ref{new:gorenstein-figure} displays the conjugate diagrams and
the complementary filtration degrees. The numbers \(42\) and \(15\)
are multiplicities: the regular representation contains \(42\)
copies of the indicated Specht module, while the polynomial
filtration can provide at most \(15\) copies at the required degree.

\FloatBarrier

\begin{proposition}[A middle orthogonal quotient]
\label{new:middle-defect}
For the ten-point cone, \(F_6\subset F_6^\perp\), and
\(F_6^\perp/F_6\) contains at least \(208\) copies of
\(S^{(4,3,2,1)}\). In particular
\[
 \dim_F(F_6^\perp/F_6)\ge159\,744,\qquad
 \delta(B)\ge23\,667\,072.
\]
\end{proposition}

\begin{proof}
The partition \(\alpha=(4,3,2,1)\) is self-conjugate. The hook-length
and hook-content formulas give
\[
 f^\alpha=\frac{10!}{4725}=768,\qquad
 m_\alpha(6)=\dim\mathbf S_{(3,2,1)}(F^5)
 =\frac{5\cdot6\cdot7\cdot4\cdot5\cdot3}{5\cdot3\cdot1\cdot3\cdot1\cdot1}
 =280.
\]
The sign pairing is nondegenerate on this self-conjugate isotypic
component. If its multiplicity in \(F_6\) is \(a\), the multiplicity
in \(F_6^\perp/F_6\) is \(768-2a\). Since \(a\le280\), this is
at least \(208\), giving the dimension bound
\[
 208\cdot768=159\,744.
\]
The single middle summand in
~\eqref{new:orthogonal-defect-identity} then yields
\[
 \delta(B)\ge\frac{13\cdot10!}{2}
             +\frac{159\,744}{2}
 =23\,667\,072.
\]
\end{proof}

\subsection{The quantitative defect interval}

For arbitrary filtration degrees, the multiplicity bounds can also
be extracted from the polynomial characters. If a permutation has
cycle lengths \(\ell_1,\ldots,\ell_s\), its cumulative polynomial
character series on \(\operatorname{Std}_n\otimes F^d\) is
\begin{equation}
 \label{new:cumulative-character}
 \sum_{r\ge0}\operatorname{Tr}
 \bigl(g\mid\operatorname{Sym}^{\le r}
       (\operatorname{Std}_n\otimes F^d)\bigr)q^r
 =\frac{(1-q)^{d-1}}{\prod_{i=1}^s(1-q^{\ell_i})^d}.
\end{equation}
Indeed the determinant on the standard representation is
\(\prod_i(1-q^{\ell_i})/(1-q)\), and cumulative degrees contribute
one additional factor \(1/(1-q)\). If \(\rho\vdash n\) has
\(m_j(\rho)\) parts of size \(j\), put
\(z_\rho=\prod_jj^{m_j(\rho)}m_j(\rho)!\). Character orthogonality
then gives
\begin{equation}
 \label{new:character-multiplicity}
 m_\alpha(r)=
 \sum_{\rho\vdash n}\frac{\chi^\alpha(\rho)}{z_\rho}
 [q^r]\frac{(1-q)^{d-1}}
 {\prod_j(1-q^j)^{d\,m_j(\rho)}}.
\end{equation}
Combining~\eqref{new:source-capacity} and
~\eqref{new:paired-multiplicities} gives the general bound
\begin{equation}
 \label{new:character-capacity}
 \dim F_r+\dim F_{\mu-1-r}
 \le
 \sum_{\alpha\vdash n} f^\alpha
 \min\{f^\alpha,\,
        m_\alpha(r)+m_{\alpha'}(\mu-1-r)\}.
\end{equation}

For the numerical bound needed here, every contributing multiplicity
can be calculated by Lemma~\ref{new:first-degree}. This avoids any
dependence on an extensive character-table computation.

\begin{theorem}[Normalization defect]
\label{new:defect-interval}
At the ten-point divisor,
\[
 \boxed{23\,738\,952\le\delta(B)\le29\,746\,080.}
\]
Consequently
\[
 20\,110\,153
 \le\operatorname{length}_B(B^{\mathrm{sn}}/B)
 \le26\,117\,281.
\]
\end{theorem}

\begin{proof}
For a partition \(\alpha\), the deficit in the paired dimensions
from~\eqref{new:character-capacity} is at least
\[
 f^\alpha
 \max\{0,f^\alpha-m_\alpha(r)-m_{\alpha'}(12-r)\}.
\]
Table~\ref{new:hook-deficit-table} lists eight nonzero contributions.
Each multiplicity is either zero below its first possible degree or
the Schur dimension in~\eqref{new:first-degree-formula}. The numerical
entries follow from
\[
 f^\alpha=\frac{10!}{\prod_{\square\in\alpha}h(\square)},
 \qquad
 \dim\mathbf S_\beta(F^5)
 =\prod_{(i,j)\in\beta}\frac{5+j-i}{h(i,j)}.
\]
The remaining deficits are nonnegative, so these eight
contributions give a lower bound for their sum.

\begin{table}[ht]
\centering
\small
\setlength{\tabcolsep}{4pt}
\begin{tabular}{@{}c c c r r r r@{}}
\toprule
\(r\)&\(\alpha\)&\(\alpha'\)&\(f^\alpha\)&
\(m_\alpha(r)\)&\(m_{\alpha'}(12-r)\)&deficit\\
\midrule
4&\((6,4)\)&\((2,2,2,2,1,1)\)&90&70&10&900\\
4&\((5,5)\)&\((2,2,2,2,2)\)&42&0&15&1,134\\
5&\((5,4,1)\)&\((3,2,2,2,1)\)&288&224&40&6,912\\
5&\((5,3,2)\)&\((3,3,2,1,1)\)&450&175&175&45,000\\
5&\((4,4,2)\)&\((3,3,2,2)\)&252&0&210&10,584\\
6&\((4,3,3)\)&\((3,3,3,1)\)&210&175&0&7,350\\
6&\((3,3,3,1)\)&\((4,3,3)\)&210&0&175&7,350\\
6&\((4,3,2,1)\)&\((4,3,2,1)\)&768&280&280&159,744\\
\bottomrule
\end{tabular}
\caption{Eight hook-formula contributions to the failure of complementary
filtration dimensions.}
\label{new:hook-deficit-table}
\end{table}

The losses are \(2\,034\) for \(r=4\), \(62\,496\) for \(r=5\),
and \(174\,444\) for \(r=6\). Together with the general bound \(10!\)
for the other pairs, they give Table~\ref{new:paired-dimension-table}.

\begin{table}[ht]
\centering
\begin{tabular}{@{}c r@{}}
\toprule
\(r\)&upper bound on \(\dim F_r+\dim F_{12-r}\)\\
\midrule
0&3,628,800\\
1&3,628,800\\
2&3,628,800\\
3&3,628,800\\
4&3,626,766\\
5&3,566,304\\
6&3,454,356\\
\bottomrule
\end{tabular}
\caption{Complementary-degree bounds for the ten-point orbit filtration.
The last row bounds \(2\dim F_6\).}
\label{new:paired-dimension-table}
\end{table}

Sum the first six paired bounds and half the last one. This gives
\[
 \sum_{r=0}^{12}\dim F_r\le23\,435\,448,
\]
and hence
\[
 \delta(B)
 =13\cdot10!-\sum_{r=0}^{12}\dim F_r
 \ge47\,174\,400-23\,435\,448
 =23\,738\,952.
\]
The upper bound is~\eqref{new:ten-delta-upper}. Subtracting
\(10!-1=3\,628\,799\) and using
~\eqref{new:seminormal-defect} gives the seminormalization bounds.
\end{proof}

This proves the defect interval in Theorem~\ref{new:quantitative}.
The conductor follows from interpolation, and the failure of
Gorenstein duality follows from the first occurrences of a few
Specht modules. Computing every intermediate \(F_r\) is unnecessary.

\Needspace{10\baselineskip}
\section{The canonical module and the alternating Gorenstein quotient}
\label{next:canonical}

Let \(F\) have characteristic zero, let \(E/F\) be an \(S_n\)-torsor,
and let
\[
 F_{-1}=0,\qquad F_0=F,\qquad F_rF_s\subset F_{r+s},
 \qquad F_\mu=E
\]
be its polynomial orbit filtration.  Thus \(F_r\) is the image of
polynomials of degree at most \(r\) in the centered coordinates of an
ordered configuration, and
\[
 R=\bigoplus_{r\geq0}t^rF_r,\qquad
 A=R/(t)=\bigoplus_{r=0}^{\mu}F_r/F_{r-1}.
\]
Write \(b=[E:F]=n!\).  Assume that the sign representation first occurs
in degree \(\mu\).  These hypotheses hold for the generic configurations
and the transverse rings considered above.

Choose a nonzero alternating element \(\Delta\in E\), and put
\[
 \lambda(a)=\operatorname{Tr}_{E/F}(a/\Delta).
\]
The pairing \(\lambda(ab)\) on \(E\) is nondegenerate and transforms
by the sign character.  In particular, \(\lambda(F_{\mu-1})=0\).
All orthogonal complements below refer to this pairing.  For convenience
put \(F_j=0\) when \(j<0\) and \(F_j=E\) when \(j\geq\mu\).

\begin{theorem}[Canonical module with its group action]
\label{next:canonical-module}
The graded \(R\)-module
\[
 K=\bigoplus_{r\geq0}t^rF_{\mu-1-r}^{\perp}
\]
satisfies \(R\subset K\subset E[t]\), and
\begin{equation}
\omega_R\simeq K(\mu-1)\otimes\operatorname{sgn}.
 \label{next:dual-eq-1}
\end{equation}
Here \(M(a)_j=M_{a+j}\) is the grading convention.  The finite length
module \(M=K/R\) has
\begin{equation}
M_r=F_{\mu-1-r}^{\perp}/F_r,\qquad
 \dim_F M=2\delta-\mu b,
 \label{next:dual-eq-2}
\end{equation}
where \(\delta=\dim_F(E[t]/R)\).  Multiplication and \(\lambda\) induce
perfect sign-valued pairings
\begin{equation}
M_r\times M_{\mu-1-r}\longrightarrow F.
 \label{next:dual-eq-3}
\end{equation}
\end{theorem}

\begin{proof}
The filtration gives \(F_r\subset F_{\mu-1-r}^{\perp}\), hence
\(R\subset K\). It also gives
\[
 \begin{gathered}
 x\in F_s,\quad a\in F_{\mu-1-r}^{\perp},\quad
 y\in F_{\mu-1-r-s}\\
 \Longrightarrow\quad xy\in F_{\mu-1-r},\quad\lambda(axy)=0.
 \end{gathered}
\]
Thus \(F_sF_{\mu-1-r}^{\perp}\subset
F_{\mu-1-r-s}^{\perp}\), which proves that \(K\) is an \(R\)-module.

Set \(S=F[t]\).  The module \(R\) is finite and torsion free over the
principal ideal domain \(S\), hence free.  Identifying its \(S\)-dual
inside \(E[t,t^{-1}]\) by the pairing \(\lambda\) gives
\[
 \operatorname{Hom}_S(R,S)_s=t^sF_{-s-1}^{\perp}.
\]
To see the degree condition, a functional represented by \(t^sa\)
sends \(t^rF_r\) into \(t^{r+s}F\). Its negative-degree values
vanish precisely when \(a\perp F_{-s-1}\). Multiplication by
\(t^\mu\) then identifies this dual with \(K\), with shift \(\mu\).
The identity \(\omega_S=S(-1)\) changes the shift to \(\mu-1\).
Finally,
\[
 \lambda(ga\,gb)=\operatorname{sgn}(g)\lambda(ab)
\]
accounts for the sign twist in~\eqref{next:dual-eq-1}.

Both filtrations equal \(E\) after degree \(\mu\), so \(M\) has finite
length.  Summing
\[
 \dim M_r=b-\dim F_{\mu-1-r}-\dim F_r
\]
for \(0\leq r<\mu\) gives \eqref{next:dual-eq-2}.  For \eqref{next:dual-eq-3}, a representative in
\(F_{\mu-1-r}^{\perp}\) pairs with one in \(F_r^\perp\).
Changing either representative by its indicated filtration subspace
does not change the value.  The annihilator of \(F_r^\perp\) in the
first space is \(F_r\), which proves perfection.
\end{proof}

Let \(\lambda_{\mathrm{top}}\) be the functional on \(A_\mu\) induced
by \(\lambda\), extended by zero to the other degrees.  Define
\[
 J=\{a\in A:\lambda_{\mathrm{top}}(ab)=0
          \text{ for every }b\in A\},\qquad G=A/J .
\]
Scaling \(\lambda_{\mathrm{top}}\) leaves its radical \(J\), and
hence the quotient \(G\), unchanged.

\begin{theorem}[The canonical alternating quotient]
\label{next:apolar-quotient}
The algebra \(G\) is Artinian Gorenstein, has socle degree \(\mu\),
and has sign socle.  There are natural identifications
\begin{equation}
J\simeq(0:_M t)(-1),\qquad
 \dim_F J=\dim_F(M/tM)\leq\dim_FM.
 \label{next:dual-eq-4}
\end{equation}
Consequently
\begin{equation}
2\delta-\mu n!\ \geq\ n!-\dim_FG.
 \label{next:dual-eq-5}
\end{equation}
\end{theorem}

\begin{proof}
The radical \(J\) is an ideal because
\(\lambda_{\mathrm{top}}((ca)b)=\lambda_{\mathrm{top}}(a(cb))\).
The induced pairing on \(G\) is perfect and has one-dimensional
sign-valued top degree, proving the Gorenstein assertion.

The obstruction to perfection on \(A\) appears when one specializes
\(R\subset K\) at \(t=0\). Since \(R\) and \(K\) are free over
\(S=F[t]\), the specialization sequence is
\begin{equation}
0\longrightarrow(0:_Mt)(-1)\longrightarrow A
 \longrightarrow K/tK\longrightarrow M/tM\longrightarrow0.
 \label{next:dual-eq-6}
\end{equation}
In degree \(r\), the middle map is the pairing
\begin{equation}
F_r/F_{r-1}\longrightarrow
 (F_{\mu-r}/F_{\mu-r-1})^*.
 \label{next:dual-eq-7}
\end{equation}
Its kernel is \(J_r\). Multiplication by \(t\) on the
finite-dimensional space \(M\) has equally dimensional kernel
and cokernel. Therefore
\[
 \dim J=\dim(0:_Mt)=\dim(M/tM)\leq\dim M,
 \qquad \dim J=\dim A-\dim G=n!-\dim G.
\]
Together with~\eqref{next:dual-eq-2}, these identities prove
\eqref{next:dual-eq-4}--\eqref{next:dual-eq-5}.
\end{proof}

\begin{corollary}[Every hook survives]
\label{next:hook-survival}
The canonical defect \(M\) and the apolar kernel \(J\) contain no
hook representations.  The alternating quotient \(G\) retains
every hook isotypic component of the regular representation, with
its entire grading.  In particular
\begin{equation}
\dim_FG\geq\sum_{k=0}^{n-1}\binom{n-1}{k}^{\!2}
          =\binom{2n-2}{n-1}.
 \label{next:dual-eq-8a}
\end{equation}
\end{corollary}

\begin{proof}
The integral exterior-covariant calculation in
Corollary~\ref{next:orbit-hooks} gives
\begin{equation}
\sum_{k=0}^{n-1}u^k
   \operatorname{Hilb}\operatorname{Hom}_{S_n}
        (S^{(n-k,1^k)},A;z)
       =\prod_{i=1}^{n-1}(1+uz^{w_i}),
 \qquad \sum_iw_i=\mu .
 \label{next:dual-eq-8b}
\end{equation}
That calculation is independent of the canonical-module argument:
super multisymmetric power sums generate the invariant exterior
algebra integrally over \(F[t]\), and the point evaluation matrix
proves their independence after inverting \(t\).

The multiplicity of the \(k\)-th hook in \(F_s\) is therefore the
number of \(k\)-subsets of \(\{1,\ldots,n-1\}\) whose weights sum
to at most \(s\).  Taking complementary subsets shows that this
number plus the multiplicity of the conjugate hook in
\(F_{\mu-1-s}\) equals \(\binom{n-1}{k}\).
The formula for \(M_s\) in \eqref{next:dual-eq-2} makes its hook multiplicity zero.
Equation \eqref{next:dual-eq-4} then gives the assertion about \(J\).  Summing the
full hook dimensions proves the bound.
\end{proof}

\subsection{The determinant inverse system at ten points}

At the generic homogeneous algebra with Hilbert function \((1,5,4)\),
let
\[
 Q=\langle q_1,\ldots,q_4\rangle\subset
       \operatorname{Sym}^2((F^5)^*)
\]
be its quadratic Macaulay inverse system.  Write
\(z_i=(z_{i1},\ldots,z_{i5})\), for \(1\leq i\leq10\), and set
\begin{equation}
P_Q(z_1,\ldots,z_{10})=
 \det\begin{bmatrix}
 1&z_{11}&\cdots&z_{15}&q_1(z_1)&\cdots&q_4(z_1)\\
 \vdots&\vdots&&\vdots&\vdots&&\vdots\\
 1&z_{10,1}&\cdots&z_{10,5}&q_1(z_{10})&\cdots&q_4(z_{10})
 \end{bmatrix}.
 \label{next:dual-eq-8}
\end{equation}

\begin{theorem}[Explicit inverse polynomial]
\label{next:quadratic-determinant}
For the generic ten-point cone, \(G\) is the apolar algebra of \(P_Q\).
The polynomial \(P_Q\) has degree \(13\) and \(45\) essential
variables.  Moreover
\begin{equation}
\dim_FG\leq3\,325\,296,\qquad
 \dim_FJ\geq303\,504 .
 \label{next:dual-eq-9}
\end{equation}
The bound on the apolar dimension, and the assertion about essential
variables, hold for every four-plane \(Q\subset\operatorname{Sym}^2
((F^5)^*)\), after extending the field if necessary.
\end{theorem}

\begin{proof}
For an ordered configuration \(p_1,\ldots,p_{10}\), the inverse
polynomial of \(\lambda_{\mathrm{top}}\), up to a nonzero scalar, is
the degree-\(13\) part of
\begin{equation}
\det\bigl(\exp(p_i\cdot z_j)\bigr)_{i,j=1}^{10}.
 \label{next:dual-eq-10}
\end{equation}
This follows either by expanding the determinant as a signed sum,
or by evaluating
\(\lambda((\sum_i p_i\cdot z_i)^{13})/13!\).
In the Cauchy--Binet expansion of \eqref{next:dual-eq-10}, the least possible total
degree is obtained by selecting the constant monomial, the five
linear monomials, and four quadratic monomials.

Let \(U\) be a matrix whose rows span the relations among the
evaluation columns of the constant and linear monomials.  Applying
\(U\) to the quadratic evaluation matrix gives the four-dimensional
quadratic quotient.  Its transpose dual, with the usual factorial
normalization of inverse systems, is precisely \(Q\).  The minors
in Cauchy--Binet are therefore the Pl\"ucker coordinates of \(Q\).
Expanding \eqref{next:dual-eq-8} by its last four columns gives the same minors, proving
the inverse-polynomial assertion.

Common translation of the ten rows changes each quadratic column
only by a combination of linear and constant columns.  Hence \(P_Q\)
depends on at most \(45\) variables.  Any further space of linear
differential annihilators is \(S_{10}\)-stable inside
\(\mathrm{Std}_{10}\otimes F^5\), and so contains a submodule
\(\mathrm{Std}_{10}\otimes Fw\) for some \(w\ne0\).
After changing the five coordinate directions, the polynomial would
then depend on only four alphabets.  A nonzero alternating polynomial
in ten rows and four alphabets has degree at least
\[
 0+4\cdot1+5\cdot2=14,
\]
because the ten exponents in each alternant must be distinct.
This contradicts degree \(13\).  The determinant in \eqref{next:dual-eq-8} is nonzero:
its ten column polynomials are linearly independent.

The dimension bound is proved in Table~\ref{next:ten-capacity} below.
It holds on the dense open set of quadratic planes arising from
generic configurations.  Each graded apolar dimension is the rank
of a catalecticant matrix whose entries are regular functions of the
Pl\"ucker coordinates, after choosing a local frame.  Rank cannot
increase on specialization.  The Grassmannian is irreducible, so
the bound extends to every \(Q\).
\end{proof}

\begin{lemma}[Essential variables of quadratic evaluation determinants]
\label{next:quadratic-essential}
Let \(Q\) be any independent \(r\)-plane of quadrics in \(d\) variables,
put \(n=1+d+r\), and choose a basis \(q_1,\ldots,q_r\).  The polynomial
\begin{equation}
\det[\,1,z_{i1},\ldots,z_{id},q_1(z_i),\ldots,q_r(z_i)\,]_{i=1}^n
 \label{next:dual-eq-8c}
\end{equation}
is nonzero, has degree \(d+2r\), and has \(d(n-1)\)
essential variables.
\end{lemma}

\begin{proof}
The column polynomials are independent, so their evaluation determinant
is nonzero.  Translation invariance removes the \(d\) diagonal
directions.  As in Theorem~\ref{next:quadratic-determinant}, any
further linear annihilator removes a whole alphabet by
\(S_n\)-equivariance.  In \(d-1\) alphabets there are only \(d-1\)
linear monomials.  An alternating polynomial on \(n=1+d+r\) rows
would therefore have degree at least
\[
 (d-1)+2(n-d)=d+2r+1.
\]
If quadratic monomials do not suffice, the minimum only increases.
This excludes the required nonzero polynomial of degree \(d+2r\).
\end{proof}

Lemma~\ref{new:first-degree} gives the two inputs to the dimension
estimate. For \(\alpha\vdash n\), the first possible degree of
\(S^\alpha\) in polynomials on \(\mathrm{Std}_n\otimes F^d\) is
\begin{equation}
a_\alpha=n-\alpha_1.
 \label{next:dual-eq-11}
\end{equation}
At that first possible degree its multiplicity is
\begin{equation}
\dim\mathbf S_{(\alpha_2,\alpha_3,\ldots)}(F^d).
 \label{next:dual-eq-12}
\end{equation}
The first-degree calculation comes from the multilinear subspace
on \(s=a_\alpha\) occupied rows,
\begin{equation}
\operatorname{Ind}_{S_{n-s}\times S_s}^{S_n}
 \bigl(\mathbf1\boxtimes(F^d)^{\otimes s}\bigr).
 \label{next:dual-eq-13}
\end{equation}
Schur--Weyl duality and Pieri's rule leave the tail partition
\((\alpha_2,\alpha_3,\ldots)\). Centering removes invariant variables,
which cannot affect a representation's first occurrence.

Perfect sign duality in \(G\) gives
\begin{equation}
G_k\simeq G_{\mu-k}^*\otimes\operatorname{sgn}.
 \label{next:dual-eq-14}
\end{equation}
Thus \(S^\alpha\) can occur only in the interval
\begin{equation}
a_\alpha\leq k\leq\mu-a_{\alpha'},
 \label{next:dual-eq-15}
\end{equation}
where \(\alpha'\) is the conjugate partition.  At ten points,
the intervals for the shapes in the following table have one or
two degrees, so \eqref{next:dual-eq-12} and \eqref{next:dual-eq-14} bound their total multiplicities.
The ambient ungraded representation \(A\) is regular.

\begin{table}[ht]
\centering
\begin{tabular}{c r r r}
\(\alpha\)&\(f^\alpha\)&multiplicity in \(G\), at most&
lost dimension for this shape\\ \hline
\((6,4)\)&90&\(70+10=80\)&900\\
\((5,5)\)&42&15&1\,134\\
\((5,4,1)\)&288&\(224+40=264\)&6\,912\\
\((5,3,2)\)&450&\(175+175=350\)&45\,000\\
\((4,4,2)\)&252&210&10\,584\\
\((4,3,3)\)&210&175&7\,350\\
\((4,3,2,1)\)&768&\(280+280=560\)&159\,744
\end{tabular}
\caption{Hook-length and hook-content bounds for the alternating
quotient.  The first six rows also contribute their distinct conjugate
partitions; the last row is self-conjugate.}
\label{next:ten-capacity}
\end{table}

The conjugate pairs and the self-conjugate last row give
\[
 \dim J\geq2(900+1\,134+6\,912+45\,000+10\,584+7\,350)
                  +159\,744=303\,504.
\]
Subtracting from \(\dim A=10!\) proves~\eqref{next:dual-eq-9}.

\Needspace{10\baselineskip}
\section{A uniform restriction on the alternating quotient}
\label{next:uniform-apolar}

Consider \(n=d+1+r\) generic points in \(F^d\), where \(r\geq0\)
and quadrics interpolate the points:
\begin{equation}
n\leq\binom{d+2}{2}.
 \label{next:dual-eq-16}
\end{equation}
Their minimal alternating degree is
\(\mu=d+2r=n-1+r\).  For a partition \(\alpha\vdash n\), define
\begin{equation}
q(\alpha)=n+1-\alpha_1-\ell(\alpha).
 \label{next:dual-eq-17}
\end{equation}
This is the number of boxes outside the first row and first column.

\begin{theorem}[Bounded interiors in the alternating quotient]
\label{next:interior-support}
Every irreducible constituent of \(G\) satisfies \(q(\alpha)\leq r\).
Consequently
\begin{equation}
\dim_FG\leq\sum_{\substack{\alpha\vdash n\\q(\alpha)\leq r}}
             (f^\alpha)^2
 \leq
 \sum_{q=0}^r
   \binom nq^2 q!\binom{2n-2q-2}{n-q-1}.
 \label{next:dual-eq-18}
\end{equation}
For fixed \(r\), this is
\begin{equation}
\dim_FG=O_r\!\left(n^{2r-1/2}4^n\right),
 \qquad \frac{\dim_FG}{n!}\longrightarrow0,\qquad
 \lim_{n\to\infty}(\dim_FG)^{1/n}=4.
 \label{next:dual-eq-19}
\end{equation}
In addition the canonical defect satisfies
\begin{equation}
2\delta-\mu n!
 \ \geq\sum_{\alpha\vdash n}
       \max\{q(\alpha)-r,0\}(f^\alpha)^2.
 \label{next:dual-eq-20}
\end{equation}
\end{theorem}

\begin{proof}
The interval~\eqref{next:dual-eq-15} is nonempty precisely when
\[
 \begin{aligned}
 a_\alpha+a_{\alpha'}\leq\mu
 &\quad\Longleftrightarrow\quad
 2n-\alpha_1-\ell(\alpha)\leq n-1+r\\
 &\quad\Longleftrightarrow\quad q(\alpha)\leq r.
 \end{aligned}
\]
As an ungraded quotient of the regular representation, \(G\)
contains at most \(f^\alpha\) copies of \(S^\alpha\). This gives
the first bound in~\eqref{next:dual-eq-18}.

Remove the first row and first column of a standard tableau of
shape \(\alpha\), leaving an interior partition \(\nu\vdash q\).
The remaining hook has \(n-q\) boxes.  Choosing the \(q\) interior
labels and standardizing the two resulting tableaux gives an
injection, hence
\[
 f^\alpha\leq\binom nq\,f^\nu f^{(a,1^c)},
 \qquad a+c=n-q.
\]
For a fixed pair consisting of an interior partition and a hook,
there is at most one original partition.  Summing squares and
discarding their compatibility conditions gives
\[
 \sum_{q(\alpha)=q}(f^\alpha)^2
 \leq\binom nq^2
       \left(\sum_{\nu\vdash q}(f^\nu)^2\right)
       \left(\sum_{c=0}^{n-q-1}\binom{n-q-1}{c}^2\right).
\]
The two sums are \(q!\) and
\(\binom{2n-2q-2}{n-q-1}\), respectively.  For fixed \(r\), the
central binomial estimate and \(\binom nq=O_r(n^q)\) give the upper
bound and the vanishing ratio in \eqref{next:dual-eq-19}.  The matching lower
exponential rate follows from Corollary~\ref{next:hook-survival}.

For the stronger assertion about \(M\), fix a shape with
\(q(\alpha)>r\).  If
\begin{equation}
r+\ell(\alpha)-1\leq s\leq n-\alpha_1-1,
 \label{next:dual-eq-21}
\end{equation}
then \(S^\alpha\) is absent from \(F_s\), and \(S^{\alpha'}\) is
absent from \(F_{\mu-1-s}\), by \eqref{next:dual-eq-11}.  The latter assertion means
that the whole regular \(\alpha\)-isotypic component belongs to
\(F_{\mu-1-s}^{\perp}\).  Therefore \(M_s\) contains this entire
component, of dimension \((f^\alpha)^2\). The number of such degrees is
\[
 (n-\alpha_1-1)-(r+\ell(\alpha)-1)+1=q(\alpha)-r.
\]
Summing over shapes and using~\eqref{next:dual-eq-2} proves
\eqref{next:dual-eq-20}.
\end{proof}

\subsection{Endpoint capacities and stronger finite bounds}

Write
\begin{equation}
m_\alpha=\dim\mathbf S_{(\alpha_2,\alpha_3,\ldots)}(F^d)
          =\prod_{(i,j)\in(\alpha_2,\alpha_3,\ldots)}
                 \frac{d+j-i}{h(i,j)}.
 \label{next:dual-eq-21a}
\end{equation}
The interval \eqref{next:dual-eq-15} has length \(r-q(\alpha)+1\).
Its one-degree and two-degree cases consequently give the following
upper bound on the total multiplicity of \(S^\alpha\) in \(G\):
\begin{equation}
c_\alpha=
 \begin{cases}
 0,&q(\alpha)>r,\\
 \min\{f^\alpha,m_\alpha,m_{\alpha'}\},&q(\alpha)=r,\\
 \min\{f^\alpha,m_\alpha+m_{\alpha'}\},&q(\alpha)=r-1,\\
 f^\alpha,&q(\alpha)\leq r-2.
 \end{cases}
 \label{next:dual-eq-21b}
\end{equation}
In the one-degree case both ends bound the same vector space.
In the two-degree case sign duality identifies the upper endpoint
with the first degree of the conjugate partition, so the two
endpoint bounds add.

Set
\begin{equation}
L_{n,d,q}=
 \sum_{\substack{\alpha\vdash n\\q(\alpha)=q}}
       f^\alpha(f^\alpha-c_\alpha),
 \qquad q\leq r.
 \label{next:dual-eq-21c}
\end{equation}
Since \(J\) is a graded submodule of \(M\) up to a shift, these losses
also give lower bounds on the corresponding isotypic dimensions of
\(M\).  Their shapes have \(q\leq r\), whereas \eqref{next:dual-eq-20} uses \(q>r\);
the contributions are disjoint.  We obtain
\begin{equation}
\begin{split}
 \dim G
 &\leq \sum_{q(\alpha)\leq r}(f^\alpha)^2
                      -\sum_{q\leq r}L_{n,d,q},\\
 2\delta-\mu n!
 &\geq \sum_{\alpha\vdash n}(q(\alpha)-r)_+(f^\alpha)^2
                      +\sum_{q\leq r}L_{n,d,q}.
\end{split}
 \label{next:dual-eq-21d}
\end{equation}
These are explicit hook-length and hook-content formulas.
For the three transverse cones they give the following compact
table.  The third column counts shapes whose summand in \eqref{next:dual-eq-21c}
is positive; conjugate shapes are counted separately.
\begin{equation}
\begin{array}{c|c|r|r}
(n,d)&q&\text{positive shapes}&L_{n,d,q}\\ \hline
(10,5)&3&7&265\,368\\
(10,5)&4&6&38\,136\\
(12,5)&5&7&20\,001\,828\\
(12,5)&6&2&265\,188\\
(15,11)&2&18&959\,228\,076\\
(15,11)&3&18&8\,922\,097\,028
\end{array}
 \label{next:dual-eq-21e}
\end{equation}
No other \(q\leq r\) contributes to this endpoint estimate.

For twelve points, the two shapes \((4,4,4)\) and \((3,3,3,3)\)
already give a short obstruction: both have Specht dimension
\(462\), their first possible degrees are \(8\) and \(9\), and their
first-degree multiplicities are \(490\) and \(175\).  The sole
possible degree for each is paired with the sole degree of the other,
so each retains at most \(175\) copies.  Each loses
\((462-175)\cdot462=132\,594\) dimensions.  The two \(q=5,6\)
rows in \eqref{next:dual-eq-21e} give
\begin{equation}
\dim G_{12}\leq458\,734\,584,\quad
 \dim J_{12}\geq20\,267\,016,\quad
 \delta_{12}\geq4\,081\,647\,108.
 \label{next:dual-eq-21f}
\end{equation}
In particular, the generic twelve-point transverse ring is
non-Gorenstein.  Its inverse polynomial is \eqref{next:dual-eq-8c} for \(d=5,r=6\),
so it has degree \(17\) and \(55\) essential variables.
The generic smoothable locus is dense in the Grassmannian of
six-planes of quadrics in five variables.  Consequently the
apolar-dimension upper bound extends to every independent
six-plane, by the same specialization argument as at ten points.

\begin{corollary}[Finite apolar and normalization bounds]
\label{next:three-apolar-bounds}
For the generic transverse cones of lengths $10,12,15$, the
alternating quotient and normalization defect satisfy
\[
\begin{array}{c|r|r}
 n&\dim G\text{ at most}&\delta\text{ at least}\\\hline
10&3\,325\,296&23\,738\,952\\
12&458\,734\,584&4\,081\,647\,108\\
15&108\,621\,259\,564&12\,356\,860\,349\,568
\end{array}
\]
At length fifteen, $\dim(K/R)\ge2\,483\,256\,443\,136$.
\end{corollary}

\begin{proof}
Apply \eqref{next:dual-eq-21d} using the endpoint sums in
\eqref{next:dual-eq-21e}. For length fifteen, the interior
sums are displayed in \eqref{next:dual-eq-22} below.
The two contributions concern disjoint isotypic components,
so their sum bounds the full canonical defect. Finally use
$2\delta-\mu n!=\dim(K/R)$ with $\mu=13,17,17$.
\end{proof}

\subsection{The fifteen-point transverse singularity}

The homogeneous \((1,11,3)\) divisor has \(n=15\), \(d=11\),
\(r=3\), and \(\mu=17\).  The negative-tangent calculation and its
attracting-chart argument establish generic regularity of the
principal component there; hence the orbit-cone description applies
with \(b=15!=1\,307\,674\,368\,000\).

The following table is obtained solely from the hook-length formula
\(f^\alpha=15!/\prod_{u\in\alpha}h(u)\).  It groups all partitions
by their number of interior boxes.
\begin{equation}
\begin{array}{c|r}
q&\displaystyle\sum_{q(\alpha)=q}(f^\alpha)^2\\ \hline
0&40\,116\,600\\
1&1\,291\,754\,518\\
2&16\,236\,508\,646\\
3&100\,934\,204\,904\\
4&323\,101\,421\,032\\
5&500\,858\,379\,450\\
6&312\,940\,289\,948\\
7&51\,622\,396\,254\\
8&649\,296\,648
\end{array}
\label{next:dual-eq-22}
\end{equation}
The entries sum to \(15!\).  The support estimate alone gives
\(\dim G\leq118\,502\,584\,668\).  Subtracting the endpoint loss
\[
 959\,228\,076+8\,922\,097\,028=9\,881\,325\,104
\]
gives the stronger bound.  Adding the same loss to the weighted
defect bound in \eqref{next:dual-eq-20} yields
\begin{equation}
\begin{split}
 \dim_FG&\leq108\,621\,259\,564,\\
 \dim_FJ&\geq1\,199\,053\,108\,436,\\
 \dim_FM&\geq2\,483\,256\,443\,136,\\
 \delta&\geq12\,356\,860\,349\,568.
\end{split}
\label{next:dual-eq-23}
\end{equation}
Thus the canonical alternating quotient loses more than \(91.69\%\)
of the regular representation.  The stronger bound on \(M\) records
the degrees through which the missing constituents persist.

The inverse polynomial also has a concrete form.  If
\(Q=\langle q_1,q_2,q_3\rangle\) is the quadratic inverse system of
the generic \((1,11,3)\) algebra, it is
\[
 \det[\,1,z_{i1},\ldots,z_{i,11},
                   q_1(z_i),q_2(z_i),q_3(z_i)\,]_{i=1}^{15},
\]
of degree \(17\), by the lowest-degree Cauchy--Binet expansion
in~\eqref{next:dual-eq-10}. Lemma~\ref{next:quadratic-essential}
gives \(11(15-1)=154\) essential variables for every independent
three-plane of quadrics. The dimension bounds
in~\eqref{next:dual-eq-23} apply to planes from the generic
smoothable configuration locus.

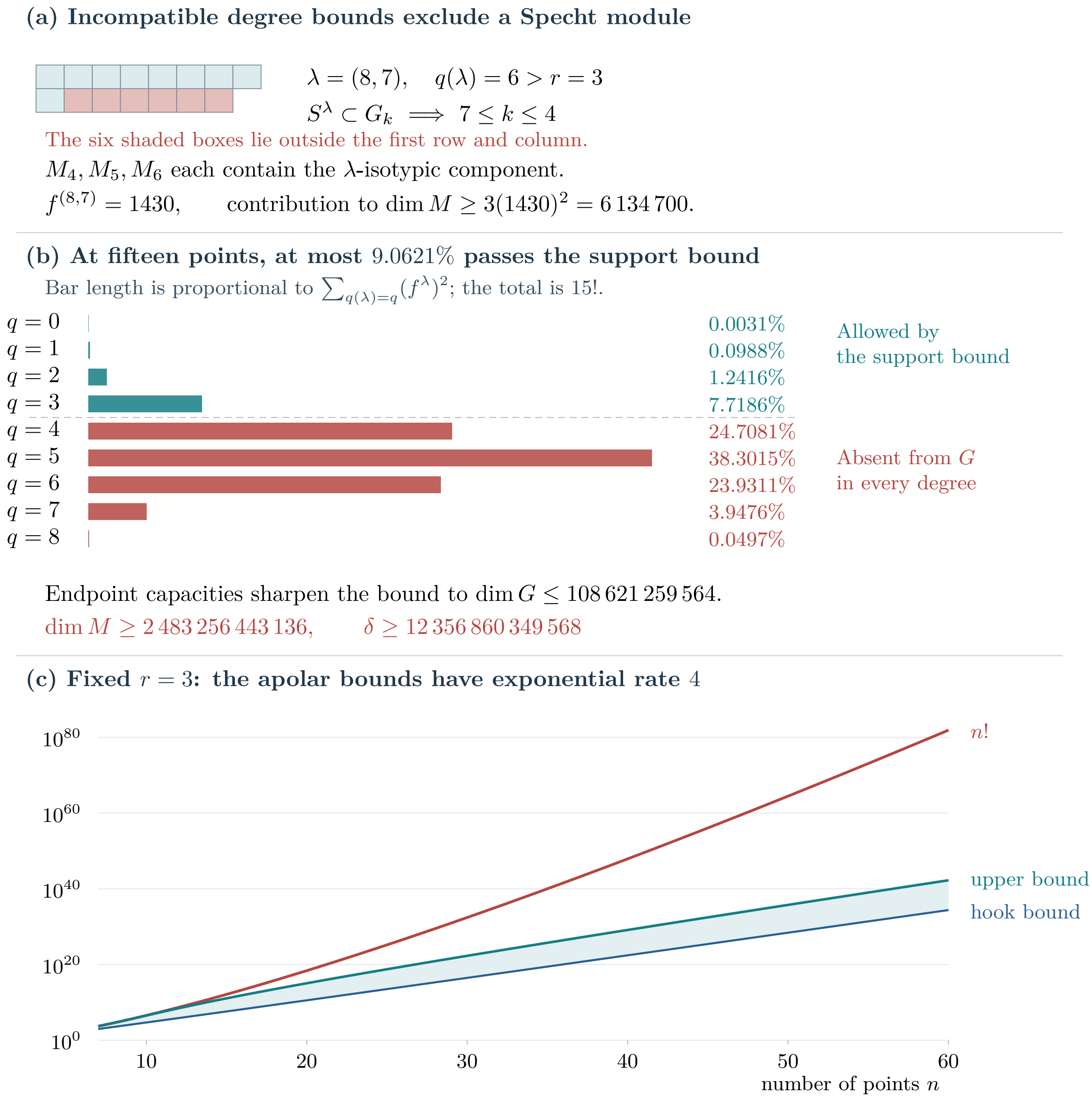
\begin{figure}[!htbp]
\centering
\begin{tikzpicture}[x=1cm,y=.84cm,font=\small]
\path[use as bounding box] (0,-.5) rectangle (16.1,20.0);
\node[anchor=west,font=\bfseries\small,Ink] at (0,19.65)
 {(a) Incompatible degree bounds exclude a Specht module};

\draw[Ink!55,fill=Teal!15,line width=.4pt] (0.3,18.8) rectangle ++(0.43,-0.43);
\draw[Ink!55,fill=Teal!15,line width=.4pt] (0.73,18.8) rectangle ++(0.43,-0.43);
\draw[Ink!55,fill=Teal!15,line width=.4pt] (1.16,18.8) rectangle ++(0.43,-0.43);
\draw[Ink!55,fill=Teal!15,line width=.4pt] (1.59,18.8) rectangle ++(0.43,-0.43);
\draw[Ink!55,fill=Teal!15,line width=.4pt] (2.02,18.8) rectangle ++(0.43,-0.43);
\draw[Ink!55,fill=Teal!15,line width=.4pt] (2.4499999999999997,18.8) rectangle ++(0.43,-0.43);
\draw[Ink!55,fill=Teal!15,line width=.4pt] (2.88,18.8) rectangle ++(0.43,-0.43);
\draw[Ink!55,fill=Teal!15,line width=.4pt] (3.3099999999999996,18.8) rectangle ++(0.43,-0.43);
\draw[Ink!55,fill=Teal!15,line width=.4pt] (0.3,18.37) rectangle ++(0.43,-0.43);
\draw[Ink!55,fill=Coral!35,line width=.4pt] (0.73,18.37) rectangle ++(0.43,-0.43);
\draw[Ink!55,fill=Coral!35,line width=.4pt] (1.16,18.37) rectangle ++(0.43,-0.43);
\draw[Ink!55,fill=Coral!35,line width=.4pt] (1.59,18.37) rectangle ++(0.43,-0.43);
\draw[Ink!55,fill=Coral!35,line width=.4pt] (2.02,18.37) rectangle ++(0.43,-0.43);
\draw[Ink!55,fill=Coral!35,line width=.4pt] (2.4499999999999997,18.37) rectangle ++(0.43,-0.43);
\draw[Ink!55,fill=Coral!35,line width=.4pt] (2.88,18.37) rectangle ++(0.43,-0.43);
\node[anchor=west] at (4.3,18.55) {$\lambda=(8,7),\quad q(\lambda)=6>r=3$};
\node[anchor=west] at (4.3,17.93)
 {$S^\lambda\subset G_k\ \Longrightarrow\ 7\le k\le4$};
\node[anchor=west,font=\footnotesize,Coral] at (.3,17.45)
 {The six shaded boxes lie outside the first row and column.};
\node[anchor=west,align=left] at (.3,16.55)
 {$M_4,M_5,M_6$ each contain the $\lambda$-isotypic component.\\[3pt]
 $f^{(8,7)}=1430$,\qquad contribution to $\dim M\ge3(1430)^2=6\,134\,700$.};
\draw[Ink!25] (0,15.75)--(16,15.75);
\node[anchor=west,font=\bfseries\small,Ink] at (0,15.3)
 {(b) At fifteen points, at most $9.0621\%$ passes the support bound};
\node[anchor=west,font=\footnotesize,Ink!90] at (.3,14.7)
 {Bar length is proportional to $\sum_{q(\lambda)=q}(f^\lambda)^2$; the total is $15!$.};

\node[anchor=east] at (.8,14.1) {$q=0$};
\fill[Teal!85] (1.1,13.95) rectangle ++(0.00069026,.30);
\node[anchor=west,font=\footnotesize,Teal] at (10.45,14.1) {0.0031\%};
\node[anchor=east] at (.8,13.61) {$q=1$};
\fill[Teal!85] (1.1,13.459999999999999) rectangle ++(0.02223,.30);
\node[anchor=west,font=\footnotesize,Teal] at (10.45,13.61) {0.0988\%};
\node[anchor=east] at (.8,13.12) {$q=2$};
\fill[Teal!85] (1.1,12.969999999999999) rectangle ++(0.27937,.30);
\node[anchor=west,font=\footnotesize,Teal] at (10.45,13.12) {1.2416\%};
\node[anchor=east] at (.8,12.629999999999999) {$q=3$};
\fill[Teal!85] (1.1,12.479999999999999) rectangle ++(1.73669,.30);
\node[anchor=west,font=\footnotesize,Teal] at (10.45,12.629999999999999) {7.7186\%};
\node[anchor=east] at (.8,12.14) {$q=4$};
\fill[Coral!85] (1.1,11.99) rectangle ++(5.55932,.30);
\node[anchor=west,font=\footnotesize,Coral] at (10.45,12.14) {24.7081\%};
\node[anchor=east] at (.8,11.649999999999999) {$q=5$};
\fill[Coral!85] (1.1,11.499999999999998) rectangle ++(8.61783,.30);
\node[anchor=west,font=\footnotesize,Coral] at (10.45,11.649999999999999) {38.3015\%};
\node[anchor=east] at (.8,11.16) {$q=6$};
\fill[Coral!85] (1.1,11.01) rectangle ++(5.38449,.30);
\node[anchor=west,font=\footnotesize,Coral] at (10.45,11.16) {23.9311\%};
\node[anchor=east] at (.8,10.67) {$q=7$};
\fill[Coral!85] (1.1,10.52) rectangle ++(0.88822,.30);
\node[anchor=west,font=\footnotesize,Coral] at (10.45,10.67) {3.9476\%};
\node[anchor=east] at (.8,10.18) {$q=8$};
\fill[Coral!85] (1.1,10.03) rectangle ++(0.01117,.30);
\node[anchor=west,font=\footnotesize,Coral] at (10.45,10.18) {0.0497\%};
\node[anchor=west,align=left,font=\footnotesize,Teal] at (12.4,13.7)
 {Allowed by\\the support bound};
\node[anchor=west,align=left,font=\footnotesize,Coral] at (12.4,11.4)
 {Absent from $G$\\in every degree};
\draw[Ink!35,densely dashed] (.2,12.385)--(11.9,12.385);
\node[anchor=west] at (.3,9.15)
 {Endpoint capacities sharpen the bound to $\dim G\le108\,621\,259\,564$.};
\node[anchor=west,Coral] at (.3,8.55)
 {$\dim M\ge2\,483\,256\,443\,136,\qquad
   \delta\ge12\,356\,860\,349\,568$};
\draw[Ink!25] (0,8.05)--(16,8.05);
\node[anchor=west,font=\bfseries\small,Ink] at (0,7.6)
 {(c) Fixed $r=3$: the apolar bounds have exponential rate $4$};

\draw[Ink!12] (1.25,1.05)--(14.25,1.05);\node[anchor=east,font=\footnotesize] at (1.13,1.05) {$10^{0}$};
\draw[Ink!12] (1.25,2.428048780487805)--(14.25,2.428048780487805);\node[anchor=east,font=\footnotesize] at (1.13,2.428048780487805) {$10^{20}$};
\draw[Ink!12] (1.25,3.8060975609756103)--(14.25,3.8060975609756103);\node[anchor=east,font=\footnotesize] at (1.13,3.8060975609756103) {$10^{40}$};
\draw[Ink!12] (1.25,5.184146341463415)--(14.25,5.184146341463415);\node[anchor=east,font=\footnotesize] at (1.13,5.184146341463415) {$10^{60}$};
\draw[Ink!12] (1.25,6.56219512195122)--(14.25,6.56219512195122);\node[anchor=east,font=\footnotesize] at (1.13,6.56219512195122) {$10^{80}$};
\draw[Ink!35] (1.9858490566037736,1.05)--++(0,-.08);\node[below,font=\footnotesize] at (1.9858490566037736,0.9500000000000001) {$ 10 $};
\draw[Ink!35] (4.438679245283019,1.05)--++(0,-.08);\node[below,font=\footnotesize] at (4.438679245283019,0.9500000000000001) {$ 20 $};
\draw[Ink!35] (6.8915094339622645,1.05)--++(0,-.08);\node[below,font=\footnotesize] at (6.8915094339622645,0.9500000000000001) {$ 30 $};
\draw[Ink!35] (9.34433962264151,1.05)--++(0,-.08);\node[below,font=\footnotesize] at (9.34433962264151,0.9500000000000001) {$ 40 $};
\draw[Ink!35] (11.797169811320755,1.05)--++(0,-.08);\node[below,font=\footnotesize] at (11.797169811320755,0.9500000000000001) {$ 50 $};
\draw[Ink!35] (14.25,1.05)--++(0,-.08);\node[below,font=\footnotesize] at (14.25,0.9500000000000001) {$ 60 $};
\fill[Teal!12] (1.2500,1.3051) -- (1.4953,1.3673) -- (1.7406,1.4329) -- (1.9858,1.5001) -- (2.2311,1.5667) -- (2.4764,1.6313) -- (2.7217,1.6938) -- (2.9670,1.7543) -- (3.2123,1.8130) -- (3.4575,1.8702) -- (3.7028,1.9260) -- (3.9481,1.9807) -- (4.1934,2.0345) -- (4.4387,2.0874) -- (4.6840,2.1395) -- (4.9292,2.1911) -- (5.1745,2.2420) -- (5.4198,2.2925) -- (5.6651,2.3424) -- (5.9104,2.3920) -- (6.1557,2.4412) -- (6.4009,2.4900) -- (6.6462,2.5385) -- (6.8915,2.5868) -- (7.1368,2.6348) -- (7.3821,2.6825) -- (7.6274,2.7300) -- (7.8726,2.7773) -- (8.1179,2.8243) -- (8.3632,2.8712) -- (8.6085,2.9180) -- (8.8538,2.9645) -- (9.0991,3.0110) -- (9.3443,3.0572) -- (9.5896,3.1034) -- (9.8349,3.1494) -- (10.0802,3.1953) -- (10.3255,3.2410) -- (10.5708,3.2867) -- (10.8160,3.3323) -- (11.0613,3.3777) -- (11.3066,3.4231) -- (11.5519,3.4684) -- (11.7972,3.5136) -- (12.0425,3.5587) -- (12.2877,3.6037) -- (12.5330,3.6487) -- (12.7783,3.6935) -- (13.0236,3.7384) -- (13.2689,3.7831) -- (13.5142,3.8278) -- (13.7594,3.8724) -- (14.0047,3.9170) -- (14.2500,3.9615) -- (14.2500,3.4193) -- (14.0047,3.3781) -- (13.7594,3.3369) -- (13.5142,3.2956) -- (13.2689,3.2544) -- (13.0236,3.2132) -- (12.7783,3.1720) -- (12.5330,3.1308) -- (12.2877,3.0896) -- (12.0425,3.0484) -- (11.7972,3.0073) -- (11.5519,2.9661) -- (11.3066,2.9249) -- (11.0613,2.8837) -- (10.8160,2.8426) -- (10.5708,2.8014) -- (10.3255,2.7603) -- (10.0802,2.7192) -- (9.8349,2.6780) -- (9.5896,2.6369) -- (9.3443,2.5958) -- (9.0991,2.5547) -- (8.8538,2.5136) -- (8.6085,2.4726) -- (8.3632,2.4315) -- (8.1179,2.3904) -- (7.8726,2.3494) -- (7.6274,2.3084) -- (7.3821,2.2674) -- (7.1368,2.2264) -- (6.8915,2.1854) -- (6.6462,2.1444) -- (6.4009,2.1035) -- (6.1557,2.0625) -- (5.9104,2.0216) -- (5.6651,1.9808) -- (5.4198,1.9399) -- (5.1745,1.8991) -- (4.9292,1.8583) -- (4.6840,1.8175) -- (4.4387,1.7768) -- (4.1934,1.7361) -- (3.9481,1.6955) -- (3.7028,1.6549) -- (3.4575,1.6144) -- (3.2123,1.5739) -- (2.9670,1.5335) -- (2.7217,1.4932) -- (2.4764,1.4530) -- (2.2311,1.4129) -- (1.9858,1.3729) -- (1.7406,1.3332) -- (1.4953,1.2936) -- (1.2500,1.2543) -- cycle;
\draw[Coral,line width=1.2pt] (1.2500,1.3051) -- (1.4953,1.3673) -- (1.7406,1.4331) -- (1.9858,1.5020) -- (2.2311,1.5737) -- (2.4764,1.6481) -- (2.7217,1.7248) -- (2.9670,1.8038) -- (3.2123,1.8849) -- (3.4575,1.9678) -- (3.7028,2.0526) -- (3.9481,2.1391) -- (4.1934,2.2272) -- (4.4387,2.3168) -- (4.6840,2.4080) -- (4.9292,2.5004) -- (5.1745,2.5943) -- (5.4198,2.6894) -- (5.6651,2.7857) -- (5.9104,2.8832) -- (6.1557,2.9818) -- (6.4009,3.0815) -- (6.6462,3.1823) -- (6.8915,3.2841) -- (7.1368,3.3868) -- (7.3821,3.4905) -- (7.6274,3.5952) -- (7.8726,3.7007) -- (8.1179,3.8071) -- (8.3632,3.9143) -- (8.6085,4.0224) -- (8.8538,4.1312) -- (9.0991,4.2408) -- (9.3443,4.3512) -- (9.5896,4.4624) -- (9.8349,4.5742) -- (10.0802,4.6867) -- (10.3255,4.8000) -- (10.5708,4.9139) -- (10.8160,5.0285) -- (11.0613,5.1437) -- (11.3066,5.2595) -- (11.5519,5.3760) -- (11.7972,5.4930) -- (12.0425,5.6107) -- (12.2877,5.7289) -- (12.5330,5.8477) -- (12.7783,5.9671) -- (13.0236,6.0870) -- (13.2689,6.2075) -- (13.5142,6.3285) -- (13.7594,6.4500) -- (14.0047,6.5720) -- (14.2500,6.6945);
\draw[Teal,line width=1.1pt] (1.2500,1.3051) -- (1.4953,1.3673) -- (1.7406,1.4329) -- (1.9858,1.5001) -- (2.2311,1.5667) -- (2.4764,1.6313) -- (2.7217,1.6938) -- (2.9670,1.7543) -- (3.2123,1.8130) -- (3.4575,1.8702) -- (3.7028,1.9260) -- (3.9481,1.9807) -- (4.1934,2.0345) -- (4.4387,2.0874) -- (4.6840,2.1395) -- (4.9292,2.1911) -- (5.1745,2.2420) -- (5.4198,2.2925) -- (5.6651,2.3424) -- (5.9104,2.3920) -- (6.1557,2.4412) -- (6.4009,2.4900) -- (6.6462,2.5385) -- (6.8915,2.5868) -- (7.1368,2.6348) -- (7.3821,2.6825) -- (7.6274,2.7300) -- (7.8726,2.7773) -- (8.1179,2.8243) -- (8.3632,2.8712) -- (8.6085,2.9180) -- (8.8538,2.9645) -- (9.0991,3.0110) -- (9.3443,3.0572) -- (9.5896,3.1034) -- (9.8349,3.1494) -- (10.0802,3.1953) -- (10.3255,3.2410) -- (10.5708,3.2867) -- (10.8160,3.3323) -- (11.0613,3.3777) -- (11.3066,3.4231) -- (11.5519,3.4684) -- (11.7972,3.5136) -- (12.0425,3.5587) -- (12.2877,3.6037) -- (12.5330,3.6487) -- (12.7783,3.6935) -- (13.0236,3.7384) -- (13.2689,3.7831) -- (13.5142,3.8278) -- (13.7594,3.8724) -- (14.0047,3.9170) -- (14.2500,3.9615);
\draw[Blue,line width=.9pt] (1.2500,1.2543) -- (1.4953,1.2936) -- (1.7406,1.3332) -- (1.9858,1.3729) -- (2.2311,1.4129) -- (2.4764,1.4530) -- (2.7217,1.4932) -- (2.9670,1.5335) -- (3.2123,1.5739) -- (3.4575,1.6144) -- (3.7028,1.6549) -- (3.9481,1.6955) -- (4.1934,1.7361) -- (4.4387,1.7768) -- (4.6840,1.8175) -- (4.9292,1.8583) -- (5.1745,1.8991) -- (5.4198,1.9399) -- (5.6651,1.9808) -- (5.9104,2.0216) -- (6.1557,2.0625) -- (6.4009,2.1035) -- (6.6462,2.1444) -- (6.8915,2.1854) -- (7.1368,2.2264) -- (7.3821,2.2674) -- (7.6274,2.3084) -- (7.8726,2.3494) -- (8.1179,2.3904) -- (8.3632,2.4315) -- (8.6085,2.4726) -- (8.8538,2.5136) -- (9.0991,2.5547) -- (9.3443,2.5958) -- (9.5896,2.6369) -- (9.8349,2.6780) -- (10.0802,2.7192) -- (10.3255,2.7603) -- (10.5708,2.8014) -- (10.8160,2.8426) -- (11.0613,2.8837) -- (11.3066,2.9249) -- (11.5519,2.9661) -- (11.7972,3.0073) -- (12.0425,3.0484) -- (12.2877,3.0896) -- (12.5330,3.1308) -- (12.7783,3.1720) -- (13.0236,3.2132) -- (13.2689,3.2544) -- (13.5142,3.2956) -- (13.7594,3.3369) -- (14.0047,3.3781) -- (14.2500,3.4193);
\node[anchor=west,Coral,font=\footnotesize] at (14.45,6.68) {$n!$};
\node[anchor=west,Teal,font=\footnotesize] at (14.45,3.96) {upper bound};
\node[anchor=west,Blue,font=\footnotesize] at (14.45,3.40) {hook bound};
\node[anchor=east,font=\footnotesize] at (14.25,.23) {number of points $n$};
\end{tikzpicture}
\caption{Restrictions on the alternating quotient of a generic orbit
algebra. In (a), first occurrence and sign duality exclude
$S^{(8,7)}$ from $G$ and place its full isotypic component in
$M_4,M_5,M_6$, where $M=K/R$. Panel (b) groups all partitions
of fifteen by $q(\lambda)$; coral bars are absent from $G$,
and endpoint multiplicities further reduce the teal contribution.
Percentages are rounded. In (c), the shaded region bounds $\dim G$
between $\binom{2n-2}{n-1}$ and
$\sum_{q(\lambda)\leq3}(f^\lambda)^2$ as $n$ varies.
Both bounds have exponential rate $4$, and the orbit algebra has
dimension $n!$.}
\label{next:apolar-figure}
\end{figure}

\Needspace{10\baselineskip}
\section{Alternating quotients in fixed ambient dimension}
\label{dim:apolar}

Let \(k\) be a characteristic-zero field and let \(n\geq2\).  Put
\[
 N_j(d)=\binom{d+j}{d},
 \qquad
 s=\min\{j:N_j(d)\geq n\},
 \qquad r=n-N_{s-1}(d),
\]
where \(N_{-1}(d)=0\).  Choose an \(r\)-plane
\(Q\subset k[z_1,\ldots,z_d]_s\), with basis \(q_1,\ldots,q_r\).
Let \(f_1,\ldots,f_n\) consist of all monomials of degree less than
\(s\), followed by this basis of \(Q\).  Define
\begin{equation}\label{dim:fixed-eq:D1}
\Delta_Q(z_1,\ldots,z_n)
   =\det(f_j(z_i))_{i,j=1}^n,
 \qquad
 G_Q=k[\partial_{ia}]/\operatorname{Ann}(\Delta_Q).
\end{equation}
Changing the basis of \(Q\) scales \(\Delta_Q\) by a nonzero constant
and leaves its apolar algebra \(G_Q\) unchanged.

The degree of this determinant is
\begin{equation}\label{dim:fixed-eq:D2}
\mu=\mu(n,d)
  =d\binom{d+s-1}{d+1}+sr
  =sn-\binom{d+s}{d+1}.
\end{equation}
These are the minimum-degree alternating determinants in \(d\)
alphabets.  For \(s\geq2\), the columns include every linear monomial.

\begin{lemma}[Evaluation determinants and alternating pairings]
\label{dim:determinant-source}
The polynomial \(\Delta_Q\) is nonzero and is invariant under
simultaneous translation of its \(n\) rows.  If \(s\geq2\), it has
\(D=d(n-1)\) essential variables.

Let \(T_0\) be the graded length-\(n\) algebra whose Macaulay inverse
system is
\begin{equation}\label{dim:fixed-eq:D3}
\bigoplus_{j<s}k[z_1,\ldots,z_d]_j\ \oplus Q.
\end{equation}
The image of the alternating multiplication pairing
\begin{equation}\label{dim:fixed-eq:D4}
\phi_{T_0}:T_0^{\otimes n}
 \longrightarrow (T_0^{\otimes n})^*\otimes\Lambda^nT_0
\end{equation}
has dimension \(\dim_kG_Q\).
If \(T_0\) arises by taking the homogeneous limit of a generic
configuration, \(G_Q\) is the canonical alternating quotient of
its Artinian orbit algebra.
\end{lemma}

\begin{proof}
The column polynomials are linearly independent; their evaluation
determinant is therefore nonzero.  A common translation acts
triangularly on all columns of degree less than \(s\), with diagonal
entries one, and changes each degree-\(s\) column by lower-degree
columns.  Thus the determinant is unchanged.

Its linear annihilators beyond the diagonal translation directions
form an \(S_n\)-submodule of \(\mathrm{Std}_n\otimes k^d\).
Any nonzero such submodule contains
\(\mathrm{Std}_n\otimes kw\) for some \(w\ne0\).
Together with the diagonal direction, this would remove one entire
alphabet.  However, the minimum alternating degree satisfies
\(\mu(n,d-1)>\mu(n,d)\) when \(n\geq d+1\): the ordered list of
monomial degrees in \(d-1\) variables is termwise at least the
corresponding list in \(d\) variables, and its \((d+1)\)-st entry
is strictly larger.  No nonzero polynomial of degree \(\mu(n,d)\)
can therefore remain.  For \(d=1\), removal of the only alphabet
is immediately impossible.

The top exterior functional on \(T_0^{\otimes n}\) is alternation
followed by the determinant in a basis of \(T_0\).
Writing its inverse polynomial in the dual basis gives \eqref{dim:fixed-eq:D1},
with the usual factorial normalization of Macaulay duality.
Its multiplication radical is the differential annihilator
of that polynomial, proving the statement about \eqref{dim:fixed-eq:D4}.

For a configuration, expand
\(\det(\exp(p_i\cdot z_j))\) by Cauchy--Binet.  Its lowest-degree
nonzero part selects all monomials of degree below \(s\), followed
by \(r\) degree-\(s\) columns.  Their Pl\"ucker coordinates are the
inverse-system plane of the homogeneous limit.  This identifies
the resulting top alternating functional with \(\Delta_Q\).
\end{proof}

\begin{theorem}[Uniform fixed-dimension loss]
\label{dim:fixed-dimension-loss}
For \(d\geq3\), uniformly over every plane \(Q\) in \eqref{dim:fixed-eq:D1},
\begin{equation}\label{dim:fixed-eq:D5}
\limsup_{n\to\infty}
 \left(\frac{\dim_kG_Q}{n!}\right)^{1/n}
 \leq
 \rho_d:=
 \frac{e^{d+1}d!}{2^d(d+1)^d}
 <1.
\end{equation}
In particular, for each fixed \(d\geq3\), there is an \(N(d)\) such
that every one of these alternating quotients has dimension
strictly less than \(n!\) for every \(n\geq N(d)\).  Their fractions
of the factorial dimension tend to zero exponentially.
The first two constants are
\begin{equation}\label{dim:fixed-eq:D6}
\rho_3=\frac{3e^4}{256}<0.640,\qquad
 \rho_4=\frac{3e^5}{1250}<0.357.
\end{equation}
No smoothability or regularity hypothesis on the homogeneous
algebra is required for this upper bound.
\end{theorem}

\begin{proof}
For any nonzero homogeneous polynomial of degree \(\mu\) in at
most \(D\) variables, its apolar algebra has a perfect pairing
between degrees \(j\) and \(\mu-j\).  Therefore
\begin{equation}\label{dim:fixed-eq:D7}
\dim(G_Q)_j\leq
 \min\left\{\binom{D+j-1}{j},
             \binom{D+\mu-j-1}{\mu-j}\right\}.
\end{equation}
Summing gives the compressed-algebra upper bound
\begin{equation}\label{dim:fixed-eq:D8}
\dim G_Q\leq C(D,\mu):=
 \begin{cases}
 2\binom{D+a}{D},&\mu=2a+1,\\[2mm]
 \binom{D+a}{D}+\binom{D+a-1}{D},&\mu=2a.
 \end{cases}
\end{equation}
The bound depends on the number of variables and the degree,
so it is uniform in \(Q\). The comparison with \(n!\) comes from
the growth of the minimum alternating degree. For fixed \(d\),
the monomial count gives
\begin{equation}\label{dim:fixed-eq:D9}
s=(d!n)^{1/d}+O_d(1),
 \qquad
 \mu=c_dn^{1+1/d}+O_d(n),
 \qquad
 c_d=\frac{d}{d+1}(d!)^{1/d}.
\end{equation}
Set \(D=d(n-1)\) and \(a=\lfloor\mu/2\rfloor\).  Since
\(a/D\to\infty\), expansion of the factorial ratio yields
\begin{equation}\label{dim:fixed-eq:D10}
\log\binom{D+a}{D}
   =D\bigl(\log(a/D)+1\bigr)+o(n)
   =n\log n+
     n\left(d+d\log\frac{c_d}{2d}\right)+o(n).
\end{equation}
Indeed,
\[
 \log\binom{D+a}{D}
 =D\log a-\log(D!)+O(D^2/a),\qquad
 D^2/a=O_d(n^{1-1/d})=o(n).
\]
Stirling's formula contributes \(O(\log D)\). The leading
\(n\log n\) term in~\eqref{dim:fixed-eq:D10} cancels the same term
in \(\log(n!)=n\log n-n+o(n)\), leaving
\begin{equation}\label{dim:fixed-eq:D11}
\log\frac{C(D,\mu)}{n!}
   =n\log\left(\frac{e^{d+1}d!}{[2(d+1)]^d}\right)+o(n).
\end{equation}
This proves the bound uniformly in \(Q\).

Finally \(\rho_3<1\), and
\begin{equation}\label{dim:fixed-eq:D12}
\frac{\rho_{d+1}}{\rho_d}
  =\frac e2\left(\frac{d+1}{d+2}\right)^{d+1}<1.
\end{equation}
For example, \(\log(1-x)\leq-x\) bounds the ratio by
\(\frac12e^{1/(d+2)}<1\).  Thus all \(d\geq3\) have
\(\rho_d<1\).
\end{proof}

An explicit finite estimate is available without asymptotic notation.
Since \(\mu/n\leq ds/(d+1)\), replacing \(D\) by \(dn\) in \eqref{dim:fixed-eq:D8}
and using the elementary factorial estimates gives
\begin{equation}\label{dim:fixed-eq:D13}
\frac{\dim G_Q}{n!}
 \leq
 2\left[
  \rho_d\left(1+\frac{2(d+1)}s\right)^d
 \right]^n.
\end{equation}
Thus, for example, choose an integer \(S_d\) satisfying
\begin{equation}\label{dim:fixed-eq:D14}
S_d\geq
 \frac{2(d+1)}{\rho_d^{-1/(2d)}-1}.
\end{equation}
For every \(n>\binom{d+S_d-1}{d}\), the bound is
\(2\rho_d^{n/2}\).  Increasing this explicit threshold to satisfy
\(n>2\log2/(-\log\rho_d)\) makes the upper bound smaller than one.
This gives an explicit sufficient threshold for strict rank loss.

\subsection{The source algebra and a sharper coefficient bound}

The point algebra imposes relations beyond those from essential
variables.  These relations give a stronger bound without computing an
individual catalecticant.

\begin{lemma}[Differential trace identity]
\label{dim:differential-trace}
Let \(T=k[x_1,\ldots,x_d]/I\) be a length-\(n\) local algebra
supported at the origin, with inverse-system basis \(f_1,\ldots,f_n\),
and put \(\Delta_T=\det(f_j(z_i))\).  For every polynomial \(g\),
\begin{equation}\label{dim:fixed-eq:D27}
\sum_{i=1}^n g(\partial_{i1},\ldots,\partial_{id})\Delta_T
   =\operatorname{Tr}_T(m_g)\Delta_T
   =n g(0)\Delta_T.
\end{equation}
In addition, \(I(\partial_{i1},\ldots,\partial_{id})\Delta_T=0\)
for each row \(i\).  Thus every positive multisymmetric power-sum
differential operator annihilates \(\Delta_T\).
\end{lemma}

\begin{proof}
Let \(C\) represent \(g(\partial)\) on the inverse system, and let
\(M=(f_j(z_i))\). The matrix \(C\) is dual to multiplication by
\(g\) on \(T\), and multilinearity gives
\[
 \sum_i g(\partial_i)\det M
 =[\varepsilon]\det\bigl(M(1+\varepsilon C)\bigr)
 =\operatorname{Tr}(C)\det M.
\]
Since \(m_{g-g(0)}\) is nilpotent,
\(\operatorname{Tr}(C)=\operatorname{Tr}(m_g)=ng(0)\).
Each element of \(I\) annihilates each inverse-system basis
polynomial, proving the row-ideal assertion.
\end{proof}

Let
\[
 h(z)=\sum_{j=0}^{s-1}\binom{d+j-1}{d-1}z^j+r z^s,
 \qquad c_j=[z^j]h(z)^{n-1}.
\]
This is the Hilbert polynomial of the point algebra \(T_0\) in
\eqref{dim:fixed-eq:D3}, and the \(c_j\) are the graded dimensions of
\(T_0^{\otimes(n-1)}\).

\begin{proposition}[A bound from the point algebra]
\label{dim:source-bound}
For every plane \(Q\) in \eqref{dim:fixed-eq:D1},
\begin{equation}\label{dim:fixed-eq:D28}
\dim G_Q\leq
 U(n,d):=\sum_{j=0}^{\mu}\min\{c_j,c_{\mu-j}\}.
\end{equation}
For each fixed \(d\geq3\), this implies the stronger estimate
\begin{equation}\label{dim:fixed-eq:D29}
\limsup_{n\to\infty}
 \left(\frac{\dim G_Q}{n!}\right)^{1/n}
 \leq\beta_d:=
 e d\inf_{a>0}\left[
 e^{ad/(2(d+1))}\int_0^1u^{d-1}e^{-au}\,du\right].
\end{equation}
In particular, the explicit constant
\begin{equation}\label{dim:fixed-eq:D30}
\eta_d:=\rho_d\left(
 1-e^{-2(d+1)}\sum_{j=0}^{d-1}\frac{[2(d+1)]^j}{j!}
 \right)
\end{equation}
satisfies \(\beta_d\leq\eta_d<\rho_d<1\).
Numerically \(\beta_3<0.630068\) and \(\beta_4<0.352215\).
\end{proposition}

\begin{proof}
Translation invariance lets us set the last row to zero without
changing the apolar algebra. Lemma~\ref{dim:differential-trace}
then gives
\[
 T_0^{\otimes(n-1)}\twoheadrightarrow G_Q,
 \qquad \dim(G_Q)_j\leq c_j.
\]
Perfect apolar duality gives the complementary bound
\[
 \dim(G_Q)_j=\dim(G_Q)_{\mu-j}\leq c_{\mu-j}.
\]
Summing the smaller of these two bounds proves~\eqref{dim:fixed-eq:D28}.

For \(0<z<1\), symmetry and positivity give
\begin{equation}\label{dim:fixed-eq:D31}
\dim G_Q\leq2\sum_{j\leq\mu/2}c_j
       \leq2z^{-\mu/2}h(z)^{n-1}.
\end{equation}
Fix \(a>0\) and set \(z=e^{-a/s}\).  Since the last partially
filled degree contributes \(O(s^{d-1})\) of the \(n\) basis
elements, a Riemann sum yields
\[
 \frac{h(e^{-a/s})}{n}
    \longrightarrow d\int_0^1u^{d-1}e^{-au}\,du,
 \qquad
 \frac{\mu}{ns}\longrightarrow\frac d{d+1}.
\]
Take \(n\)-th roots of \eqref{dim:fixed-eq:D31}, divide by \((n!)^{1/n}\), and
use \((n!)^{1/n}\sim n/e\).  Taking the infimum over the fixed
parameter \(a\) gives \eqref{dim:fixed-eq:D29}.  At \(a=2(d+1)\), the elementary
integer incomplete-gamma formula evaluates the integral and gives
\eqref{dim:fixed-eq:D30}.  The minimum in \eqref{dim:fixed-eq:D29} is characterized uniquely by
\[
 \frac{\int_0^1u^d e^{-au}\,du}
      {\int_0^1u^{d-1}e^{-au}\,du}
     =\frac{d}{2(d+1)}.
\]
The left side decreases strictly from \(d/(d+1)\) to zero;
its derivative is minus the variance of \(u\) for the indicated
positive probability density.  This also gives a direct way to
evaluate the displayed numerical upper bounds.
\end{proof}

For explicit finite consequences, take the full fat point
\(T_0=k[x_1,\ldots,x_d]/\mathfrak m^{s+1}\).  The integer
bound \eqref{dim:fixed-eq:D28} gives
\begin{equation}\label{dim:fixed-eq:D32}
\begin{array}{c|c|c|c}
 (n,d)&\text{point Hilbert function}&\mu& U(n,d)/n!\\ \hline
 (70,4)&(1,4,10,20,35)&224&<12181/10^6\\
 (120,3)&(1,3,6,10,15,21,28,36)
     &630&<208011/10^6
 \end{array}
\end{equation}
These rational inequalities are obtained by the finite coefficient
sums in \eqref{dim:fixed-eq:D28}, with integer arithmetic.  The resulting alternating
pairing loses more than \(98.7819\%\) and \(79.1989\%\), respectively,
of its generic factorial rank.  Both fat points are homogeneous
limits of generic configurations and therefore lie on the principal
Hilbert component.  The inequalities concern the rank of the
pairing. Corollary~\ref{dim:pairing-rank-loss} relates this rank
to the scheme-theoretic cover fibres.

\subsection{Planar equality}

\begin{theorem}[Planar equality]
\label{dim:planar-equality}
For \(d=2\), every plane \(Q\) in \eqref{dim:fixed-eq:D1} has
\begin{equation}\label{dim:fixed-eq:D15}
\dim_kG_Q=n!.
\end{equation}
More generally, for every configuration of \(n\) distinct points in
the affine plane, its filtered permutation-orbit algebra \(R\)
has a Gorenstein Artinian reduction \(A=R/(t)\), and its canonical
alternating quotient is all of \(A\).  Its sign socle has degree
equal to the sum of the degrees of a filtered point-algebra basis.
The same equality holds for configurations on the affine line.
\end{theorem}

\begin{proof}
We apply Haiman's theorem that the planar isospectral Hilbert
scheme is Cohen--Macaulay and Gorenstein.  Since
\(\operatorname{Hilb}^n(\mathbb A^2)\) is smooth, its Procesi algebra
\(\mathcal P\) is finite locally free of rank \(n!\).
Haiman's alternating multiplication pairing gives a perfect
pairing of algebra bundles
\begin{equation}\label{dim:fixed-eq:D16}
\mathcal P\otimes\mathcal P
    \longrightarrow\det\mathcal T,
\end{equation}
where \(\mathcal T\) is the rank-\(n\) tautological algebra.
These are Theorem 3.1 and Propositions 3.7.2 and 3.7.4 of
\cite{gg-HaimanPolygraph}.
They apply to every homogeneous ideal defined by \eqref{dim:fixed-eq:D3}, because
every length-\(n\) planar scheme is a point of that Hilbert scheme.
Specializing the perfect pairing preserves its rank \(n!\).
Lemma~\ref{dim:determinant-source} identifies that rank with
\(\dim G_Q\), proving \eqref{dim:fixed-eq:D15}.

For the assertion about an arbitrary configuration \(Z\), take the
flat dilation family \(tZ\), including its homogeneous limit at
zero, and pull \(\mathcal P\) back to \(\mathbb A^1\).
The resulting algebra is free over \(k[t]\).  Its generic algebra
is the coordinate algebra \(E\) of the ordered permutation orbit,
so torsion-freeness embeds it into \(E[t,t^{-1}]\).
The closed embedding of the isospectral family into
\(\mathbb A^1\times(\mathbb A^2)^n\) shows that its image is
precisely
\begin{equation}\label{dim:fixed-eq:D17}
k[t,tv_{i1},tv_{i2}]=R.
\end{equation}
Thus the orbit filtration describes the pulled-back algebra
also at \(t=0\).

The pairing \eqref{dim:fixed-eq:D16} remains perfect after this base change and after
setting \(t=0\).  A filtered basis of the point algebra, with
degrees \(w_0=0,w_1,\ldots,w_{n-1}\), gives a homogeneous generator
of its determinant line in degree \(\sum_iw_i\).
The resulting pairing on \(A\) is sign-valued and homogeneous in
that total degree.  Thus \(A\) is Gorenstein with the asserted
sign socle, and its top alternating pairing has no radical.
Embedding a line configuration into the plane gives the final
assertion; its orbit algebra is unchanged.
\end{proof}

Together, Theorems~\ref{dim:planar-equality}
and~\ref{dim:fixed-dimension-loss} determine the ambient boundary:
\(\dim G_Q=n!\) for every planar source, and
\(\dim G_Q/n!\) tends to zero exponentially in each fixed
ambient dimension \(d\geq3\).

\subsection{The universal pairing and its specialization}

On the reduced principal component \(H_{n,d}\), let \(\mathcal T\)
be the tautological algebra.  The kernel description in Haiman's
Proposition 3.7.2 extends to this setting:
\begin{equation}\label{dim:fixed-eq:D18}
\phi:\mathcal T^{\otimes n}\longrightarrow
       (\mathcal T^{\otimes n})^*\otimes\det\mathcal T,\qquad
 \operatorname{im}\phi=\mathcal P=\rho_*\mathcal O_{X_{n,d}}.
\end{equation}
Both sides agree on the dense reduced-configuration locus.
For the kernel statement, alternation detects a function at an
ordered tuple of distinct points, while it vanishes at repeated
tuples.  Since \(\mathcal T^{\otimes n}\) and the target of \(\phi\)
are locally free over the reduced integral base, identities
checked on that dense locus hold everywhere.  This is the same
kernel argument as in the planar proof.

\begin{corollary}[Rank loss of the alternating multiplication pairing]
\label{dim:pairing-rank-loss}
Let \(I\in H_{n,d}\) be a homogeneous point with the maximal point
Hilbert function used in \eqref{dim:fixed-eq:D3}.  Then
\begin{equation}\label{dim:fixed-eq:D19}
\operatorname{rank}\phi(I)=\dim G_Q\leq C(d(n-1),\mu(n,d)).
\end{equation}
For fixed \(d\geq3\), this rank is an exponentially vanishing
fraction of the generic rank \(n!\), uniformly in such points.

The induced bilinear pairing on the scheme-theoretic
fibre algebra \(\mathcal P(I)\) has radical of dimension at least
\begin{equation}\label{dim:fixed-eq:D20}
n!-\dim G_Q.
\end{equation}
\end{corollary}

\begin{proof}
The first assertion is Lemma~\ref{dim:determinant-source}.
The quotient map \(\mathcal T^{\otimes n}\to\mathcal P\) remains
surjective on fibres.  The specialized pairing on \(\mathcal P(I)\)
therefore has the same image rank as \(\phi(I)\), and its quotient
by the pairing radical is \(G_Q\).

The finite module \(\mathcal P\) has generic rank \(n!\).
Its fibre dimension at every point is at least this rank:
a minimal set of local generators still generates after passage
to the fraction field.  Thus
\[
 \dim\operatorname{rad}\bigl(\mathcal P(I)\bigr)
   =\dim\mathcal P(I)-\dim G_Q
   \geq n!-\dim G_Q.
\]
Here the radical is the radical of the indicated bilinear
pairing, rather than the nilradical of the fibre algebra.
\end{proof}

The two fibre constructions in~\eqref{dim:fixed-eq:D18} are related by
\[
 \mathcal P\otimes k(I)=\mathcal P(I)\twoheadrightarrow
 \operatorname{im}\phi(I)=G_Q.
\]
The kernel is the pairing radical measured in the corollary.
Nonflatness of \(\mathcal P\) in dimensions three and four follows
from the punctual-divisor arguments; the present calculation
describes the pairing on its fibres.

\subsection{Hooks and Plancherel multiplicities}

For a generic configuration, let \(w_1,\ldots,w_{n-1}\) be the
positive degrees of a filtered point-algebra basis.  The integral
exterior-covariant theorem is independent of their particular
values:
\begin{equation}\label{dim:fixed-eq:D21}
\sum_{k=0}^{n-1}u^k
  \operatorname{Hilb}
   \operatorname{Hom}_{S_n}(S^{(n-k,1^k)},A;z)
       =\prod_{i=1}^{n-1}(1+uz^{w_i}).
\end{equation}
Complementary subsets, as in Corollary~\ref{next:hook-survival},
show that every hook isotypic component survives in \(G_Q\).  In particular
\begin{equation}\label{dim:fixed-eq:D22}
\binom{2n-2}{n-1}\leq\dim G_Q.
\end{equation}

Let \(m_\alpha\) be the ungraded multiplicity of \(S^\alpha\) in
\(G_Q\), let \(f^\alpha=\dim S^\alpha\), and equip partitions of
\(n\) with Plancherel measure
\(\mathbb P_n(\alpha)=(f^\alpha)^2/n!\).  Since \(G_Q\) is a quotient
of the regular orbit algebra,
\begin{equation}\label{dim:fixed-eq:D23}
0\leq\frac{m_\alpha}{f^\alpha}\leq1,\qquad
 \mathbb E_{\mathbb P_n}
       \left[\frac{m_\alpha}{f^\alpha}\right]
       =\frac{\dim G_Q}{n!}.
\end{equation}
Hence, for every fixed \(\varepsilon>0\),
\begin{equation}\label{dim:fixed-eq:D24}
\mathbb P_n\!\left\{
    \alpha:m_\alpha\geq\varepsilon f^\alpha\right\}
 \leq\varepsilon^{-1}(\rho_d+o(1))^n
\end{equation}
in each fixed dimension \(d\geq3\).
Thus the partitions whose multiplicities retain a fixed fraction
of their regular multiplicities have exponentially small
Plancherel probability.

\subsection{All minimum-degree alternants and polarization}

Write \(S=k[z_{ia}:1\leq i\leq n,\ 1\leq a\leq d]\).
The space of diagonal harmonics is
\[
 \mathrm{DH}_{n,d}
  =\{P\in S:q(\partial)P=0
            \text{ for every }q\in S^{S_n}_{+}\}.
\]
Here \(S^{S_n}_{+}\) denotes the invariant polynomials with zero
constant term.  Let \(\mathcal A_{n,d}^{\min}\) be the space of
alternating polynomials of minimum total degree \(\mu(n,d)\),
and let
\[
 \mathcal D_{n,d}^{\min}
  =\operatorname{span}_k
     \{\partial^\alpha P:
           P\in\mathcal A_{n,d}^{\min},\ \alpha\geq0\}.
\]

\begin{theorem}[The entire minimum alternating component]
\label{dim:all-minimum-alternants}
For every fixed \(d\geq3\),
\begin{equation}\label{dim:fixed-eq:D33}
\mathcal D_{n,d}^{\min}\subseteq\mathrm{DH}_{n,d},\qquad
 \limsup_{n\to\infty}
  \left(\frac{\dim\mathcal D_{n,d}^{\min}}{n!}\right)^{1/n}
       \leq\beta_d<1.
\end{equation}
The same upper bound holds with \(n!\) in the denominator
replaced by \(\dim\mathrm{DH}_{n,d}\).
Moreover \(\mathcal D_{n,d}^{\min}\) is stable under every
generalized polarization operator
\begin{equation}\label{dim:fixed-eq:D34}
E_{a,f}=\sum_{i=1}^n z_{ia}f(\partial_{i1},\ldots,\partial_{id}),
\end{equation}
where \(f\) is homogeneous of positive degree.  Thus adjoining all
these operators to differentiation does not enlarge this space.
\end{theorem}

\begin{proof}
Alternating monomials are determinants of distinct monomials in
\(d\) variables.  A basis in the minimum total degree uses every
monomial of degree less than \(s\), together with \(r\) distinct
monomials of degree \(s\).  Consequently, with
\(M_s=\binom{d+s-1}{d-1}\),
\begin{equation}\label{dim:fixed-eq:D35}
\dim\mathcal A_{n,d}^{\min}=\binom{M_s}{r}
      \leq2^{M_s}
      =\exp\bigl(O_d(n^{1-1/d})\bigr).
\end{equation}
Each basis determinant is one of the sources \eqref{dim:fixed-eq:D1}, so the
dimension of its derivative span is bounded by \eqref{dim:fixed-eq:D28}.
The dimension of the sum of these spans is at most the sum of
their dimensions.  Multiplication by the subexponential factor
in \eqref{dim:fixed-eq:D35} preserves the \(n\)-th-root upper bound \(\beta_d\).

For \(P\in\mathcal A_{n,d}^{\min}\) and a homogeneous positive-degree
invariant \(q\),
\[
 q(\partial)P\in S^{\sgn}_{\mu-\deg q}=0.
\]
Invariant differential operators commute with all partial
derivatives, so the same annihilation holds throughout
\(\mathcal D_{n,d}^{\min}\). This proves the harmonic containment.

The comparison with harmonic dimension uses
\(\dim\mathrm{DH}_{n,d}\geq n!\).
Choose any configuration of \(n\) distinct points and its
ordered permutation orbit.  For a homogeneous invariant \(q\),
the polynomial \(q-q(p)\) vanishes on the orbit.  Hence the ideal
of positive invariants lies in the leading-form ideal of that
orbit.  The diagonal coinvariant algebra therefore surjects
onto the associated graded orbit algebra, of length \(n!\).
The constant-coefficient differential pairing identifies its
graded dual with diagonal harmonics, giving the assertion.

Finally, \(E_{a,f}\) is \(S_n\)-equivariant and has degree
\(1-\deg f\). Hence
\[
 E_{a,f}\mathcal A_{n,d}^{\min}
 \subseteq\mathcal A_{n,d}^{\min}\quad(\deg f=1),
 \qquad
 E_{a,f}\mathcal A_{n,d}^{\min}=0\quad(\deg f>1).
\]
To pass from generators to their derivatives, use the commutator
\begin{equation}\label{dim:fixed-eq:D36}
[E_{a,f},\partial_{jc}]
       =-\delta_{ac}f(\partial_{j1},\ldots,\partial_{jd}).
\end{equation}
Commuting \(E_{a,f}\) through an arbitrary product of partial
derivatives proves stability of their entire span.
\end{proof}

At the two full fat-point lengths in \eqref{dim:fixed-eq:D32}, the minimum alternating
component has dimension one, because \(r=M_s\).  Thus those
finite bounds apply to the entire module generated by all
minimum-degree alternants and all operators \eqref{dim:fixed-eq:D34}:
\[
 \begin{aligned}
 \frac{\dim\mathcal D_{70,4}^{\min}}
      {\dim\mathrm{DH}_{70,4}}
 &\leq\frac{\dim\mathcal D_{70,4}^{\min}}{70!}
 <\frac{12181}{10^6},\\
 \frac{\dim\mathcal D_{120,3}^{\min}}
      {\dim\mathrm{DH}_{120,3}}
 &\leq\frac{\dim\mathcal D_{120,3}^{\min}}{120!}
 <\frac{208011}{10^6}.
 \end{aligned}
\]
These are finite harmonic-generation obstructions, in addition
to the uniform asymptotic theorem.

The same mechanism applies to a growing band above the minimum.
For an integer \(e\geq0\), let \(\mathcal A_{n,d}^{\leq\mu+e}\)
be the sum of the homogeneous alternating spaces of degrees at
most \(\mu+e\), and let \(\mathcal D_{n,d}^{\leq\mu+e}\) be
its full span of partial derivatives.

\begin{theorem}[A band above the minimum degree]
\label{dim:near-minimum-band}
Fix \(d\geq3\), and let \(e_n=o(n^{1+1/d})\) be nonnegative
integers.  Then
\begin{equation}\label{dim:fixed-eq:D37}
\limsup_{n\to\infty}
 \left(
  \frac{\dim\mathcal D_{n,d}^{\leq\mu+e_n}}{n!}
 \right)^{1/n}\leq\rho_d<1.
\end{equation}
This derivative space is stable under every operator \eqref{dim:fixed-eq:D34}.
Its intersection with \(\mathrm{DH}_{n,d}\), and the space
generated by harmonic alternants in the same band, satisfy
the same bound.  The whole nonminimal band is not assumed
to be harmonic.
\end{theorem}

\begin{proof}
We bound the number of alternating generators and then the
derivatives of each generator. Put \(M_j=\binom{d+j-1}{d-1}\).
The generating function for alternating monomial determinants is
\begin{equation}\label{dim:fixed-eq:D38}
\prod_{j\geq0}(1+u z^j)^{M_j}.
\end{equation}
Its coefficient of \(u^n z^k\) is the dimension of the
degree-\(k\) alternating space: one chooses \(n\) distinct
monomials and sums their degrees.  Let \(B(n,d,e)\) be the
sum of these dimensions for \(k\leq\mu+e\). The catalecticant
bound~\eqref{dim:fixed-eq:D8} gives
\[
 \dim\mathcal D_{n,d}^{\leq\mu+e}
 \leq B(n,d,e)\,C(dn,\mu+e).
\]
It remains to show that the first factor is subexponential and
that the second has the same rate as at degree \(\mu\).
For \(0<z<1\) and \(u>0\), positivity gives
\[
 B(n,d,e)\leq
 u^{-n}z^{-\mu-e}
       \prod_{j\geq0}(1+u z^j)^{M_j}.
\]
Set \(z=e^{-a/s}\), \(u=e^a\), with \(1\leq a\leq s\),
and write \(\gamma=a/s\).  The identity
\(\mu-ns=-\sum_{j<s}(s-j)M_j\) cancels the contributions
from the fully occupied degrees below \(s\).  We obtain
\begin{equation}\label{dim:fixed-eq:D39}
\begin{split}
 \log B(n,d,e)\leq {}&\frac{ae}{s}
 +\sum_{j<s}M_j\log(1+e^{-\gamma(s-j)})\\
 &+\sum_{j\geq s}M_j\log(1+e^{-\gamma(j-s)}).
 \end{split}
\end{equation}
Using \(\log(1+x)\leq x\) and
\(M_{s+\ell}=O_d((s+\ell)^{d-1})\), the last two sums
are at most
\begin{equation}\label{dim:fixed-eq:D40}
C_d\left(\frac{s^{d-1}}\gamma+\gamma^{-d}\right)
 \leq C'_d\frac n a.
\end{equation}
The first estimate follows by summing a polynomial times
\(e^{-\gamma\ell}\); it is uniform for \(0<\gamma\leq1\).
The second uses \(a\geq1\) and \(n\asymp_d s^d\).

For \(e>0\), choose \(a=\min\{s,\sqrt{ns/e}\}\), and take
\(a=s\) when \(e=0\).  Since \(e=o(ns)\), this is an
admissible choice for all sufficiently large \(n\).  Equations
\eqref{dim:fixed-eq:D39}--\eqref{dim:fixed-eq:D40} give
\begin{equation}\label{dim:fixed-eq:D41}
\log B(n,d,e)
   =O_d\left(s^{d-1}+\sqrt{ne/s}\right)=o(n).
\end{equation}
Since \(e=o(ns)\), the degree \(\mu+e\) has the leading
asymptotic of \(\mu\). Equation~\eqref{dim:fixed-eq:D11}, with
\(D=dn\), therefore gives
\[
 \log\frac{\dim\mathcal D_{n,d}^{\leq\mu+e}}{n!}
 \leq\underbrace{\log B(n,d,e)}_{o(n)}
       +\log\frac{C(dn,\mu+e)}{n!}
 =n\log\rho_d+o(n).
\]
Taking \(n\)-th roots proves~\eqref{dim:fixed-eq:D37}.

Every \(E_{a,f}\) in \eqref{dim:fixed-eq:D34} preserves the alternating band,
because it commutes with \(S_n\) and does not increase total
degree.  The commutator \eqref{dim:fixed-eq:D36} then proves stability of its
derivative span.  Taking subspaces gives the harmonic
assertions.
\end{proof}

\subsection{The two new divisor sources}

For the four-dimensional point algebra with Hilbert function
\((1,4,10,5)\), the determinant has \(20\) rows, degree \(39\),
and \(76\) essential variables. For the three-dimensional point
algebra with Hilbert function \((1,3,6,10,15,21,28,12)\), it has
\(96\) rows, degree \(462\), and \(285\) essential variables.  Their hook multiplicity
polynomials are respectively
\begin{equation}\label{dim:fixed-eq:D25}
(1+uz)^4(1+uz^2)^{10}(1+uz^3)^5
\end{equation}
and
\begin{equation}\label{dim:fixed-eq:D26}
(1+uz)^3(1+uz^2)^6(1+uz^3)^{10}
 (1+uz^4)^{15}(1+uz^5)^{21}(1+uz^6)^{28}
 (1+uz^7)^{12}.
\end{equation}
At these lengths, the two dimension estimates are
\[
 \begin{array}{c|cc}
 (n,d)&C(d(n-1),\mu)/n!&U(n,d)/n!\\\hline
 (20,4)&\approx37.02&\approx4.799\\
 (96,3)&\approx651.91&\approx5.678
 \end{array}
\]
Both estimates exceed one at these particular points. The
finite rank-loss examples are the fat points in
\eqref{dim:fixed-eq:D32}; the estimates become exponentially
small uniformly as \(n\) grows in each fixed dimension.

\paragraph{Historical comparison.}
Haiman explicitly records the \(n=4,d=3\) determinant with
twenty independent derivatives in place of \(4!=24\), and credits
Tesler with the general full-span dimension
\(\binom{2d}{d}\), in Section 5.2 of \cite{gg-HaimanPolygraph}.
Theorems~\ref{dim:fixed-dimension-loss}
and~\ref{dim:all-minimum-alternants} extend the rank-loss
phenomenon to every minimum-degree determinant and to their
combined derivative span, with a uniform exponential bound in
each fixed \(d\geq3\). The planar equality is a consequence
of Haiman's perfect pairing, as proved above.

\subsection{A positive relative band of alternating degrees}
\label{dim:relative-band}

The loss persists when the permitted excess is a fixed positive
fraction of the minimum degree.  Retain the notation
\(\mu=\mu(n,d)\), \(s\), \(\rho_d\), and
\(\mathcal D_{n,d}^{\leq M}\) from the preceding results, and put
\[
 R_d(a)=d\int_0^\infty u^{d-1}
             \log\bigl(1+e^{-a|1-u|}\bigr)\,du
 \qquad(a>0).
\]

\begin{theorem}[A fixed relative alternating band]
\label{dim:relative-alternating-band}
For every fixed \(d\geq3\), \(\varepsilon\geq0\), and \(a>0\),
\begin{equation}
\limsup_{n\to\infty}
 \left(
  \frac{\dim\mathcal D_{n,d}^{\leq
             \lfloor(1+\varepsilon)\mu\rfloor}}{n!}
 \right)^{1/n}
 \leq
 \rho_d(1+\varepsilon)^d
 \exp\left(\frac{ad\varepsilon}{d+1}+R_d(a)\right).
 \label{dim:relative-band-bound}
\end{equation}
The infimum of the right-hand side over \(a>0\) is therefore
also an upper bound.  For every fixed \(d\geq3\), some
\(\varepsilon_d>0\) makes this bound strictly smaller than one.

In dimension three the explicit estimates are
\begin{equation}
\begin{aligned}
 \limsup_{n\to\infty}
 \left(
  \frac{\dim\mathcal D_{n,3}^{\leq
                    \lfloor(126/125)\mu\rfloor}}{n!}
 \right)^{1/n}&<\frac{24}{25},\\[2mm]
 \limsup_{n\to\infty}
 \left(
  \frac{\dim\mathcal D_{n,3}^{\leq
                    \lfloor(101/100)\mu\rfloor}}{n!}
 \right)^{1/n}&<\frac{49}{50}.
 \end{aligned}
 \label{dim:three-relative-band}
\end{equation}
These spaces are stable under every generalized polarization
\(E_{a,f}\) in \eqref{dim:fixed-eq:D34}.  Whenever the right-hand side of \eqref{dim:relative-band-bound} is
smaller than one, differentiation together with all those
polarizations, applied to every alternant in the band, generates
an exponentially vanishing fraction of the dimension of
\(\mathrm{DH}_{n,d}\).  In particular this collection cannot
span the full harmonic space for large \(n\).
The conclusion also holds for the harmonic alternants in the
band and their generated subspace.
\end{theorem}

\begin{figure}[!htbp]
\centering
\resizebox{\linewidth}{!}{%
\begin{tikzpicture}[x=.92cm,y=.73cm,every node/.style={font=\small}]
\node[anchor=west,font=\bfseries,Ink] at (0,16.8)
 {(a) Upper bounds for the alternating-pairing dimension};
\node[anchor=west,font=\footnotesize] at (0,16.1)
 {Curves: $C(d(n-1),\mu)/n!$.  Points: stronger bounds from the point algebra.};

\node[anchor=west,font=\footnotesize,Ink!90] at (0.2,15.35) {log scale};
\fill[Teal!5] (1.65,10.35) rectangle (14.0,14.53065);
\draw[Ink!13] (1.65,10.65968)--(14.0,10.65968);\node[anchor=east,font=\footnotesize] at (1.48,10.65968) {$10^{-25}$};
\draw[Ink!13] (1.65,11.43387)--(14.0,11.43387);\node[anchor=east,font=\footnotesize] at (1.48,11.43387) {$10^{-20}$};
\draw[Ink!13] (1.65,12.20806)--(14.0,12.20806);\node[anchor=east,font=\footnotesize] at (1.48,12.20806) {$10^{-15}$};
\draw[Ink!13] (1.65,12.98226)--(14.0,12.98226);\node[anchor=east,font=\footnotesize] at (1.48,12.98226) {$10^{-10}$};
\draw[Ink!13] (1.65,13.75645)--(14.0,13.75645);\node[anchor=east,font=\footnotesize] at (1.48,13.75645) {$10^{-5}$};
\draw[Ink!13] (1.65,14.53065)--(14.0,14.53065);\node[anchor=east,font=\footnotesize] at (1.48,14.53065) {$10^{0}$};
\node[below,font=\footnotesize] at (1.65000,10.2) {$ 20 $};
\node[below,font=\footnotesize] at (3.85536,10.2) {$ 70 $};
\node[below,font=\footnotesize] at (6.06071,10.2) {$ 120 $};
\node[below,font=\footnotesize] at (8.70714,10.2) {$ 180 $};
\node[below,font=\footnotesize] at (11.35357,10.2) {$ 240 $};
\node[below,font=\footnotesize] at (14.00000,10.2) {$ 300 $};
\draw[Blue,line width=1.1pt] (1.65000,14.71770) -- (1.73821,14.77298) -- (1.82643,14.81266) -- (1.91464,14.83880) -- (2.00286,14.85292) -- (2.09107,14.85621) -- (2.17929,14.84962) -- (2.26750,14.83391) -- (2.35571,14.85481) -- (2.44393,14.89288) -- (2.53214,14.92186) -- (2.62036,14.94246) -- (2.70857,14.95530) -- (2.79679,14.96089) -- (2.88500,14.95968) -- (2.97321,14.95206) -- (3.06143,14.93836) -- (3.14964,14.91889) -- (3.23786,14.89391) -- (3.32607,14.92649) -- (3.41429,14.95290) -- (3.50250,14.97347) -- (3.59071,14.98853) -- (3.67893,14.99833) -- (3.76714,15.00313) -- (3.85536,15.00314) -- (3.94357,14.99856) -- (4.03179,14.98957) -- (4.12000,14.97633) -- (4.20821,14.95901) -- (4.29643,14.93774) -- (4.38464,14.91265) -- (4.47286,14.88387) -- (4.56107,14.90759) -- (4.64929,14.92713) -- (4.73750,14.94264) -- (4.82571,14.95426) -- (4.91393,14.96214) -- (5.00214,14.96639) -- (5.09036,14.96714) -- (5.17857,14.96450) -- (5.26679,14.95856) -- (5.35500,14.94943) -- (5.44321,14.93719) -- (5.53143,14.92193) -- (5.61964,14.90372) -- (5.70786,14.88265) -- (5.79607,14.85877) -- (5.88429,14.83216) -- (5.97250,14.80288) -- (6.06071,14.77100) -- (6.14893,14.78725) -- (6.23714,14.80053) -- (6.32536,14.81093) -- (6.41357,14.81852) -- (6.50179,14.82338) -- (6.59000,14.82557) -- (6.67821,14.82516) -- (6.76643,14.82221) -- (6.85464,14.81678) -- (6.94286,14.80891) -- (7.03107,14.79867) -- (7.11929,14.78611) -- (7.20750,14.77127) -- (7.29571,14.75419) -- (7.38393,14.73493) -- (7.47214,14.71352) -- (7.56036,14.69001) -- (7.64857,14.66443) -- (7.73679,14.63682) -- (7.82500,14.60721) -- (7.91321,14.57564) -- (8.00143,14.54215) -- (8.08964,14.52601) -- (8.17786,14.53486) -- (8.26607,14.54156) -- (8.35429,14.54615) -- (8.44250,14.54868) -- (8.53071,14.54918) -- (8.61893,14.54770) -- (8.70714,14.54427) -- (8.79536,14.53892) -- (8.88357,14.53168) -- (8.97179,14.52260) -- (9.06000,14.51169) -- (9.14821,14.49900) -- (9.23643,14.48454) -- (9.32464,14.46835) -- (9.41286,14.45045) -- (9.50107,14.43088) -- (9.58929,14.40964) -- (9.67750,14.38677) -- (9.76571,14.36230) -- (9.85393,14.33624) -- (9.94214,14.30861) -- (10.03036,14.27945) -- (10.11857,14.24876) -- (10.20679,14.21657) -- (10.29500,14.18289) -- (10.38321,14.14776) -- (10.47143,14.11118) -- (10.55964,14.11568) -- (10.64786,14.11855) -- (10.73607,14.11982) -- (10.82429,14.11950) -- (10.91250,14.11762) -- (11.00071,14.11420) -- (11.08893,14.10927) -- (11.17714,14.10284) -- (11.26536,14.09493) -- (11.35357,14.08557) -- (11.44179,14.07478) -- (11.53000,14.06257) -- (11.61821,14.04896) -- (11.70643,14.03397) -- (11.79464,14.01762) -- (11.88286,13.99992) -- (11.97107,13.98089) -- (12.05929,13.96054) -- (12.14750,13.93890) -- (12.23571,13.91597) -- (12.32393,13.89177) -- (12.41214,13.86632) -- (12.50036,13.83962) -- (12.58857,13.81170) -- (12.67679,13.78256) -- (12.76500,13.75222) -- (12.85321,13.72070) -- (12.94143,13.68800) -- (13.02964,13.65413) -- (13.11786,13.61911) -- (13.20607,13.58295) -- (13.29429,13.54566) -- (13.38250,13.50726) -- (13.47071,13.50708) -- (13.55893,13.50565) -- (13.64714,13.50297) -- (13.73536,13.49907) -- (13.82357,13.49395) -- (13.91179,13.48763) -- (14.00000,13.48012);
\draw[Purple,line width=1.1pt] (1.65000,14.77351) -- (1.73821,14.79982) -- (1.82643,14.81271) -- (1.91464,14.81360) -- (2.00286,14.80359) -- (2.09107,14.78357) -- (2.17929,14.75432) -- (2.26750,14.71645) -- (2.35571,14.73455) -- (2.44393,14.77544) -- (2.53214,14.80741) -- (2.62036,14.83115) -- (2.70857,14.84724) -- (2.79679,14.85618) -- (2.88500,14.85839) -- (2.97321,14.85424) -- (3.06143,14.84408) -- (3.14964,14.82819) -- (3.23786,14.80684) -- (3.32607,14.78026) -- (3.41429,14.74868) -- (3.50250,14.71230) -- (3.59071,14.67129) -- (3.67893,14.62583) -- (3.76714,14.57606) -- (3.85536,14.52213) -- (3.94357,14.54711) -- (4.03179,14.56723) -- (4.12000,14.58271) -- (4.20821,14.59374) -- (4.29643,14.60050) -- (4.38464,14.60314) -- (4.47286,14.60182) -- (4.56107,14.59666) -- (4.64929,14.58780) -- (4.73750,14.57536) -- (4.82571,14.55943) -- (4.91393,14.54014) -- (5.00214,14.51756) -- (5.09036,14.49180) -- (5.17857,14.46293) -- (5.26679,14.43104) -- (5.35500,14.39620) -- (5.44321,14.35848) -- (5.53143,14.31795) -- (5.61964,14.27468) -- (5.70786,14.22872) -- (5.79607,14.18014) -- (5.88429,14.12899) -- (5.97250,14.07533) -- (6.06071,14.01920) -- (6.14893,13.96065) -- (6.23714,13.89973) -- (6.32536,13.83649) -- (6.41357,13.84419) -- (6.50179,13.84915) -- (6.59000,13.85146) -- (6.67821,13.85116) -- (6.76643,13.84833) -- (6.85464,13.84301) -- (6.94286,13.83527) -- (7.03107,13.82516) -- (7.11929,13.81272) -- (7.20750,13.79801) -- (7.29571,13.78106) -- (7.38393,13.76193) -- (7.47214,13.74065) -- (7.56036,13.71726) -- (7.64857,13.69181) -- (7.73679,13.66433) -- (7.82500,13.63485) -- (7.91321,13.60341) -- (8.00143,13.57004) -- (8.08964,13.53477) -- (8.17786,13.49764) -- (8.26607,13.45867) -- (8.35429,13.41790) -- (8.44250,13.37534) -- (8.53071,13.33103) -- (8.61893,13.28498) -- (8.70714,13.23724) -- (8.79536,13.18782) -- (8.88357,13.13674) -- (8.97179,13.08402) -- (9.06000,13.02970) -- (9.14821,12.97378) -- (9.23643,12.91630) -- (9.32464,12.85726) -- (9.41286,12.79670) -- (9.50107,12.73463) -- (9.58929,12.67107) -- (9.67750,12.60603) -- (9.76571,12.53953) -- (9.85393,12.47160) -- (9.94214,12.40225) -- (10.03036,12.33149) -- (10.11857,12.32494) -- (10.20679,12.31675) -- (10.29500,12.30694) -- (10.38321,12.29553) -- (10.47143,12.28255) -- (10.55964,12.26803) -- (10.64786,12.25197) -- (10.73607,12.23441) -- (10.82429,12.21536) -- (10.91250,12.19484) -- (11.00071,12.17288) -- (11.08893,12.14949) -- (11.17714,12.12468) -- (11.26536,12.09848) -- (11.35357,12.07091) -- (11.44179,12.04197) -- (11.53000,12.01169) -- (11.61821,11.98009) -- (11.70643,11.94717) -- (11.79464,11.91295) -- (11.88286,11.87745) -- (11.97107,11.84068) -- (12.05929,11.80266) -- (12.14750,11.76340) -- (12.23571,11.72291) -- (12.32393,11.68120) -- (12.41214,11.63829) -- (12.50036,11.59419) -- (12.58857,11.54892) -- (12.67679,11.50248) -- (12.76500,11.45488) -- (12.85321,11.40614) -- (12.94143,11.35627) -- (13.02964,11.30528) -- (13.11786,11.25318) -- (13.20607,11.19998) -- (13.29429,11.14569) -- (13.38250,11.09032) -- (13.47071,11.03388) -- (13.55893,10.97638) -- (13.64714,10.91783) -- (13.73536,10.85824) -- (13.82357,10.79762) -- (13.91179,10.73597) -- (14.00000,10.67331);
\draw[Ink,line width=.85pt] (1.65,14.53065)--(14.0,14.53065);
\node[anchor=west,font=\footnotesize,fill=white,inner sep=2pt] at (9.3,14.80065) {$d=2$: $\dim G/n!=1$};
\node[anchor=east,Blue,font=\footnotesize] at (14.0,13.16084) {$d=3$: upper bound};
\node[anchor=east,Purple,font=\footnotesize] at (14.0,11.02361) {$d=4$: upper bound};
\fill[Purple] (3.85536,14.23423) circle (2.7pt);
\node[anchor=west,align=left,font=\footnotesize,Purple] at (2.1,9.05) {Point algebra: $n=70$, $d=4$\\$\dim G/n!<0.012181$};
\draw[Purple!80,thin] (3.85536,14.11423)--(3.7000,9.5000);
\fill[Blue] (6.06071,14.42506) circle (2.7pt);
\node[anchor=west,align=left,font=\footnotesize,Blue] at (8.1,9.05) {Point algebra: $n=120$, $d=3$\\$\dim G/n!<0.208011$};
\draw[Blue!80,thin] (6.06071,14.30506)--(9.7000,9.5000);
\draw[Ink!30] (0,8.45)--(15.1,8.45);
\node[anchor=west,font=\bfseries,Ink] at (0,7.9)
 {(b) Bounds for the entire alternating degree band in dimension three};
\node[anchor=west,font=\footnotesize] at (0,7.2)
 {Dimension three: upper bound for $\limsup(\dim\mathcal D/n!)^{1/n}$.};

\fill[Teal!6] (1.65,2.35) rectangle (14.0,4.7500);
\draw[Ink!15] (1.65,2.3500)--(14.0,2.3500);\node[anchor=east,font=\footnotesize] at (1.48,2.3500) {$ 0.6 $};
\draw[Ink!15] (1.65,3.5500)--(14.0,3.5500);\node[anchor=east,font=\footnotesize] at (1.48,3.5500) {$ 0.8 $};
\draw[Ink!15] (1.65,4.7500)--(14.0,4.7500);\node[anchor=east,font=\footnotesize] at (1.48,4.7500) {$ 1.0 $};
\draw[Ink!15] (1.65,5.9500)--(14.0,5.9500);\node[anchor=east,font=\footnotesize] at (1.48,5.9500) {$ 1.2 $};
\draw[Ink!65,line width=.8pt] (1.65,4.75000)--(14.0,4.75000);
\draw[Teal,line width=1.25pt] (1.65000,2.58893) -- (1.71175,2.74429) -- (1.77350,2.81118) -- (1.83525,2.86356) -- (1.89700,2.90844) -- (1.95875,2.94853) -- (2.02050,2.98524) -- (2.08225,3.01938) -- (2.14400,3.05150) -- (2.20575,3.08196) -- (2.26750,3.11105) -- (2.32925,3.13895) -- (2.39100,3.16585) -- (2.45275,3.19185) -- (2.51450,3.21707) -- (2.57625,3.24159) -- (2.63800,3.26548) -- (2.69975,3.28879) -- (2.76150,3.31159) -- (2.82325,3.33390) -- (2.88500,3.35578) -- (2.94675,3.37726) -- (3.00850,3.39836) -- (3.07025,3.41911) -- (3.13200,3.43953) -- (3.19375,3.45965) -- (3.25550,3.47949) -- (3.31725,3.49905) -- (3.37900,3.51836) -- (3.44075,3.53743) -- (3.50250,3.55628) -- (3.56425,3.57491) -- (3.62600,3.59334) -- (3.68775,3.61157) -- (3.74950,3.62962) -- (3.81125,3.64749) -- (3.87300,3.66519) -- (3.93475,3.68274) -- (3.99650,3.70013) -- (4.05825,3.71737) -- (4.12000,3.73447) -- (4.18175,3.75144) -- (4.24350,3.76828) -- (4.30525,3.78499) -- (4.36700,3.80158) -- (4.42875,3.81805) -- (4.49050,3.83441) -- (4.55225,3.85067) -- (4.61400,3.86682) -- (4.67575,3.88287) -- (4.73750,3.89882) -- (4.79925,3.91468) -- (4.86100,3.93045) -- (4.92275,3.94613) -- (4.98450,3.96172) -- (5.04625,3.97723) -- (5.10800,3.99267) -- (5.16975,4.00802) -- (5.23150,4.02330) -- (5.29325,4.03850) -- (5.35500,4.05364) -- (5.41675,4.06870) -- (5.47850,4.08370) -- (5.54025,4.09863) -- (5.60200,4.11350) -- (5.66375,4.12830) -- (5.72550,4.14305) -- (5.78725,4.15774) -- (5.84900,4.17237) -- (5.91075,4.18694) -- (5.97250,4.20146) -- (6.03425,4.21593) -- (6.09600,4.23034) -- (6.15775,4.24471) -- (6.21950,4.25903) -- (6.28125,4.27329) -- (6.34300,4.28751) -- (6.40475,4.30169) -- (6.46650,4.31582) -- (6.52825,4.32991) -- (6.59000,4.34395) -- (6.65175,4.35796) -- (6.71350,4.37192) -- (6.77525,4.38584) -- (6.83700,4.39972) -- (6.89875,4.41357) -- (6.96050,4.42738) -- (7.02225,4.44115) -- (7.08400,4.45489) -- (7.14575,4.46859) -- (7.20750,4.48225) -- (7.26925,4.49589) -- (7.33100,4.50949) -- (7.39275,4.52306) -- (7.45450,4.53659) -- (7.51625,4.55010) -- (7.57800,4.56358) -- (7.63975,4.57702) -- (7.70150,4.59044) -- (7.76325,4.60383) -- (7.82500,4.61719) -- (7.88675,4.63053) -- (7.94850,4.64383) -- (8.01025,4.65711) -- (8.07200,4.67037) -- (8.13375,4.68360) -- (8.19550,4.69681) -- (8.25725,4.70999) -- (8.31900,4.72314) -- (8.38075,4.73628) -- (8.44250,4.74939) -- (8.50425,4.76248) -- (8.56600,4.77555) -- (8.62775,4.78859) -- (8.68950,4.80161) -- (8.75125,4.81462) -- (8.81300,4.82760) -- (8.87475,4.84056) -- (8.93650,4.85351) -- (8.99825,4.86643) -- (9.06000,4.87933) -- (9.12175,4.89222) -- (9.18350,4.90509) -- (9.24525,4.91793) -- (9.30700,4.93077) -- (9.36875,4.94358) -- (9.43050,4.95638) -- (9.49225,4.96916) -- (9.55400,4.98192) -- (9.61575,4.99467) -- (9.67750,5.00740) -- (9.73925,5.02012) -- (9.80100,5.03282) -- (9.86275,5.04550) -- (9.92450,5.05817) -- (9.98625,5.07083) -- (10.04800,5.08347) -- (10.10975,5.09610) -- (10.17150,5.10872) -- (10.23325,5.12132) -- (10.29500,5.13390) -- (10.35675,5.14648) -- (10.41850,5.15904) -- (10.48025,5.17159) -- (10.54200,5.18413) -- (10.60375,5.19665) -- (10.66550,5.20916) -- (10.72725,5.22166) -- (10.78900,5.23415) -- (10.85075,5.24663) -- (10.91250,5.25910) -- (10.97425,5.27155) -- (11.03600,5.28400) -- (11.09775,5.29643) -- (11.15950,5.30886) -- (11.22125,5.32127) -- (11.28300,5.33367) -- (11.34475,5.34607) -- (11.40650,5.35845) -- (11.46825,5.37083) -- (11.53000,5.38319) -- (11.59175,5.39555) -- (11.65350,5.40790) -- (11.71525,5.42023) -- (11.77700,5.43256) -- (11.83875,5.44488) -- (11.90050,5.45720) -- (11.96225,5.46950) -- (12.02400,5.48180) -- (12.08575,5.49409) -- (12.14750,5.50637) -- (12.20925,5.51864) -- (12.27100,5.53091) -- (12.33275,5.54316) -- (12.39450,5.55542) -- (12.45625,5.56766) -- (12.51800,5.57990) -- (12.57975,5.59213) -- (12.64150,5.60435) -- (12.70325,5.61657) -- (12.76500,5.62878) -- (12.82675,5.64098) -- (12.88850,5.65318) -- (12.95025,5.66537) -- (13.01200,5.67756) -- (13.07375,5.68974) -- (13.13550,5.70191) -- (13.19725,5.71408) -- (13.25900,5.72624) -- (13.32075,5.73840) -- (13.38250,5.75055) -- (13.44425,5.76270) -- (13.50600,5.77484) -- (13.56775,5.78698) -- (13.62950,5.79911) -- (13.69125,5.81124) -- (13.75300,5.82336) -- (13.81475,5.83548) -- (13.87650,5.84759) -- (13.93825,5.85970) -- (14.00000,5.87181);
\node[below,font=\footnotesize] at (1.65000,2.23) {$ 0\% $};
\node[below,font=\footnotesize] at (4.73750,2.23) {$ 0.5\% $};
\node[below,font=\footnotesize] at (7.82500,2.23) {$ 1\% $};
\node[below,font=\footnotesize] at (10.91250,2.23) {$ 1.5\% $};
\node[below,font=\footnotesize] at (14.00000,2.23) {$ 2\% $};
\fill[Teal] (7.82500,4.61740) circle (2.8pt);
\draw[Teal!65,densely dashed] (7.82500,2.35)--(7.82500,4.61740);
\node[anchor=west,align=left,font=\footnotesize,Teal,fill=white,inner sep=3pt] at (8.2,3.65) {$1\%$ above the minimum\\has rate $<49/50$.};
\node[anchor=west,font=\footnotesize] at (1.65,1.20) {excess degree divided by the minimum alternating degree};
\node[anchor=west,font=\footnotesize,Teal] at (1.65,.30)
 {Below $1$, these bounds prove exponential loss for the entire derivative span.};
\end{tikzpicture}}
\caption{Bounds for alternating quotients and degree bands.
Panel (a) plots the finite catalecticant bound; the marked fat points
use~\eqref{dim:fixed-eq:D28}. The planar line gives equality.
Panel (b) bounds the derivative span of every alternant of degree
at most $\lfloor(1+\varepsilon)\mu\rfloor$. The curve uses
$R_3(a)\leq31/(6a)+1231/(108a^3)$ with
$a=\sqrt{62/(9\varepsilon)}$ for $\varepsilon>0$, and its limit
at zero. The marked $1\%$ estimate uses $a=26$ and the rational
inequality proving the bound $49/50$.
Values below one prove a deficit; values above one give an
upper bound greater than the comparison dimension. The derivative
spaces are stable under positive-order generalized polarizations.}
\label{dim:pairing-band-figure}
\end{figure}
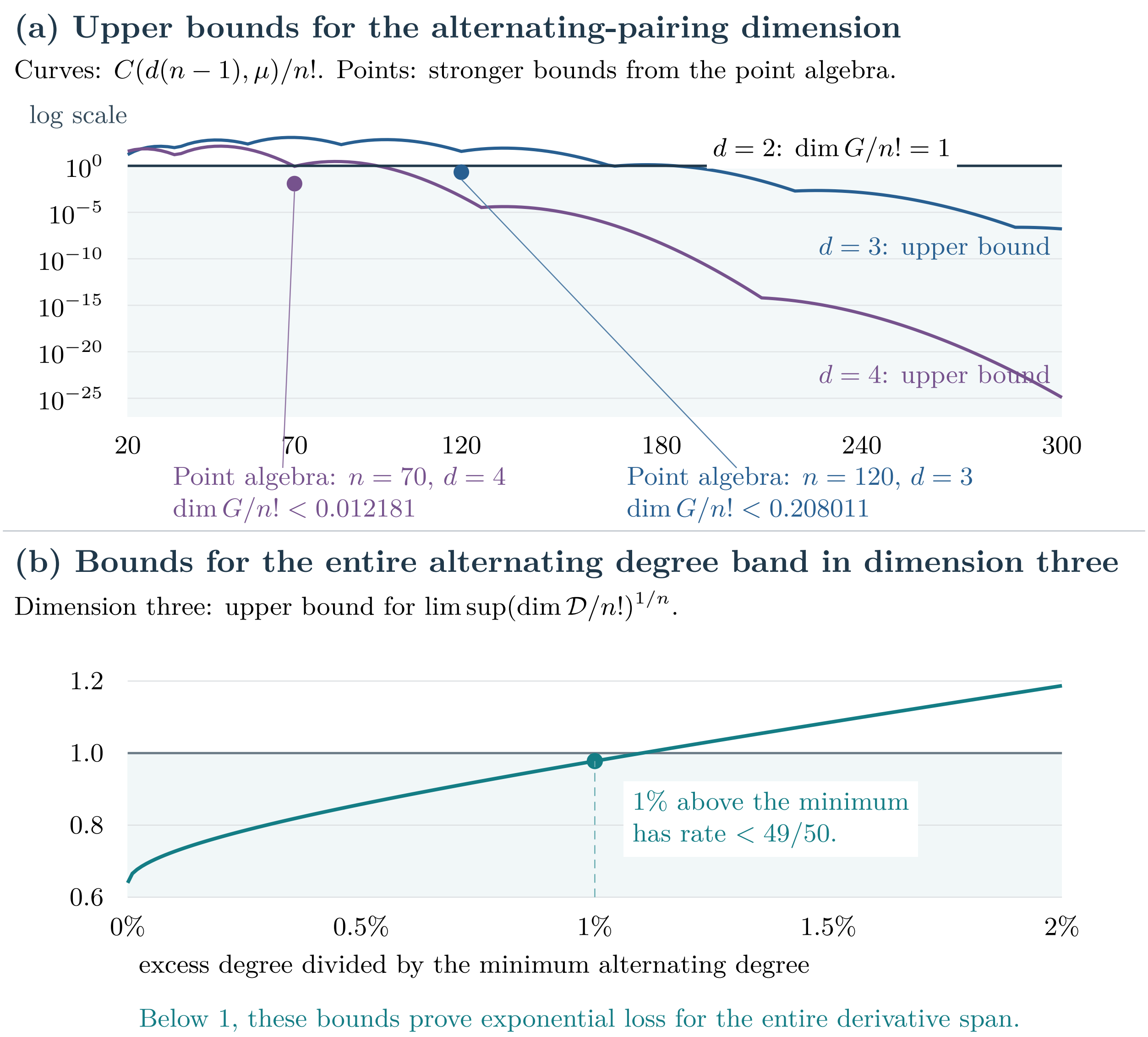

\begin{proof}
As in Theorem~\ref{dim:near-minimum-band}, the estimate has the form
\[
 \dim\mathcal D_{n,d}^{\leq\mu+e}
 \leq B(n,d,e)\,C(dn,\mu+e),
\]
where \(B(n,d,e)\) counts monomial determinants in the band.
A fixed relative excess allows this number to grow exponentially;
we compare its rate with the apolar bound.

Write \(M_k=\binom{d+k-1}{d-1}\). Evaluating the generating
function~\eqref{dim:fixed-eq:D38} at \(u=e^a\), \(z=e^{-a/s}\) gives
\begin{equation}
\log B(n,d,e)\leq
 \frac{ae}{s}
 +\sum_{k\geq0}M_k
       \log\bigl(1+e^{-a|1-k/s|}\bigr).
 \label{dim:relative-band-count}
\end{equation}
Indeed,
\(\log(1+e^{a(1-k/s)})
 =a(1-k/s)_++\log(1+e^{-a|1-k/s|})\),
and
\[
 \sum_{k<s}aM_k(1-k/s)=a(n-\mu/s).
\]
Thus the terms for occupied degrees below \(s\) cancel the
coefficient-extraction factor. The identity also covers a
partly occupied top degree.

For fixed \(a>0\), the monomial-count asymptotics and a Riemann
sum imply
\begin{equation}
\frac1n\sum_{k\geq0}M_k
       \log\bigl(1+e^{-a|1-k/s|}\bigr)
       \longrightarrow R_d(a).
 \label{dim:relative-band-riemann}
\end{equation}
To justify the unbounded interval, use
\(M_k=O_d((1+k)^{d-1})\) and
\(\log(1+e^{-a|1-u|})\leq e^{-a|1-u|}\).
The tails are then bounded by a polynomial times a decaying
exponential, uniformly for large \(s\).  On each bounded
interval, \(M_k=s^{d-1}(k/s)^{d-1}/(d-1)!+O_d(s^{d-2})\)
and \(n=s^d/d!+O_d(s^{d-1})\), which give \eqref{dim:relative-band-riemann}.
Moreover \(\mu/(ns)\to d/(d+1)\).
With \(e=\lfloor\varepsilon\mu\rfloor\), \eqref{dim:relative-band-count} yields
\[
 \limsup_{n\to\infty}\frac1n\log B(n,d,e)
       \leq\frac{ad\varepsilon}{d+1}+R_d(a).
\]

For the second factor, the calculation
\eqref{dim:fixed-eq:D10}--\eqref{dim:fixed-eq:D11}, with \(dn\)
variables and degree \(\lfloor(1+\varepsilon)\mu\rfloor\), gives
\[
 \lim_{n\to\infty}
 \left(\frac{C(dn,\lfloor(1+\varepsilon)\mu\rfloor)}{n!}
 \right)^{1/n}
       =\rho_d(1+\varepsilon)^d.
\]
Multiplying these two rate bounds proves~\eqref{dim:relative-band-bound}.
Using all \(dn\) variables makes the estimate apply to every
alternant in the band.

The elementary estimate
\[
 R_d(a)\leq
 2d\int_0^\infty(1+v)^{d-1}e^{-av}\,dv
 =2d\sum_{j=0}^{d-1}\binom{d-1}{j}\frac{j!}{a^{j+1}}
\]
shows that \(R_d(a)\to0\) as \(a\to\infty\).
Since \(\rho_d<1\), first choose a fixed \(a\) sufficiently
large and then choose \(\varepsilon_d>0\) sufficiently small.

For the explicit three-dimensional estimates, splitting the
integral at one and extending the lower half to the full
positive half-line gives
\[
 R_3(a)\leq
 6\int_0^\infty(1+v^2)\log(1+e^{-av})\,dv
 \leq\frac6a+\frac{12}{a^3}.
\]
Take \(a=32\) and \(\varepsilon=1/125\).  Since
\(\rho_3<16/25\) and \(\log(1+x)\leq x\), the logarithmic
correction to \(\rho_3\) in \eqref{dim:relative-band-bound} is at most
\[
 \frac{27}{125}+\frac3{16}+\frac3{8192}
 <\frac{152}{375}<\log\frac32.
\]
The last inequality follows by retaining the first two terms
of
\(\log(3/2)=2\sum_{j\geq0}(1/5)^{2j+1}/(2j+1)\).
This proves the first assertion in \eqref{dim:three-relative-band}.

For the wider band, the inequality
\(\log(1+x)\leq x-x^2/2+x^3/3\) for \(0\leq x\leq1\)
improves the same integral estimate to
\[
 R_3(a)\leq\frac{31}{6a}+\frac{1231}{108a^3}.
\]
Taking \(a=26\), \(\varepsilon=1/100\), the logarithmic
correction is at most
\[
 \frac9{40}+\frac{31}{156}+\frac{1231}{1898208}
 <2\left(\frac{17}{81}+\frac{(17/81)^3}{3}\right)
 <\log\frac{49}{32}.
\]
Therefore the bound in \eqref{dim:relative-band-bound} is strictly smaller than
\((16/25)(49/32)=49/50\).

The commutator \eqref{dim:fixed-eq:D36} proves polarization stability for every
degree band: these operators preserve alternating polynomials
and do not increase their total degree.  Finally
\(\dim\mathrm{DH}_{n,d}\geq n!\), as established in
Theorem~\ref{dim:all-minimum-alternants}, so the same estimates
hold relative to the full harmonic dimension.  The estimate bounds the full derivative space, so it also
bounds the subspace generated by harmonic inputs in the band.
\end{proof}

\section{Projective dimension of marked polygraphs}
\label{new:polygraph-homology}

Write
\[
 A_\ell=\C[Z(n,\ell;d)],\qquad
 S=\C[x_{ia}:1\le i\le n,\ 1\le a\le d],
 \qquad J=S\cdot S^{\sgn}.
\]
Fix a coordinate alphabet $P_j=\C[x_{1j},\ldots,x_{nj}]$.
A divisor whose image has codimension $c$ produces Koszul homology
through degree $c-1$. We transfer these classes through the Rees
algebra to marked polygraphs, where they determine projective
dimension over $P_j$. The argument applies to modules of arbitrary
rank over $P_j$.

\subsection{Blowup charts and alternating blocks}

Haiman's blowup description, valid in every ambient dimension, gives
\begin{equation}\label{fh:j:rees}
 X_{n,d}=\operatorname{Proj}\bigoplus_{r\ge0}J^rT^r
 \qquad
 \text{\cite[Proposition~3.4.2 and Section~5.2]{gg-HaimanPolygraph}.}
\end{equation}
Here $T$ records the Rees degree. For $a\in J\setminus\{0\}$, the
corresponding affine chart has coordinate ring
\begin{equation}\label{new:rees-colimit}
 S[J/a]=\bigcup_{r\ge0}\frac{J^r}{a^r}
       \cong\varinjlim_r(J^r,\times a).
\end{equation}
The transition map sends $f\in J^r$ to $af\in J^{r+1}$, so the
fraction represented by $f$ remains $f/a^r$.

\begin{lemma}[Alternating blocks]\label{new:alternating-blocks}
For every $r\ge1$, $J^r$ is a direct summand of $A_{nr}$ as an
$S$-module. If $\ell'\ge\ell$, then $A_\ell$ is a direct summand
of $A_{\ell'}$ as an $S$-module.
\end{lemma}

\begin{proof}
Group the $nr$ marks into $r$ blocks of size $n$. Let
$G=(S_n)^r$ permute the marks within these blocks, and let
$e\in\C[G]$ be the idempotent for the tensor product of their sign
characters. Evaluation at $z_{b,i}=x_i$ induces
\begin{equation}\label{fh:j:alternating-summand}
 eA_{nr}\xrightarrow[\ \mathrm{ev}\ ]{\ \sim\ }J^r.
\end{equation}
For an alternating function $F$, let $\Gamma_f$ be the component
graph indexed by $f=(f_1,\ldots,f_r)$. Restriction gives
\[
 F|_{\Gamma_f}=
 \begin{cases}
 0,&\text{some }f_b\text{ repeats a label},\\
 \bigl(\prod_b\sgn(f_b)\bigr)\mathrm{ev}(F),
    &f_b\in S_n\text{ for every }b.
 \end{cases}
\]
In the first case, a transposition fixes the component pointwise
and negates $F$. Thus $\mathrm{ev}(F)=0$ forces $F$ to vanish on
every component, and reducedness gives injectivity.

Every alternating polynomial in one block of marks is a linear
combination of alternants
$\det(m_b(z_a))_{a,b=1}^n$, where the $m_b$ are monomials in $d$
variables. Alternating the monomial basis proves this assertion
directly. Evaluation therefore sends the alternating part of a
marked polynomial to an $S$-linear combination of products of $r$
such determinants. Conversely, these products generate $J^r$:
the same monomial alternation spans $S^{\sgn}$. This proves
surjectivity in \eqref{fh:j:alternating-summand}. The decomposition
$A_{nr}=eA_{nr}\oplus(1-e)A_{nr}$ gives the splitting.

For the second statement, the morphism forgetting the last
$\ell'-\ell$ marks has a section obtained by setting every extra
mark equal to $x_1$. The two induced algebra maps compose to the
identity on $A_\ell$.
\end{proof}

\subsection{Invariant Tor at every punctual divisor}
\label{d3:invariant-punctual}

Alternating charts express invariant Koszul homology in terms of the
trace pairing of the tautological algebra. This description also applies
when the principal Hilbert component is singular along the punctual divisor.

\begin{proposition}[Exterior covariants on an alternating chart]
\label{d3:chart-exterior}
Let \(k\) have characteristic zero, let
\(S=k[x_{ia}:1\le i\le n,\ 1\le a\le d]\), and let \(J\) be its
alternating ideal. For a nonzero monomial alternant
\[
 a=\det(m_b(x_i))_{1\le i,b\le n},
 \qquad C=S[J/a],\qquad C_0=C^{S_n},
\]
there is an isomorphism of exterior \(C_0\)-algebras
\begin{equation}\label{d3:chart-exterior-formula}
 (C\otimes_k\Lambda(k^n))^{S_n}
 \cong \Lambda_{C_0}(\gamma_1,\ldots,\gamma_n),
 \qquad \gamma_b=\sum_i m_b(x_i)\theta_i.
\end{equation}
In particular, the invariant exterior coefficient algebra is free over
\(C_0\), including on the exceptional divisor.

The Cramer coefficients make the span of \(m_1,\ldots,m_n\) into a
finite free point algebra \(\mathcal T\) of rank \(n\) over \(C_0\).
Locally choose a trace-free basis \(1,f_1,\ldots,f_{n-1}\).
For a selected coordinate alphabet \(j\), put
\[
 \bar x_j=\frac1n\sum_i x_{ij},
 \qquad \beta_i=x_{ij}-\bar x_j,\qquad \beta\in\mathcal T.
\]
Then the invariant coordinate Koszul DGA has the integral description
\begin{equation}\label{d3:trace-koszul}
 K(x_{2j}-x_{1j},\ldots,x_{nj}-x_{1j};C)^{S_n}
 \cong
 K_{C_0}\bigl(\Tr(\beta f_1),\ldots,\Tr(\beta f_{n-1})\bigr).
\end{equation}
\end{proposition}

\begin{proof}
Every ratio of two alternating polynomials with denominator \(a\) is
invariant. Thus \(C\) is generated over \(C_0\) by the coordinates
\(x_{ia}\). The exactness of averaging and the super multisymmetric power-sum
argument show that the invariant exterior algebra is generated over
\(C_0\) by the odd sums \(\sum_i f(x_i)\theta_i\). The even power sums
already belong to \(C_0\). Cramer's rule gives
\[
 f(x_i)=\sum_b c_b(f)m_b(x_i),\qquad
 c_b(f)=
 \frac{\det(m_1,\ldots,m_{b-1},f,m_{b+1},\ldots,m_n)}a\in C_0.
\]
These identities hold in \(C\), since they hold in its fraction field,
and they express each odd generator as
\[
 \sum_i f(x_i)\theta_i=\sum_b c_b(f)\gamma_b.
\]
This proves surjectivity in \eqref{d3:chart-exterior-formula} on the
whole chart. Over the fraction field, the evaluation matrix has
determinant \(a\ne0\), so all its exterior powers are invertible.
The kernel is therefore torsion; as a submodule of a free module
over the domain \(C_0\), it is zero.

The same coefficients reduce products of the \(m_b\), constructing
\(\mathcal T\). On the dense reduced locus its trace is
\(\Tr(f)=\sum_i f(x_i)\), so this identity holds over the whole chart.
Since \(\Tr(1)=n\) is a unit, the unit splits off and a trace-free
complement is locally free. Quotienting the permutation exterior
algebra by \(\gamma_0=\sum_i\theta_i\) gives the standard representation.
The differential \(d\theta_i=\beta_i\) kills \(\gamma_0\), and
\[
 d\gamma_{f_b}=\sum_i\beta_i f_b(x_i)=\Tr(\beta f_b).
\]
This proves \eqref{d3:trace-koszul}. Cramer's rule is the step that
keeps the exterior generators defined on the exceptional divisor.
\end{proof}

\begin{corollary}[The global hook covariant bundles]
\label{d3:hook-bundles}
Let \(\rho:X_{n,d}\to H_{n,d}\) be the reduced isospectral cover of
the principal component, let \(\mathcal T\) be its rank-\(n\)
tautological point algebra, and put
\[
 \mathcal T_0=\ker(\Tr:\mathcal T\to\mathcal O_{H_{n,d}}),
 \qquad \mathcal P=\rho_*\mathcal O_{X_{n,d}}.
\]
For every \(0\le k\le n-1\), there is a canonical isomorphism
\begin{equation}\label{d3:global-hook-identity}
 \bigl(\mathcal P\otimes S^{(n-k,1^k)}\bigr)^{S_n}
 \cong\Lambda^k\mathcal T_0.
\end{equation}
Every hook covariant is therefore locally free of rank
\(\binom{n-1}{k}\), in every ambient dimension. If the finite cover
fails to be flat over the principal component, that failure occurs
in its nonhook isotypic summands.
\end{corollary}

\begin{proof}
Evaluation at the ordered points defines the canonical map
\[
 \mathcal T\longrightarrow
 \bigl(\mathcal P\otimes k^n\bigr)^{S_n},
 \qquad f\longmapsto\sum_i f(x_i)\theta_i.
\]
Its exterior powers induce the algebra map of
Proposition~\ref{d3:chart-exterior} on every alternating chart.
That proposition proves the local isomorphisms. Their construction
is independent of the chosen monomial basis, so they glue.
Splitting off the unit gives \(\mathcal T_0\), and
\(\Lambda^k\operatorname{Std}_n=S^{(n-k,1^k)}\) gives
\eqref{d3:global-hook-identity}. The finite characteristic-zero
isotypic decomposition of \(\mathcal P\) then proves the last assertion.
\end{proof}

\begin{theorem}[All-degree invariant obstruction]
\label{d3:punctual-invariant-tor}
Let \(n\geq3\), and suppose the principal component \(H_{n,d}\) contains a divisor whose
general point is supported at one point. Let \(\xi\) be the generic
point of its reduced preimage in \(X_{n,d}\), and put
\[
 B=\mathcal O_{X_{n,d},\xi},\qquad
 P_j=k[x_{1j},\ldots,x_{nj}],\qquad
 \mathfrak q=(x_{ij}-x_{1j}:2\le i\le n).
\]
Then
\begin{equation}\label{d3:local-invariant-all-degrees}
 \Tor_i^{(P_j)_{\mathfrak q}}
 \bigl(B,\kappa(\mathfrak q)\bigr)^{S_n}\ne0
 \qquad(0\le i\le n-2).
\end{equation}
There is an integer \(L\) such that, for every \(\ell\ge L\),
the polygraph satisfies
\begin{equation}\label{d3:marked-invariant-all-degrees}
 \Tor_i^{(P_j)_{\mathfrak q}}
 \bigl((A_\ell)_{\mathfrak q},\kappa(\mathfrak q)\bigr)^{S_n}\ne0
 \qquad(1\le i\le n-2).
\end{equation}
Its global and local projective and flat dimensions over the selected
coordinate alphabet are \(n-2\). The same dimensions hold
for every sufficiently high power \(J^r\).
\end{theorem}

\begin{proof}
The group \(S_n\) fixes \(\xi\) and its residue field pointwise.
The punctual locus is a proper closed subset of the principal
component, so a divisor contained in it is one of its irreducible
components. The connected translation group preserves each such
component. Its support consequently varies over all of \(\mathbb A^d\),
and the alphabet image of \(\xi\) is the generic point of the diagonal
defined by \(\mathfrak q\).
The invariant local ring \(H=B^{S_n}\) is a one-dimensional Noetherian
local domain, as is \(B\). Choose a monomial-alternant chart containing
\(\xi\); these charts cover the blowup because the monomial alternants
generate \(J\). The preceding proposition localizes to \(B,H\).
Invariant denominators are cofinal, by replacing a nonvanishing
denominator \(s\) with \(\prod_{\sigma\in S_n}\sigma(s)\).

The parameter Koszul complex of \((P_j)_{\mathfrak q}\) resolves
its residue field. Averaging and \eqref{d3:trace-koszul} give
\[
 \Tor_i^{(P_j)_{\mathfrak q}}(B,\kappa(\mathfrak q))^{S_n}
 \cong H_iK_H(g_1,\ldots,g_{n-1}),
 \qquad g_b=\Tr(\beta f_b).
\]
All \(g_b\) vanish on the punctual divisor. Their ideal is nonzero:
expressing \(\beta\) in the trace-free basis shows that it contains
\(\Tr(\beta^2)=\sum_i(x_{ij}-\bar x_j)^2\), a nonzero function on the
reduced locus. Hence this ideal is primary to the maximal ideal of \(H\).

Renumber so that \(g_1\ne0\). Because \(H\) is a domain, its
Koszul factor resolves \(H/g_1H\), giving
\[
 K_H(g_1,\ldots,g_{n-1})
 \simeq K_{H/g_1H}(g_2,\ldots,g_{n-1}).
\]
The coefficient ring \(H/g_1H\) is nonzero and Artinian. Adjoining
one remaining entry \(g\) gives the exact sequences
\[
 0\longrightarrow H_i(K')/gH_i(K')
 \longrightarrow H_i(K'\otimes K(g))
 \longrightarrow (0:_{H_{i-1}(K')}g)\longrightarrow0.
\]
For any nonzero homology module \(M\), finite length gives
\[
 M/gM\ne0,\qquad (0:_M g)\ne0:
\]
the first assertion is Nakayama's lemma, and the second follows
from nilpotence of \(g\). Each new Koszul generator therefore
extends the interval of nonzero homology by one. Induction proves
\eqref{d3:local-invariant-all-degrees}.

On the chosen chart the Rees colimit is equivariant in the form
\[
 C=\varinjlim_r
 \bigl(J^r\otimes\sgn^r,\times a\bigr),
 \qquad f\longmapsto f/a^r.
\]
Invariant norm denominators lift each local invariant class to an
invariant chart class. Koszul homology and averaging commute with
the filtered colimit. A representative of a nonzero colimit class
remains nonzero at every later stage. Taking the maximum of the
finitely many starting stages gives \(r_0\) for all degrees
\(1,\ldots,n-2\), with
\[
 \bigl(\Tor_i^{(P_j)_{\mathfrak q}}
 (J^r_{\mathfrak q},\kappa(\mathfrak q))\otimes\sgn^r\bigr)^{S_n}
 \ne0\qquad(r\ge r_0).
\]
Alternating separately in \(r\) blocks of \(n\) marks gives the
equivariant direct summand \(J^r\otimes\sgn^r\) of \(A_{nr}\).
Its label sign twist agrees with the Rees twist. The alternating
projector transfers the classes to the polygraph. Forgetting
extra marks has the equivariant retraction that sets each extra mark
equal to one retained mark, so \(L=nr_0\) proves
\eqref{d3:marked-invariant-all-degrees}.

Centering the selected alphabet splits off a translation variable
and embeds the remaining polygraph module in a free module over a
polynomial ring in \(n-1\) variables. Dimension shifting gives
\(\operatorname{pd}\le(n-1)-1=n-2\), for modules of arbitrary rank.
The top Tor class attains this bound for both projective and flat
dimension. Centering \(J^r\) by the invariant barycenter gives the
same upper bound. The centered decompositions are written out in
the proof of Corollary~\ref{new:seed-pd}.
\end{proof}

\begin{corollary}[Separated blocks and uniform marks]
\label{d3:uniform-mixed}
For a finite collection of punctual-divisor seeds \(m\), choose their
mark constants \(L_m\) in their minimal ambient dimensions. The same
constants work in every larger ambient dimension. At a partition
collision with nonsingleton blocks \(m_1,\ldots,m_b\), if
\(\ell\ge\sum_\nu L_{m_\nu}\), then
\[
 \Tor_i^{(P_j)_{\mathfrak q_\pi}}
 \bigl((A_\ell)_{\mathfrak q_\pi},\kappa(\mathfrak q_\pi)\bigr)^{
 \prod_\nu S_{m_\nu}}\ne0
 \quad\left(1\le i\le\sum_\nu(m_\nu-2)\right).
\]
The local projective and flat dimensions are
\(\sum_\nu(m_\nu-2)\), which is also a global lower bound.
\end{corollary}

\begin{proof}
Ambient projection and its zero section split the polygraph maps,
and coordinate permutations treat every selected alphabet.
After inverting cross-block alphabet differences, the polygraph is
a finite product indexed by assignments of marks to blocks.
Assign at least \(L_{m_\nu}\) marks to block \(\nu\), and apply the
equivariant Koszul K\"unneth formula to that factor. Degree zero
contains a nonzero unit in each factor, so every stated total degree
occurs. Once centered differences kill Koszul homology, each
inverted cross-block difference is a nonzero difference of independent
support parameters. Localization therefore preserves these classes.
The centered upper bounds for every factor and every mark assignment
give the matching local upper bound.
\end{proof}

\begin{proposition}[Principal trace ideals and nonzero products]
\label{d3:principal-trace}
In Theorem~\ref{d3:punctual-invariant-tor}, suppose the trace ideal
\((\Tr(\beta f_b))\subset H\) is principal, generated by \(g\).
Let \(K\) be the coordinate Koszul complex over the localized chart.
Then
\[
 K^{S_n}\cong K(g;H)\otimes_H
 \Lambda_H(\eta_1,\ldots,\eta_{n-2}),
 \qquad
 H_*(K)^{S_n}\cong H/(g)\otimes_H
 \Lambda_H(\eta_1,\ldots,\eta_{n-2}).
\]
The invariant DGA is formal, and the product of the \(n-2\) displayed
degree-one classes is nonzero. Disjoint alternating mark blocks
transfer an exterior subalgebra to the polygraph Tor algebra.
\end{proposition}

\begin{proof}
The row \((g_b/g)\) generates the unit ideal of the local ring \(H\).
An invertible change of odd generators therefore changes the
differential row to \((g,0,\ldots,0)\), giving the displayed DGA.
Its projection modulo \(g\), killing the single nonclosed generator,
is a multiplicative quasi-isomorphism. The disjoint-mark construction
in the proof of Corollary~\ref{next:polygraph-exterior} transfers
these products and preserves their nonzero top class.
\end{proof}

The principal-ideal condition holds when \(H\) is a discrete valuation
ring. The all-degree nonvanishing theorem does not require it.

\subsection{The three-dimensional interval}
\label{d3:interval}

The explicit collision calculation extends to every length
\(96\le m\le105\). In each case the point Hilbert function is
\[
 (1,3,6,10,15,21,28,m-84),
\]
and the collision differential has the maximum rank \(3m-4\).
The following nonzero maximal minors use the prime \(2147483647\).
The certificates retain the integer configurations, all evaluation and
Jacobian matrices, the selected rows and columns, and the four
translation and dilation kernel directions.
\[
\begin{array}{c|r|r@{\qquad}c|r|r}
 m&\operatorname{rank}&\det\bmod p&
 m&\operatorname{rank}&\det\bmod p\\ \hline
 96&284&1915860956&101&299&799176677\\
 97&287&496682692 &102&302&1502039734\\
 98&290&370051734 &103&305&1982918366\\
 99&293&1690310144&104&308&2138877458\\
100&296&1532159477&105&311&842809162
\end{array}
\]
All denominators in the evaluation inversions and Grassmannian
coordinates are units modulo \(p\), so these minors prove the
characteristic-zero ranks.

\begin{corollary}\label{d3:exact-seed-pd}
For every \(96\le m\le105\), there is an integer \(L_m\) such
that, for \(d\ge3\), every coordinate alphabet \(P_j\), and
\(\ell\ge L_m\),
\[
 \operatorname{pd}_{P_j}k[Z(m,\ell;d)]
 =\operatorname{fd}_{P_j}k[Z(m,\ell;d)]=m-2.
\]
At the full coordinate collision, invariant Tor is nonzero in every
degree \(1,\ldots,m-2\). Thus the dimensions attained by these
seeds range from \(94\) through \(103\).
\end{corollary}

\begin{proof}
The maximum collision rank gives a full punctual divisor.
Apply Theorem~\ref{d3:punctual-invariant-tor} and ambient projection
with its zero section.
\end{proof}

\begin{corollary}\label{d3:all-size-bound}
Put \(L=\max_{96\le m\le105}L_m\). For \(n\ge1056\) and
\(\ell\ge L\lceil n/105\rceil\), in every ambient dimension \(d\ge3\),
\[
 \operatorname{pd}_{P_j}k[Z(n,\ell;d)]
 \ge n-2\left\lceil\frac n{105}\right\rceil.
\]
There is a partition collision at which this value is the local
projective and flat dimension. Under its product of block symmetric
groups, invariant Tor is nonzero in every positive degree up to this
value.
\end{corollary}

\begin{proof}
For \(b=\lceil n/105\rceil\), one has \(96b\le n\le105b\) throughout
the stated range. Therefore \(n\) is a sum of \(b\) integers in
\([96,105]\). Apply Corollary~\ref{d3:uniform-mixed}.
\end{proof}

\subsection{The first homogeneous punctual divisors}
\label{d3:contraction-threshold}

Lemma~\ref{dim:automatic-maximal-hilbert} forces the general point
of every homogeneous punctual divisor to have maximal point Hilbert
function. Here such a divisor has general point a translate of a
homogeneous scheme supported at one point.

\begin{proposition}[The first homogeneous punctual divisors]
\label{d3:homogeneous-optimal}\label{d3:uniform-optimal}
Among lengths \(n\ge3\), the first length at which \(H_{n,3}\)
has a homogeneous punctual divisor is \(96\). The corresponding
first length for \(H_{n,4}\) is \(20\).
The familiar double-point divisor occurs at \(n=2\) in both
dimensions.
\end{proposition}

\begin{proof}
By Lemma~\ref{dim:automatic-maximal-hilbert}, the general point
of every homogeneous punctual divisor has maximal point Hilbert
function. Write
\[
 L_s=\binom{d+s-1}{d},\qquad N_s=\binom{d+s-1}{d-1},
 \qquad n=L_s+r,\quad1\le r\le N_s.
\]
Its homogeneous ideal has no elements in degrees below \(s\),
has a codimension-\(r\) subspace of relations in degree \(s\),
and contains every degree above \(s\). Hence the whole homogeneous
locus is \(\operatorname{Gr}(r,N_s)\). Varying the support adds
\(d\) parameters. A necessary condition for such a divisor is
therefore
\[
 r(N_s-r)\ge nd-d-1.
\]
The condition therefore applies to every homogeneous punctual
divisor, including the smoothable locus described in
Corollary~\ref{dim:entire-smoothable-locus}.

In dimension three this becomes
\[
 L_s=\binom{s+2}{3},\qquad N_s=\binom{s+2}{2},
 \qquad
 r(N_s-3-r)-3L_s+4\ge0.
\]
For \(s=2,\ldots,6\), the integral maxima are
\(-6,-14,-20,-20,-8\), respectively.
For \(s=1\), the only solution is \(r=1\), namely the familiar
double-point divisor \(n=2\). For \(s=7\), the inequality becomes
\[
 r(33-r)\ge248,
\]
whose integral solutions are \(12\le r\le21\), giving
\(96\le n\le105\). The displayed maximal minors attain the bound
throughout this interval.

In dimension four, \(s=1\) again permits only \(n=2\).
For \(s=2\), the condition is \(r(6-r)\ge15\), whose maximum
left-hand side is nine. For \(s=3\), it is
\[
 r(16-r)\ge55,
\]
with integral solutions \(5\le r\le11\), namely
\(20\le n\le26\). Theorem~\ref{dim:four-divisors} establishes
existence throughout this interval. These inequalities exclude
every smaller length \(n\ge3\) in the respective dimensions.
\end{proof}

The minimality assertion concerns homogeneous punctual divisors;
it leaves open the first lengths for nonhomogeneous punctual
families and other collision multiplicities.

\begin{proposition}[A fixed unequal-weight contraction]
\label{d3:fixed-weight-obstruction}
For a fixed positive diagonal contraction with unequal weights,
the translated limits of reduced configurations cannot sweep out
a divisor, in any dimension \(d\ge2\) and length \(n\ge2\).
This assertion keeps the weighted coordinate flag fixed.
\end{proposition}

\begin{proof}
The map taking a configuration to its initial ideal is a morphism
on finitely many constructible rank strata. Translation and the
chosen one-parameter dilation preserve each initial ideal. Their
orbits have dimension \(d+1\) on every configuration with at least
two distinct points. Thus a proper rank stratum, of dimension at
most \(nd-1\), gives a translated limit locus of dimension at most
\((nd-1)-(d+1)+d=nd-2\). Only the dense generic stratum could
produce a divisor.

For unequal positive weights, choose \(w_i<w_j\). In addition to
translations and weighted dilation, the transformation
\(x_j\mapsto x_j+c x_i\) preserves the weighted initial ideal,
because it induces the identity on the associated graded polynomial
ring. These \(d+2\) directions are independent on a dense open
configuration set. Its translated limit locus therefore has
dimension at most \(nd-2\); the proper strata were already handled.
\end{proof}

If the weighted flag is allowed to vary, its parameters can offset
the linear shear directions in the last proof. That broader
construction is not excluded by
Proposition~\ref{d3:fixed-weight-obstruction}.

\subsection{A divisor forces the maximal local dimension}

For a partition $\pi$ of $[n]$ into $s$ blocks, let
\[
 \mathfrak q_\pi
 =(x_{ij}-x_{kj}:i,k\text{ belong to the same block of }\pi)
 \subset P_j.
\]
This is a prime of height $c=n-s$. Its residue field
$\kappa(\mathfrak q_\pi)$ is the rational function field of the
independent coordinate values at the $s$ supports.

\begin{theorem}[Divisorial collision]\label{new:divisor-pd}
Let $\mathfrak q=\mathfrak q_\pi$ have height $c\ge2$. Suppose
$X_{n,d}$ has a codimension-one point $\xi$ whose image in
$\Spec(P_j)$ is the generic point of $V(\mathfrak q)$.
There exists $r_0$ such that, for every $r\ge r_0$,
\[
 \operatorname{pd}_{(P_j)_{\mathfrak q}}(J^r)_{\mathfrak q}
 =\operatorname{fd}_{(P_j)_{\mathfrak q}}(J^r)_{\mathfrak q}
 =c-1.
\]
Moreover, for every $\ell\ge nr_0$,
\begin{equation}\label{new:polygraph-pd-local}
 \operatorname{pd}_{(P_j)_{\mathfrak q}}(A_\ell)_{\mathfrak q}
 =\operatorname{fd}_{(P_j)_{\mathfrak q}}(A_\ell)_{\mathfrak q}
 =c-1,
\end{equation}
and
\begin{equation}\label{new:all-positive-tor}
 \Tor_i^{(P_j)_{\mathfrak q}}
 \bigl((A_\ell)_{\mathfrak q},\kappa(\mathfrak q)\bigr)\ne0
 \qquad(1\le i\le c-1).
\end{equation}
The same nonvanishing holds for $(J^r)_{\mathfrak q}$ whenever
$r\ge r_0$.
\end{theorem}

\begin{proof}
Set $R=(P_j)_{\mathfrak q}$ and $B=\mathcal O_{X_{n,d},\xi}$.
The blowup is integral, so $B$ is a one-dimensional Noetherian
local domain. The map $P_j\to B$ is injective, and the preimage
of the maximal ideal of $B$ is $\mathfrak q$. Hence
$\mathfrak qB$ is a nonzero proper ideal. In a one-dimensional
local domain every such ideal is primary to the maximal ideal.

Choose a regular system of parameters $u_1,\ldots,u_c$ of $R$.
The Koszul complex $K(u;R)$ resolves $\kappa(\mathfrak q)$.
Since $u_1$ is nonzero in the domain $B$, the coefficient complex
satisfies
\[
 K(u;B)\simeq K(u_2,\ldots,u_c;B/u_1B).
\]
The module $B/u_1B$ is nonzero and has finite length over $B$.
Its nonzero socle is annihilated by the remaining $u_i$, so
\[
 \Tor_{c-1}^{R}(B,\kappa(\mathfrak q))\ne0.
\]
The finite-length induction in the proof of
Theorem~\ref{d3:punctual-invariant-tor} applies to
$K(u_2,\ldots,u_c;B/u_1B)$ without taking invariants. It gives
\[
 H_i(u;B)\ne0\qquad(0\le i\le c-1).
\]

Choose a chart $S[J/a]$ containing $\xi$. Localization is exact
and takes its Koszul complex to $K(u;B)$. Thus, for each
$1\le i\le c-1$,
\[
 \Tor_i^R(S[J/a]_{\mathfrak q},\kappa(\mathfrak q))\ne0.
\]
Using \eqref{new:rees-colimit} and commutation of Tor with filtered
colimits gives
\[
 \Tor_i^R(S[J/a]_{\mathfrak q},\kappa(\mathfrak q))
 \cong
 \varinjlim_r
 \Tor_i^R((J^r)_{\mathfrak q},\kappa(\mathfrak q)).
\]
A representative of a nonzero colimit class remains nonzero at
every later stage; otherwise its colimit image would vanish.
Taking the maximum of the finitely many starting indices proves
nonvanishing in all these degrees for every $r\ge r_0$.

For the upper bound, restriction to the component graphs embeds
\begin{equation}\label{new:polygraph-free-embedding}
 A_\ell\hookrightarrow
 \prod_{f:[\ell]\to[n]}S.
\end{equation}
The target is a finite direct sum of copies of $S$, and
$S_{\mathfrak q}$ is a free $R$-module, with basis the monomials
in the other coordinate alphabets. Both $(A_\ell)_{\mathfrak q}$
and $(J^r)_{\mathfrak q}$ are therefore submodules of free
$R$-modules.

A regular local ring of dimension $c$ has global dimension $c$.
For a submodule $M$ of a free module $F$, dimension shifting gives
\[
 \Ext_R^c(M,N)\cong\Ext_R^{c+1}(F/M,N)=0
 \quad\text{for every }N,
 \qquad \operatorname{pd}_R M\le c-1.
\]
This applies to arbitrary $R$-modules $M$.
The nonzero Tor in degree $c-1$ proves equality for both
projective and flat dimension of $(J^r)_{\mathfrak q}$.
Lemma~\ref{new:alternating-blocks} transfers all the Tor classes
of $J^{r_0}$ to $A_{nr_0}$ and then to every $A_\ell$ with
$\ell\ge nr_0$. The same upper bound proves
\eqref{new:polygraph-pd-local} and \eqref{new:all-positive-tor}.
The standard homological facts used here hold for arbitrary
modules \cite[Tags 00O7, 0BNF and 062F]{new:StacksHomological}.
\end{proof}

\begin{corollary}[Projective dimensions eight, ten and thirteen]
\label{new:seed-pd}
For $n\ge10$ and $d\ge5$, at a coordinate collision with one
block of size ten and all other blocks singletons, every
sufficiently large $A_\ell$ has projective and flat dimension
eight over $(P_j)_{\mathfrak q}$, with nonzero Tor in every
degree $1,\ldots,8$. For $n\ge12$ and $d\ge5$, a block
of size twelve gives dimension ten and nonzero Tor in
degrees $1,\ldots,10$. For $n\ge15$ and $d\ge11$, the corresponding
statement for one block of size fifteen has dimension thirteen
and nonzero Tor in every degree $1,\ldots,13$.

At the seed lengths the global dimensions have the same values:
for all sufficiently large $\ell$,
\begin{align}
 \operatorname{pd}_{P_j}\C[Z(10,\ell;d)]
 &=\operatorname{fd}_{P_j}\C[Z(10,\ell;d)]=8
 &&(d\ge5),\label{new:global-pd-eight}\\
 \operatorname{pd}_{P_j}\C[Z(12,\ell;d)]
 &=\operatorname{fd}_{P_j}\C[Z(12,\ell;d)]=10
 &&(d\ge5),\label{next:global-pd-ten}\\
 \operatorname{pd}_{P_j}\C[Z(15,\ell;d)]
 &=\operatorname{fd}_{P_j}\C[Z(15,\ell;d)]=13
 &&(d\ge11).\label{new:global-pd-thirteen}
\end{align}
At these respective lengths, the same global equalities hold
for $J^r$ for all sufficiently large $r$.
All statements hold for every coordinate alphabet, with one
common sufficiently large threshold if desired.
\end{corollary}

\begin{proof}
The ten-point divisor of Lemma~\ref{new:divisor} and its
graph-bundle and separated-support extensions dominate the
specified coordinate collision locus: the common support
coordinate and the coordinates of the simple supports vary
independently. Its prime has height nine. The twelve-point
divisor of Theorem~\ref{next:twelve-smoothability} gives height
eleven by the same argument. The fifteen-point
divisor of Section~\ref{new:fifteen} gives the same conclusion
with height fourteen. Apply Theorem~\ref{new:divisor-pd}.

For the global upper bound at length $n$, replace the selected
coordinates by
\[
 x_{1j},\quad x_{ij}-x_{1j}\ (i\ge2),\quad
 z_{aj}-x_{1j}\ (1\le a\le\ell).
\]
All component equations are independent of $x_{1j}$. Thus
$A_\ell=\C[x_{1j}]\otimes_\C M_\ell$, where $M_\ell$ is the
centered polygraph ring over the polynomial base
\[
 P_j^0=\C[x_{2j}-x_{1j},\ldots,x_{nj}-x_{1j}].
\]
The centered version of \eqref{new:polygraph-free-embedding}
embeds $M_\ell$ into a free $P_j^0$-module. Since $P_j^0$ has
global dimension $n-1$, its submodule has projective dimension
at most $n-2$. Polynomial extension gives the same upper bound
over $P_j$. At $n=10,12,15$, the localized lower bounds
attain it.

For the alternating powers, use the invariant barycenter
$p_j=n^{-1}\sum_i x_{ij}$ and the centered coordinates
$y_{ij}=x_{ij}-p_j$. Then $S=\C[p_j]\otimes_\C S^0$, where
$S_n$ acts trivially on $p_j$. Consequently
$S^{\sgn}=\C[p_j]\otimes_\C(S^0)^{\sgn}$ and
$J=\C[p_j]\otimes_\C J^0$, with
$J^0=S^0(S^0)^{\sgn}$. Thus
$J^r=\C[p_j]\otimes_\C(J^0)^r$.
The ideal $(J^0)^r\subset S^0$ embeds into a free module over
the centered coordinate alphabet, a polynomial ring in $n-1$
variables. Dimension shifting gives projective dimension at
most $n-2$, and polynomial extension preserves that upper bound.
The localized lower bounds for $J^r$ in
Theorem~\ref{new:divisor-pd} attain it at $n=10,12,15$.
Finally, take the maximum of the thresholds over the finitely
many coordinate alphabets.
\end{proof}

\subsection{Separated collisions add their homological dimensions}

\begin{proposition}[Factorization by mark allocations]
\label{new:polygraph-clusters}
Let $\pi=\{B_1,\ldots,B_s\}$ be a partition of $[n]$, with
$m_b=|B_b|$. Let $T_\pi\subset\Spec(P_j)$ be the open set on
which every cross-block difference $x_{ij}-x_{kj}$ is invertible.
For $\alpha:[\ell]\to[s]$, put
$\ell_b=|\alpha^{-1}(b)|$. There is an isomorphism of algebras
\begin{equation}\label{new:cluster-algebra}
 A_\ell\otimes_{P_j}\mathcal O(T_\pi)
 \cong
 \prod_{\alpha:[\ell]\to[s]}
 \left(\bigotimes_{b=1}^s\C[Z(m_b,\ell_b;d)]\right)
 \otimes_{P_j}\mathcal O(T_\pi).
\end{equation}
The tensor factors use disjoint label and mark variables.
Consequently, at $\mathfrak q=\mathfrak q_\pi$,
\begin{equation}\label{new:cluster-upper}
 \operatorname{pd}_{(P_j)_{\mathfrak q}}(A_\ell)_{\mathfrak q}
 \le\sum_{b=1}^s\max(m_b-2,0)
 \qquad\text{for every }\ell.
\end{equation}
\end{proposition}

\begin{proof}
Over $\mathcal O(T_\pi)$ the polynomials
\[
 f_b(U)=\prod_{i\in B_b}(U-x_{ij})
\]
are pairwise comaximal, because their resultants are products
of the inverted cross-block differences. Each selected marked
coordinate satisfies $\prod_bf_b(z_{aj})=0$ in $A_\ell$.
The Chinese remainder idempotents therefore specify, for each
mark, the block containing its assigned label. Multiplying these
idempotents gives the factors indexed by $\alpha$.

In the factor for $\alpha$, the component graphs correspond to
the assignments respecting that allocation. The choices
within distinct blocks are independent. Their union is the
product of the corresponding cluster polygraphs; these products
are reduced over $\C$. This proves \eqref{new:cluster-algebra}.
A block with no marks contributes its ordinary label-coordinate
polynomial ring.

By the centered argument in Corollary~\ref{new:seed-pd}, a block
of size $m\ge2$ has projective dimension at most $m-2$ over its
coordinate alphabet for every number of marks. A singleton
block contributes zero. Tensoring projective resolutions over
$\C$ gives a projective resolution over the tensor product of
the independent bases: tensor products of their projective
modules remain projective, and exactness follows because
$\C$ is a field. Localizing and taking the finite product in
\eqref{new:cluster-algebra} proves \eqref{new:cluster-upper}.
\end{proof}

\Needspace{25\baselineskip}
\begin{theorem}[Additive growth on collision strata]
\label{new:mixed-pd}
There exist positive integers $L_{10},L_{12},L_{15}$ with the
following property. Suppose $\pi$ has $a$ blocks of size ten, $b$ blocks of size twelve,
$c$ blocks of size fifteen, and any number of singleton blocks. Assume $d\ge5$,
and assume $d\ge11$ if $c>0$. Then, for every
\begin{equation}\label{new:uniform-mark-threshold}
 \ell\ge aL_{10}+bL_{12}+cL_{15},
\end{equation}
\begin{equation}\label{new:mixed-pd-formula}
 \operatorname{pd}_{(P_j)_{\mathfrak q_\pi}}
       (A_\ell)_{\mathfrak q_\pi}
 =\operatorname{fd}_{(P_j)_{\mathfrak q_\pi}}
       (A_\ell)_{\mathfrak q_\pi}
 =8a+10b+13c.
\end{equation}
Every Tor group against $\kappa(\mathfrak q_\pi)$ in degrees
$1,\ldots,8a+10b+13c$ is nonzero. In particular, for $d\ge5$,
\begin{equation}\label{new:global-linear-pd}
 \operatorname{fd}_{P_j} A_\ell
 \ge8\lfloor n/10\rfloor
 \qquad\bigl(\ell\ge L_{10}\lfloor n/10\rfloor\bigr),
\end{equation}
and the same lower bound holds for projective dimension.
For $d\ge5$, with $c=0$ unless $d\ge11$, each $a,b,c$ with $10a+12b+15c\le n$ gives the
global lower bound $8a+10b+13c$ whenever
$\ell\ge aL_{10}+bL_{12}+cL_{15}$. In particular the eventual lower
bound is $\max\{8a+10b+13c:10a+12b+15c\le n\}$.

The constants are independent of $n$, $d$ and the selected
alphabet. Thus one ten-point block gives local dimension
eight for every $n\ge10$, $d\ge5$ and $\ell\ge L_{10}$;
one twelve-point block gives ten for every $n\ge12$,
$d\ge5$ and $\ell\ge L_{12}$; one fifteen-point block gives thirteen for every $n\ge15$,
$d\ge11$ and $\ell\ge L_{15}$.
\end{theorem}

\begin{proof}
Choose $L_{10}$ from the seed $(n,d)=(10,5)$ in
Corollary~\ref{new:seed-pd}, so all eight positive Tor degrees
are nonzero whenever $\ell\ge L_{10}$. Choose $L_{12}$ and $L_{15}$ in the
same way from $(12,5)$ and $(15,11)$. Coordinate permutation makes these
choices uniform over the alphabets of each seed.

For $d\ge d_0$, choose $d_0$ ambient coordinates containing
the selected alphabet. Projection to them gives
$Z(m,\ell;d)\to Z(m,\ell;d_0)$, and setting all omitted
coordinates of both labels and marks equal to zero gives a
section. Thus the smaller polygraph ring is a direct summand
over the selected alphabet. Every seed Tor class persists
in every larger ambient dimension, with the same mark
threshold. The centered upper bound still gives $m-2$.

For every $\ell$ satisfying \eqref{new:uniform-mark-threshold},
choose an allocation $\alpha$ giving at least $L_{10}$ marks
to each ten-point block, at least $L_{12}$ to each twelve-point
block, and at least $L_{15}$ to each
fifteen-point block. Extra marks can be assigned to any block.

Write $s_b$ for the selected support coordinate in block $b$
and $u_b$ for its $m_b-1$ coordinate differences. Centering that
block writes its ring as $\C[s_b]\otimes_\C M_b$. Let
\[
 H_{b,i}=H_i(K(u_b;M_b)).
\]
Each $H_{b,i}$ is a $\C$-vector space, since the differences
annihilate Koszul homology. Corollary~\ref{new:seed-pd} gives
$H_{b,i}\ne0$ for all $0\le i\le m_b-2$ in a nonsingleton
block; nonvanishing after adjoining the independent support
coordinate implies nonvanishing before that field extension.
For a singleton, $H_{b,0}\ne0$ and there are no difference
parameters.

The Koszul complex on all the block differences is the tensor
product of the block complexes, so over $\C$ its homology is
\[
 \bigoplus_{\sum_b i_b=i}\ \bigotimes_b H_{b,i_b}.
\]
The support parameters are independent of these centered complexes.
After localization at $\mathfrak q_\pi$, the homology becomes
\[
 \C(s_1,\ldots,s_s)\otimes_\C
 \bigoplus_{\sum_b i_b=i}\ \bigotimes_b H_{b,i_b}.
\]
Indeed, the block differences act by zero on homology, while
an inverted cross-block difference acts as $s_b-s_{b'}$.
Every element of $P_j\setminus\mathfrak q_\pi$ likewise becomes
a nonzero function of the support parameters. Localization is
therefore a field extension on homology and preserves each
nonzero tensor class.

The K\"unneth formula now gives nonzero Tor in every degree
between zero and $8a+10b+13c$ for this allocation. It is a direct
summand of the polygraph algebra by
\eqref{new:cluster-algebra}. The upper bound
\eqref{new:cluster-upper} proves \eqref{new:mixed-pd-formula}.
Taking $a=\lfloor n/10\rfloor$, $b=c=0$, and the remaining labels
as singletons proves \eqref{new:global-linear-pd}, because
localization cannot increase flat or projective dimension.
The mixed bound follows by choosing the maximizing partition.
\end{proof}

\subsection{An invariant class in internal degree twenty-one}

The last nonzero graded piece of the normalization defect produces
the highest local Tor class. The connecting homomorphism gives
a formula valid for any filtered finite field extension.

\begin{proposition}[The last defect layer]\label{new:defect-tor}
Let $E/F$ be a finite field extension and let
\[
 \mathcal R=\bigoplus_{r\ge0}t^rV^r\subset E[t],
 \qquad 1\in V,\qquad V^h=E,\quad V^{h-1}\ne E,
\]
where $V^r$ denotes the span of products of $r$ elements of
$V$. Suppose $h\ge2$ and $u_i=t\beta_i\in\mathcal R_1$
for $1\le i\le c$, with $\beta_1\ne0$. Give each exterior
generator of $K(u;\mathcal R)$ internal degree one.
Then there is a canonical isomorphism
\begin{equation}\label{new:last-layer-tor}
 H_{c-1}(u;\mathcal R)_{h+c-1}
 \cong
 (E/V^{h-1})\otimes_F\bigwedge^c F^c.
\end{equation}
It is equivariant for group actions fixing $F$ and $t$,
preserving $V$, and acting on the Koszul generator space $F^c$
so that $e_i\mapsto u_i$ is equivariant.
For $w\in t^{h-1}E$ representing a nonzero defect class, the
corresponding cycle is
\begin{equation}\label{new:explicit-defect-cycle}
 z(w)=\sum_{i=1}^c(-1)^{i-1}u_iw\,
 e_1\wedge\cdots\widehat e_i\cdots\wedge e_c.
\end{equation}
It remains nonzero after localization at the cone vertex.
\end{proposition}

\begin{proof}
Put $C=E[t]$ and $D=C/\mathcal R$. Since $1\in V$ and
$V^h=E$, the last nonzero degree of $D$ is $h-1$, where
$D_{h-1}=t^{h-1}(E/V^{h-1})$. Every $u_i$ annihilates that
layer. Thus the top Koszul homology in total degree $h+c-1$ is
\[
 H_c(u;D)_{h+c-1}
 =D_{h-1}\otimes_F\bigwedge^c F^c.
\]
Because $\beta_1$ is a unit in $E$, a change of Koszul basis
over $E$ transforms the sequence $(t\beta_1,\ldots,t\beta_c)$
into $(t,0,\ldots,0)$. It follows that $H_i(u;C)$ is concentrated in internal degree
$i$. Since $h\ge2$, the long exact sequence for
$0\to\mathcal R\to C\to D\to0$ has the segment
\[
 0\longrightarrow H_c(u;D)_{h+c-1}
 \xrightarrow{\ \partial\ }
 H_{c-1}(u;\mathcal R)_{h+c-1}\longrightarrow0.
\]
This proves \eqref{new:last-layer-tor}. The degree separation is
what prevents the defect class from being killed by normalization
homology. The connecting map is equivariant.

Its explicit formula is \eqref{new:explicit-defect-cycle}.
Each coefficient $u_iw$ lies in $\mathcal R_h$, and $d^2=0$
shows that it is a cycle. If it were the boundary of
$v\,e_1\wedge\cdots\wedge e_c$, then
$u_i(w-v)=0$ for every $i$. Since $u_1$ is a nonzerodivisor
on $C$, this would force $w=v$, contradicting
$w\notin\mathcal R$.

Finally, $E[t]=E[u_1]$ is finite over $F[u_1]$ because
$t=u_1/\beta_1$. Its $F[u_1]$-submodule $\mathcal R$ is
therefore finite as well. Hence
$\mathcal R/u_1\mathcal R$ is finite-dimensional over $F$.
The Koszul homology is supported at the cone vertex, so
localizing there preserves the nonzero class.
\end{proof}

\begin{corollary}[The invariant top Tor class]
\label{new:degree-twenty-one}
For the ten-point cone of Theorem~\ref{new:local-cone}, let
$u_i=x_{i+1,j}-x_{1j}$, $1\le i\le9$, be the differences in
one coordinate alphabet. Then
\begin{equation}\label{new:degree-twenty-one-invariant}
 \dim_F H_8(u;\mathcal R)_{21}^{S_{10}}=1.
\end{equation}
This class survives in the Koszul homology of the
local ring at the divisor. After adjoining the moving support
coordinates, it is a nonzero class in
$\Tor_8^{(P_j)_{\mathfrak q}}
(B,\kappa(\mathfrak q))$, where $\mathfrak q$ is the tenfold
coordinate-collision prime and $B$ is the local ring of the
ordered cover over the moving divisor.
\end{corollary}

\begin{proof}
Theorem~\ref{new:orbit-conductor} gives $h=13$.
The field $E$ is an $S_{10}$-torsor algebra over $F$, so its
sign representation occurs once. The filtration piece
$V^{12}$ contains no sign representation, since the least
degree of an alternant in ten blocks of five variables is
$0+5+2\cdot4=13$. The differences span the standard
nine-dimensional representation $U$ of $S_{10}$, and
$\bigwedge^9U$ is its sign character. Formula
\eqref{new:last-layer-tor} consequently identifies the invariant
subspace with the unique sign line in $E/V^{12}$, proving
\eqref{new:degree-twenty-one-invariant}.

Explicitly, let
$\Delta(x)=\det(m_b(x_i))_{i,b=1}^{10}$ be a nonzero
homogeneous alternating polynomial formed from the constant
monomial, five linear monomials and four quadratic monomials.
It has degree thirteen. Then
$w=\Delta(x)/t$ has internal degree twelve and lies outside
$\mathcal R$. Formula \eqref{new:explicit-defect-cycle}
is the required invariant cycle: its coefficients have degree
thirteen and its exterior factor has degree eight.
Theorem~\ref{new:local-cone} identifies the local ring
with the localization of $\mathcal R$ at its vertex.
Adjoining the independent support coordinates is a faithful
field extension. The nine differences are a regular system of
parameters of $(P_j)_{\mathfrak q}$, so the resulting Koszul
complex computes the stated Tor group.
\end{proof}

\begin{corollary}[Invariant classes in the marked polygraphs]
\label{new:equivariant-transfer}
After increasing the constant $L_{10}$ in
Theorem~\ref{new:mixed-pd} if necessary, the following holds
for every $n\ge10$, $d\ge5$ and $\ell\ge L_{10}$. Let
$\mathfrak q$ be the coordinate-collision prime with one
ten-point block and all remaining blocks singletons. Let
$G\simeq S_{10}$ permute the labels in that block. Then
\begin{equation}\label{new:invariant-polygraph-tor}
 \Tor_8^{(P_j)_{\mathfrak q}}
 \bigl((A_\ell)_{\mathfrak q},\kappa(\mathfrak q)\bigr)^G
 \ne0.
\end{equation}
At length ten, for every $d\ge5$ and all sufficiently large
$r$, the label representation
\[
 \Tor_8^{(P_j)_{\mathfrak q}}
 \bigl((J^r)_{\mathfrak q},\kappa(\mathfrak q)\bigr)
\]
contains the character $\sgn^r$.
\end{corollary}

\begin{proof}
First take $(n,d)=(10,5)$ and $G=S_{10}$, and put
$R=(P_j)_{\mathfrak q}$ and $\kappa=\kappa(\mathfrak q)$.
The group preserves $\mathfrak q$ and acts trivially on
$\kappa$, so these Tor groups have their natural
$\kappa$-linear $G$-actions.
Choose a nonzero degree-thirteen alternant $a$ and write $JB=gB$;
the blowup makes this ideal invertible. The normalization valuation
satisfies
\[
 v(a)=v(JB)=v(g)=13,\qquad a/g\in B,\qquad v(a/g)=0.
\]
The maximal ideal of $B$ is the inverse image of the normalization
maximal ideal, so $a/g$ is a unit. Hence $JB=aB$: the chart
$S[J/a]$ contains the generic divisor and is $G$-stable because
$\sigma(a)=\sgn(\sigma)a$ for $\sigma\in G$.

Set $Q=S[J/a]_{\mathfrak q}$, and let $\mathfrak p$ be the
prime whose localization is $B$. It is $G$-stable. Invariant
denominators are cofinal for this localization: if
$s\in Q\setminus\mathfrak p$, then
\[
 N(s)=\prod_{g\in G}g(s)\in Q^G\setminus\mathfrak p
\]
is a multiple of $s$. Koszul homology commutes with
localization. By Corollary~\ref{new:degree-twenty-one},
$\Tor_8^R(B,\kappa)$ has a nonzero invariant class $z$. Writing
it with a norm denominator gives
\[
 z=\frac{h}{N(s)}
   =\frac{\operatorname{Av}_G(h)}{N(s)},
 \qquad
 \operatorname{Av}_G(h)=\frac1{|G|}\sum_{\sigma\in G}\sigma h.
\]
Its averaged numerator is a class on $Q$ whose localization is
nonzero. Consequently
\[
 \Tor_8^R(Q,\kappa)^G\ne0.
\]

The equivariant version of \eqref{new:rees-colimit} is
\begin{equation}\label{new:equivariant-rees}
 Q\cong\varinjlim_r
 \bigl((J^r)_{\mathfrak q}\otimes\sgn^r,\times a\bigr).
\end{equation}
Here the sign twist is forced by division by $a^r$:
\[
 \sigma(f/a^r)=\sgn(\sigma)^r\,\sigma(f)/a^r.
\]
Tor and averaging commute with the filtered colimit. A representative
of a nonzero colimit class stays nonzero at every later stage,
so for all sufficiently large $r$,
\[
 \bigl(
 \Tor_8^R((J^r)_{\mathfrak q},\kappa)
 \otimes\sgn^r
 \bigr)^G\ne0.
\]
This proves the assertion about the character of the
alternating-power Tor group.

For the label action, the evaluation isomorphism in
\eqref{fh:j:alternating-summand} satisfies
\[
 \mathrm{ev}(g f)=\sgn(g)^r\,g\bigl(\mathrm{ev}(f)\bigr).
\]
The sign arises by reordering the marks in each of the $r$
blocks back to the identity assignment. Thus its equivariant
form is $eA_{10r}\cong J^r\otimes\sgn^r$.
The matching twists therefore carry the invariant chart class to
invariant polygraph Tor. For $\ell\ge10r_0$, forgetting extra
marks has the $G$-equivariant algebra section that sets each
extra mark equal to a retained mark. The induced splitting preserves
the class for every such $\ell$. Increase
$L_{10}$ to include this seed threshold.

Setting additional ambient coordinates equal to zero gives
$G$-equivariant retractions both for the polygraph rings and
for the alternating powers. For the latter, specialization
sends $J_d^r$ onto $J_5^r$ and is split by inclusion; an
alternant specializes to an alternant, and every
five-coordinate alternant already belongs to $J_d$.
Thus the assertions persist for every $d\ge5$.

Finally, for arbitrary $n\ge10$, use the factor in
\eqref{new:cluster-algebra} which assigns every mark to the
ten-point block. It is $G$-stable and is the tensor product
of the ten-point polygraph ring with polynomial rings for
the singleton labels. Centering the ten-point block by its
barycenter makes its support coordinate $G$-invariant.
The localization argument in Theorem~\ref{new:mixed-pd}
is then a faithful extension by $G$-fixed support parameters.
It preserves the invariant Tor class, proving
\eqref{new:invariant-polygraph-tor} with the same $L_{10}$
for all $n$ and $d$ in the stated ranges.
\end{proof}

\begin{figure}[!htbp]
\centering
\begin{tikzpicture}[x=1cm,y=1cm,font=\small]
\path[use as bounding box] (0,-1.1) rectangle (16.2,9.6);
\node[anchor=west,font=\bfseries\small,Ink] at (0,9.2) {(a) The ten-point divisor has eight excess fibre dimensions};
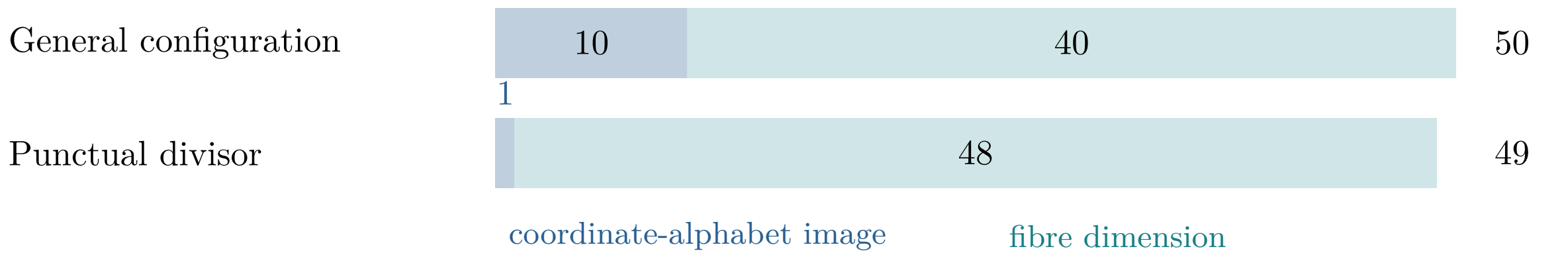
\node[anchor=west] at (0,8.45) {General configuration};
\fill[Blue!30] (5.0,8.1) rectangle (6.92,8.8);
\fill[Teal!20] (6.92,8.1) rectangle (14.6,8.8);
\node at (5.96,8.45) {$10$};\node at (10.76,8.45) {$40$};
\node[anchor=west] at (14.85,8.45) {$50$};
\node[anchor=west] at (0,7.35) {Punctual divisor};
\fill[Blue!30] (5.0,7.0) rectangle (5.192,7.7);
\fill[Teal!20] (5.192,7.0) rectangle (14.408,7.7);
\node[above,Blue] at (5.10,7.7) {$1$};\node at (9.8,7.35) {$48$};
\node[anchor=west] at (14.85,7.35) {$49$};
\node[anchor=west,font=\footnotesize,Blue] at (5,6.52) {coordinate-alphabet image};
\node[anchor=west,font=\footnotesize,Teal] at (10.0,6.52) {fibre dimension};
\draw[Ink!25] (0,6.05)--(16,6.05);
\node[anchor=west,font=\bfseries\small,Ink] at (0,5.64) {(b) The last alternating function gives top Koszul homology};
\node[anchor=west] at (.2,4.90) {$0\longrightarrow R\longrightarrow E[t]\longrightarrow M\longrightarrow0$};
\node[anchor=west] at (8.25,4.90) {$M_{12}=E/F_{12}\supset\sgn$,\quad $M_{13}=0$};
\node[anchor=west] at (.2,4.20) {$q_i=t(v_{i1}-v_{10,1}),\quad 1\le i\le9$};
\node[anchor=west] at (8.25,4.20) {$\bigwedge^9\operatorname{Std}_{10}\simeq\sgn$};
\draw[Ink!30,fill=Teal!4,rounded corners=3pt] (.2,2.78) rectangle (15.8,3.78);
\node at (8,3.28)
 {$H_8(q;R)_{21}\simeq
 H_9(q;M)_{21}\simeq
 M_{12}\otimes\bigwedge^9\operatorname{Std}_{10}\supset\mathbf1$};
\node[anchor=west,font=\footnotesize,text width=15.5cm,align=left] at (.2,2.29)
 {The normalization complex is exact in this degree. The missing sign class therefore survives in degree eight.};
\draw[Ink!25] (0,1.81)--(16,1.81);
\node[anchor=west,font=\bfseries\small,Ink] at (0,1.39) {(c) Matching sign twists preserve the invariant Tor class};
\node[font=\footnotesize] at (2.3,.61) {$\displaystyle\varinjlim_r(J^r\otimes\sgn^r,\times a)$};
\draw[-{Latex[length=2mm]},Ink!60] (4.5,.61)--(5.35,.61);
\node[font=\footnotesize] at (8.1,.61) {$eA_{10r}\simeq J^r\otimes\sgn^r$};
\draw[-{Latex[length=2mm]},Ink!60] (10.8,.61)--(11.65,.61);
\node[font=\footnotesize] at (13.8,.61) {$\Tor_8(A_\ell,\kappa)^{S_{10}}\ne0$};
\node[anchor=west,font=\footnotesize,text width=15.5cm,align=left] at (.2,-.18)
 {The left colimit is the stable chart $S[J/a]$. Both sign twists agree, so the invariant class transfers to the polygraph.};
\node[anchor=west,font=\footnotesize,text width=15.5cm,align=left] at (.2,-.85)
 {One threshold $L_{10}$ gives local projective dimension eight for every $n\ge10$, $d\ge5$ and $\ell\ge L_{10}$.};
\end{tikzpicture}
\caption{From the punctual divisor to invariant Tor.
Panel (a) compares dimensions under projection to a coordinate alphabet.
In (b), the equality $21=12+9$ places the final normalization defect
at the top of the nine-generator Koszul complex. Its sign character
cancels that of the top exterior power, giving a nonzero invariant class.
Panel (c) transfers this class through the Rees colimit and the alternating
mark idempotent $e$. Tor is taken over the localized coordinate alphabet,
with residue field $\kappa$. The mark threshold is uniform in $n\geq10$
and $d\geq5$.}
\label{new:tor-figure}
\end{figure}

Figure~\ref{new:tor-figure} relates the divisorial collision to
Koszul homology and the transfer through Rees powers. The internal
degree $21$ belongs to the transverse cone. The Rees exponent $r$
is determined by clearing the degree-zero coefficient denominators
in the blowup chart, and the alternating-block construction then
uses $nr$ marks. Thus the mark thresholds $L_{10},L_{12},L_{15}$
are uniform in the stated ranges of $n$ and $d$, with existence
proved by this finite-stage argument.

\FloatBarrier

\subsection{The entire invariant Koszul algebra}
\label{next:invariant-koszul}

The invariant class constructed above belongs to the algebra of
invariant Tor classes. Its multiplication is determined by the polynomial
filtration of the degree-$n$ point algebra.

Let $E/F$ be the $S_n$-torsor algebra of an ordered configuration of
$n$ distinct points. Write $v_{ia}\in E$ for its centered coordinates,
and put
\[
 R=F[t,tv_{ia}]\subset E[t],\qquad
 T=E^{S_{\{2,\ldots,n\}}}.
\]
The finite \'{e}tale point algebra $T$ has degree $n$ over $F$.
Give it the filtration $T_{\le r}$ by evaluations of polynomials of
degree at most $r$ in the coordinates of one point. Choose a filtered
basis
\[
 1,f_1,\ldots,f_{n-1},
 \qquad \operatorname{Tr}_{T/F}(f_i)=0,\qquad
 \deg f_i=w_i.
\]
Here a filtered basis means that the elements with $w_i\le r$,
together with $1$, form a basis of $T_{\le r}$.

Fix an alphabet $j$, write $\beta=v_{1j}\in T$, and suppose
\begin{equation}\label{next:trace-square}
 \operatorname{Tr}_{T/F}(\beta^2)\ne0.
\end{equation}
This holds for a generic configuration. We may choose
$f_1=\beta$ and $w_1=1$. Let
\[
 K=K\bigl(t(v_{2j}-v_{1j}),\ldots,
          t(v_{nj}-v_{1j});R\bigr).
\]
We give a Koszul exterior generator homological degree one and
internal degree one.

\begin{theorem}[Invariant Koszul algebra]
\label{next:invariant-koszul-algebra}
There is an isomorphism of differential graded $F[t]$-algebras
\begin{equation}\label{next:invariant-dga}
 K^{S_n}\cong
 K(t^2;F[t])\otimes_F
 \Lambda_F(\eta_2,\ldots,\eta_{n-1}),
 \qquad |\eta_i|=(1,w_i+1).
\end{equation}
The generators $\eta_i$ are closed. Consequently
\begin{equation}\label{next:invariant-tor-algebra}
 H_*(K)^{S_n}\cong
 F[t]/(t^2)\otimes_F
 \Lambda_F(\eta_2,\ldots,\eta_{n-1}).
\end{equation}
This identifies the full algebra, including its internal grading.
The invariant DGA is formal.
\end{theorem}

\begin{proof}
Begin with the permutation representation and its odd generators
$\theta_1,\ldots,\theta_n$. The super multisymmetric algebra
\[
 \bigl(F[t,X_{ia}]\otimes
       \Lambda(\theta_1,\ldots,\theta_n)\bigr)^{S_n}
\]
is generated by the even and odd power sums
\begin{equation}\label{next:super-power-sums}
 p_\alpha=\sum_iX_i^\alpha,\qquad
 \psi_\alpha=\sum_iX_i^\alpha\theta_i
 \qquad(\alpha\in\N^d).
\end{equation}
This is the one-fermionic-alphabet instance of the power-sum
generation theorem \cite[Lemma~15 and Theorem~17]{next:OrellanaZabrocki}.
One can also prove generation by symmetrizing a monomial and inducting on the number
of occupied row indices. Multiplying its one-row sums gives the
desired distinct-row orbit sum, up to a nonzero integer stabilizer
factor, together with terms occupying fewer rows. In each collision,
the polynomial factors multiply; a collision of two odd factors
vanishes. Induction therefore proves generation. A monomial whose
orbit average vanishes contributes nothing. The construction
preserves polynomial degree throughout.

The evaluation $X_{ia}\mapsto tv_{ia}$ is an equivariant
surjection onto $R$. Averaging over the finite group is exact,
so it remains surjective on the invariant exterior algebras.
Under evaluation, each even power sum belongs to $F[t]$:
\[
 p_\alpha\longmapsto
 t^{|\alpha|}\operatorname{Tr}_{T/F}(v_1^\alpha).
\]
Define
\[
 \gamma_0=\sum_i\theta_i,\qquad
 \gamma_a=t^{w_a}\sum_i f_a(v_i)\theta_i
 \quad(1\le a\le n-1).
\]
Set $f_0=1$ and $w_0=0$. The filtered-basis property gives
\[
 v_1^\alpha=\sum_{w_a\le r}c_a f_a,
 \qquad
 \psi_\alpha\longmapsto
 \sum_{w_a\le r}c_a t^{r-w_a}\gamma_a
 \quad(r=|\alpha|).
\]
The nonnegative exponents keep this expression over $F[t]$.
Hence
\begin{equation}\label{next:integral-exterior-surjection}
 F[t]\otimes\Lambda(\gamma_0,\ldots,\gamma_{n-1})
 \longrightarrow
 \bigl(R\otimes\Lambda(F^n)\bigr)^{S_n}
\end{equation}
is surjective before inverting $t$.

After inverting $t$ and passing to a splitting algebra, the
coefficient matrix of the $\gamma_a$ is the evaluation matrix of
a basis of the finite \'{e}tale point algebra. Its determinant and
all its exterior powers are invertible. Thus the kernel of
\eqref{next:integral-exterior-surjection} is $t$-torsion; the source
is free over $F[t]$, so the kernel is zero. The filtered expansion
above is what extends this isomorphism across $t=0$.

Since the configuration is centered, the differential
$d\theta_i=tv_{ij}$ kills $\gamma_0$. Quotienting by the
trivial exterior direction $\gamma_0$ gives the standard
representation. This quotient commutes with invariants because
$F^n={\bf1}\oplus\operatorname{Std}_n$ in characteristic zero.
In the basis $e_i=\theta_i-\theta_1$, $2\le i\le n$,
the differential is $de_i=t(v_{ij}-v_{1j})$, so the quotient
is the stated Koszul complex. Moreover
\begin{equation}\label{next:trace-differential}
 d\gamma_a=
 \operatorname{Tr}_{T/F}(\beta f_a)t^{w_a+1}.
\end{equation}

Put $c=\operatorname{Tr}_{T/F}(\beta^2)$ and set
\[
 \epsilon=c^{-1}\gamma_1,\qquad
 \eta_a=\gamma_a-
 \frac{\operatorname{Tr}_{T/F}(\beta f_a)}{c}
 t^{w_a-1}\gamma_1
 \quad(a\ge2).
\]
Every $w_a\ge1$, so this is an invertible triangular change
of generators over $F[t]$. It preserves internal degree.
Equation~\eqref{next:trace-differential} gives
$d\epsilon=t^2$ and $d\eta_a=0$, proving
\eqref{next:invariant-dga}.

The projection
\[
 \epsilon\longmapsto0,\qquad
 F[t]\longrightarrow F[t]/(t^2),\qquad
 \eta_a\longmapsto\eta_a
\]
is a multiplicative quasi-isomorphism onto
\eqref{next:invariant-tor-algebra} with zero differential.
This proves formality, and exactness of invariants gives
$H(K^{S_n})=H(K)^{S_n}$.
\end{proof}

\begin{corollary}[All hook components]\label{next:orbit-hooks}
Write $R=\bigoplus_{r\ge0}t^rF_r$ and
$A=R/(t)=\bigoplus_r F_r/F_{r-1}$. For $0\le k\le n-1$,
\begin{align}
 \sum_{r\ge0}[F_r:S^{(n-k,1^k)}]\,z^r
 &=\frac{e_k(z^{w_1},\ldots,z^{w_{n-1}})}{1-z},
 \label{next:filtered-hooks}\\
 \sum_{r\ge0}[A_r:S^{(n-k,1^k)}]\,z^r
 &=e_k(z^{w_1},\ldots,z^{w_{n-1}}).
 \label{next:artinian-hooks}
\end{align}
These formulas require only the filtered point basis; they
do not require the trace-square hypothesis
\eqref{next:trace-square}.
\end{corollary}

\begin{proof}
The integral exterior-algebra description preceding
\eqref{next:trace-differential} computes every invariant
coefficient space
$(R\otimes\Lambda^k\operatorname{Std}_n)^{S_n}$.
Subtracting the exterior contribution $k$ from the internal
degree gives the numerator in \eqref{next:filtered-hooks}.
The representation
$\Lambda^k\operatorname{Std}_n=S^{(n-k,1^k)}$ is self-dual,
so these invariant dimensions are its multiplicities.
Passing to $R/(t)$ removes the factor $(1-z)^{-1}$.
\end{proof}

For a point Hilbert function $(1,d,r)$, a filtered basis has
degrees $0,1^d,2^r$. Removing the selected linear coordinate
in Theorem~\ref{next:invariant-koszul-algebra} gives
\begin{equation}\label{next:invariant-series}
 \sum_{i,m}\dim_F H_i(K)_m^{S_n}\,s^iz^m
 =(1+z)(1+sz^2)^{d-1}(1+sz^3)^r.
\end{equation}
In particular, the three divisorial seeds give
\[
\begin{array}{c|c|c|c}
 n&H&\text{invariant Tor series}&
       \text{top invariant degrees}\\ \hline
 10&(1,5,4)&(1+z)(1+sz^2)^4(1+sz^3)^4&20,\ 21\\
 12&(1,5,6)&(1+z)(1+sz^2)^4(1+sz^3)^6&26,\ 27\\
 15&(1,11,3)&(1+z)(1+sz^2)^{10}(1+sz^3)^3&29,\ 30
\end{array}
\]
At ten points, the invariant group in each homological degree
$i$ has dimension $2\binom8i$. The eight invariant classes in
degree one have a nonzero product in internal degree twenty;
multiplying that product by $t$ gives the invariant line in
degree twenty-one. The hook series of the Artinian orbit
algebra is $(1+uz)^5(1+uz^2)^4$, with $u$ recording the leg
length of the hook.

\subsection{Exterior products in the polygraphs}
\label{next:polygraph-products}

Separate blocks of marks make multiplication compatible with
alternating evaluation. This retains the nonzero exterior products
when the classes pass from a blowup chart to a marked polygraph.

\begin{corollary}[Maximal nonzero products of invariant Tor classes]
\label{next:polygraph-exterior}
For each of the ten-, twelve- and fifteen-point seeds, and for
every sufficiently large number of marks, the Tor algebra of
the polygraph over the localized coordinate alphabet
contains an exterior algebra over $\kappa(\mathfrak q)$ on
respectively $8$, $10$ and $13$ invariant degree-one generators.
Their top product is nonzero.

There are uniform seed mark constants $C_{10},C_{12},C_{15}$
such that a collision stratum with $a$ ten-point blocks,
$b$ twelve-point blocks and $c$ fifteen-point blocks has an
exterior subalgebra on
\[
 8a+10b+13c
\]
invariant generators whenever
\[
 \ell\ge aC_{10}+bC_{12}+cC_{15}.
\]
The invariance group is
$(S_{10})^a\times(S_{12})^b\times(S_{15})^c$.
The ambient dimension is at least five, and at least eleven
when a fifteen-point block occurs. These products attain the
projective dimension on the stated collision strata.
\end{corollary}

\begin{proof}
At a seed of length $n$, the local orbit-cone ring $B$
has $n-2$ invariant degree-one Koszul classes with nonzero
product by Theorem~\ref{next:invariant-koszul-algebra}.
Use the stable alternating chart and invariant denominator
norms from Corollary~\ref{new:equivariant-transfer}.
A single common invariant denominator, which is a unit at
$B$, lifts all these classes to the chart. Their product
remains nonzero after multiplication by its $(n-2)$nd power.
The equivariant Rees colimit puts the finitely many classes
at stages $J^{r_1},\ldots,J^{r_{n-2}}$, with the twists
$\sgn^{r_i}$.

Give the $i$th class its own block of $nr_i$ marks. The
alternating-block evaluation maps of
\eqref{fh:j:alternating-summand} identify it with an invariant
degree-one Tor class in that polygraph. For disjoint mark
sets, there is a natural multiplication map
\begin{equation}\label{next:disjoint-mark-product}
 A_{nr_i}\otimes_S A_{nr_j}
 \longrightarrow A_{n(r_i+r_j)}
\end{equation}
induced by using their separate mark variables and common
label variables. With evaluation at the identity assignments and normalized
alternating projectors, disjoint mark sets give
\[
 \mathrm{ev}(F_iF_j)=\mathrm{ev}(F_i)\mathrm{ev}(F_j),
 \qquad J^{r_i}J^{r_j}=J^{r_i+r_j}.
\]
Thus \eqref{next:disjoint-mark-product} commutes with ideal
multiplication and with the Koszul product. Keeping the mark sets
disjoint is what places the product in the alternating summand
for their union.

Put $L=n\sum_i r_i$. In $A_L$, the product of all the
transferred classes belongs to the alternating summand for
all the disjoint blocks. Its evaluated image has the
nonzero local product just constructed. If it were a
boundary in the full polygraph complex, applying the
alternating projector, evaluation and localization would
make that local product a boundary. Therefore it is nonzero.

Write the transferred generators as $\xi_1,\ldots,\xi_{n-2}$
and their nonzero product as $\Omega$. Graded commutativity gives
$\xi_i^2=0$. For a homogeneous relation and a fixed $k$-subset $I_0$,
\[
 \sum_{|I|=k}c_I\xi_I=0
 \quad\Longrightarrow\quad
 \left(\sum_{|I|=k}c_I\xi_I\right)\xi_{I_0^c}
   =\pm c_{I_0}\Omega=0
 \quad\Longrightarrow\quad c_{I_0}=0.
\]
Here $c_I\in\kappa(\mathfrak q)$ and $\xi_I$ is the product in
increasing index order. Every other term has a repeated odd factor.
Hence the exterior algebra embeds over the collision residue field.

Forgetting extra marks has an algebra section, so the
nonzero product persists for all $\ell\ge L$. Ambient
coordinate projection with its zero section preserves
multiplication and the label action. The cluster
factorization and the faithful Koszul K\"unneth argument
then preserve the product of the top classes from the
different blocks. Choosing enough marks in each seed
block proves the uniform bound and the stated mixed
product length.
\end{proof}

The local description
\eqref{next:invariant-tor-algebra} identifies the
invariant Tor algebra of the orbit cone. The transfer
establishes the displayed exterior subalgebra of the
polygraph Tor algebra. The constants $C_{10},C_{12},C_{15}$
give uniform existence thresholds for the required marks.

\begin{figure}[!htbp]
\centering
\begin{tikzpicture}[x=1cm,y=1cm,font=\small,line cap=round,line join=round]
\path[use as bounding box] (0,0) rectangle (16.2,16.35);
\node[anchor=west,font=\bfseries,Ink] at (0,16.05)
 {(a) The point filtration determines the invariant Koszul generators};
\node at (4.0,15.28) {constant};
\node[Blue] at (7.0,15.28) {selected linear};
\node[Teal] at (10.5,15.28) {four other linear};
\node[Amber] at (14.35,15.28) {four quadratic};
\draw[Ink!35] (0,14.92)--(16.2,14.92);
\node[anchor=west] at (0,14.47) {point basis};
\node[anchor=west] at (0,13.76) {point degree};
\node[anchor=west] at (0,13.05) {Koszul image};
\node at (4.0,14.47) {$1$};
\node at (4.0,13.76) {$0$};
\node at (4.0,13.05) {$\gamma_0=0$};
\node[Blue] at (7.0,14.47) {$\beta$};
\node[Blue] at (7.0,13.76) {$1$};
\node[Blue] at (7.0,13.05) {$d\epsilon=t^2$};
\node[Teal] at (10.5,14.47) {$f_2,\ldots,f_5$};
\node[Teal] at (10.5,13.76) {$1,1,1,1$};
\node[Teal] at (10.5,13.05) {$d\eta_a=0$};
\node[Amber] at (14.35,14.47) {$q_1,\ldots,q_4$};
\node[Amber] at (14.35,13.76) {$2,2,2,2$};
\node[Amber] at (14.35,13.05) {$d\zeta_b=0$};
\draw[Ink!35] (0,12.61)--(16.2,12.61);
\node[Ink] at (8.1,11.94)
 {$H_*(K)^{S_{10}}\cong F[t]/(t^2)\otimes
   \Lambda(\eta_1,\ldots,\eta_4,\zeta_1,\ldots,\zeta_4)$};
\node[Blue] at (3.2,11.25) {$|t|=(0,1)$};
\node[Teal] at (8.0,11.25) {$|\eta_a|=(1,2)$};
\node[Amber] at (12.9,11.25) {$|\zeta_b|=(1,3)$};
\draw[Ink!25] (0,10.78)--(16.2,10.78);
\node[anchor=west,font=\bfseries,Ink] at (0,10.30)
 {(b) The top exterior product};
\node[anchor=west,font=\bfseries,Ink] at (8.2,10.30)
 {(c) Homological and internal degrees};
\fill[Teal] (0.65,9.52) circle (2.7pt);
\node[Teal] at (0.65,9.02) {$\eta_1$};
\fill[Teal] (1.53,9.52) circle (2.7pt);
\node[Teal] at (1.53,9.02) {$\eta_2$};
\fill[Teal] (2.41,9.52) circle (2.7pt);
\node[Teal] at (2.41,9.02) {$\eta_3$};
\fill[Teal] (3.29,9.52) circle (2.7pt);
\node[Teal] at (3.29,9.02) {$\eta_4$};
\fill[Amber] (4.17,9.52) circle (2.7pt);
\node[Amber] at (4.17,9.02) {$\zeta_1$};
\fill[Amber] (5.05,9.52) circle (2.7pt);
\node[Amber] at (5.05,9.02) {$\zeta_2$};
\fill[Amber] (5.93,9.52) circle (2.7pt);
\node[Amber] at (5.93,9.02) {$\zeta_3$};
\fill[Amber] (6.81,9.52) circle (2.7pt);
\node[Amber] at (6.81,9.02) {$\zeta_4$};
\node at (3.75,8.29)
 {$\Omega=\eta_1\eta_2\eta_3\eta_4\zeta_1\zeta_2\zeta_3\zeta_4\ne0$};
\node[font=\Large] at (1.65,7.37) {$\Omega$};
\node[font=\Large,Blue] at (5.65,7.37) {$t\Omega$};
\draw[-{Latex[length=2.2mm]},Blue,line width=.85pt]
 (2.35,7.37)--(4.75,7.37);
\node[Blue,above] at (3.55,7.47) {$\times t$};
\node at (1.65,6.78) {$(8,20)$};
\node[Blue] at (5.65,6.78) {$(8,21)$};
\node[align=left,anchor=west] at (0,5.72)
 {$|\Omega|=4(1,2)+4(1,3)=(8,20)$,\\[6pt]
  $|t\Omega|=(0,1)+(8,20)=(8,21)$.};
\draw[Ink!25] (0,4.76)--(7.6,4.76);
\node[align=left,anchor=north west,text width=7.4cm] at (0,4.40)
 {Each subset of the eight generators gives a nonzero exterior product.
  Disjoint alternating mark blocks transfer these products to the
  polygraph Tor algebra.};
\node at (11.25,9.65) {$k=j-2i$};
\node at (8.25,9.05) {$i$};
\node[font=\footnotesize] at (15.14,9.05) {$\dim H_i^{S_{10}}$};
\node at (9.00,9.05) {$0$};
\node at (9.90,9.05) {$1$};
\node at (10.80,9.05) {$2$};
\node at (11.70,9.05) {$3$};
\node at (12.60,9.05) {$4$};
\node at (13.50,9.05) {$5$};
\node at (8.25,8.43) {$0$};
\path[fill=Teal!6,draw=Ink!22,line width=.35pt] (8.585,8.135) rectangle (9.415,8.725);
\node[Ink] at (9.00,8.43) {$1$};
\path[fill=Teal!6,draw=Ink!22,line width=.35pt] (9.485,8.135) rectangle (10.315,8.725);
\node[Ink] at (9.90,8.43) {$1$};
\path[fill=white,draw=Ink!22,line width=.35pt] (10.385,8.135) rectangle (11.215,8.725);
\path[fill=white,draw=Ink!22,line width=.35pt] (11.285,8.135) rectangle (12.115,8.725);
\path[fill=white,draw=Ink!22,line width=.35pt] (12.185,8.135) rectangle (13.015,8.725);
\path[fill=white,draw=Ink!22,line width=.35pt] (13.085,8.135) rectangle (13.915,8.725);
\node at (15.14,8.43) {$2$};
\node at (8.25,7.77) {$1$};
\path[fill=Teal!7,draw=Ink!22,line width=.35pt] (8.585,7.475) rectangle (9.415,8.065);
\node[Ink] at (9.00,7.77) {$4$};
\path[fill=Teal!9,draw=Ink!22,line width=.35pt] (9.485,7.475) rectangle (10.315,8.065);
\node[Ink] at (9.90,7.77) {$8$};
\path[fill=Teal!7,draw=Ink!22,line width=.35pt] (10.385,7.475) rectangle (11.215,8.065);
\node[Ink] at (10.80,7.77) {$4$};
\path[fill=white,draw=Ink!22,line width=.35pt] (11.285,7.475) rectangle (12.115,8.065);
\path[fill=white,draw=Ink!22,line width=.35pt] (12.185,7.475) rectangle (13.015,8.065);
\path[fill=white,draw=Ink!22,line width=.35pt] (13.085,7.475) rectangle (13.915,8.065);
\node at (15.14,7.77) {$16$};
\node at (8.25,7.11) {$2$};
\path[fill=Teal!8,draw=Ink!22,line width=.35pt] (8.585,6.815) rectangle (9.415,7.405);
\node[Ink] at (9.00,7.11) {$6$};
\path[fill=Teal!17,draw=Ink!22,line width=.35pt] (9.485,6.815) rectangle (10.315,7.405);
\node[Ink] at (9.90,7.11) {$22$};
\path[fill=Teal!17,draw=Ink!22,line width=.35pt] (10.385,6.815) rectangle (11.215,7.405);
\node[Ink] at (10.80,7.11) {$22$};
\path[fill=Teal!8,draw=Ink!22,line width=.35pt] (11.285,6.815) rectangle (12.115,7.405);
\node[Ink] at (11.70,7.11) {$6$};
\path[fill=white,draw=Ink!22,line width=.35pt] (12.185,6.815) rectangle (13.015,7.405);
\path[fill=white,draw=Ink!22,line width=.35pt] (13.085,6.815) rectangle (13.915,7.405);
\node at (15.14,7.11) {$56$};
\node at (8.25,6.45) {$3$};
\path[fill=Teal!7,draw=Ink!22,line width=.35pt] (8.585,6.155) rectangle (9.415,6.745);
\node[Ink] at (9.00,6.45) {$4$};
\path[fill=Teal!21,draw=Ink!22,line width=.35pt] (9.485,6.155) rectangle (10.315,6.745);
\node[Ink] at (9.90,6.45) {$28$};
\path[fill=Teal!32,draw=Ink!22,line width=.35pt] (10.385,6.155) rectangle (11.215,6.745);
\node[Ink] at (10.80,6.45) {$48$};
\path[fill=Teal!21,draw=Ink!22,line width=.35pt] (11.285,6.155) rectangle (12.115,6.745);
\node[Ink] at (11.70,6.45) {$28$};
\path[fill=Teal!7,draw=Ink!22,line width=.35pt] (12.185,6.155) rectangle (13.015,6.745);
\node[Ink] at (12.60,6.45) {$4$};
\path[fill=white,draw=Ink!22,line width=.35pt] (13.085,6.155) rectangle (13.915,6.745);
\node at (15.14,6.45) {$112$};
\node at (8.25,5.79) {$4$};
\path[fill=Teal!6,draw=Ink!22,line width=.35pt] (8.585,5.495) rectangle (9.415,6.085);
\node[Ink] at (9.00,5.79) {$1$};
\path[fill=Teal!14,draw=Ink!22,line width=.35pt] (9.485,5.495) rectangle (10.315,6.085);
\node[Ink] at (9.90,5.79) {$17$};
\path[fill=Teal!34,draw=Ink!22,line width=.35pt] (10.385,5.495) rectangle (11.215,6.085);
\node[Ink] at (10.80,5.79) {$52$};
\path[fill=Teal!34,draw=Ink!22,line width=.35pt] (11.285,5.495) rectangle (12.115,6.085);
\node[Ink] at (11.70,5.79) {$52$};
\path[fill=Teal!14,draw=Ink!22,line width=.35pt] (12.185,5.495) rectangle (13.015,6.085);
\node[Ink] at (12.60,5.79) {$17$};
\path[fill=Teal!6,draw=Ink!22,line width=.35pt] (13.085,5.495) rectangle (13.915,6.085);
\node[Ink] at (13.50,5.79) {$1$};
\node at (15.14,5.79) {$140$};
\node at (8.25,5.13) {$5$};
\path[fill=white,draw=Ink!22,line width=.35pt] (8.585,4.835) rectangle (9.415,5.425);
\path[fill=Teal!7,draw=Ink!22,line width=.35pt] (9.485,4.835) rectangle (10.315,5.425);
\node[Ink] at (9.90,5.13) {$4$};
\path[fill=Teal!21,draw=Ink!22,line width=.35pt] (10.385,4.835) rectangle (11.215,5.425);
\node[Ink] at (10.80,5.13) {$28$};
\path[fill=Teal!32,draw=Ink!22,line width=.35pt] (11.285,4.835) rectangle (12.115,5.425);
\node[Ink] at (11.70,5.13) {$48$};
\path[fill=Teal!21,draw=Ink!22,line width=.35pt] (12.185,4.835) rectangle (13.015,5.425);
\node[Ink] at (12.60,5.13) {$28$};
\path[fill=Teal!7,draw=Ink!22,line width=.35pt] (13.085,4.835) rectangle (13.915,5.425);
\node[Ink] at (13.50,5.13) {$4$};
\node at (15.14,5.13) {$112$};
\node at (8.25,4.47) {$6$};
\path[fill=white,draw=Ink!22,line width=.35pt] (8.585,4.175) rectangle (9.415,4.765);
\path[fill=white,draw=Ink!22,line width=.35pt] (9.485,4.175) rectangle (10.315,4.765);
\path[fill=Teal!8,draw=Ink!22,line width=.35pt] (10.385,4.175) rectangle (11.215,4.765);
\node[Ink] at (10.80,4.47) {$6$};
\path[fill=Teal!17,draw=Ink!22,line width=.35pt] (11.285,4.175) rectangle (12.115,4.765);
\node[Ink] at (11.70,4.47) {$22$};
\path[fill=Teal!17,draw=Ink!22,line width=.35pt] (12.185,4.175) rectangle (13.015,4.765);
\node[Ink] at (12.60,4.47) {$22$};
\path[fill=Teal!8,draw=Ink!22,line width=.35pt] (13.085,4.175) rectangle (13.915,4.765);
\node[Ink] at (13.50,4.47) {$6$};
\node at (15.14,4.47) {$56$};
\node at (8.25,3.81) {$7$};
\path[fill=white,draw=Ink!22,line width=.35pt] (8.585,3.515) rectangle (9.415,4.105);
\path[fill=white,draw=Ink!22,line width=.35pt] (9.485,3.515) rectangle (10.315,4.105);
\path[fill=white,draw=Ink!22,line width=.35pt] (10.385,3.515) rectangle (11.215,4.105);
\path[fill=Teal!7,draw=Ink!22,line width=.35pt] (11.285,3.515) rectangle (12.115,4.105);
\node[Ink] at (11.70,3.81) {$4$};
\path[fill=Teal!9,draw=Ink!22,line width=.35pt] (12.185,3.515) rectangle (13.015,4.105);
\node[Ink] at (12.60,3.81) {$8$};
\path[fill=Teal!7,draw=Ink!22,line width=.35pt] (13.085,3.515) rectangle (13.915,4.105);
\node[Ink] at (13.50,3.81) {$4$};
\node at (15.14,3.81) {$16$};
\node at (8.25,3.15) {$8$};
\path[fill=white,draw=Ink!22,line width=.35pt] (8.585,2.855) rectangle (9.415,3.445);
\path[fill=white,draw=Ink!22,line width=.35pt] (9.485,2.855) rectangle (10.315,3.445);
\path[fill=white,draw=Ink!22,line width=.35pt] (10.385,2.855) rectangle (11.215,3.445);
\path[fill=white,draw=Ink!22,line width=.35pt] (11.285,2.855) rectangle (12.115,3.445);
\path[fill=Teal!6,draw=Blue,line width=.85pt] (12.185,2.855) rectangle (13.015,3.445);
\node[Ink] at (12.60,3.15) {$1$};
\path[fill=Teal!6,draw=Blue,line width=.85pt] (13.085,2.855) rectangle (13.915,3.445);
\node[Ink] at (13.50,3.15) {$1$};
\node at (15.14,3.15) {$2$};
\draw[Ink!35] (14.08,8.78)--(14.08,2.82);
\node[font=\footnotesize,align=center] at (11.95,2.38)
 {entry $=[s^iz^{2i+k}](1+z)(1+sz^2)^4(1+sz^3)^4$};
\draw[Ink!25] (0,1.81)--(16.2,1.81);
\node[align=center] at (8.1,1.14)
 {$\Lambda_{\kappa(\mathfrak q)}(\xi_1,\ldots,\xi_8)
   \hookrightarrow\operatorname{Tor}_*(A_\ell,\kappa(\mathfrak q))^{S_{10}}$};
\node[font=\footnotesize,Ink!90] at (8.1,.49)
 {For sufficiently many marks; Tor is over $(P_j)_{\mathfrak q}$ and $A_\ell$ is localized at $\mathfrak q$.};
\end{tikzpicture}
\caption{The invariant Koszul algebra at the generic ten-point divisor.
The point Hilbert function $(1,5,4)$ gives eight closed generators after
removing the constant exterior direction and imposing $d\epsilon=t^2$.
Their bidegrees give the nonzero classes $\Omega$ and $t\Omega$ in
$(8,20)$ and $(8,21)$. In (c), row $i$ records homological degree,
and column $k$ records internal degree $j=2i+k$; blank cells are zero.
Row sums are $2\binom8i$. The displayed internal degrees belong to the
local orbit cone. For sufficiently many marks, disjoint alternating
mark blocks give the indicated exterior subalgebra of polygraph Tor.}
\label{next:invariant-tor-figure}
\end{figure}

\subsection{Persistent failure of integral closure}
\begin{proposition}\label{new:integral-powers}
For every pair $(n,d)$ in the ranges of Theorem~\ref{fh:j:main},
all sufficiently large powers of the alternating ideal satisfy
\eqref{fh:j:ideal}. More generally this holds whenever the ordered
cover has a nonseminormal alternating chart.
\end{proposition}
\begin{proof}
Choose a nonseminormal affine chart $A=S[J/a]$. The
square-and-cube criterion gives
\[
 u\in\operatorname{Frac}(S)\setminus A,
 \qquad u^2,u^3\in A.
\]
Since $u$ is integral over $S[1/a]$ and this ring is normal,
$u=f/a^r$. Increasing $r$ clears the two denominators, giving
$f^2\in J^{2r}$ and $f^3\in J^{3r}$. Hence, for every $m\geq0$,
\[
 fa^m\notin J^{r+m},\qquad
 (fa^m)^2\in J^{2(r+m)},\qquad
 (fa^m)^3\in J^{3(r+m)}.
\]
Thus $fa^m\in\overline{J^{r+m}}\setminus J^{r+m}$.
Since $J\subset I_{ij}$ and each power of a coordinate ideal
$I_{ij}$ is integrally closed, \eqref{fh:j:ideal} follows.
\end{proof}
\Needspace{10\baselineskip}
\section{The smallest length for divisorial failure}
\label{new:small}
We prove Theorem~\ref{new:sharp} by bounding punctual strata.
Suppose a punctual family has length $\ell$, embedding dimension
$r$ and dimension $a$ with its support fixed. The graph-bundle
description gives its moving-support dimension in $\A^d$ as at most
$r+a+\ell(d-r)$. Thus its codimension in the principal component
is at least
\begin{equation}\label{new:codim}
 c=r(\ell-1)-a.
\end{equation}
Codimensions add for disjoint support factors.

For $\ell\le8$, Jelisiejew--Keneshlou prove
$a\le(r-1)(\ell-1)$ \cite[Theorem~1.2]{new:JK}. The same bound
holds for $r\le3$ and $\ell\le11$ by their Theorem~1.3, together
with the classical one- and two-dimensional cases. Their
Theorem~5.7 gives the bound for strata with $H(2)\le2$.
The square-zero case has $a=0$ and can also be handled directly.
Hence these strata have codimension at least $\ell-1$.

At length nine, the remaining Hilbert functions and dimension
bounds are
\[
\begin{array}{c|c|c}
 H & a\text{ at most}& c\text{ at least}\\\hline
 (1,4,4)&24&8\\
 (1,5,3)&36&4\\
 (1,4,3,1)&24&8
\end{array}
\]
The first two rows are Grassmannian dimensions. For the last,
a homogeneous inverse system is specified by a cubic $F$ up to
scalar and a three-dimensional quadratic space containing its first
partials. If $F$ has $t$ essential variables, where $1\le t\le3$,
the parameter count is
\[
 t(4-t)+\binom{t+2}{3}-1+(3-t)(10-3)\le17.
\]
The fibre of the associated-graded map has dimension seven by
\cite[Proposition~5.11]{new:JK}, giving $a\le24$.
The length bound and Macaulay growth leave no additional Hilbert
functions with $r\ge4$ and $H(2)\ge3$.

It follows that every boundary divisor at total length at most
nine has one length-two local factor and all other factors simple.
The ordered length-two cover is
\[
 X_{2,d}=\operatorname{Bl}_{\Delta}(\A^d\times\A^d),
\]
which is smooth. Its centered factor is the total space of
$\mathcal O_{\mathbb P^{d-1}}(-1)$. The disjoint-support
decomposition is an \'etale local product, so the ordered cover is
smooth over the generic point of every such divisor. The distinct
configuration locus is also smooth. Therefore $X_{n,d}$ satisfies
$R_1$ for every $n\le9$ and every $d$.

For length ten and $d\le4$, every local factor has embedding
dimension at most four. Proper local factors are covered by the
preceding argument. At length ten, $r\le3$ and $H(2)\le2$ are
covered by the same theorems. The remaining bounds are
\begin{equation}\label{new:ten-table}
\begin{array}{c|c|c}
H&a\text{ at most}&c\text{ at least}\\\hline
(1,4,5)&25&11\\
(1,4,4,1)&27&9\\
(1,4,3,2)&27&9\\
(1,4,3,1,1)&34&2
\end{array}
\end{equation}
These counts also occur in \cite[Proposition~5.1]{new:Michalek};
we give their ingredients. The first is
$\dim\Gr(5,10)=25$. For $(1,4,4,1)$, the homogeneous inverse-system
count has maximum $21$ and the associated-graded fibre has dimension
six. For $(1,4,3,2)$, the homogeneous locus has components of
dimensions $13$ and $11$ by \cite[Proposition~2.6]{new:GKS};
the associated-graded fibre has dimension $2(10-3)=14$ by
\cite[Proposition~5.11]{new:JK}.

For $(1,4,3,1,1)$, the highest inverse-system pieces are generated
by $L^4$ and $L^3$. The quadratic part is a three-dimensional
subspace containing $L^2$. Its parameter count is
$\dim\mathbb P^3+\dim\Gr(2,9)=17$. The homogeneous ideal has seven
quadratic minimal generators and at most three additional cubic
generators: the quotient by its quadrics has degree-three dimension
at most four, by Macaulay growth, while the prescribed quotient has
degree-three dimension one. A filtered lift permits two higher
coefficients for each quadratic generator and one for each cubic
generator, giving at most another $7\cdot2+3=17$ parameters.
Generators of degree at least four have no higher tails.
This proves the final bound in \eqref{new:ten-table}.

All punctual length-ten strata in ambient dimension at most four
therefore have codimension at least two. The only boundary divisors
again consist of one double point and simple points, and their
ordered covers are smooth. Thus $X_{10,d}$ satisfies $R_1$ for
$d\le4$. Theorem~\ref{fh:j:main} gives divisorial nonnormality
for $n=10,d\ge5$, completing the proof of
Theorem~\ref{new:sharp}.

\Needspace{10\baselineskip}
\section{A multiplication chart for every ambient dimension}
\label{sec:all-dimensional-chart}

Let $V=\C^{n-1}$. Commutative unital algebra structures on
$\C1\oplus V$ have multiplication
\[
 e_i e_j=a_{ij}1+\sum_k b_{ij}^ke_k,
 \qquad a_{ij}=a_{ji},\quad b_{ij}^k=b_{ji}^k,
\]
subject to associativity. Impose
$\Tr(L_{e_i})=\sum_jb_{ij}^j=0$, and let $C_n$ be the reduced
closure of the semisimple laws. This is the centered principal
border-basis chart with basis $1,x_1,\ldots,x_{n-1}$.
Its reduced ordered cover $\widetilde C_n$ has dimension $(n-1)^2$.
The dilation
\[
 (a_{ij},b_{ij}^k,(p_\nu)_i)
 \longmapsto(t^2a_{ij},tb_{ij}^k,t(p_\nu)_i)
\]
contracts it to the vertex $0$, the square-zero algebra
$\C[x_1,\ldots,x_{n-1}]/(x_1,\ldots,x_{n-1})^2$
with its unique ordered support.

\begin{theorem}\label{fh:j:cone}
Let $n\geq3$, $d\geq1$ and $q=\max(d,n-1)$. For every closed
point $x\in X_{n,d}$, there is $c\in\widetilde C_n$ such that
\begin{equation}\label{fh:j:completion}
 \widehat{\mathcal O}_{X_{n,d},x}[[u_1,\ldots,u_{n(q-d)}]]
 \cong
 \widehat{\mathcal O}_{\widetilde C_n,c}
 [[v_1,\ldots,v_{n(q-n+1)+n-1}]].
\end{equation}
Consequently, at fixed $n$, normality or Cohen--Macaulayness of
$X_{n,d}$ for every $d$ is equivalent to the corresponding property
of $\mathcal O_{\widetilde C_n,0}$.
\end{theorem}
\begin{proof}
An extra coordinate is an element of the rank-$n$ quotient algebra.
Its coefficients in a quotient basis and its ordered evaluations
lift the graph bundle to the ordered cover, as in
Section~\ref{new:ten}. The formulas agree on distinct configurations
and hence on their reduced closures. At the embedding $x'$ with
zero extra coordinates,
\[
 \widehat{\mathcal O}_{X_{n,q},x'}
 \cong\widehat{\mathcal O}_{X_{n,d},x}
 [[u_1,\ldots,u_{n(q-d)}]].
\]
After renumbering coordinates, choose a minimal algebra-generating
set of coordinate images $a_1,\ldots,a_r$. Minimality makes
$1,a_1,\ldots,a_r$ linearly independent. Extend these elements to
a basis by polynomials $h_{r+1}(a),\ldots,h_{n-1}(a)$, and write
the remaining coordinate images as $g_j(a)$. The triangular
polynomial automorphism fixing $x_1,\ldots,x_r$ and defined by
\[
 x_j\longmapsto
 \begin{cases}
 x_j-g_j(x_1,\ldots,x_r)+h_j(x_1,\ldots,x_r),&r<j<n,\\
 x_j-g_j(x_1,\ldots,x_r),&n\leq j\leq q
 \end{cases}
\]
puts the point in the basis chart with zero extra coordinate images.
It acts on each ordered point. Splitting off translation gives the chart
\[
 \A^{n-1}\times\widetilde C_n\times\A^{n(q-n+1)}.
\]
This proves \eqref{fh:j:completion}; both sides have dimension
$nq=(n-1)^2+n(q-n+1)+n-1$.

The normal and Cohen--Macaulay loci of $\widetilde C_n$ are open
and dilation invariant. An invariant open set containing $0$
meets every contracting orbit at a nonzero parameter and therefore
contains the whole cone. Completion and adjoining formal variables
preserve and reflect these properties for excellent rings, so
\eqref{fh:j:completion} proves the assertion. The converse uses
the displayed chart in $X_{n,n-1}$.
\end{proof}

\begin{corollary}\label{fh:j:vertex}
The vertex of $\widetilde C_n$ has the maximum Cohen--Macaulay
defect among all local rings of all $X_{n,d}$ at length $n$.
For $n\geq10$, the vertex is nonnormal, nonseminormal, non-Gorenstein,
and not Du Bois. In particular, it is not a local complete intersection.
\end{corollary}
\begin{proof}
On a reduced equidimensional affine scheme, Cohen--Macaulay defect
is upper semicontinuous at closed points: embed in affine space
and use upper semicontinuity of projective dimension. Dilation
therefore maximizes the defect at $0$.
Equation~\eqref{fh:j:completion} preserves it, and localization
cannot increase its maximum.

For $n\geq10$, Theorems~\ref{fh:j:main} and
\ref{new:not-gorenstein} give nonnormal, nonseminormal and
non-Gorenstein points of $X_{n,5}$. The graph and formal local
comparisons transfer these failures to $\widetilde C_n$.
Their loci are closed and dilation invariant, so they contain $0$.
Du Bois local rings are seminormal, and local complete intersections
are Gorenstein, giving the remaining assertions.
\end{proof}

\begin{corollary}[The sublinear defect criterion]
\label{fh:j:defect-growth}
Put
\[
 \delta_n=\dim\mathcal O_{\widetilde C_n,0}
 -\operatorname{depth}\mathcal O_{\widetilde C_n,0}
 \quad(n\geq3),\qquad\delta_1=\delta_2=0.
\]
Then
\[
 \delta_{n+m}\geq\delta_n+\delta_m,\qquad
 \lim_{n\to\infty}\frac{\delta_n}{n}
 =\sup_{n\geq1}\frac{\delta_n}{n}\in[0,+\infty].
\]
All the reduced isospectral Hilbert schemes $X_{n,d}$ are
Cohen--Macaulay if and only if $\delta_{n_j}=o(n_j)$ along some
unbounded sequence. If $\delta_N>0$, then
$\delta_n\geq\lfloor n/N\rfloor\delta_N$ for every $n$.
In particular, $\delta_9\geq1$ and
$\delta_n\geq\lfloor n/9\rfloor$ for every $n$.
\end{corollary}
\begin{proof}
The length-one cover and the ordered length-two blowup are smooth
(Section~\ref{new:small}). Use graph coordinates to place points
attaining $\delta_n,\delta_m$ in one ambient space, then translate
their supports apart. The support idempotents lift uniquely in
infinitesimal deformations, so the local ordered covers split as
products. Their distinct-configuration loci also factor, which
identifies the reduced principal closures. For their completed
local rings $B,B'$, dimension and depth add:
\[
 \operatorname{cmd}(B\widehat\otimes_{\C}B')
 =\operatorname{cmd}B+\operatorname{cmd}B'
 =\delta_n+\delta_m\leq\delta_{n+m},
\]
where $\operatorname{cmd}B=\dim B-\operatorname{depth}B$.
Iteration, with reduced points added, gives
$\delta_n\geq\lfloor n/N\rfloor\delta_N$. Hence
\[
 \liminf_{n\to\infty}\frac{\delta_n}{n}
 \geq\frac{\delta_N}{N}\quad\text{for every }N,
\]
proving the limit formula. The same bound makes sublinear growth
along an unbounded sequence equivalent to $\delta_N=0$ for every
$N$. Theorem~\ref{fh:j:cone} gives the stated criterion.

Lee proves that $H_{9,8}$ is not locally Cohen--Macaulay
\cite[Theorem~A]{dim:Lee}. Since $H_{9,8}=X_{9,8}/S_9$,
Cohen--Macaulayness of $X_{9,8}$ would descend: a parameter
sequence is regular upstairs, and Reynolds averaging descends
its regularity to the invariant ring. Thus $X_{9,8}$ is not
Cohen--Macaulay, so $\delta_9\geq1$. The preceding bound gives
$\delta_n\geq\lfloor n/9\rfloor$.
\end{proof}

\appendix
\Needspace{10\baselineskip}
\section{An explicit rational point and its tangent representation}
The main proof uses the published smooth point and the radial
valuation argument. The following computation gives a rational
witness and independently verifies the primitive normal character.
All matrices are constructed from the eleven displayed quadrics;
the accompanying rational-arithmetic script retains their rows and ranks.
\begin{lemma}[The normal dilation weight]\label{new:weight}
Let $H^0$ be the centered principal component and let
$D^0=\Gr(4,15)\subset H^0$ be its homogeneous $(1,5,4)$ divisor.
The divisor is fixed pointwise by coordinate dilation. At its
generic point the normal line has dilation weight $1$.
\end{lemma}
\begin{proof}
We use the rational witness of \cite[Appendix~A]{new:Michalek},
obtained from the smoothing in \cite[Proposition~4.2]{new:GKS}.
In $\C[a,b,c,d,e]$ its ideal has the eleven quadratic generators
\[
\begin{gathered}
ae-a^2,\quad 3b^2-ab-10ac,\quad bc-3ad,\\
3bd-8a^2+ad,\quad 3be-a^2+4ab,\quad c^2-3a^2,\\
cd-a^2+2ab,\quad ce+ac,\quad 3d^2+ab+4ac,\\
3de-ac+2ad,\quad 3e^2+a^2-2ad.
\end{gathered}
\]
The quadratic quotient has basis $a^2,ab,ac,ad$. Multiplication
of these eleven relations by the five variables has rank $35$,
so they generate all cubics and have twenty independent linear
syzygies. The tangent space is $\Hom(I,\C[a,b,c,d,e]/I)$.
Its weight-zero part has dimension at most $11\cdot4=44$.
For weight one, the twenty linear syzygies impose a matrix of
rank $49$ on the $55$ possible linear images of the relations.
For weight two they impose a matrix of rank $11$ on the eleven
possible scalar images. Consequently the total tangent dimension
is at most $44+6=50$. Rational row reduction gives these ranks; the accompanying
script constructs the matrices from the displayed generators.

The explicit smoothing places this point in the fifty-dimensional
principal component. The tangent bound is therefore attained,
and the point is smooth, with tangent representation
$\C^{44}_{(0)}\oplus\C^6_{(1)}$. Translation by the trace
barycenter splits $H_{10,5}\simeq\A^5\times H^0$ and removes
five independent weight-one directions. Thus the normal line to
$D^0$ has weight one at this point. On the smooth locus of the
irreducible fixed divisor its character is locally constant, so
the generic normal line has the same weight.
\end{proof}

\section*{Computational certificates}
\label{sec:certificates}
The ancillary files reconstruct the matrices used in the finite rank
calculations from the configurations and formulas in the paper.
Finite rank calculations use prime fields; coefficient and band estimates
use integer or rational arithmetic. Since the matrix coefficients are rational
and every denominator used in row reduction is invertible at the chosen
prime, a nonzero minor modulo that prime gives the same rank lower bound
in characteristic zero. The geometric dimension bounds and the displayed
translation and dilation vectors give the matching upper bounds where
equality is required.

\begin{center}
\small
\begin{tabularx}{\textwidth}{@{}>{\raggedright\arraybackslash}p{.47\textwidth}>{\raggedright\arraybackslash}X@{}}
\toprule
\textbf{Ancillary file}&\textbf{Calculation}\\
\midrule
\texttt{verify\_ten\_point\_witness.py}&
The ten-point smoothing and tangent ranks $49$ and $11$.\\
\texttt{verify\_dimension\_three.py}&
Collision differentials at lengths $96,\ldots,105$.\\
\texttt{verify\_dimension\_four.py}&
Collision differentials at lengths $20,\ldots,26$ and the negative tangent at length $20$.\\
\texttt{verify\_twenty\_obstruction.py}&
The rank-$154$ quadratic obstruction used in the bound $h\leq2^{30}$.\\
\texttt{verify\_tangents\_twelve\_fifteen.py}&
The negative tangents at lengths $12$ and $15$.\\
\texttt{verify\_source\_bounds.py}&
Finite coefficient bounds and the fixed-dimension exponential estimates.\\
\texttt{verify\_band\_bounds.py}&
The rational inequalities in the alternating degree-band estimates.\\
\bottomrule
\end{tabularx}
\end{center}
The ancillary README gives the commands and software requirements.

\end{document}